\documentclass[reqno,11pt,a4paper]{amsart}
\usepackage{amsmath}
\usepackage{amssymb}
\usepackage{graphicx}
\usepackage{mathrsfs}
\usepackage{amssymb}
\usepackage{amsfonts,bm,bbm}
\usepackage{amsthm}
\usepackage[colorlinks,linkcolor=blue,anchorcolor=blue,citecolor=blue,urlcolor=purple]{hyperref}
\usepackage{fullpage}
\usepackage{color}
\usepackage{enumitem}
\usepackage{underscore}
\usepackage{cite}
\usepackage{doi}
\DeclareFontFamily{U}{mathx}{}
\DeclareFontShape{U}{mathx}{m}{n}{<-> mathx10}{}
\DeclareSymbolFont{mathx}{U}{mathx}{m}{n}
\DeclareMathAccent{\widecheck}{0}{mathx}{"71}

\usepackage[open,openlevel=2]{bookmark}
\numberwithin{equation}{section}
\newtheorem{Thm}{Theorem}[section]
\newtheorem{Def}[Thm]{Definition}
\newtheorem{Lem}[Thm]{Lemma}
\newtheorem{Coro}[Thm]{Corollary}

\theoremstyle{remark}
\newtheorem{remark}{Remark}[section]

\newcommand{\R}{\mathbb{R}}
\newcommand{\T}{\mathbb{T}}

\renewcommand{\S}{\mathbb{S}}

\newcommand{\vertiii}[1]{{\left\vert\kern-0.25ex\left\vert\kern-0.25ex\left\vert #1 \right\vert\kern-0.25ex\right\vert\kern-0.25ex\right\vert}}

\usepackage{tikz,pgfplots}
\usepackage{tikz-cd}
\usepackage{tikz-3dplot}
\usepackage{amsmath,amssymb,amsthm,amsfonts,bm}
\usepackage{color}
\usepackage{tcolorbox}
\usepackage{pdfpages}
\usepackage{textcomp}
\usetikzlibrary{patterns}
\usetikzlibrary{shapes.geometric}
\usetikzlibrary{arrows.meta,arrows}
\usetikzlibrary{decorations.pathreplacing,decorations.pathmorphing}
\usetikzlibrary{hobby}
\usetikzlibrary{math}
\usepgfplotslibrary{fillbetween}
\usepackage{graphicx}
\usepackage{multicol}
\usepackage{float}
\usepackage{subfigure}

\allowdisplaybreaks

\makeatletter
\def\@tocline#1#2#3#4#5#6#7{\relax
	\ifnum #1>\c@tocdepth
	\else
	\par \addpenalty\@secpenalty\addvspace{#2}
	\begingroup \hyphenpenalty\@M
	\@ifempty{#4}{
		\@tempdima\csname r@tocindent\number#1\endcsname\relax
	}{
		\@tempdima#4\relax
	}
	\parindent\z@ \leftskip#3\relax \advance\leftskip\@tempdima\relax
	\rightskip\@pnumwidth plus4em \parfillskip-\@pnumwidth
	#5\leavevmode\hskip-\@tempdima
	\ifcase #1
	\or\or \hskip 1em \or \hskip 2em \else \hskip 3em \fi
	#6\nobreak\relax
	\hfill\hbox to\@pnumwidth{\@tocpagenum{#7}}\par
	\nobreak
	\endgroup
	\fi}
\makeatother

\begin{document}

\title[$L^p_{v}L^\infty_{x}$ solution of the Boltzmann Equation]
{Global Solutions in $L^p_{v} L^\infty_{x}$ for the Boltzmann Equation in Bounded Domains}

\author{Dingqun Deng}

\address{Graduate School of Engineering Science, Akita University, Japan }
\email{dingqun.deng@s.akita-u.ac.jp}

\author{Jong-in Kim}

\address{Department of Mathematics, Pohang University of Science and Technology, South Korea }
\email{kimjim@postech.ac.kr}

\author{Donghyun Lee}
\address{Department of Mathematics, Pohang University of Science and Technology, South Korea }
\email{donglee@postech.ac.kr}

\begin{abstract}	
	The existence theory for solutions to the Boltzmann equation in bounded domains has primarily been developed  within uniformly bounded function classes, such as $L^{\infty}_{x,v}$, as in \cite{Duan2017a, Duan2019b, GuoDecay, Guo2010b}. In this paper, we investigate solutions in relaxed function spaces $L^{p}_{v}L^\infty_{x}$ for the initial-boundary value problem of the Boltzmann equation in bounded domains. We consider the case of hard potential under diffuse reflection boundary conditions and assume cutoff model. For large initial data in a weighted $L^{p}_{v}L^\infty_{x}$ space with small relative entropy, we construct unique global-in-time mild solution that converge exponentially to the global Maxwellian. A pointwise estimate for the gain term, bounded in terms of $L^p_v$ and $L^2_v$ norms, is essential to prove our main results. Relative to \cite{Gualdani2017}, our work provides an alternative perspective on convergence to equilibrium in the presence of boundary conditions. 	
\end{abstract}

\maketitle
\tableofcontents

\vspace{0.2cm}
\section{Introduction}
\subsection{Boltzmann equation}
We consider a rarefied gas in a domain $\Omega$ whose boundaries are kept at constant temperature. We assume that $\Omega$ is a bounded open set in $\R^3$ with a $C^1$ boundary $\partial \Omega$.  
Our governed equation is the Boltzmann equation 
\begin{align} \label{Boltzmanneq}
	\partial_t F+v\cdot\nabla_xF =Q(F,F), \quad (t,x,v) \in [0, \infty) \times \Omega \times \R^3 
\end{align}
with initial data
\begin{align} \label{initialdatacond}
	F(0,x,v) = F_0(x,v)
\end{align}
and with some boundary condition. Here, $F=F(t,x,v)$ means the density distribution function of gas particles with position $x\in \Omega$ and velocity $v \in \R^3$ at time $t\ge 0$. The collision operator $Q$ has the bilinear form
\begin{align*}
	Q(F_1,F_2) = \int_{\R^3}\int_{\S^2}B(v-u,\omega)[F_1(v')F_2(u')-F_1(v)F_2(u)]d\omega du.
\end{align*}
Here, the post-collision velocity pair $(v',u')$ and the pre-collision velocity pair $(v,u)$ satisfy the relation
\begin{align} \label{postcollisionrep}
	v' = \frac{v+u}{2}+\frac{|v-u|}{2}\omega, \quad u'=\frac{v+u}{2}-\frac{|v-u|}{2}\omega
\end{align}
with $\omega \in \S^2$, according to the conservation of momentum and energy of two particles before and after the collision
\begin{align} \label{momentumenergyconserv}
	v+u=v'+u', \quad       |v|^2+|u|^2=|v'|^2+|u'|^2.
\end{align}
The collision kernel $B(v-u,\omega)$ is composed of the relative velocity $|v-u|$ and the angular part $\cos\theta: = \frac{(v-u)\cdot \omega}{|v-u|}$. Throughout this paper, we assume that it has the form 
\begin{align*}
	B(v-u,\omega) = b(\cos \theta)|v-u|^\gamma, \quad 0\le b(\cos \theta) \le C_b,
\end{align*}
where $0\le \gamma \le 1$ and $0\le \theta \le \frac{\pi}{2} $, which represents the hard potential model with angular cutoff. For this assumption, we can write the collision operator $Q$ as
\begin{align*}
	Q(F_1,F_2) &= \int_{\R^3}\int_{\S^2}B(v-u,\omega)F_1(v')F_2(u')d\omega du -\int_{\R^3}\int_{\S^2}B(v-u,\omega)F_1(v)F_2(u)d\omega du\\
	& =: Q^+(F_1,F_2) - Q^-(F_1,F_2),
\end{align*}
where $ Q^+(F_1,F_2)$ and $Q^-(F_1,F_2)$ mean the gain term and the loss term, respectively.
Note that the Boltzmann equation \eqref{Boltzmanneq} has the normalized global Maxwellian 
\begin{align*}
	\mu(v) := \frac{1}{(2\pi)^{3/2}}e^{-\frac{|v|^2}{2}}.
\end{align*}

\bigskip

\subsection{Domain and Boundary condition}
Throughout this paper, we assume $\Omega: = \{x \in \R^3 : \xi(x) <0 \}$ is connected and bounded, where $\xi$ is a $C^1$ function. Suppose that $\nabla_x \xi(x) \not= 0$ at the boundary $\partial \Omega=\{x: \xi(x) =0\}$, and the outward normal vector at $x \in \partial \Omega$ is given by $n(x) = \frac{\nabla_x \xi(x)}{|\nabla_x \xi(x)|}$. We denote the phase boundary in the space $\Omega \times \mathbb{R}^3$ by $\gamma:=\partial \Omega \times \mathbb{R}^3$. We decompose $\gamma$ into the outgoing set $\gamma_+$, the grazing set $\gamma_0$, and the incoming set $\gamma_-$ :
\begin{align*}
	&\gamma_+=\{(x,v) \in \partial \Omega \times \mathbb{R}^3 : n(x) \cdot v>0 \},\\
	&\gamma_0=\{(x,v) \in \partial \Omega \times \mathbb{R}^3 : n(x) \cdot v=0 \},\\
	&\gamma_-=\{(x,v) \in \partial \Omega \times \mathbb{R}^3 : n(x) \cdot v<0 \}.
\end{align*}
 In this paper, we consider the diffuse reflection boundary condition
\begin{align} \label{Diffuseboundarycond}
	F(t,x,v)|_{\gamma_-} = c_\mu \mu(v) \int_{n(x) \cdot v'>0}F(t,x,v')\{n(x) \cdot v'\}dv',
\end{align}
where the constant $c_\mu$ is given by
\begin{align} \label{diffusemassconvconst}
	c_\mu \int_{n(x) \cdot v'>0}\mu(v')\{n(x) \cdot v'\}dv'=1.
\end{align}
Note that $\mu(v)$ satisfies the boundary condition \eqref{Diffuseboundarycond}. For the diffuse reflection boundary condition \eqref{Diffuseboundarycond}, it holds that for $t>0$,
\begin{align*}
	\int_{\Omega \times \R^3} F(t,x,v)dxdv = \int_{\Omega \times \R^3} F_0(x,v)dxdv
\end{align*}
for any solution to the equation \eqref{Boltzmanneq} with the conditions \eqref{initialdatacond}, \eqref{Diffuseboundarycond}, which means that the mass is conserved over time.

\subsection{Back-time cycle } \label{Backtimecycle}

Given $(t,x,v) \in [0,\infty) \times \Omega \times \mathbb{R}^3$, the backward characteristic $[X(s;t,x,v),V(s;t,x,v)]$ for the Boltzmann equation \eqref{Boltzmanneq} is determined by
\begin{align} \label{Backwardcharacteristic}
	\frac{dX(s;t,x,v)}{ds} = V(s;t,x,v), \quad \frac{dV(s;t,x,v)}{ds} = 0
\end{align}
with the condition $[X(t;t,x,v),V(t;t,x,v)] = [x,v]$.\\
\begin{Def} \label{Backwardexit}
	(Backward exit time) For $(x,v) \in \bar{\Omega} \times \mathbb{R}^3$ with $v \not=0$, we define its backward exit time $t_\mathbf{b}(x,v)>0$ to be the last moment at which the back-time characteristic straight line $[X(s;t,x,v),V(s;t,x,v)]$ remains in the interior of $\Omega$ : 
	\begin{align*}
		t_\mathbf{b}(x,v) = \sup\{\tau>0 : X(s;t,x,v) \in \Omega \ \text{for all }t-\tau \le s \le t \}.
	\end{align*}
 For any $(x,v)$, we use $t_\mathbf{b}(x,v)$ whenever it is well-defined. We have $x-t_\mathbf{b}v\in \partial \Omega$ and $\xi(x-t_\mathbf{b}v) = 0$. We also define
	\begin{align*}
		x_\mathbf{b}(x,v):=x-t_\mathbf{b}v \in \partial \Omega.
	\end{align*}
\end{Def}
For each $x \in \partial \Omega$, we introduce the velocity space for the outgoing particles
\begin{align*}
	\mathcal{V}(x):= \{v\in \mathbb{R}^3 : n(x) \cdot v >0\}
\end{align*}
with the probability measure $d\sigma=d\sigma(x)$ by
\begin{align*}
	d\sigma(x) := c_\mu \mu(v') \{n(x) \cdot v'\}dv'.
\end{align*}
Fix any point $(x,v) \notin \gamma_0 \cup \gamma_-$, and let $(t_0,x_0,v_0) = (t,x,v)$ and $k \ge 1$. For $v_{k+1} \in \mathcal{V}_{k+1}:=\{v\in \mathbb{R}^3 : n(x_{k+1}) \cdot v >0\}$, we define the back-time cycle as
\begin{equation} \label{Backcycle}
	\begin{aligned}
	&X_{\mathbf{cl}}(s;t,x,v) = \sum_k \mathbf{1}_{[t_{k+1},t_k)}(s) \{x_k-v_k(t_k-s)\}, \\
	&V_{\mathbf{cl}}(s;t,x,v) = \sum_k \mathbf{1}_{[t_{k+1},t_k)}(s) v_k,
\end{aligned}
\end{equation}
with
\begin{align*}
	(t_{k+1},x_{k+1},v_{k+1}) = (t_k-t_\mathbf{b}(x_k,v_k), x_\mathbf{b}(x_k,v_k),v_{k+1}).
\end{align*}
We note that each $v_j$ is an independent variable, and $t_k, x_k$ depend on $t_j, x_j, v_j$ for $1\le j \le k-1$. However, the phase space $\mathcal{V}_j$ depends on $(t,x,v,v_1,v_2,\cdots, v_{j-1})$.\\
For $k\ge 2$, we define the iterated integral 
\begin{align} \label{iteratedintegral}
	\int_{\prod_{j=1}^{k-1}\mathcal{V}_j}\prod_{j=1}^{k-1}d\sigma_j \equiv \int_{\mathcal{V}_1} \int_{\mathcal{V}_2} \cdots\left\{ \int_{\mathcal{V}_{k-1}}d\sigma_{k-1} \right\} \cdots d\sigma_2 d\sigma_1,
\end{align}
where $d\sigma_j:=d\sigma(x_j)=c_\mu \mu(v_j) \{n(x_j) \cdot v_j\}dv_j$.

\bigskip

\subsection{Notations}

\begin{enumerate}
\item We denote the Japanese bracket $\langle v \rangle = (1+|v|^2)^{1/2}$. 
\item As a convention, we denote the following function spaces for $p \in [1,\infty]$,
\begin{align} \label{fntspace1}
	L^p_{x,v}=L^{p}(\Omega \times \mathbb{R}^3_v), \quad L^p_x=L^p(\Omega), \quad L^p_v=L^p(\mathbb{R}^3_v).
\end{align}
\item  We define the mixed function space $L^p_vL^\infty_x$ with the following norm:
\begin{align*}
	\|f\|_{L^p_vL^\infty_x}:= \left(\int_{\R^3} \Big|\sup_{x\in \Omega}|f(x,v)|\Big|^p dv \right)^{1/p}.
\end{align*} 
\item We define the space $L^\infty_xL^p_v(\gamma)=L^\infty_xL^p_v(\partial\Omega\times \mathbb{R}^3)$ with the norm
\begin{align*}
	\|f\|_{L^\infty_xL^p_v(\gamma)}:=\sup_{x\in \partial\Omega}\left[\left(\int_{\{n(x) \cdot v > 0\} \cup \{n(x) \cdot v < 0\}}\left|f(x,v)\right|^{p}|n(x)\cdot v|dv\right)^{1/p}\right].
\end{align*}
Usually, We define the norm on the boundary $\|\cdot\|_{L^\infty_xL^p_v(\gamma_{\pm})}$ by
\begin{align*}
	\|f\|_{L^\infty_xL^p_v(\gamma_{\pm})}:=\sup_{x\in \partial\Omega}\left[\left(\int_{n(x) \cdot v \gtrless 0}\left|f(x,v)\right|^{p}|n(x)\cdot v|dv\right)^{1/p}\right].
\end{align*}
\item We define the space $L^2_\gamma = L^2(\partial\Omega\times \mathbb{R}^3)$ with the norm
\begin{align*}
	\|f\|_{L^2_\gamma}:= \left[\int_{\gamma}|f(x,v)|^2|n(x) \cdot v| dS(x)dv\right]^{1/2}.
\end{align*}
Usually, we denote the norm on the boundary $\gamma_+$ in $L^2$ by
\begin{align} \label{L2bnorm}
	\|f\|_{L_{\gamma_+}^2}:=\left[\int_{\gamma_+}|f(x,v)|^2\{n(x) \cdot v\} dS(x)dv\right]^{1/2}.
\end{align}
\item For a positive Lebesgue measurable function $m$ on $\R^3$, we define the weighted function space $L^p_vL^\infty_x(m)$ given by the norm
\begin{align} \label{mnorm}
	\|f\|_{L^p_vL^\infty_x(m)} := \|mf\|_{L^p_vL^\infty_x}= \left(\int_{\R^3} \left|\sup_{x\in \Omega}|m(v)f(x,v)|\right|^p dv \right)^{1/p}.
\end{align}
\item We set $(f,g)_{L^2_v} = \int_{\mathbb{R}^3} f(v)g(v)dv$ the inner product in $L^2(\mathbb{R}^3)$\\
\item $p'$ is the exponent conjugate of $p$.\\
\item If not specifically mentioned, $C_a$ or $C(a)$ is the generic positive constant depending on $a$, while $C_0, C_1, C_2, \cdots$ denote some specific positive constants.\\
\item $\partial_{x_i}$ and $\partial_{v_i}$ mean the partial derivative with respect to $x_i$ and $v_i$, respectively. We also abbreviate $\partial_{ij}=\partial_{x_i}\partial_{x_j}$.
\end{enumerate}

\subsection{Perturbation theory near Maxwellian}
Setting $F(t,x,v) = \mu(v) + \mu^{1/2}(v) f(t,x,v)$, the Boltzmann equation \eqref{Boltzmanneq} becomes
\begin{align} \label{FPBER}
	\partial_t f+v \cdot \nabla_x f +Lf = \Gamma(f,f)
\end{align} 
and the diffuse reflection boundary condition \eqref{Diffuseboundarycond} can be written as
\begin{align} \label{PDRBC}
	f(t,x,v)|_{\gamma_-} = c_\mu \mu^{1/2}(v) \int_{n(x) \cdot v'>0 }f(t,x,v')\mu^{1/2}(v')\{n(x) \cdot v'\}dv'.
\end{align}
Here, the operator $L$ is a linear operator of the form
\begin{align*}
	Lf = -\mu^{-1/2} \left\{Q(\mu, \mu^{1/2}f)+Q(\mu^{1/2}f,\mu)\right\} = \nu(v) f -Kf,
\end{align*}
with the collision frequency $\nu(v) = \int_{\R^3}\int_{\S^2} B(v-u,\omega)\mu(u) d\omega du \sim (1+|v|)^\gamma$. Note that $\nu(v)$ has the greatest lower bound, denoted by $\nu_0$. And the operator $K=K_1-K_2$ is defined by
\begin{align} \label{operatorK1}
	(K_1f)(t,x,v) &= \int_{\R^3} \int_{\S^2} B(v-u,\omega) \mu^{1/2}(v)\mu^{1/2}(u)f(t,x,u)d\omega du,\\
	(K_2f)(t,x,v) &= \int_{\R^3} \int_{\S^2} B(v-u,\omega) \mu^{1/2}(u)\mu^{1/2}(u')f(t,x,v')d\omega du\nonumber\\ \label{operatorK2}
	& \quad + \int_{\R^3} \int_{\S^2} B(v-u,\omega) \mu^{1/2}(u)\mu^{1/2}(v')f(t,x,u')d\omega du.
\end{align}
It is well-known that the operator $L$ has a kernel
\begin{align*}
	\text{Ker}(L) = \text{span}\left\{\mu^{1/2},v_1\mu^{1/2},v_2\mu^{1/2},v_3\mu^{1/2}, \frac{|v|^2-3}{\sqrt{6}}\mu^{1/2} \right\}.
\end{align*}
The nonlinear term $\Gamma(f,f) = \Gamma^+(f,f) - \Gamma^-(f,f)$ is given by
\begin{align*}
	\Gamma^+(f,f) = \mu^{-1/2}Q^+(\mu^{1/2}f,\mu^{1/2}f) , \quad \Gamma^-(f,f)=\mu^{-1/2}Q^-(\mu^{1/2}f,\mu^{1/2}f).
\end{align*}
Denoting
\begin{align} \label{Rfdefinition}
	R(f)(t,x,v) = \int_{\R^3} \int_{\S^2} B(v-u,\omega) \left[\mu(u)+\mu^{1/2}(u)f(t,x,u)\right]d\omega du,
\end{align}
the equation \eqref{FPBER} can be rewritten as
\begin{align} \label{FPBE}
	\partial_t f +v\cdot \nabla_x f + R(f)f =Kf+ \Gamma^+(f,f).
\end{align}
As previously mentioned, for the diffuse reflection boundary condition \eqref{Diffuseboundarycond}, the mass is conserved.
Therefore, we may assume 
\begin{align} \label{massconserv}
	\int_{\Omega \times \mathbb{R}^3} f(t,x,v) \mu^{1/2}(v)dxdv = 0 \quad \text{for all } t\ge 0.
\end{align}
by imposing initial data $F_0$ such that
\begin{align} \label{initialmass}
	\int_{\Omega \times \mathbb{R}^3} F_0(x,v) dxdv = 	\int_{\Omega \times \mathbb{R}^3} \mu(v) dxdv.
\end{align}

\subsection{Main results}

We introduce the weight function
\begin{align} \label{weightfnt}
	w(v) = \left(1+|v|^2\right)^{\beta/2}.
\end{align}
By considering a perturbation of the form $\mu+\mu^{1/2}f$ near the Maxwellian $\mu$ and setting $h(t,x,v)= w(v) f(t,x,v)$, the Boltzmann equation \eqref{Boltzmanneq} can be written as the perturbed Boltzmann equation for $h$ with $R(f)$:
\begin{align} \label{WPBE2}
	\partial_t h + v\cdot \nabla_x h + R(f)h = K_wh +w\Gamma^+\left(f,f\right),
\end{align}
with initial data $h_0$ and the diffuse reflection boundary condition for $h$, where the weighted operator $K_w$ is defined by
\begin{align} \label{woperK}
	K_w h := wK\left(\frac{h}{w}\right).
\end{align}
Applying Duhamel principle to the equation \eqref{WPBE2}, we obtain the mild form for $h$
\begin{equation} \label{mildsolh}
	\begin{aligned}
	h(t,x,v) &= S_{G_f}(t,0)h_0(x,v)+\int_0^t (S_{G_f}(t,s)K_wh(s))(x,v)ds\\
	& \quad + \int_0^t (S_{G_f}(t,s)w\Gamma\left(f, f\right)(s))(x,v)ds,
\end{aligned}
	\end{equation}
where $S_{G_f}(t,s)$ is the solution operator for solutions to the equation
\begin{align*}
	\partial_t h +v\cdot \nabla_x h  +R(f)h=0
\end{align*}
with the diffuse reflection boundary condition for $h$.\\
Before achieving our main goal, we need to demonstrate the small perturbation problem to the Boltzmann equation \eqref{Boltzmanneq} with the smallness of relative entropy $\mathcal{E}(F)$ defined by
\begin{align*}
	\mathcal{E}(F)= \int_{\Omega \times \mathbb{R}^3} \left(\frac{F}{\mu}\log\frac{F}{\mu}-\frac{F}{\mu}+1\right) \mu dxdv.
\end{align*}

\begin{Thm}[Small perturbation problem] \label{Mainresult1}
	Let $w(v) = (1+|v|^2)^{\beta /2}$. Let $p$ satisfy the condition 
	\begin{align*}
		p> \begin{cases}
 4, & \text{if}\quad 0\le \gamma \le \frac{3}{4}, \\
 \frac{5}{2-\gamma}, & \text{if}\quad \frac{3}{4}<\gamma \le 1,
\end{cases}
	\end{align*}
 a and $\beta>\max\Big\{\frac{3p-6}{2p}+2\gamma, \beta(p,\gamma) \Big\}$, where $\beta(p,\gamma)$ is given by
	\begin{align*}
	\beta(p,\gamma) := \begin{cases}
 \frac{3p-5}{p}, & \text{if} \quad 0\le \gamma \le \frac{3}{4},p>4\ \text{or}\ \frac{3}{4}<\gamma \le 1, p>\frac{1}{1-\gamma},  \\
 \frac{7p-12}{2p}, & \text{if}\quad \frac{3}{4}<\gamma \le 1, \frac{5}{2-\gamma}< p \le \frac{1}{1-\gamma}.
\end{cases}
	\end{align*}
Assume that $F_0(x,v) = \mu(v) + \mu^{1/2}(v) f_0(x,v) \ge 0$ satisfying the mass conservation $\int_{\Omega \times \R^3}f_0 \mu^{1/2}dxdv =0$. Then there is  $\eta_0  \ll  1$ so that there exists a constant $\epsilon_0 = \epsilon_0(\eta_0) >0$, depending only on $\eta_0$, such that if
	\begin{align*}
		\|wf_0\|_{L^p_vL^\infty_x} \le \eta_0, \quad \mathcal{E}(F_0)\le \epsilon_0,
	\end{align*}
	then there exists a unique solution $F(t,x,v) =\mu(v) + \mu^{1/2}(v) f(t,x,v) \ge 0 $ to the Boltzmann equation \eqref{Boltzmanneq} with the diffuse reflection boundary condition \eqref{Diffuseboundarycond} and initial data $F(0,x,v)=F_0(x,v)$. Moreover, there exist $C_0>0$ and $ \lambda_0>0$ such that
	\begin{align*}
		\|wf(t)\|_{L^p_v L^\infty_x} \le C_{0} e^{-\lambda_0 t}\|wf_0\|_{L^p_v L^\infty_x}
	\end{align*}
	for all $t \ge 0$.
	\end{Thm}

\begin{remark}
	Note that $\eta_0$ depends only on $p$, $\gamma$, $\beta$, and $\mathcal{T}$, where $\mathcal{T}$ is a large constant determined by \eqref{assump1}. 
\end{remark}

\begin{remark}
	In Theorem \ref{Mainresult1}, the condition on $p$ is specified in \eqref{pcondition}, and it comes from the derivation of the $L^p$ estimate for the gain term $\Gamma^+$. (See Corollary \ref{LpGamma+est}.) 
\end{remark}

\begin{remark}
	In Theorem \ref{Mainresult1}, the condition on $\beta$ is determined by \eqref{betacondition1} and \eqref{betacondL2}. The condition \eqref{betacondition1} essentially arises from the application of the pointwise estimate \eqref{Gamma+estimatecombine} for the gain term $\Gamma^+$ in the proof of Lemma \ref{Rfest1}. (See \eqref{Rf26}.) On the other hand, the condition \eqref{betacondL2} originates from the step \eqref{betacondL2} to derive the estimate \eqref{Gamma+L2}. (See Lemma \ref{nonlinear L^2 decay}.)
\end{remark}

\bigskip

The main goal of the paper is to prove the global existence of the solution to the Boltzmann equation \eqref{Boltzmanneq} with a large oscillation initial data and smallness for relative entropy near the global maxwellian $\mu(v)$.
\begin{Thm}[Large amplitude problem] \label{Mainresult2}
	Let $w(v) = (1+|v|^2)^{\beta /2}$. Let $p$ satisfy the condition 
	\begin{align*}
		p> \begin{cases}
 4, & \text{if}\quad 0\le \gamma \le \frac{3}{4}, \\
 \frac{5}{2-\gamma}, & \text{if}\quad \frac{3}{4}<\gamma \le 1,
\end{cases}
	\end{align*}
	 and $\beta>\max\Big\{\frac{3p-6}{2p}+2\gamma, \beta(p,\gamma) \Big\}$, where $\beta(p,\gamma)$ is given by
	\begin{align*}
	\beta(p,\gamma) := \begin{cases}
 \frac{3p-5}{p}, & \text{if} \quad 0\le \gamma \le \frac{3}{4},p>4\ \text{or}\ \frac{3}{4}<\gamma \le 1, p>\frac{1}{1-\gamma},  \\
 \frac{7p-12}{2p}, & \text{if}\quad \frac{3}{4}<\gamma \le 1, \frac{5}{2-\gamma}< p \le \frac{1}{1-\gamma}.
\end{cases}
	\end{align*} 
	Assume that $F_0(x,v) = \mu(v) + \mu^{1/2}(v) f_0(x,v) \ge 0$ satisfying the mass conservation $\int_{\Omega \times \R^3} f_0\mu^{1/2}dxdv=0$. For any $M_0 \ge 1$, there exists a constant $\tilde{\epsilon}_0>0$, depending only on $\eta_0$ and $M_0$, such that if
	\begin{align} \label{mainconditionLA}
		\|wf_0\|_{L^p_vL^\infty_x} \le M_0, \quad \mathcal{E}(F_0)\le \tilde{\epsilon}_0,
	\end{align}
	then there exists a unique solution $F(t,x,v) =\mu(v) + \mu^{1/2}(v) f(t,x,v) \ge 0 $ to the Boltzmann equation \eqref{Boltzmanneq} with the diffuse reflection boundary condition \eqref{Diffuseboundarycond} and initial data $F(0,x,v)=F_0(x,v)$ satisfying
	\begin{align*}
		\|wf(t)\|_{L^p_vL^\infty_x} \le  \tilde{C}_L M_0^5\exp\left\{ \frac{2\tilde{C}_L}{\nu_0}M_0^5 \right\}e^{-\lambda_L t}
	\end{align*}
	for all $t \ge 0$, where $\tilde{C}_L \ge 1$ is a generic constant, $\lambda_L:= \min\{\lambda_0, \frac{\nu_0}{16}\}>0$, where $\lambda_0$ is the same constant in Theorem \ref{Mainresult1}.
\end{Thm}

\bigskip

\begin{remark}
	In Theorem \ref{Mainresult2}, the conditions on $p$ and $\beta$ are determined for similar reasons as in Theorem \ref{Mainresult1}.	
\end{remark}

\begin{remark} \label{noconflict}
	We note that in \eqref{mainconditionLA}, the large amplitude initial condition in $L^p_vL^\infty_x$ is not in conflict with the smallness condition for initial relative entropy $\mathcal{E}(F_0)$. We now provide an example to support this observation. Let $p\ge1$ and let
	\begin{align} \label{exampleinitialcond}
		F_0(x,v) = \mu(v)+\mu^{1/2}(v)f_0(x,v)\quad \text{with}\quad f_0(x,v)=\frac{\phi(x)-1}{w(v)}\mu^{1/2}(v),
	\end{align}
	where $w(v) = (1+|v|^2)^{\beta/2}$ and $\phi(x)$ is chosen later, depending only on $x$. First of all, we check that $F_0$ satisfies the positivity and the mass conservation. For positivity, we have
	\begin{align*}
		F_0(x,v)=\frac{\mu(v)}{w(v)} \left(w(v)-1+\phi(x)\right) \ge 0
	\end{align*}
	provided that
	\begin{align} \label{xfntcondition1}
		\phi(x) \ge 0 \quad \text{for all}\  x \in \Omega.
	\end{align}
	Also, if $\phi(x)$ satisfies the following condition
	\begin{align} \label{xfntcondition2}
	 \int_{\Omega} (\phi(x)-1)dx=0,	
	\end{align}
	the mass conservation \eqref{massconserv} is guaranteed. Next, we consider the large amplitude initial condition. Let $M_0>0$ be finite but arbitrarily large. From the large amplitude condition in \eqref{mainconditionLA}, the function $\phi(x)$ need to satisfy the following condition
	\begin{align} \label{xfntcondition3}
		M_0 = \|wf_0\|_{L^p_vL^\infty_x} = \sup_{x\in \Omega}|\phi(x)-1|\left(\int_{\R^3}\mu^{p/2}(v)dv\right)^{1/p}=C_{p} \sup_{x\in \Omega}|\phi(x)-1|,
	\end{align}
	where $C_{p}$ is a constant depending on $p$. The only remaining requirement for $F_0$ is the smallness condition for the initial relative entropy $\mathcal{E}(F_0)$. From \eqref{exampleinitialcond}, the initial relative entropy $\mathcal{E}(F_0)$ becomes
	\begin{align*}
		\mathcal{E}(F_0) =\int_{\Omega}\int_{\R^3}\mu(v) \left(1+\frac{\phi(x)-1}{w(v)}\right) \log\left(1+\frac{\phi(x)-1}{w(v)} \right)dxdv.
	\end{align*}
	Thanks to the convexity of $\Psi(s):=s\log s$ on $s>0$, we obtain
	\begin{align*}
		\Psi\left(1+\frac{\phi(x)-1}{w(v)} \right) &= \Psi\left(\left(1-\frac{1}{w(v)}\right) \cdot 1+\frac{1}{w(v)} \phi(x) \right)\\ 
		&\le \left(1 -\frac{1}{w(v)}\right) \Psi(1) + \frac{1}{w(v)} \Psi(\phi(x))\\
		& = \frac{1}{w(v)}\phi(x) \log \phi(x),
	\end{align*}
	where we have used $\Psi(1)=0$.
	From the above inequality, the initial relative entropy is bounded by
	\begin{align*}
		\mathcal{E}(F_0) \le \int_{\Omega}\phi(x) \log \phi(x)dx \int_{\R^3} \frac{\mu(v)}{w(v)}dv \lesssim \int_{\Omega}\phi(x) \log \phi(x)dx.
	\end{align*}
	If it holds that
	\begin{align} \label{xfntcondition4}
		\int_{\Omega}\phi(x) \log \phi(x)dx \ll 1,	\end{align}
then the initial relative entropy $\mathcal{E}(F_0)$ can be sufficiently small. Gathering the conditions \eqref{xfntcondition1}, \eqref{xfntcondition2}, \eqref{xfntcondition3}, \eqref{xfntcondition4}, $\phi(x)$ is required to satisfy the following conditions:
\begin{equation} \label{xfntcombined}
\begin{aligned}
	&\phi(x) \ge 0 \quad \text{for all} \ x\in \Omega; \quad  \int_{\Omega} (\phi(x)-1)dx=0;\\
	& \sup_{x\in \Omega}|\phi(x)-1| = C_p M_0; \quad \int_{\Omega}\phi(x) \log \phi(x)dx \ll 1.
\end{aligned}
\end{equation}
Thus, the initial data $F_0$ given by \eqref{exampleinitialcond} with $\phi(x)$ satisfying \eqref{xfntcombined} has a large amplitude in $L^p_vL^\infty_x$, but its relative entropy can be sufficiently small.\\
\indent  Note that initial data $F_0(x,v)$ with a large amplitude in the $L^p_vL^\infty_x$ norm may contain locally the vacuum. However, after sufficient time has passed, the initial vacuum can disappear.	
\end{remark}

\bigskip

\begin{figure}[H]
\begin{tikzpicture}
  \begin{axis}[
    axis lines=middle,
    xlabel={$t$},
    ylabel={$\mathcal{E}(F(t))$},
    xmin=0, xmax=10,
    ymin=0, ymax=2,
    thick,
    width=12cm,
    height=5cm,
    xtick=\empty,
    ytick=\empty,
    clip=false,
    enlargelimits
  ]

    \addplot[
      blue, 
      thick, 
      smooth, 
      tension=0.7
    ] coordinates {
      (0,0.5)
      (1,0.48)
      (2,0.46)
      (3,0.42)
      (4,0.36)
      (5,0.27)
      (6.3,0.15)
    };

    \addplot[
      red, 
      thick, 
      smooth, 
      tension=0.7
    ] coordinates {
      (6.3,0.15)
      (8,0.10)
      (9,0.06)
      (10,0.03)
      (11,0.02)
    };

    \draw[dashed] (axis cs:0,0.5) -- (axis cs:0,0);

    \node[anchor=east] at (axis cs:-0.1,0.52) {$\tilde{\epsilon}_0$};

    \draw[dashed] (axis cs:6.3,0) -- (axis cs:6.3,0.8);

    \node[below] at (axis cs:6.3,0) {$T_0$};

    \node[anchor=east] at (axis cs:-0.1,0.8) {$\epsilon_0$};

    \draw[dashed] (axis cs:0,0.8) -- (axis cs:6.3,.8);

    \node[below, font=\bfseries] at (axis cs:3.1,0) {\small large amplitude problem};

    \node[below, anchor=west, font=\bfseries] at (axis cs:6.6,-0.2) {\small small perturbation problem};

  \end{axis}
\end{tikzpicture}
\caption{The relative entropy $\mathcal{E}(F_0)$}
 \label{figrelative}
\end{figure}
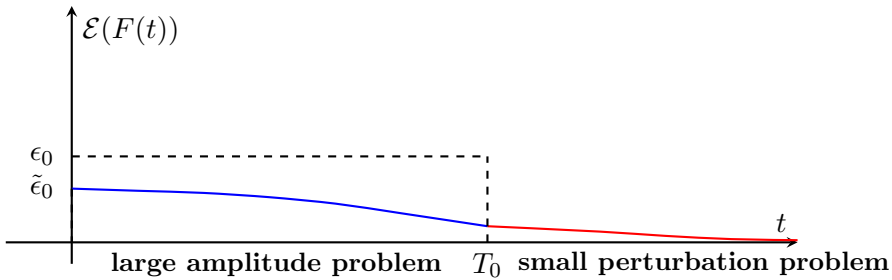

\indent There are several significant results concerning the spatially homogeneous Boltzmann equation. In Arkeryd \cite{Arkeryd1972, arkeryd1988}, the author established the well-posedness and stability for solutions in $L^1$. Later, Mischler-Wennberg \cite{mischler1999} constructed global-in time solutions for any initial data with finite mass and energy. Additionally, Mouhot \cite{mouhot2006b} proved convergence to equilibrium whose exponential-decay rate is determined by the spectral gap for the linearized operator.   \\
\indent The inhomogeneous Boltzmann equation has also been studied. For whole space, the global existence of renormalized solutions  with general initial data $F_0(x,v) \ge 0$ satisying bounded physical quanitities has been proven by DiPerna-Lions \cite{Diperna1989}.  Later, for domains with general boundary condition, Hamdache \cite{Hamdache1992} established the global existence of weak solutions. However, the uniqueness for those solutions is an open problem. On the other hand, the convergence to equilibrium is one of the key issues studied for the Boltzmann equation. Regarding such issue, Desvillettes-Villani \cite{Desvillettes2005} demonstrated that under the assumption of the existence of a global solution satisfying some a priori high-order Sobolev bounds and a Gaussian lower bound, the solution converges algebraically to the global equilibrium. Following this result, Gualdani-Mischler-Mouhot \cite{Gualdani2017} enhanced the convergence result on the cutoff hard sphere model by showing that the solution converges exponentially to equilibrium under suitable high-order a priori assumptions. A detailed comparison of our result with the well-known work \cite{Gualdani2017} will be provided later. \\
\indent Meanwhile, the well-posedness and asymptotic behavior have been studied for the perturbation theory near Maxwellians. When perturbation $f_0 = \frac{F_0-\mu}{\sqrt{\mu}}$ is sufficiently small in some high-order Sobolev spaces, Ukai \cite{Ukai1974} constructed the global-in-times solutions to the Boltzmann equation in a periodic box. Subsequently, in high-order Sobolev spaces, the global well-posedness theory and the convergence to the equilibrium was developed by Guo \cite{Guo2002}, \cite{Guo2003a},\cite{Guo2012} in periodic box. We mention the other contributions \cite{Guo2012a}, \cite{Strain2007} in this framework.\\
\indent Unfortunately, for general bounded domains, high-order regularity of solutions to the Boltzmann equation is generally not guaranteed. This has been observed in \cite{Guo2016, Kim2011}. Instead, Guo \cite{GuoDecay} proposed an $L^2$-$L^\infty$ bootstrap argument to construct solutions with low regularity under physical boundary conditions, assuming small perturbations. Motivated by this study, low regularity approach in weighted $L^\infty$ space has been further developed. In particular, under the specular reflection boundary condition, Kim-Lee \cite{Kim2017} extended the setting from real analytic uniformly convex domains to general $C^3$ uniformly convex bounded domains. Subsequently, Ko-Kim-Lee \cite{Ko2023} addressed the small perturbation problem for certain non-convex domains. Additionally, initial-boundary value problems  with polynomial tails under various boundary conditions were studied in \cite{Briant2017,Briant2016}. We also refer to \cite{Kim2018, Cao2019} for related works on initial-boundary value problems.\\
\indent On the other hand, $L^\infty$ theory for the Boltzmann equation was extended to allow sufficiently large amplitude initial data with a smallness condition in $L^p$ class. In Duan-Huang-Wang-Yang \cite{Duan2017a}, the authors first replaced the $L^\infty$ smallness condition for initial data with the $L^p$ type smallness condition in a periodic box or whole space. Later on, research on the large amplitude problems has also been applied to problems with boundary condition. We refer Duan-Wang \cite{Duan2019b} for diffuse reflection boundary conditions and Duan-Ko-Lee \cite{Ko2023a} for specular reflection boundary conditions. We also mention \cite{Duan2019, Cao2022a, Guanfa2019, Jong2025, Ko2022, Ko2025,Jiang2025, Ko2023BGK} for additional works on large amplitude problems. Nearly all of these results concern the Boltzmann equation in uniformly bounded function classes, i.e. in weighted $L^\infty_{x,v}$ spaces. In this paper, we study the Boltzmann boundary condition problem in pointwise unbounded function class $L^p_vL^\infty_x(\langle v \rangle^\beta \mu^{-1/2})$. It is meaningful in that we address the most general function classes among the known results on large amplitude problems with boundary conditions.\\
\indent In the view of the renowned work \cite{Gualdani2017}, our results provide an alternative perspective. Unlike our results, the paper \cite{Gualdani2017} addressed the small data problem for solutions with polynomial tails on $\T^3$. If we consider our result on the small data problem in the torus setting, it is possible that from the result in the paper \cite{Gualdani2017}, we obtain the exponential decay of solutions in the  $L^1_vL^\infty_x(\langle v\rangle^\beta)$ sense, since the initial data $F_0-\mu$ in the space $L^p_vL^\infty_x(\langle v\rangle^\beta\mu^{-1/2})$ belongs to $L^1_vL^\infty_x(\langle v\rangle^\beta)$. However, the situation is different when the norm $\left\|F_0-\mu\right\|$ becomes large. In \cite{Gualdani2017}, under the high-order a priori assumption, the quantity $\|F-\mu\|_{L^1_vL^\infty_x(\langle v\rangle^\beta)}$ is controlled by relative entropy. In \cite{Desvillettes2005}, the authors  demonstrated that relative entropy decreases over time, provided that the high-order a priori assumption holds, and thus the paper \cite{Gualdani2017} treated the case where the $L^1_vL^\infty_x(\langle v \rangle^\beta)$ norm is small. On the other hand, as mentioned earlier, the smoothness for solutions cannot be assumed in bounded domains, unlike in periodic domains. This means that small relative entropy does not necessarily imply smallness of the norm in bounded domains. Therefore, even if the relative entropy is small, the initial data in $L^p_vL^\infty_x(\langle v \rangle^\beta \mu^{-1/2})$ can be large. (See remark \ref{noconflict}.) Based on this observation, we establish  convergence to equilibrium in the $L^p_vL^\infty_x(\langle v\rangle^\beta\mu^{-1/2})$ sense.\\
\indent The key idea of this paper is to use a pointwise estimate for the gain term $\Gamma^+(f,g)$ with the weight $w$, which is controlled in terms of $L^p_v$ and $L^2_v$ norms. The pointwise estimate contains a suitable polynomial decay in the velocity variable, and such decay play a crucial role in achieving our goal in the $L^p_vL^\infty_x$ framework.\\ 
\indent We outline potential directions for the future research originating from our work. In this paper, we consider weighted $L^p_vL^\infty_x$ spaces for solutions to the perturbed Boltzmann equation. As the future work, we hope to extend this result to general mixed Lebesgue spaces $L^p_vL^q_x$ with $1\le p,q<\infty$. However, it is hard to deal with our problem in the framework because the collision operator $Q(f,f)$ does not have favorable structures on $L^q_x$ spaces, in other words, only the $L^\infty$ space among $L^q$ spaces have an algebraic property: $\|fg\| \lesssim \|f\|\|g\|$. Turning to another topic, we may address the large amplitude problem for solutions in polynomial tail class, which is larger than exponential tail class. It is very interesting problem, but in order to solve the problem with polynomial class, it is important to control the collision operator $Q(f,f)$ in function space $L^p_vL^\infty_x(\langle v\rangle^\beta)$. As in Lemma \ref{pointwiseGamma+estimate}, it is difficult to establish a pointwise estimate for the collision term $Q^+(f,f)$ because the estimate has no sufficient decaying factor with respect to velocity variables which is compared to the nonlinear term $\Gamma^+(f,f)$. Rather, it generally has the velocity growth like $\nu(v)$ as presented in \cite{Gualdani2017}. We may probably need a new approach to solve the large amplitude problem for the Boltzmann equation with polynomial perturbation.\\
\indent We now briefly explain our strategy to achieve the main goal. To control initial data with large amplitude in the $L^p_vL^\infty_x$ norm, we impose the smallness condition for the initial relative entropy $\mathcal{E}(F_0)$. Since it is difficult to obtain the exponential time decay in the large amplitude problem, we instead derive convergence to equilibrium by utilizing the exponential time decay established in the small perturbation problem.\\
\indent First of all, we consider the small perturbation problem. In this problem, we face the challenge of controlling the nonlinear term. At this stage, we may apply the Minkowski's integration inequality to handle the $L^p_vL^\infty_x$ norm of the time integration for the nonlinear term, in particular the loss term:
\begin{align*}
	\left\|\int_0^t e^{-\nu(v)(t-s)} w\Gamma^-(f,f)(s)ds\right\| _{L^p_vL^\infty_x} 
	& \lesssim  \left(\int_{\R_v^3}\sup_{x\in \Omega}\left|\int_0^t e^{-\nu(v)(t-s)} \nu(v)|wf(s)|\|wf(s)\|_{L^p_vL^\infty_x}ds\right|^p dv \right)^{\frac{1}{p}}\\
	& \nleq \int_0^t e^{-\nu_0 (t-s)}\|h^{(m)}(s) \|_{L^p_vL^\infty_x}^2ds.
\end{align*}
However, we cannot apply the argument mentioned above due to the presence of the factor $\nu(v)$ in $w\Gamma^-(f,f)$. Thus, as an alternative approach, we will consider the term $R(f)f$ formed by combining the term $\nu(v)f$ and the loss term $\Gamma^-(f,f)$. (See the equation \eqref{FPBE}.) By considering the term $R(f)$ instead of the collision frequency $\nu(v)$, we need to establish the time decay based on this term. Therefore, by obtaining a lower bound for the term $R(f)$, such as $\nu(v)/2$, we will achieve the desired goal. At this stage, we introduce the a priori assumption $\|wf(t)\|_{L^p_vL^\infty_x} \le \bar{\eta}\ll 1$ and impose the smallness condition for relative entropy. (See subsection \ref{apriorisamll}.) Based on the a priori assumption, we can also derive the nonlinear $L^p_vL^\infty_x$ estimate for the Boltzmann equation, and we demonstrate the small perturbation problem with the smallness condition for initial relative entropy, $\mathcal{E}(F_0)\le \epsilon_0.$\\
\indent Next, we consider the process combining the small perturbation regime and the large amplitude regime. Thanks to the smallness condition for initial relative entropy, the initial large amplitude $\|wf_0\|_{L^p_vL^\infty_x}$ gradually decreases and can be smaller than the small amplitude $\eta_0$, as mentioned in the small perturbation problem(Theorem \ref{Mainresult1}). We set the time point as $T_0$, after which we can attach the small amplitude regime. (See figure \ref{figamplitude}.) Moreover, thanks to Lemma \ref{relativedecrease}, relative entropy is non-decreasing, and if initial relative entropy $\mathcal{E}(F_0)$ is chosen to be smaller than $\epsilon_0$, relative entropy at time $T_0$ is also smaller than $\epsilon_0$. (See figure \ref{figrelative}.) This does not cause any issues in connecting the small amplitude regime. Hence, we complete the large amplitude problem in the $L^p_vL^\infty_x$ framework.

\begin{figure}[H]
\begin{tikzpicture}
  \begin{axis}[
    axis lines=middle,
    xlabel={$t$},
    ylabel={$\|wf(t)\|_{L^p_vL^\infty_x}$},
    xmin=0, xmax=8,
    ymin=0, ymax=6,
    thick,
    width=12cm,
    height=6cm,
    xtick=\empty,
    ytick=\empty,
    clip=false,
    enlargelimits
  ]
    
    \addplot[
      blue, 
      thick, 
      smooth
    ] coordinates {
      (0,4)
      (0.3,4.2)
      (0.5,3.8)
      (1,4.1)
      (1.4,4.1)
      (1.6,4.6)
      (2,4.5)
      (2.5,5)
      (3,3.51)
      (3.5,3.3)
      (3.6,3)
      (4,1.5)
      (5,1)
    };

    \addplot[
      red, 
      thick, 
      smooth
    ] coordinates {
      (5,1)
      (5.2,0.9)
      (5.4,1.1)
      (5.7,1.3)
      (6,1.4)
      (6.3,1.2)
      (6.6,1.3)
      (7,1)
      (7.5,0.5)
      (9,0.2)
    };

    \draw[dashed] (axis cs:0,5) -- (axis cs:5,5);
    \draw[dashed] (axis cs:0,1.) -- (axis cs:5,1.);
    \draw[dashed] (axis cs:5,0) -- (axis cs:5,5);
   
    \node[anchor=east] at (axis cs:-0.1,5) {$\bar{M}$};
    \node[anchor=east] at (axis cs:-0.1,4) {$M_0$};
    \node[anchor=east] at (axis cs:-0.1,1.) {$\eta_0$};

    \node[below] at (axis cs:5.0,-0.1) {$T_0$};

    \node[below, font=\bfseries] at (axis cs:2.5,0) {\small large amplitude problem};
    \node[below right, font=\bfseries] at (axis cs:5.3,0) {\small small perturbation problem};
  \end{axis}
  \end{tikzpicture}
  \caption{Amplitude for the norm $\|wf(t)\|_{L^p_vL^\infty_x}$}
 \label{figamplitude}
\end{figure}
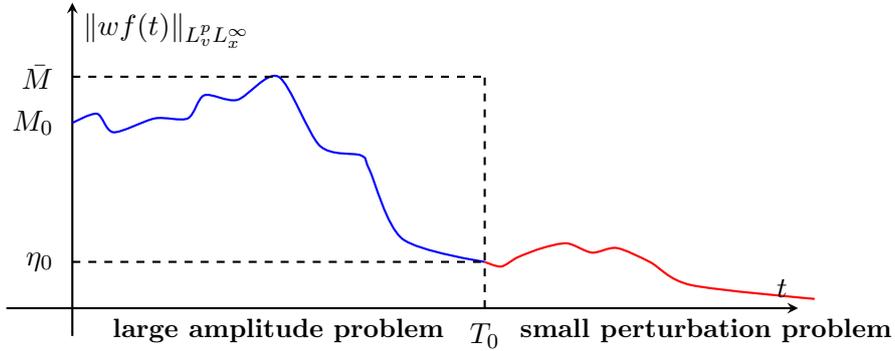

\bigskip

\subsection{Organization of the paper} In Section \ref{Preliminaries}, we present estimates for operators associated with the full perturbed Boltzmann equation \eqref{FPBER}. In particular, we derive a pointwise estimate for $w\Gamma^+(f,f)$. We also deal with relative entropy and its applications.
In Section \ref{smallperturbationproblem}, we focus on the small perturbation problem. We first deal with the elementary equation associated with the equation \eqref{FPBER}, followed by a study of the corresponding solution operators and their estimates. Under the a priori assumption in the small perturbation problem, we use these estimates to establish the global solution to the Boltzmann equation with initial data with small amplitude in $L^p_vL^\infty_x$ and small relative entropy.
In Section \ref{Largeamplitudeproblem} as a main part, given a large oscillation initial datum in $L^p_vL^\infty_x$, we solve the global existence of the solution to the Boltzmann equation. To achieve our main aim, we will introduce estimates to derive a main goal under the a priori assumption and then from a Gr\"onwall type we apply the global-in time existence of a solution given by the small amplitude initial data to prove our main goal. 
In Section \ref{Appendix}, we provide an appendix for the proof of Lemma \ref{coer}.

\bigskip

\section{Preliminaries} \label{Preliminaries}

\subsection{Operators and the related estimates}
We first consider the operator $K$ and the related estimate. Recall the operator $K$ can be decomposed into $K_1$ and $K_2$ as \eqref{operatorK1} and \eqref{operatorK2}. In \cite{Glassey1996}, it is well-known that the operator $K$ can be written as 
\begin{align*}
	(Kf)(t,x,v) = \int_{\R^3}k(v,u)f(t,x,u)du,
\end{align*} 
where the kernel $k$ is symmetric. Now we introduce a useful estimate for the operator $K$.
\bigskip

\begin{Lem} \cite{Deng2018} \label{KESTIMATE}
For $0< \epsilon<1$, we have
\begin{align*}
	|k(v,v')| \lesssim \left(\frac{1}{|v-v'|}+|v-v'|^\gamma\right)\frac{1}{(1+|v|+|v'|)^{1-\gamma}}e^{-\frac{1-\epsilon}{8}|v-v'|^2-\frac{1-\epsilon}{8}\frac{||v|^2-|v'|^2|^2}{|v-v'|^2}}
\end{align*}
for any $v, v' \in \R^3$ with $v \not=v'$.
Let $1\le q<3$ and $\ell \in \R$. Then there exists $C_{\beta,q}>0$ such that
	\begin{align} \label{pkestimate}
		\int_{\R^3} \left|k(v,v') \frac{\langle v \rangle^\ell}{\langle v' \rangle^\ell}\right|^qdv' \le \frac{C_{\beta,q}}{(1+|v|)^{q(1-\gamma)+1}}
	\end{align}
	for any $v \in \R^3$.
\end{Lem}

\bigskip

From now on, we consider the key estimate to deal with the nonlinear term $\Gamma^+(f,f)$. Before presenting the pointwise estimate for $w\Gamma^+(f,f)$, we need to compute the Jacobian determinant $\left|\frac{du'}{du}\right|$ used to
derive this estimate. Recall \eqref{postcollisionrep}. For each $\omega$ and fixed $v$, we make the change of variables $u \mapsto u'$. This change of variables is well-defined on the set $\{\theta : \cos \theta >0\}$, and from a direct computation, we can derive
\begin{align} \label{Jacobianuu'}
	\left|\frac{du'}{du}\right| = \left|\frac{1}{2}I+\frac{1}{2} \frac{(v-u)}{|v-u|} \otimes \omega \right| = \frac{1}{8}\left(1+\frac{(v-u)}{|v-u|}\cdot \omega \right) = \frac{1+\cos \theta}{8} \ge \frac{1}{8}
\end{align}
for $\theta \in [0,\frac{\pi}{2}]$.

We now present the core estimate, namely the pointwise estimate for $w\Gamma^+(f,f)$. This will be derived using the Carleman's representation from \cite{Gamba2009}. 

\bigskip

\begin{Lem} \label{pointwiseGamma+estimate} Let $0 \le \gamma \le 1$. Then the following hold:
\begin{enumerate}[label=(\arabic*)]
  \item  If $0\le \gamma \le 3/4$ and $p> 4$, then there exists a constant $C_{\beta,p,\gamma}>0$ such that
	\begin{align*}
		\left|w(v) \Gamma^+(f,f)(v) \right| \le \frac{C_{\beta,p,\gamma} \|wf\|_{L^p_v}}{(1+|v|)^{1/p'}}\left(\int_{\R^3}(1+|\eta|)^{\frac{2p-4}{p}}|f(\eta)|^2 d\eta \right)^{1/2}
	\end{align*}
	for all $v\in \R^3$, where $p'$ is the exponent conjugate of $p.$\\
  \item If $3/4< \gamma \le 1$, we divide it into two cases:
  \begin{enumerate}[label=(2-{\arabic*})]
    \item If $\frac{2}{2\gamma-1}< p \le \frac{1}{1-\gamma}$, then there exists a constant $C_{\beta,p,\gamma}>0$ such that
\begin{align*}
	\left|w(v) \Gamma^+(f,f)(v) \right| \le \frac{C_{\beta,p,\gamma}\|wf\|_{L^p_v}}{(1+|v|)^{\frac{2-p'\gamma}{p'}}}\left(\int_{\R^3} (1+|\eta|)^{\frac{3p-6}{p} } |f(\eta)|^{2}   d\eta \right)^{1/2}
\end{align*}
for all $v\in \R^3$, where $p'$ is the exponent conjugate of $p.$\\
    \item If $p > \frac{1}{1-\gamma}$, then there exists a constant $C_{\beta,p,\gamma}>0$ such that
\begin{align*}
		\left|w(v) \Gamma^+(f,f)(v) \right| \le \frac{C_{\beta,p,\gamma} \|wf\|_{L^p_v}}{(1+|v|)^{1/p'}}\left(\int_{\R^3}(1+|\eta|)^{\frac{2p-4}{p}}|f(\eta)|^2 d\eta \right)^{1/2}
	\end{align*}
	for all $v\in \R^3$, where $p'$ is the exponent conjugate of $p.$
  \end{enumerate}
\end{enumerate}
\end{Lem}
\begin{proof}
	From the energy conservation law \eqref{momentumenergyconserv}, we get $|v|^2 \le |u'|^2+|v'|^2$, which implies that
	\begin{align*}
		\text{either} \quad |v|^2 \le 2|u'|^2 \quad \text{or} \quad |v|^2 \le 2|v'|^2.
	\end{align*}
	Thus it holds that either $w(v) \le 2^{\beta/2} w(u')$ or $w(v) \le 2^{\beta/2} w(v')$. Then we have
	\begin{equation} \label{gamma1}
	\begin{aligned}
		\left|w(v) \Gamma^+(f,f)(v) \right| &\le 2^{\beta /2}\int_{\R^3}\int_{\S^2}B(v-u, \omega) \mu^{1/2}(u)\left|w(u')f(u')f(v') \right|d\omega du\\
		&\quad +2^{\beta /2}\int_{\R^3}\int_{\S^2}B(v-u, \omega) \mu^{1/2}(u)\left|f(u')w(v')f(v') \right|d\omega du.
\end{aligned}
\end{equation}
By the rotation, one can make an interchange of $v'$ and $u'$, and change the second term in \eqref{gamma1} to the same form as the first term in \eqref{gamma1}. Thus it suffices to estimate the first term in \eqref{gamma1}. We denote the first term in \eqref{gamma1} by $I_1$.\\
\indent We use the H\"older inequality to obtain
\begin{align*}
	I_1 &\le C_\beta \left( \int_{\R^3}\int_{\S^2}b(\cos \theta) \left|w(u')f(u')\right|^pd\omega du\right)^{1/p}\left(\int_{\R^3}\int_{\S^2}b(\cos \theta)|v-u|^{p'\gamma} \left|\mu^{1/2}(u)f(v') \right|^{p'}d\omega du \right)^{1/p'},
\end{align*}
where $p'$ is the exponent conjugate of $p$. By making the change of variables $u \mapsto u'$, we have
\begin{align*}
	\left(\int_{\R^3}\int_{\S^2}b(\cos \theta) \left|w(u')f(u')\right|^pd\omega du\right)^{1/p} &= \left(\int_{\R^3}\int_{\S^2}b(\cos \theta) \left|w(u')f(u')\right|^p\left|\frac{du}{du'}\right|d\omega du'\right)^{1/p}\\
	& \le C_p \|wf\|_{L^p_v},
\end{align*}
where we have used \eqref{Jacobianuu'}.
As in \cite{Gamba2009}, by using the delta function, we rewrite
\begin{align*}
	&\int_{\R^3}\int_{\S^2}b(\cos \theta)|v-u|^{p'\gamma} \left|\mu^{1/2}(u)f(v') \right|^{p'}d\omega du\\
	& = \int_{\R^3}\int_{\R^3}b\left(\frac{(v-u)}{|v-u|}\cdot \frac{k}{|k|}\right)|v-u|^{p'\gamma} \left|\mu^{1/2}(u)f(v') \right|^{p'} \delta\left(\frac{|k|^2-1}{2}\right) dk du,
\end{align*}
where $v' = \frac{v+u}{2} +\frac{|v-u|}{2}k$. We make a change of variables $u\mapsto z:=v-u$ to obtain
\begin{align} \label{zzzz1}
	I_1 \le C_{p,\beta} \|wf\|_{L^p_v} \left(\int_{\R^3}\int_{\R^3} |z|^{p'\gamma} \mu^{p'/2}(v-z)|f(v+z_-)|^{p'}  \delta\left(\frac{|k|^2-1}{2}\right)dk dz \right)^{1/p'},
\end{align}
where $z_- = -\frac{z}{2}+\frac{|z|}{2}k$ and $b$ is bounded by $C_b$. We make a change of variables $k \mapsto z_-$ with the Jacobian determinant $\left(\frac{|z|}{2}\right)^3$, and we have
\begin{align*}
	k = \frac{2z_-+z}{|z|} \quad \text{and} \quad \frac{|k|^2-1}{2} = \frac{|2z_-+z|^2-|z|^2}{2|z|^2} = \frac{2z_- \cdot (z_-+z)}{|z|^2}.
\end{align*}
Thus the integral in \eqref{zzzz1} can be written as
\begin{align*}
	J:=\int_{\R^3}\int_{\R^3} \frac{2^3}{|z|^3}|z|^{p'\gamma} \mu^{p'/2}(v-z)|f(v+z_-)|^{p'}  \delta\left(\frac{2z_- \cdot (z_-+z)}{|z|^2}\right)dz_- dz.
\end{align*}
We set $y:= -z_{-}-z$ and $\eta := v+z_-$, and then $|z|= |y+\eta-v|$ and $ \delta\left(\frac{2z_- \cdot (z_-+z)}{|z|^2}\right) = \frac{|y+\eta-v|^2}{2} \delta((\eta-v) \cdot y)$. Note that for any test function $\varphi$,
\begin{align*}
	\int_{\R^3} \delta((\eta-v) \cdot y) \varphi(y)dy = \frac{1}{|\eta-v|}\int_{(\eta-v) \cdot y=0} \varphi(y)d\pi_y ,
\end{align*}
where $d\pi_y$ is the Lebesgue measure on the hyperplane $\{y:(\eta-v)\cdot y=0\}$. It follows that
\begin{align*}
	J = 2^3\int_{\R^3}\frac{|f(\eta)|^{p'}}{|\eta-v|}\int_{(\eta-v)\cdot y=0} |y+\eta-v|^{p'\gamma-1} \mu^{p'/2}(\eta+y) d\pi_y d\eta.
\end{align*}
Note that $|y+\eta-v|^2 = |y|^2+|\eta-v|^2$ on $\{y:(\eta-v)\cdot y=0\}$.
\newline
$\mathbf{Case\ 1:}$ $\left(|y|^2+|\eta-v|^2\right)^{\frac{p'\gamma-1}{2}} \le |y|^{p'\gamma-1}$ if $p'\gamma-1 \le 0$.\\
Note that $p'\gamma-1 \le 0$ is equivalent to $p \ge \frac{1}{1-\gamma}$. In this case, we have
\begin{align*}
	I_1 \le C_{p,\beta} \|wf\|_{L^p_v} \left(\int_{\R^3}  \frac{|f(\eta)|^{p'}}{|\eta-v|}  \left(\int_{(\eta-v)\cdot y=0}|y|^{p'\gamma-1} e^{-\frac{p'|\eta+y|^2}{4}} d\pi_y\right)  d\eta \right)^{1/p'}.
\end{align*}
Here, since $-1\le p'\gamma-1<0$, we obtain
\begin{align*}
	\int_{(\eta-v)\cdot y=0}|y|^{p'\gamma-1} e^{-\frac{p'|\eta+y|^2}{4}} d\pi_y \le \int_{|y|\le 1} |y|^{p'\gamma-1}d\pi_y + \int_{|y|\ge 1}  e^{-\frac{p'|\eta+y|^2}{4}} d\pi_y \le C_{p,\gamma}
\end{align*}
for some constant $C_{p,\gamma}>0$. It follows that
\begin{align*}
	I_1 \le C_{\beta,p,\gamma}\|wf\|_{L^p_v} \left(\int_{\R^3}  \frac{|f(\eta)|^{p'}}{|\eta-v|}d\eta \right)^{1/p'}.
\end{align*}
Since $p'<2$, we use the Cauchy-Schwarz inequality to derive that
\begin{align*}
	\int_{\R^3}  \frac{|f(\eta)|^{p'}}{|\eta-v|}d\eta \le \left(\int_{\R^3} \frac{1}{|\eta-v|^q} \frac{1}{(1+|\eta|)^2}d \eta  \right)^{1/q} \left(\int_{\R^3}(1+|\eta|)^{\frac{2}{q}\cdot\frac{2}{p'}}|f(\eta)|^2 d\eta\right)^{p'/2},
\end{align*}
where $q$ satisfies $\frac{1}{q} +\frac{p'}{2}=1$ i.e. $q =\frac{2p-2}{p-2}$.\\
Since $p>4$, which is equivalent to $q<3$, it holds that
\begin{align*}
	\int_{\R^3} \frac{1}{|\eta-v|^q} \frac{1}{(1+|\eta|)^2}d \eta \le \frac{C_p}{(1+|v|)^q},
\end{align*}
which implies that
\begin{align*}
	I_1 \le \frac{C_{\beta,p,\gamma}\|wf\|_{L^p_v}}{(1+|v|)^{1/p'}}\left(\int_{\R^3}(1+|\eta|)^{\frac{2p-4}{p}}|f(\eta)|^2 d\eta\right)^{1/2}.
\end{align*}
This estimate can cover the case when $0\le \gamma \le 3/4, p>4$ or $3/4 < \gamma \le 1, p>\frac{1}{1-\gamma}$.
\newline
$\mathbf{Case\ 2:}$ $\left(|\eta-v|^2+|y|^2\right)^{\frac{p'\gamma-1}{2}} \le C_{p,\gamma}\left(|\eta-v|^{p'\gamma-1}+ |y|^{p'\gamma-1}\right)$ if $p'\gamma-1 \ge 0$.\\
Note that $p'\gamma-1 \ge 0$ is equivalent to $p \le \frac{1}{1-\gamma}$ and this case holds when $3/4 < \gamma\le  1$. In the case, we can derive
\begin{align*}
	I_1 &\le C_{\beta,p,\gamma} \|wf\|_{L^p_v} \left(\int_{\R^3}  \frac{|f(\eta)|^{p'}}{|\eta-v|^{2-p'\gamma}}  \left(\int_{(\eta-v)\cdot y=0} e^{-\frac{p'|\eta+y|^2}{4}} d\pi_y\right)  d\eta \right)^{1/p'}\\
	& \quad + C_{\beta,p,\gamma} \|wf\|_{L^p_v}\left(\int_{\R^3}  \frac{|f(\eta)|^{p'}}{|\eta-v|}  \left(\int_{(\eta-v)\cdot y=0}|y|^{p'\gamma-1} e^{-\frac{p'|\eta+y|^2}{4}} d\pi_y\right)  d\eta \right)^{1/p'}\\
	& \le  C_{\beta,p,\gamma} \|wf\|_{L^p_v} \left(\int_{\R^3}  \frac{|f(\eta)|^{p'}}{|\eta-v|^{2-p'\gamma}}   d\eta \right)^{1/p'} + C_{\beta,p,\gamma} \|wf\|_{L^p_v}\left(\int_{\R^3}  \frac{|f(\eta)|^{p'}}{|\eta-v|}   d\eta \right)^{1/p'},
	\end{align*}
	provided that $p>\frac{2}{2-\gamma}$, which is equivalent to $2-p'\gamma>0$. Since $p'\gamma -1\ge 0 $, we have
	\begin{align*}
	I_1 \le  C_{\beta,p,\gamma} \|wf\|_{L^p_v} \left(\int_{\R^3}  \frac{|f(\eta)|^{p'}}{|\eta-v|^{2-p'\gamma}}   d\eta \right)^{1/p'}.
\end{align*}
Since $p'<2$, we use the Cauchy-Schwarz inequality to derive that
\begin{align*}
	\int_{\R^3}  \frac{|f(\eta)|^{p'}}{|\eta-v|^{2-p'\gamma}}   d\eta \le \left(\int_{\R^3}  \frac{1}{|\eta-v|^{(2-p'\gamma)q}}\frac{1}{(1+|\eta|)^{3}}   d\eta \right)^{1/q}\left(\int_{\R^3} (1+|\eta|)^{\frac{3}{q}\cdot \frac{2}{p'} } |f(\eta)|^{2}   d\eta \right)^{p'/2},
\end{align*}
where $q$ satisfies $\frac{1}{q} +\frac{p'}{2}=1$ i.e. $q =\frac{2p-2}{p-2}$. Then it holds that 
\begin{align*}
	\int_{\R^3}  \frac{1}{|\eta-v|^{(2-p'\gamma)q}}\frac{1}{(1+|\eta|)^{3}}   d\eta \le \frac{C_{p,\gamma}}{(1+|v|)^{(2-p'\gamma)q}},
\end{align*}
provided that $p>\frac{2}{2\gamma-1}$, which means $(2-p'\gamma)q<3$. It follows that
\begin{align*}
	I_1 &\le \frac{C_{\beta,p,\gamma}\|wf\|_{L^p_v}}{(1+|v|)^{\frac{2-p'\gamma}{p'}}}\left(\int_{\R^3} (1+|\eta|)^{\frac{3}{q}\cdot \frac{2}{p'} } |f(\eta)|^{2}   d\eta \right)^{1/2}\\
	& = \frac{C_{\beta,p,\gamma}\|wf\|_{L^p_v}}{(1+|v|)^{\frac{2-p'\gamma}{p'}}}\left(\int_{\R^3} (1+|\eta|)^{\frac{3p-6}{p} } |f(\eta)|^{2}   d\eta \right)^{1/2}.
\end{align*}
This estimate can cover the case when $3/4 < \gamma \le 1, \frac{2}{2\gamma -1}<p\le \frac{1}{1-\gamma}$.
\end{proof}

\bigskip

We present the following corollary to this lemma. From Lemma \ref{pointwiseGamma+estimate}, we can derive the weighted $L^p_v$ estimate for the gain term $\Gamma^+$.

\bigskip
\begin{Coro} \label{LpGamma+est}		Let $0 \le \gamma \le 1$. Let $p$ satisfy the condition
\begin{equation} \label{pcondition}
	\begin{aligned}
	p> \begin{cases}
 4, & \text{if}\quad 0\le \gamma \le \frac{3}{4}, \\
 \frac{5}{2-\gamma}, & \text{if}\quad \frac{3}{4}<\gamma \le 1
\end{cases}
\end{aligned}
\end{equation} 
and $\beta$ satisfy the condition
\begin{equation} \label{betacondition}
	\begin{aligned}
	\beta> \begin{cases}
 \frac{5p-10}{2p}, & \text{if} \quad 0\le \gamma \le \frac{3}{4},p>4\ \text{or}\ \frac{3}{4}<\gamma \le 1, p>\frac{1}{1-\gamma},  \\
 \frac{3p-6}{p}, & \text{if}\quad \frac{3}{4}<\gamma \le 1, \frac{5}{2-\gamma}< p \le \frac{1}{1-\gamma}.
\end{cases}
\end{aligned}
\end{equation}
Then there exists a constant $C_{\beta,p,\gamma}>0$ such that
	\begin{align*}
		\left\|w\Gamma^+(f,f)\right\|_{L^p_v} \le C_{\beta,p,\gamma} \|wf\|_{L^p_v}^2.
	\end{align*}
\end{Coro}
\begin{proof}
	We consider the estimates in Lemma \ref{pointwiseGamma+estimate}.
	If $0\le \gamma\le 3/4$ and $p>4$, we have
	\begin{align*}
		\left|w(v) \Gamma^+(f,f)(v) \right| \le \frac{C_{\beta,p,\gamma} \|wf\|_{L^p_v}}{(1+|v|)^{1/p'}}\left(\int_{\R^3}(1+|\eta|)^{\frac{2p-4}{p}}|f(\eta)|^2 d\eta \right)^{1/2}
	\end{align*}
	for all $v\in \R^3$.
	Here, we use the H\"older inequality to derive
	\begin{align*}
		\left(\int_{\R^3}(1+|\eta|)^{\frac{2p-4}{p}}|f(\eta)|^2 d\eta \right)^{1/2} &\le \left(\int_{\R^3}\left|\frac{(1+|\eta|)^{\frac{2p-4}{p}}}{w^2(\eta)}\right|^{\frac{p}{p-2}} d\eta\right)^{\frac{p-2}{2p}}\|wf\|_{L^p_v}\\
		& \le C_{\beta,p}\|wf\|_{L^p_v}
	\end{align*}
	since $\beta > \frac{5p-10}{2p}$. Thus we obtain
	\begin{align*}
		\|w\Gamma^+(f,f)\|_{L^p_v} \le C_{\beta,p,\gamma}\|wf\|_{L^p_v}^2
	\end{align*}
	because $p>4$.\\
	When $3/4<\gamma \le 1$ and $p > \frac{1}{1-\gamma}$, it is similar to the previous case.\\
	If $3/4 < \gamma \le 1$ and $\frac{5}{2-\gamma}< p \le \frac{1}{1-\gamma}$, we have
	\begin{align*}
	\left|w(v) \Gamma^+(f,f)(v) \right| \le \frac{C_{\beta,p,\gamma}\|wf\|_{L^p_v}}{(1+|v|)^{\frac{2-p'\gamma}{p'}}}\left(\int_{\R^3} (1+|\eta|)^{\frac{3p-6}{p} } |f(\eta)|^{2}   d\eta \right)^{1/2}
\end{align*}
for all $v\in \R^3$. Here, we use the H\"older inequality to obtain
	\begin{align*}
		\left(\int_{\R^3}(1+|\eta|)^{\frac{2p-4}{p}}|f(\eta)|^2 d\eta \right)^{1/2} &\le \left(\int_{\R^3}\left|\frac{(1+|\eta|)^{\frac{3p-6}{p}}}{w^2(\eta)}\right|^{\frac{p}{p-2}} d\eta\right)^{\frac{p-2}{2p}}\|wf\|_{L^p_v}\\
		& \le C_{\beta,p}\|wf\|_{L^p_v}
	\end{align*}
	since $\beta > \frac{3p-6}{p}$. Thus we can derive
	\begin{align*}
		\|w\Gamma^+(f,f)\|_{L^p_v} \le C_{\beta,p,\gamma}\|wf\|_{L^p_v}^2
	\end{align*}
	because $p>\frac{5}{2-\gamma}$. We conclude the proof.\\

\end{proof}

\bigskip

We define the constant $e(p,\gamma)$, depending on $p$ and $\gamma$, as follows:
\begin{equation} \label{expcstp}
	\begin{aligned}
	e(p,\gamma)= \begin{cases}
 \frac{2p-4}{p}, & \text{if}\quad 0\le \gamma \le \frac{3}{4}, p>4\ \text{or}\ \frac{3}{4}<\gamma \le 1, p>\frac{1}{1-\gamma} \\
 \frac{3p-6}{p}, & \text{if}\quad \frac{3}{4}<\gamma \le 1, \frac{2}{2\gamma-1}< p \le \frac{1}{1-\gamma}.
\end{cases}
\end{aligned}
\end{equation}
By using the constant $e(p,\gamma)$, we can combine estimates for all cases in Lemma \ref{pointwiseGamma+estimate} and it can be restated as follows:
\begin{align} \label{Gamma+estimatecombine}
	\left|w(v) \Gamma^+(f,f)(v) \right| \le \frac{C_{\beta,p,\gamma} \|wf\|_{L^p_v}}{(1+|v|)^{\min\left\{\frac{1}{p'},\frac{2-p'\gamma}{p'}\right\}}}\left(\int_{\R^3}(1+|\eta|)^{e(p,\gamma)}|f(\eta)|^2 d\eta \right)^{1/2}
\end{align}
for all $v\in \R^3$, where $C_{\beta,p,\gamma}>0$ is a constant depending on $\beta$, $p$, and $\gamma$. In what follows, we will primarily make use of the form \eqref{Gamma+estimatecombine} of the pointwise estimate for $w\Gamma^+(f,f)$.
 
\bigskip

\subsection{Relative entropy and its application}
Recall the definition for relative entropy $\mathcal{E}(F)$
\begin{align*}
	\mathcal{E}(F)= \int_{\Omega \times \mathbb{R}^3} \left(\frac{F}{\mu}\log\frac{F}{\mu}-\frac{F}{\mu}+1\right) \mu dxdv.
\end{align*}
We introduce some lemmas about relative entropy. The following lemma provides the global-in-time a priori estimate of relative entropy. This lemma means that relative entropy is non-increasing over time. Later on, we will employ it as a bridge between the small perturbation regime and the large amplitude problem.
\bigskip
\begin{Lem} \label{relativedecrease} \cite{Jong2025}
	Assume $F$ satisfies the Boltzmann equation \eqref{Boltzmanneq} with the diffuse reflection boundary condition \eqref{Diffuseboundarycond} and initial data $F_0$. Then
	\begin{align*}
		\mathcal{E}(F(t)) \le \mathcal{E}(F_0) \quad \text{for all } t\ge 0.
	\end{align*}
\end{Lem}

\bigskip

 The following lemma states initial relative entropy provides the bound for the $L^1$ and $L^2$ norms of $F-\mu$ on the different domains.

\bigskip
\begin{Lem} \label{L1L2cont} \cite{Guo2010b}
	Assume $F$ satisfies the Boltzmann equation \eqref{Boltzmanneq}. We have
	\begin{align*}
		\int_{\Omega \times \mathbb{R}^3} \frac{1}{4\mu}|F-\mu|^2 \mathbf{1}_{\{|F-\mu|\le \mu\}} dxdv + \int_{\Omega \times \mathbb{R}^3} \frac{1}{4}|F-\mu| \mathbf{1}_{\{|F-\mu| > \mu\}} dxdv \le \mathcal{E}(F_0)
	\end{align*}
	for all $t \ge 0$.
	Moreover, if we write $F = \mu + \mu^{1/2} f$, then
	\begin{align*}
		\int_{\Omega \times \mathbb{R}^3} \frac{1}{4}|f|^2 \mathbf{1}_{\{|f|\le \sqrt{\mu}\}} dxdv + \int_{\Omega \times \mathbb{R}^3} \frac{\sqrt{\mu}}{4}|f| \mathbf{1}_{\{|f| > \sqrt{\mu}\}} dxdv \le \mathcal{E}(F_0)
	\end{align*}
	for all $t \ge 0$.
\end{Lem}

\section{Small perturbation problem in $L^p_vL^\infty_x$ flamework} \label{smallperturbationproblem}
In this section, our aim is to prove the small perturbation problem in the $L^p_vL^\infty_x$ space. Recall the velocity weight function 
\begin{align} \label{weightfnt}
	w = w(v): = (1+|v|^2)^{\beta/2}
\end{align}
with $\beta > 0$.
Setting $h(t,x,v)=w(v)f(t,x,v)$, the equation \eqref{FPBER} with weight $w(v)$ becomes
\begin{align} \label{WFPBE}
	\partial_t h + v\cdot \nabla_x h + \nu(v) h -K_wh= w\Gamma(f,f),
\end{align}
where
\begin{align}\label{weightedK}
	K_wh=wK\left(\frac{h}{w}\right),
\end{align}
and the diffuse reflection boundary condition \eqref{PDRBC} with weight $w(v)$ can be written as
\begin{align} \label{WPDRBC}
	h(t,x,v)|_{\gamma_-} = \frac{1}{\tilde{w}(v)} \int_{\mathcal{V}(x) }h(t,x,v')\tilde{w}(v')d\sigma(x),
\end{align}
where $\mathcal{V}(x) = \{v' \in \R^3 : n(x) \cdot v' >0\}$, $\tilde{w}(v)=\frac{1}{w(v)\mu^{1/2}(v)}$, and $d\sigma(x) = c_\mu \mu(v')\{n(x) \cdot v'\}dv'$.
Note that from \eqref{diffusemassconvconst}, $d\sigma(x)$ is the probability measure, in other words, it satisfies
\begin{align*}
	 \int_{\mathcal{V}(x)}d\sigma(x) =1.
\end{align*}

\subsection{$L^p_vL^\infty_x$ estimate for Damped transport equation}In this subsection, before we treat the nonlinear equation \eqref{WFPBE}, we consider the following elementary equation:
\begin{align} \label{WDBE1}
	\{\partial_t+v\cdot \nabla_x +\nu(v)\}h=0
\end{align}
with initial data $h_0$ and the diffuse reflection boundary condition \eqref{WPDRBC}. We will study the solution operator for the equation \eqref{WDBE1} with the diffuse reflection boundary condition and the related estimates.\\
Before  dealing with the solution operator for the equation \eqref{WDBE1}, we present a simple lemma related to the back-time cycle and the iterated integral.
Recall subsection \ref{Backtimecycle}, especially \eqref{Backcycle} and \eqref{iteratedintegral}. The following lemma means the set in the phase space $\prod_{j=1}^{k-1} \mathcal{V}_j$ not reaching $t=0$ after bouncing $k$ times is small when $k$ is large.
\begin{Lem} \label{Lsmall} \cite{GuoDecay}
	For any $\epsilon > 0$, there exists $k_0(\epsilon, T_0)$ such that for $k \ge k_0$, for all $(t,x,v)\in [0,T_0]\times \bar{\Omega}\times \R^3$,
	\begin{align} \label{Lsmallineq}
		\int_{\prod_{j=1}^{k-1}\mathcal{V}_j}\mathbf{1}_{\{t_k(t,x,v,v_1,v_2,\cdots ,v_{k-1})>0\}} \prod_{j=1}^{k-1}d\sigma_j \le \epsilon.	
	\end{align}
	Furthermore, for $T_0$ sufficiently large, there exist constant $\tilde{C}_1, \tilde{C}_2>0$, independent of $T_0$, such that for $k=\tilde{C}_1T_0^{5/4}$,
	\begin{align*}
		\int_{\prod_{j=1}^{k-1}}\mathbf{1}_{\{t_k(t,x,v,v_1,v_2,\cdots,v_{k-1})>0\}} \prod_{j=1}^{k-1}d\sigma_j \le \left(\frac{1}{2}\right)^{\tilde{C}_2T_0^{5/4}}.
	\end{align*}	
\end{Lem}

\bigskip

The following lemma present the mild formulation for the inhomogeneous form of the equation \eqref{WDBE1} with the diffuse reflection boundary condition \eqref{WPDRBC}. Recall the definition of the iterated integral \eqref{iteratedintegral}.
\begin{Lem} \label{RepresentationforDiffuse1} \cite{GuoDecay}
	Assume that $h,q \in L^p_vL^\infty_x$ satisfy $\{\partial_t+v\cdot \nabla_x +\nu(v)\}h=q(t,x,v)$ with the diffuse reflection boundary condition \eqref{WPDRBC}. Then for any $0\le s \le t$, for almost every $x,v$, if $t_1(t,x,v) \le s$,
	\begin{align*}
		h(t,x,v) = e^{-\nu(v)(t-s)}h(s,x-v(t-s),v)+\int_s^te^{-\nu(v)(t-\tau)}q(\tau, x-v(t-\tau),v)d\tau;
	\end{align*}
	If $t_1(t,x,v)>s$, then for $k\ge 2$,
	\begin{align*}
		h(t,x,v) &= \int_{t_1}^t e^{-\nu(v)(t-\tau)}q(\tau, x-v(t-\tau),v)d\tau\\
		&\quad +\frac{e^{-\nu(v)(t-t_1)}}{\tilde{w}(v)} \sum_{l=1}^{k-1} \int_{\prod_{j=1}^{k-1}\mathcal{V}_j}\mathbf{1}_{\{t_{l+1}\le s <t_l\}} h(s,x_l-v_l(t_l-s),v_l)d\Sigma_l(s)\\
		&\quad +\frac{e^{-\nu(v)(t-t_1)}}{\tilde{w}(v)} \sum_{l=1}^{k-1} \int_{s}^{t_l}\int_{\prod_{j=1}^{k-1}\mathcal{V}_j} \mathbf{1}_{\{t_{l+1}\le s <t_l\}}q(\tau, x_l-v_l(t_l-\tau),v_l)d\Sigma_l(\tau)d\tau\\
		&\quad +\frac{e^{-\nu(v)(t-t_1)}}{\tilde{w}(v)} \sum_{l=1}^{k-1} \int_{t_{l+1}}^{t_l}\int_{\prod_{j=1}^{k-1}\mathcal{V}_j} \mathbf{1}_{\{t_{l+1}>s\}}q(\tau, x_l-v_l(t_l-\tau),v_l) d\Sigma_l(\tau)d\tau\\
		&\quad +\frac{e^{-\nu(v)(t-t_1)}}{\tilde{w}(v)} \int_{\prod_{j=1}^{k-1}\mathcal{V}_j} \mathbf{1}_{\{t_k>s\}}h(t_k,x_k,v_{k-1})d\Sigma_{k-1}(t_k),
	\end{align*}
	where 
	\begin{align*}
		d\Sigma_{l}(s) = \left\{\prod_{j=l+1}^{k-1}d\sigma_j \right\}\left\{e^{-\nu(v_l)(t_l-s)} \tilde{w}(v_l)d\sigma_l \right\}\prod_{j=1}^{l-1} \left\{ e^{-\nu(v_j)(t_j-t_{j+1})}d\sigma_j \right\}.
	\end{align*}
\end{Lem}

\bigskip
We now prove the existence and uniqueness for the equation \eqref{WDBE1} with initial data $h_0$ and the boundary condition \eqref{WPDRBC}. Before proving it, we introduce some lemmas that aid in achieving our goal in this subsection. The below lemma presents an explicit formula for a solution to the equation \eqref{WDBE1} with the in-flow boundary condition.
\begin{Lem}\label{inflowref1} \cite{GuoDecay}
	Let $p \ge 1$ and $h_0 \in L^p_vL^\infty_x$, $g(t,x,v) \in L^\infty_xL^p_v(\gamma_-)$. Let $h(t,x,v)$ be the solution to the following equation: 
	\begin{align*}
		&\{ \partial_t +v\cdot \nabla_x+\nu(v) \}h =0, \quad h(0,x,v)=h_0(x,v), \quad h(t,x,v)|_{\gamma_-}=g(t,x,v).
	\end{align*}
	Then for $(t,x,v) \in (0,\infty) \times \bar{\Omega} \times \R^3$ with $(x,v) \notin \gamma_0 \cup \gamma_-$,
	\begin{align*}
		h(t,x,v) &= \mathbf{1}_{t-t_\mathbf{b}<0}e^{-\nu_0 t} h_0(x-vt,v)\\
		& \quad +\mathbf{1}_{t-t_\mathbf{b}\ge0} e^{-\nu_0 t_\mathbf{b}} g(t-t_\mathbf{b}, x_\mathbf{b},v),
	\end{align*}
	where $t_\mathbf{b}$ and $x_\mathbf{b}$ is defined in Definition \ref{Backwardexit}.
\end{Lem}
\begin{proof}
	It is similar to \cite[Lemma 12]{GuoDecay}, so we omit the proof.
\end{proof}
As in \cite[Lemma14]{GuoDecay}, we use a similar argument to prove Lemma \ref{DampLpDecayorgin}. However, unlike in \cite{GuoDecay}, we consider the $L^p_vL^\infty_x$ space instead of the $L^\infty_{x,v}$ space. Thus, we propose a variation of the statement in \cite[Lemma14]{GuoDecay} on the diffuse reflection boundary condition. 
Let $p \ge 1$. Let $h_0\in L^p_vL^\infty_x$ and $\epsilon>0$. We consider a solution $h(t)$ to the equation 
\begin{align*}
	&\{ \partial_t +v\cdot \nabla_x+\nu(v) \}h =0, \quad h(0)=h_0,\\
	& \quad h(t,x,v)|_{\gamma_-}=(1-\epsilon) \mathcal{M} h|_{\gamma_+}:=(1-\epsilon)(c_\mu \mu(v))^{1/p}\int_{n(x)\cdot v'>0} h(t,x,v')\frac{1}{(c_\mu \mu(v'))^{1/p}} d\sigma(x).
\end{align*}
Fix $\epsilon >0$. Let $h_{\gamma_+}^{(0)}\equiv 0$. We construct the following iteration system for $k=0,1,2,\cdots$ :
\begin{align*}
	&\{ \partial_t +v\cdot \nabla_x+\nu(v)\}h^{(k+1)} =0, \quad h^{(k+1)}(0,x,v)=h_0,\\
	& \quad h^{(k+1)}_{\gamma_-}=(1-\epsilon) \mathcal{M} h^{(k)}_{\gamma_+}:=(1-\epsilon)(c_\mu \mu(v))^{1/p}\int_{n(x)\cdot v'>0} h^{(k)}_{\gamma_+}\frac{1}{(c_\mu \mu(v'))^{1/p}} d\sigma(x).
\end{align*}
We claim that $\{h^{(k)}\}$ and $\{h_\gamma^{(k)}\}$ are Cauchy sequences. Taking difference, we get
\begin{equation} \label{itersys10}
	\begin{aligned}
	&\{ \partial_t +v\cdot \nabla_x+\nu(v) \}\{h^{(k+1)}-h^{(k)}\} =0, \quad \{h^{(k+1)}-h^{(k)}\}(0,x,v)=0,\\
	& \quad h^{(k+1)}_{\gamma_-}-h^{(k)}_{\gamma_-}=(1-\epsilon) \mathcal{M}\{ h^{(k)}_{\gamma_+}-h^{(k-1)}_{\gamma_+}\}=(1-\epsilon)(c_\mu \mu(v))^{1/p}\int_{n(x)\cdot v'>0} \{ h^{(k)}_{\gamma_+}-h^{(k-1)}_{\gamma_+}\}\frac{1}{(c_\mu \mu(v'))^{1/p}} d\sigma(x).
\end{aligned}
\end{equation}
By Lemma \ref{inflowref1}, we have
\begin{align} \label{differepresen}
	\{h^{(k+1)}-h^{(k)}\}(t,x,v) = \mathbf{1}_{t-t_\mathbf{b} \ge 0} (1-\epsilon)\mathcal{M}\{ h^{(k)}_{\gamma_+}-h^{(k-1)}_{\gamma_+}\}(t-t_\mathbf{b},x_\mathbf{b},v)
\end{align}
for $(t,x,v) \in (0,\infty) \times \bar{\Omega} \times \R^3$ with $(x,v) \notin \gamma_0 \cup \gamma_-$.\\
Taking the norm $\|\cdot\|_{L^\infty_xL^p_v (\gamma_+)}$ to \eqref{differepresen} for $(t,x,v)\in (0,\infty) \times \gamma_+$,
\begin{align*}
	\|\{h^{(k+1)}_{\gamma_+}-h^{(k)}_{\gamma_+}\}(t)\|_{L^\infty_xL^p_v (\gamma_+)} &\le (1-\epsilon)\|\mathcal{M}\{ h^{(k)}_{\gamma_+}-h^{(k-1)}_{\gamma_+}\}(t-t_\mathbf{b}) \|_{L^\infty_xL^p_v (\gamma_+)}\\
	&\le (1-\epsilon) \sup_{x\in \partial\Omega}\Biggl[ \left(\int_{n(x)\cdot v>0}c_\mu \mu(v)|n(x) \cdot v|dv\right)^{1/p}\\
	& \quad \times   \int_{n(x)\cdot v'>0} \left| \{h^{(k)}_{\gamma_+}-h^{(k-1)}_{\gamma_+}\}(t-t_\mathbf{b},x,v')\right|(c_\mu \mu(v'))^{1-\frac{1}{p}} |n(x) \cdot v'|dv' \Biggr] \\
	& = (1-\epsilon) \sup_{x\in \partial\Omega}\Biggl[ \int_{n(x)\cdot v'>0} \left| \{h^{(k)}_{\gamma_+}-h^{(k-1)}_{\gamma_+}\}(t-t_\mathbf{b},x,v')\right|(c_\mu \mu(v'))^{1-\frac{1}{p}} |n(x) \cdot v'|dv' \Biggr] \\
	& \le (1-\epsilon)\sup_{x\in \partial\Omega}\Biggl[\left(\int_{n(x)\cdot v'>0} \left|\{h^{(k)}_{\gamma_+}-h^{(k-1)}_{\gamma_+}\}(t-t_\mathbf{b},x,v')\right|^p |n(x)\cdot v'|dv'\right)^{1/p}\\
	&\quad \times \left(\int_{n(x)\cdot v'>0}c_\mu \mu(v') |n(x)\cdot v'|dv'\right)^{1/p'}\Biggr]\\
	& \le (1-\epsilon) \|\{h^{(k)}_{\gamma_+}-h^{(k-1)}_{\gamma_+}\}(t-t_\mathbf{b})\|_{L^\infty_xL^p_v (\gamma_+)},
\end{align*}
where we have used the H\"older inequality and
\begin{align*}
	\int_{n(x)\cdot v>0}c_\mu \mu(v)|n(x) \cdot v|dv =\int_{n(x)\cdot v<0}c_\mu \mu(v)|n(x) \cdot v|dv =1.
\end{align*}
It follows that
\begin{align} \label{midob}
	\sup_t\|\{h^{(k+1)}_{\gamma_+}-h^{(k)}_{\gamma_+}\}(t)\|_{L^\infty_xL^p_v (\gamma_+)} \le (1-\epsilon)\sup_t \|\{h^{(k)}_{\gamma_+}-h^{(k-1)}_{\gamma_+}\}(t)\|_{L^\infty_xL^p_v (\gamma_+)},
\end{align}
which by induction implies that
\begin{align*}
	\sup_t\|\{h^{(k+1)}_{\gamma_+}-h^{(k)}_{\gamma_+}\}(t)\|_{L^\infty_xL^p_v (\gamma_+)} \le (1-\epsilon)^k\sup_t \|h^{(1)}_{\gamma_+}(t)\|_{L^\infty_xL^p_v (\gamma_+)}.
\end{align*}
Thus, $\{h^{(k)}_{\gamma_+}\}$ is a Cauchy sequence in $L^\infty_xL^p_v (\gamma_+)$. Similarly, taking the norm $\|\cdot\|_{L^\infty_xL^p_v (\gamma_-)}$ to the boundary condition in \eqref{itersys10}, we can show that $\{h^{(k)}_{\gamma_-}\}$ is a Cauchy sequence in $L^\infty_xL^p_v (\gamma_-)$.\\
Now, all that remains is to prove that $\{h^{(k)}\}$ is a Cauchy sequence in $L^p_vL^\infty_x$. Taking the norm $\|\cdot\|_{L^p_vL^\infty_x}$ to \eqref{differepresen}, we get
\begin{align*}
	\|\{h^{(k+1)}-h^{(k)}\}(t)\|_{L^p_vL^\infty_x}  
	&\le (1-\epsilon) \left(\int_{\R^3}c_\mu \mu(v)dv\right)^{1/p} \\
	& \quad \times   \sup_{x\in \partial\Omega}\Biggl[\int_{n(x)\cdot v'>0} \left| \{h^{(k)}_{\gamma_+}-h^{(k-1)}_{\gamma_+}\}(t-t_\mathbf{b},x,v')\right|(c_\mu \mu(v'))^{1-\frac{1}{p}} |n(x) \cdot v'|dv' \Biggr] \\
	& \le (1-\epsilon)(2\pi)^{\frac{1}{2p}}\sup_{x\in \partial\Omega}\Biggl[\left(\int_{n(x)\cdot v'>0} \left|\{h^{(k)}_{\gamma_+}-h^{(k-1)}_{\gamma_+}\}(t-t_\mathbf{b},x,v')\right|^p |n(x)\cdot v'|dv'\right)^{1/p}\\
	&\quad \times \left(\int_{n(x)\cdot v'>0}c_\mu \mu(v') |n(x)\cdot v'|dv'\right)^{1/p'}\Biggr]\\
	& \le (1-\epsilon) (2\pi)^{\frac{1}{2p}}\|\{h^{(k)}_{\gamma_+}-h^{(k-1)}_{\gamma_+}\}(t-t_\mathbf{b})\|_{L^\infty_xL^p_v (\gamma_+)}.
\end{align*}
where we have used the H\"older inequality and
\begin{align*}
	\int_{\R^3}c_\mu \mu(v)dv= (2\pi)^{1/2}, \qquad \int_{n(x)\cdot v>0}c_\mu \mu(v)|n(x) \cdot v|dv =1.
\end{align*}
It follows from \eqref{midob} that
\begin{align*}
	\sup_t\|\{h^{(k+1)}-h^{(k)}\}(t)\|_{L^p_v L^\infty_x} \le (1-\epsilon)^k(2\pi)^{\frac{1}{2p}} \sup_t \|h^{(1)}_{\gamma_+}(t)\|_{L^\infty_xL^p_v (\gamma_+)}.
\end{align*}
Hence, $\{h^{(k)}\}$ is a Cauchy sequence in $L^p_vL^\infty_x$. Taking $k \rightarrow \infty$, we conclude the following result:

\begin{Lem} \label{Existencesystem1}
Let $p\ge 1$ and $h_0 \in L^p_vL^\infty_x$. Then for any $\epsilon>0$, there exists $h(t) \in L^p_vL^\infty_x$ and $h_\gamma \in L^\infty_xL^p_v(\gamma)$ solving
\begin{align*}
	&\{ \partial_t +v\cdot \nabla_x+\nu(v) \}h =0, \quad h(0)=h_0,\\
	& \quad h_{\gamma_-}(x,v)=(1-\epsilon)(c_\mu \mu(v))^{1/p}\int_{n(x)\cdot v'>0} h_{\gamma_+}(x,v')\frac{1}{(c_\mu \mu(v'))^{1/p}} d\sigma(x).
\end{align*}
\end{Lem}

We denote by $S_{G_\nu}(t,0)h_0$ the solution operator to the equation $\eqref{WDBE1}$ with initial data $h_0$ and the diffuse reflection boundary condition \eqref{WPDRBC}. We now introduce an useful exponential time-decay to $S_{G_\nu}(t,0)h_0$ in $L^p_vL^\infty_x$.

\bigskip
\begin{Lem} \label{DampLpDecayorgin}
	Let $p>2$, $\beta > \frac{2p-3}{p}$ and $h_0 \in L^p_vL_x^\infty$. There exists a unique solution $h(t) = S_{G_{\nu}}(t,0)h_0  \in L^p_vL^\infty_x$ to  \eqref{WDBE1} with initial data $h_0$ and the diffuse reflection boundary condition \eqref{WPDRBC}. Moreover, for all $0< \tilde{\nu}_0 < \nu_0$, there exists $C_{p,\beta}>0$, depending on $p$ and $\beta$, such that 
	\begin{align*}
		\sup_{t\ge 0}\left\{e^{\tilde{\nu}_0 t}\|S_{G_{\nu}}(t,0)h_0\|_{L^p_vL_x^\infty}\right\} \le C_{p,\beta}\|h_0\|_{L^p_vL_x^\infty}.
	\end{align*}
\end{Lem}
\begin{proof}
	We first show the uniqueness of solution. Assume that there exists two solutions $h=wf,\tilde{h}=w\tilde{f}$ in $L^p_vL^\infty_x$. Since $\beta > \frac{2p-3}{p}$, $\|f\|_{L^2_{x,v}}\le \|wf\|_{L^p_vL^\infty_x} \left\|\frac{1}{w^2}\right\|_{L^{\frac{p}{p-2}}_v}< \infty$ and $\|f\|_{L^2_{\gamma}}\le \|wf\|_{L^p_v L^\infty_x}\left\| \frac{|v|}{w^2} \right\|_{L^{\frac{p}{p-2}}_v} < \infty$, and thus $f$, $\tilde{f}$ are in $L^2_{x,v}$ and $f|_{\gamma},\tilde{f}|_{\gamma}$ are in $L^2_{\gamma}$.
	By the energy estimate in $L^2$, we conclude that a solution is unique.\\
	Let $h(t,x,v) =w(v)f(t,x,v)$.  
	Given any $m \ge 1$, we construct a solution to
	\begin{align} \label{qq1}
		\{\partial_t + v\cdot \nabla_x +  \nu(v)\}h^{(m)}=0,
	\end{align}
	with boundary and initial condition
	\begin{equation} \label{qq2}
	\begin{aligned}
		&h^{(m)}(t,x,v) |_{\gamma_-}= \left(1-\frac{1}{m}\right)\frac{1}{\tilde{w}(v)}\int_{n(x)\cdot v'>0} \left[h^{(m)}(t,x,v')\right]\tilde{w}(v')d\sigma(x),\\
		&h^{(m)}(0,x,v) = h_0 \mathbf{1}_{\{|v|\le m\}}.
	\end{aligned}
	\end{equation}
	Setting $\tilde{h}^{(m)}(t,x,v) = (c_\mu \mu(v))^{1/p}\tilde{w}(v) h^{(m)}(t,x,v)$, where $p\ge 1$, the equation \eqref{qq1} and the condition \eqref{qq2} become
	\begin{equation} \label{qq8}
	\begin{aligned}
		&\{\partial_t + v\cdot \nabla_x + \nu(v)\}\tilde{h}^{(m)}=0,\\
		&\tilde{h}^{(m)}(t,x,v)  |_{\gamma_-}= \left(1-\frac{1}{m}\right)(c_\mu \mu(v))^{1/p}\int_{n(x)\cdot v'>0} \tilde{h}^{(m)}(t,x,v')\frac{1}{(c_\mu \mu(v'))^{1/p}} d\sigma(x),\\
		& \tilde{h}^{(m)}(0,x,v) = \tilde{h}_0 \mathbf{1}_{\{|v|\le m\}}
	\end{aligned}
	\end{equation}
	Note that
	\begin{align*}
		\|\tilde{h}^{(m)}(0)\|_{L^p_vL^\infty_{x}} \le C_{m}\|h_0\|_{L^p_vL^\infty_{x}} < \infty.
	\end{align*}
	By Lemma \ref{Existencesystem1}, there exists a solution $\tilde{h}^{(m)}(t,x,v) \in L^p_vL^\infty_{x}$ to the equation \eqref{qq8}, so that we have constructed $h^{(m)}$, which is bounded.\\ 
	\indent We show the uniform $L^p_vL^\infty_{x}$ bound for $h^{(m)}$. Assume that $T_0$ is sufficiently large. We consider the case $0\le t \le T_0$.\\
	\indent If $t_1(t,x,v) \le 0$, we know  
	\begin{align*}
		\left(S_{G_\nu}(t,0)h_0\right)(x,v) = e^{-\nu(v) t}h_0(x-vt,v),
	\end{align*}
	and we deduce
	\begin{align*} 
		\|S_{G_\nu}(t,0)h_0\|_{L^p_vL^\infty_{x}} \le e^{-\nu_0 t}\|h_0\|_{L^p_vL^\infty_{x}} \quad \text{for all }0\le t \le T_0.
	\end{align*}
	\indent We consider the case $t_1(t,x,v) > 0$. By Lemma \ref{RepresentationforDiffuse1}, we get
	\begin{align*}
		\left|h^{(m)}(t,x,v)\right| &\le \frac{e^{-\nu_0(t-t_1)}}{\tilde{w}(v)}\sum_{l=1}^{k-1}\int_{\prod_{j=1}^{k-1}\mathcal{V}_j} \mathbf{1}_{\{t_{l+1}\le 0 < t_l\}} \left|h^{(m)}(0,x_l-v_lt_l,v_l)\right|d\Sigma_l(0)\\
		&\quad +\frac{e^{-\nu_0(t-t_1)}}{\tilde{w}(v)}\int_{\prod_{j=1}^{k-1}\mathcal{V}_j} \mathbf{1}_{\{t_k>0\}} \left|h^{(m)}(t_k,x_k,v_{k-1})\right|d\Sigma_{k-1}(t_k)\\
		&=: I_1+I_2
	\end{align*}
	First of all, let us consider $I_2$. Using the boundary condition
	\begin{align*}
		h^{(m)}(t_k,x_k,v_{k-1}) = \left(1-\frac{1}{m}\right)\frac{1}{\tilde{w}(v_{k-1})}\int_{\mathcal{V}_k} h^{(m)}(t_k,x_k,v_k)\tilde{w}(v_k)d\sigma_k
	\end{align*}
and the fact $h^{(m)}(t_k,x_k,v_k) = \mathbf{1}_{\{t_{k+1} \le 0 <t_k\}}e^{-\nu_0t_k}h^{(m)}(0,x_k-v_kt_k,v_k) +\mathbf{1}_{\{t_{k+1}>0\}}h^{(m)}(t_k,x_k,v_k)$,
\begin{align*}
	I_2 &\le \frac{e^{-\nu_0(t-t_1)}}{\tilde{w}(v)}\int_{\prod_{j=1}^{k}\mathcal{V}_j} \mathbf{1}_{\{t_{k+1} \le 0 <t_k\}}\left|h^{(m)}(0,x_k-v_kt_k,v_k)\right|d\Sigma_{k}(0)\\
	&\quad + \frac{e^{-\nu_0(t-t_1)}}{\tilde{w}(v)}\int_{\prod_{j=1}^{k}\mathcal{V}_j} \mathbf{1}_{\{t_{k+1}>0\}} \left|h^{(m)}(t_k,x_k,v_k)\right|d\Sigma_{k}(t_k)\\
	&=:J_1+J_2.
\end{align*}
We know that the exponential in $d\Sigma_l(s)$ is bounded by $e^{-\nu_0(t_1-s)}$. By Lemma \ref{Lsmall}, we can choose $\tilde{C}_1,\tilde{C}_2>0$ such that for $k=\tilde{C}_1T_0^{\frac{5}{4}}$
\begin{align} \label{qq6}
	\int_{\prod_{j=1}^{k-1}\mathcal{V}_j} \mathbf{1}_{\{t_k(t,x,v,v_1,v_2,\cdots,v_{k-1})>0\}}\prod_{j=1}^{k-1}d\sigma_j \le \left(\frac{1}{2}\right)^{\tilde{C}_2T_0^{\frac{5}{4}}}.
\end{align}
Using \eqref{qq6}, we obtain
\begin{align*}
	J_2 &\le \frac{e^{-\nu_0(t-t_k)}}{\tilde{w}(v)}\left\|h^{(m)}(t_k)\mathbf{1}_{\{t_1>0\}}\right\|_{L^p_vL^\infty_x} \int_{\prod_{j=1}^{k-1}\mathcal{V}_j} \mathbf{1}_{\{t_k>0\}}\left(\int_{\R^3}\left|\frac{|v_k|}{w(v_k)}\mu^{1/2}(v_k) \right|^{p'} dv_k\right)^{1/p'}\prod_{j=1}^{k-1}d\sigma_j\\
	&\le  \frac{C_{p,\beta}}{\tilde{w}(v)}\sup_{0\le s \le t \le T_0} \left\{e^{-\nu_0(t-s)}\left\|h^{(m)}(s) \mathbf{1}_{\{t_1>0\}}\right\|_{L^p_vL^\infty_{x}}\right\}\left(\int_{\prod_{j=1}^{k-1}\mathcal{V}_j}\mathbf{1}_{\{t_{k}>0\}}\prod_{j=1}^{k-1}d\sigma_j\right)\\
	&\le \frac{C_{p,\beta}}{\tilde{w}(v)} \left(\frac{1}{2}\right)^{\tilde{C}_2T_0^{\frac{5}{4}}}\sup_{0\le s \le t \le T_0} \left\{e^{-\nu_0(t-s)}\left\|h^{(m)}(s) \mathbf{1}_{\{t_1>0\}}\right\|_{L^p_vL^\infty_{x}}\right\},
\end{align*}
which implies that
\begin{align*}
	\|J_2\|_{L^p_vL^\infty_x} \le C_{p,\beta}\left(\frac{1}{2}\right)^{\tilde{C}_2T_0^{\frac{5}{4}}}\sup_{0\le s \le t \le T_0} \left\{e^{-\nu_0(t-s)}\left\|h^{(m)}(s) \mathbf{1}_{\{t_1>0\}}\right\|_{L^p_vL^\infty_{x}}\right\}.
\end{align*}
Let us consider $I_1$ and $J_1$. By inserting $\int_{\mathcal{V}_k}d\sigma_k=1$ into $I_1$, we get
\begin{align*}
	I_1+J_1 &=  \frac{e^{-\nu_0(t-t_1)}}{\tilde{w}(v)}\sum_{l=1}^{k}\int_{\prod_{j=1}^{k}\mathcal{V}_j} \mathbf{1}_{\{t_{l+1}\le 0 < t_l\}} \left|h^{(m)}(0,x_l-v_lt_l,v_l)\right|d\Sigma_l(0).
\end{align*}
Now, we fix $l$ and consider the $l$-th term
\begin{align*}
	&\int_{\prod_{j=1}^{k}\mathcal{V}_j} \mathbf{1}_{\{t_{l+1}\le 0 < t_l\}} \left|h^{(m)}(0,x_l-v_lt_l,v_l)\right|d\Sigma_l(0)\\
	&\le e^{-\nu_0 t_1}\left\|h^{(m)}(0)\right\|_{L^p_vL^\infty_{x}}\int_{\prod_{j=1}^{l-1}\mathcal{V}_j} \left(\int_{\R^3}\left|\frac{|v_l|}{w(v_l)}\mu^{1/2}(v_l) \right|^{p'} dv_l\right)^{1/p'}\left\{\prod_{j=1}^{l-1} d\sigma_j\right\}\\
	&\le C_{p,\beta}e^{-\nu_0 t_1}\left\|h^{(m)}(0)\right\|_{L^p_vL^\infty_{x}},
\end{align*}
Summing $1\le l \le k$, it follows that
\begin{align*}
	I_1+J_1 \le C_{p,\beta}T_0^{\frac{5}{4}} \frac{e^{-\nu_0t}}{\tilde{w}(v)}\left\|h^{(m)}(0)\right\|_{L^p_vL^\infty_{x}},
\end{align*}
which implies that
\begin{align*}
	\|I_1+J_1\|_{L^p_vL^\infty_x} \le C_{p,\beta}T_0^{\frac{5}{4}} e^{-\nu_0t}\left\|h^{(m)}(0)\right\|_{L^p_vL^\infty_{x}}.
\end{align*}
Gathering $I_1$, $J_1$, and $J_2$, we deduce that for $0 \le t \le T_0$,
\begin{align*}
	e^{\nu_0 t}\left\|h^{(m)}(t)\mathbf{1}_{\{t_1>0\}}\right\|_{L^p_vL^\infty_x} &\le C_{p,\beta} \left(\frac{1}{2}\right)^{\tilde{C}_2T_0^{\frac{5}{4}}}\sup_{0\le s \le T_0} \left\{e^{\nu_0 s}\left\|h^{(m)}(s) \mathbf{1}_{\{t_1>0\}}\right\|_{L^p_vL^\infty_{x}}\right\}\\
	& \quad +C_{p,\beta} T_0^{\frac{5}{4}}\left\|h^{(m)}(0)\right\|_{L^p_vL^\infty_{x}}.
\end{align*}
Choosing sufficiently large $T_0 > 0$ such that $C_{p,\beta} \left(\frac{1}{2}\right)^{\tilde{C}_2T_0^{\frac{5}{4}}} \le \frac{1}{2}$,
\begin{align} \label{qq4}
	\sup_{0\le t \le T_0}\left\{e^{\nu_0 t}\left\|h^{(m)}(t)\mathbf{1}_{\{t_1>0\}}\right\|_{L^p_vL^\infty_{x}}\right\} \le C_{p,\beta} T_0^{\frac{5}{4}}\left\|h^{(m)}(0)\right\|_{L^p_vL^\infty_{x}} \le C_{p,\beta} T_0^{\frac{5}{4}}\left\|h_0\right\|_{L^p_vL^\infty_x}.
\end{align}
From now on, we extend the exponential decay to all time $t>0$. Letting $t=T_0$ in \eqref{qq4} and choosing sufficiently large $T_0>0$, for all $\tilde{\nu}_0 < \nu_0$,
\begin{align} \label{qq7}
	\left\|h^{(m)}\left(T_0\right)\right\|_{L^p_vL^\infty_x}\le C_{p,\beta} T_0^{\frac{5}{4}}e^{-\nu_0T_0}\left\|h_0\right\|_{L^p_vL^\infty_x} \le e^{-\tilde{\nu}_0 T_0}\left\|h_0\right\|_{L^p_vL^\infty_x},
\end{align}
and applying repeatedly the process \eqref{qq7}, we can derive for $l\ge 1$
\begin{align*}
	\left\|h^{(m)}(lT_0)\right\|_{L^p_vL^\infty_x} \le e^{-\tilde{\nu}_0 T_0} \left\|h^{(m)}\left((l-1)T_0\right)\right\|_{L^p_vL^\infty_x} \le e^{-l\tilde{\nu}_0 T_0}\left\|h_0\right\|_{L^p_vL^\infty_x}.
\end{align*}
Thus, for $lT_0 \le t \le (l+1)T_0$ with $l\ge 1$, we deduce that
\begin{align*}
	\left\|h^{(m)}\left(t\right)\right\|_{L^p_vL^\infty_x}
	&\le C_{T_0}e^{-l\tilde{\nu}_0 T_0}\left\|h_0\right\|_{L^p_vL^\infty_x}\\ 
	&\le C_{T_0}e^{- \tilde{\nu}_0 t}e^{ -\tilde{\nu}_0 T_0 }\left\|h_0\right\|_{L^p_vL^\infty_x}\\ 
	&\le C_{T_0} e^{- \tilde{\nu}_0 t}\left\|h_0\right\|_{L^p_vL^\infty_x}.
\end{align*}
since $0 \le t-lT_0 \le T_0$. Hence $\left(h^{(m)}\right)$ is uniformly bounded, and the sequence has weak* limit in $L^p_vL^\infty_x$. Letting $m \rightarrow \infty$, we conclude the existence of a solution and the exponential decay for the solution.
\end{proof}

\bigskip

\subsection{$L^p_vL^\infty_x$ estimate for the equation with $R(\varphi)$}
 We first define the function $R(\varphi)$ by
 \begin{align*}
 	R(\varphi)(t,x,v) = \int_{\R^3} \int_{\S^2} B(v-u,\omega) \left[\mu(v)+\mu^{1/2}(u)\varphi(t,x,u)\right]d\omega du.
 \end{align*} 
In this subsection, we consider the following equation with $R(\varphi)$:
\begin{align} \label{WDBE}
	\partial_t h+v\cdot \nabla_xh +R(\varphi)h=0
\end{align}
with initial data $h_0$ and the diffuse reflection boundary condition \eqref{WPDRBC}, where $\varphi = \varphi(t,x,v)$ is a given function satisfying
\begin{align} \label{varphicond}
	\mu(v) + \mu^{1/2}(v)\varphi(t,x,v)\ge 0 ,\quad  \|\varphi(t)\|_{L^p_vL^\infty_x} < \infty.
\end{align} 
We will study the solution operator for the equation \eqref{WDBE} with the diffuse reflection boundary condition and the related estimates.\\
\indent The following lemma present the mild formulation for the inhomogeneous form of the equation \eqref{WDBE} with the diffuse reflection boundary condition \eqref{WPDRBC}. Recall the definition of the iterated integral \eqref{iteratedintegral}. Before introducing the following lemma, we denote
\begin{align} \label{defIfnt}
	I^\varphi(t,s): = \exp \left\{ -\int_s^t R(\varphi)(\tau,X_{\mathbf{cl}}(\tau),V_{\mathbf{cl}}(\tau))d\tau \right\}.
\end{align}
\begin{Lem} \label{RepresentationforDiffuse} 
	Let $\varphi(t,x,v)$ satisfy the condition \eqref{varphicond}. Assume that $h,q \in L^p_vL^\infty_x$ satisfy $\partial_th+v\cdot \nabla_xh +R(\varphi)h=q(t,x,v)$ with the diffuse reflection boundary condition \eqref{WPDRBC}. Then for any $0\le s \le t$, for almost every $x,v$, if $t_1(t,x,v) \le s$,
	\begin{align*}
		h(t,x,v) = I^\varphi(t,s)h(s,x-v(t-s),v)+\int_s^tI^\varphi(t,\tau)q(\tau, x-v(t-\tau),v)d\tau;
	\end{align*}
	If $t_1(t,x,v)>s$, then for $k\ge 2$,
	\begin{align*}
		h(t,x,v) &= \int_{t_1}^tI^\varphi(t,\tau)q(\tau, x-v(t-\tau),v)d\tau\\
		&\quad +\frac{I^\varphi(t,t_1)}{\tilde{w}(v)}\sum_{l=1}^{k-1} \int_{\prod_{j=1}^{k-1}\mathcal{V}_j}\mathbf{1}_{\{t_{l+1}\le s <t_l\}} h(s,x_l-v_l(t_l-s),v_l)d\Sigma_l^\varphi(s)\\
		&\quad +\frac{I^\varphi(t,t_1)}{\tilde{w}(v)} \sum_{l=1}^{k-1} \int_{s}^{t_l}\int_{\prod_{j=1}^{k-1}\mathcal{V}_j} \mathbf{1}_{\{t_{l+1}\le s <t_l\}}q(\tau, x_l-v_l(t_l-\tau),v_l)d\Sigma_l^\varphi(\tau)d\tau\\
		&\quad +\frac{I^\varphi(t,t_1)}{\tilde{w}(v)} \sum_{l=1}^{k-1} \int_{t_{l+1}}^{t_l}\int_{\prod_{j=1}^{k-1}\mathcal{V}_j} \mathbf{1}_{\{t_{l+1}>s\}}q(\tau, x_l-v_l(t_l-\tau),v_l) d\Sigma_l^\varphi(\tau)d\tau\\
		&\quad +\frac{I^\varphi(t,t_1)}{\tilde{w}(v)} \int_{\prod_{j=1}^{k-1}\mathcal{V}_j} \mathbf{1}_{\{t_k>s\}}h(t_k,x_k,v_{k-1})d\Sigma_{k-1}^\varphi(t_k),
	\end{align*}
	where 
	\begin{align*}
		d\Sigma_{l}^\varphi(s) = \left\{\prod_{j=l+1}^{k-1}d\sigma_j \right\}\left\{I^\varphi(t_l,s) \tilde{w}(v_l)d\sigma_l \right\}\prod_{j=1}^{l-1} \left\{ I^\varphi(t_j,t_{j+1})d\sigma_j \right\}.
	\end{align*}
\end{Lem}
\begin{proof}
	It is a variant of Lemma \ref{RepresentationforDiffuse1}. Thus we omit the proof of this lemma.
\end{proof}

\bigskip
From now on, our goal in this subsection is to prove the existence and uniqueness for the equation \eqref{WDBE} with the boundary condition \eqref{WPDRBC}. The lemmas presented below are modifications of those from the previous subsection only in that $\nu(v)$ is replaced by $R(\varphi)$. Thus we omit the proof of those lemmas. (See Lemma \ref{inflowref1}, Lemma \ref{Existencesystem1}, and Lemma \ref{DampLpDecayorgin}.)
\begin{Lem}\label{inflowref} 
	Let $p \ge 1$ and $h_0 \in L^p_vL^\infty_x$, $g(t,x,v) \in L^\infty_xL^p_v(\gamma_-)$. Assume that $\varphi$ satisfies the condition \eqref{varphicond}. Let $h(t,x,v)$ be the solution to the following equation: 
	\begin{align*}
		\partial_t h+v\cdot \nabla_xh +R(\varphi)h=0, \quad h(0,x,v)=h_0(x,v), \quad h(t,x,v)|_{\gamma_-}=g(t,x,v).
	\end{align*}
	Then for $(t,x,v) \in (0,\infty) \times \bar{\Omega} \times \R^3$ with $(x,v) \notin \gamma_0 \cup \gamma_-$,
	\begin{align*}
		h(t,x,v) &= \mathbf{1}_{t-t_\mathbf{b}<0}\exp \left\{ -\int_0^t R(\varphi)(\tau,x-v(t-\tau),v)d\tau \right\} h_0(x-vt,v)\\
		& \quad +\mathbf{1}_{t-t_\mathbf{b}\ge0} \exp \left\{ -\int_{t-t_\mathbf{b}}^t R(\varphi)(\tau,x-v(t-\tau),v)d\tau \right\} g(t-t_\mathbf{b}, x_\mathbf{b},v),
	\end{align*}
	where $t_\mathbf{b}$ and $x_\mathbf{b}$ is defined in Definition \ref{Backwardexit}.
\end{Lem}
\begin{proof}
	It is similar to Lemma \ref{inflowref1}.
\end{proof}

\bigskip

\begin{Lem} \label{Existencesystem}
Let $p\ge 1$ and $h_0 \in L^p_vL^\infty_x$. Assume that $\varphi$ satisfies the condition \eqref{varphicond}. Then for any $\epsilon>0$, there exists $h(t) \in L^p_vL^\infty_x$ and $h_\gamma \in L^\infty_xL^p_v(\gamma)$ solving
\begin{align*}
	&\{ \partial_t +v\cdot \nabla_x+R(\varphi) \}h =0, \quad h(0)=h_0,\\
	& \quad h_{\gamma_-}(x,v)=(1-\epsilon)(c_\mu \mu(v))^{1/p}\int_{n(x)\cdot v'>0} h_{\gamma_+}(x,v')\frac{1}{(c_\mu \mu(v'))^{1/p}} d\sigma(x).
\end{align*}
\end{Lem}
\begin{proof}
	It is similar to the proof of Lemma \ref{Existencesystem1}.
\end{proof}
\bigskip
We denote by $S_{G_{\varphi}}(t,0)h_0$ the solution operator to the equation \eqref{WDBE} with the initial data $h_0$ and the diffuse reflection boundary condition \eqref{WPDRBC}. We now present the existence and uniqueness for the equation \eqref{WDBE} in $L^p_vL^\infty_x$ space.

\bigskip
\begin{Lem} \label{Dampexistence}
	Let $p>2$, $\beta > \frac{2p-3}{p}$ and $h_0 \in L^p_vL_x^\infty$. Assume that $\varphi$ satisfies the condition \eqref{varphicond} and $\mathcal{T}$ is sufficiently large. There exists a unique solution $h(t) = S_{G_{\varphi}}(t,0)h_0  \in L^p_vL^\infty_x$ to  \eqref{WDBE} with initial data $h_0$ and the diffuse reflection boundary condition \eqref{WPDRBC}. Moreover, there is a constant $C_1>0$, depending on $p$, $\beta$, $\gamma$, but independent of $\mathcal{T}$, such that
	\begin{align*}
		\|S_{G_{\varphi}}(t,0)h_0\|_{L^p_vL_x^\infty} \le C_1 \mathcal{T}^{\frac{5}{4}}\|h_0\|_{L^p_vL_x^\infty}
	\end{align*}
	for all $0 \le t \le \mathcal{T}$.
\end{Lem}
\begin{proof}
	It is similar to the proof of Lemma \ref{DampLpDecayorgin}
\end{proof}

\bigskip

\begin{Lem} \label{DampLpdecay}
	Let $p>2$, $\beta > \frac{2p-3}{p}$ and $h_0 \in L^p_vL_x^\infty$. Assume that $\varphi$ satisfies the condition \eqref{varphicond} and $\mathcal{T}$ is sufficiently large. Suppose that 
	\begin{align} \label{DampLpdecay1}
		R(\varphi)(t,x,v) \ge \frac{1}{2}\nu(v)
	\end{align}
	for all $(t,x,v) \in [0,\infty) \times \Omega \times \R^3$.
	Then it holds that
	\begin{align} \label{DampLpdecay2}
		\|S_{G_{\varphi}}(t,0)h_0\|_{L^p_vL_x^\infty} \le C_1 \mathcal{T}^{\frac{5}{4}} e^{-\frac{\nu_0}{2}t}\|h_0\|_{L^p_vL_x^\infty}
	\end{align}
	for $0 \le t \le \mathcal{T}$.
	Furthermore, there is a constant $C_{p,\gamma,\mathcal{T}} \ge 1$ such that
	\begin{align} \label{DampLpdecay3}
		\|S_{G_{\varphi}}(t,0)h_0\|_{L^p_vL_x^\infty} \le C_{p,\gamma,\mathcal{T}} e^{-\frac{\nu_0}{4}t}\|h_0\|_{L^p_vL_x^\infty}
	\end{align}
	for all $t \ge 0$.
\end{Lem}
\begin{proof}
	Recall the definition \eqref{defIfnt}. By the assumption \eqref{DampLpdecay1}, we have
\begin{align*}
	I^{\varphi}(t,0) \le e^{-\frac{\nu_0}{2}t}. 
\end{align*}
From Lemma \ref{RepresentationforDiffuse}, the remainder of the proof is similar to that of Lemma \ref{DampLpDecayorgin}.
\end{proof}

\bigskip
\subsection{Local-in-time existence}

Our goal in this subsection is to prove the local-in-time existence in $L^p_vL^\infty_x$ to the Boltzmann equation \eqref{Boltzmanneq} with the diffuse reflection boundary condition \eqref{Diffuseboundarycond}. 
\begin{Lem} \label{Localexistence}
	Let $p$ satisfy the condition \eqref{pcondition}, $\beta$ satisfy the condition \eqref{betacondition}, and $\mathcal{T}$ be sufficiently large.  Assume that $F_0(x,v) = \mu(v) + \mu^{1/2}(v) f_0(x,v) \ge 0$ and $ \|wf_0\|_{L^p_vL^\infty_x}<\infty$. Then there exists a positive time $\hat{t}_0:= \left(
	\hat{C}_\mathcal{T}\left[1+\|wf_0\|_{L^p_vL^\infty_x}\right]\right)^{-1}$ such that there is a unique solution $F(t,x,v) =\mu(v) + \mu^{1/2}(v) f(t,x,v) \ge 0 $ on $[0,\hat{t}_0]$ to the Boltzmann equation \eqref{Boltzmanneq} with initial data $F(0,x,v)=F_0(x,v)$ and the diffuse reflection boundary condition \eqref{Diffuseboundarycond} satisfying
	\begin{align*}
		\sup_{0\le t \le \hat{t}_0}\|wf(t)\|_{L^p_vL^\infty_x} \le 2C_1 \mathcal{T}^{\frac{5}{4}}\|wf_0\|_{L^p_vL^\infty_x},
	\end{align*}
	where $\hat{C}_\mathcal{T}$ is a positive constant depending only on $\mathcal{T}$, and $C_1$ is a constant defined in Lemma \ref{Dampexistence}
\end{Lem}
\begin{proof}
	To consider the local existence of solutions for the Boltzmann equation \eqref{Boltzmanneq}, we start from the iteration system that for $m=0,1,2,\cdots$,
	\begin{align} \label{iterativesys1}
		\partial_t F^{(m+1)} + v \cdot \nabla_x F^{(m+1)} + F^{(m+1)} \int_{\R^3} \int_{\S^2} B(v-u, \omega) F^{(m)}(t,x,u)d\omega du= Q^+(F^{(m)},F^{(m)})
	\end{align}
	with the boundary condition
	\begin{align*}
		F^{(m+1)}(t,x,v)|_{\gamma_-} = c_\mu \mu(v) \int_{n(x) \cdot v' >0} F^{(m+1)}(t,x,v') \{n(x) \cdot v'\}dv'
	\end{align*}
	and with initial data $F^{(m+1)}(t,x,v) |_{t=0} =F_0(x,v) \ge 0$, where we set $F^{(0)}(t,x,v) = \mu(v)$. We denote
	\begin{align*}
		f^{(m+1)}(t,x,v)=\frac{F^{(m+1)}(t,x,v)-\mu(v)}{\mu^{1/2}(v)},
	\end{align*}
	and 
	\begin{align*}
		h^{(m+1)}(t,x,v) = w(v) f^{(m+1)}(t,x,v).
	\end{align*}
	Then the iterative system \eqref{iterativesys1} can be written as, for $m=0,1,2,\cdots$,
	\begin{align} \label{iterativesys2}
		\partial_t h^{(m+1)} + v \cdot \nabla_x h^{(m+1)} + R(f^{(m)})h^{(m+1)} =K_wh^{(m)}+  w\Gamma^+(f^{(m)},f^{(m)})	
	\end{align}
	with the diffuse reflection boundary condition
	\begin{align} \label{iterativebdry1}
		h^{(m+1)}|_{\gamma_-} = \frac{1}{\tilde{w}(v)} \int_{\mathcal{V}(x) }h^{(m+1)}(t,x,v')\tilde{w}(v')d\sigma(x)
	\end{align} 
	and with initial data $h^{(m+1)}(t,x,v) |_{t=0} =h_0(x,v)$ and $h^{(0)}(t,x,v) \equiv 0$.\\
	\newline
	$\mathbf{(Uniform\ boundedness\ and\ Positivity)}$\\
	We shall use the induction argument on $m=0,1,2,\cdots$. Now, we show the uniform bound and positivity for solutions to the system \eqref{iterativesys1}:
	\begin{align} \label{uniformbound}
		\|h^{(m)}(t)\|_{L^p_vL^\infty_x} \le 2C_1 \mathcal{T}^{\frac{5}{4}} \|h_0\|_{L^p_vL^\infty_x}, \quad F^{(m)}(t,x,v) \ge 0 
	\end{align} 
	for all $m$ and $0 \le t \le \hat{t}_1$, where $\hat{t}_1$ is determined later. \\
	We denote the function $I^{(m)}(t,0)$ by
	\begin{align*}
		I^{(m)}(t,s) = \exp\left\{-\int_s^t R(f^{(m)})(\tau,X_{\mathbf{cl}}(\tau),V_{\mathbf{cl}}(\tau))d\tau\right\}. 
	\end{align*}
	for $m=0,1,2,\cdots$.
	First of all, for $m=0$, by Lemma \ref{Dampexistence}, there exists a solution operator $S_{G^{(1)}}(t,0)$ to the system \eqref{iterativesys2} and \eqref{iterativebdry1} such that
	\begin{align} \label{local1}
		\|h^{(1)}(t)\|_{L^p_vL^\infty_x} \le \|S_{G^{(1)}}(t,0)h_0\|_{L^p_vL^\infty_x} \le C_1 \mathcal{T}^{\frac{5}{4}}\|h_0\|_{L^p_vL^\infty_x}
	\end{align}
	for all $0 \le t \le \mathcal{T}$. By a variant of Lemma \ref{RepresentationforDiffuse},  we also have
	\begin{align*} 
		F^{(1)}(t,x,v) &= \mathbf{1}_{\{t_1 \le 0\}} I^{(0)}(t,0) F_0(x-vt,v) \nonumber\\
		& \quad + \mathbf{1}_{\{t_1 > 0\}} c_\mu \mu(v) I^{(0)}(t,t_1) \sum_{l=1}^{k-1} \int_{\prod_{j=1}^{k-1}\mathcal{V}_j}\mathbf{1}_{\{t_{l+1}\le 0 <t_l\}} F_0(x_l-v_lt_l,v_l)d\Sigma_l^{(1)}(0) \nonumber\\
		& \quad +\mathbf{1}_{\{t_1 > 0\}} c_\mu \mu(v) I^{(0)}(t,t_1)  \int_{\prod_{j=1}^{k-1}\mathcal{V}_j}\mathbf{1}_{\{t_k >0\}} F^{(1)}(t_k, x_k,v_{k-1})d\Sigma_{k-1}^{(1)}(t_k),
	\end{align*}
	where we denote
	\begin{align*} 
		d\Sigma_l^{(m)}(s):= \left\{\prod_{j=l+1}^{k-1} d\sigma_j \right\} \left\{ I^{(m)}(t_l,s)\{n(x_l) \cdot v_l \}dv_l \right\} \prod_{j=1}^{l-1} \left\{ I^{(m)}(t_j,t_{j+1})d\sigma_j \right\}.
	\end{align*}
	Note that from \eqref{local1},
	\begin{align*}
		\int_{\mathcal{V}_{k-1}}\left|F^{(1)}(t_k,x_k,v_{k-1})\right|\{n(x_k)\cdot v_{k-1}\}dv_{k-1} &\le  \int_{\mathcal{V}_{k-1}}\left(\mu(v) + \mu^{1/2}(v) |f^{(1)}(t_k,x_k,v_{k-1})|\right)|v_{k-1}|dv_{k-1}\\ &\le C_p\left(1+\sup_{0\le t \le  \mathcal{T}} \|h^{(1)}(t)\|_{L^p_vL^\infty_x} \right)\\
		& \le C_p \mathcal{T}^{\frac{5}{4}} \left(1+\|h_0\|_{L^p_vL^\infty_x}\right).
	\end{align*}
	This implies that when $k$ is sufficiently large,
	\begin{align*}
		F^{(1)}(t,x,v) &\ge -C\int_{\prod_{j=1}^{k-1}\mathcal{V}_j}\mathbf{1}_{\{t_{k}>0\}} \left|F^{(m+1)}(t_k,x_k, v_{k-1})\right|\{n(x_k) \cdot v_{k-1}\}dv_{k-1}\prod_{j=1}^{k-2}d\sigma_j	\\
		&\ge -C_p \mathcal{T}^{\frac{5}{4}}\left(1+\|h_0\|_{L^p_vL^\infty_x}\right)\int_{\prod_{j=1}^{k-2}\mathcal{V}_j}\mathbf{1}_{\{t_{k-1}>0\}} \prod_{j=1}^{k-2}d\sigma_j \\
		&\ge-C_p \mathcal{T}^{\frac{5}{4}}\left(1+\|h_0\|_{L^p_vL^\infty_x}\right) \epsilon,
	\end{align*}
	where we have used the estimate \eqref{Lsmallineq}. Since $\epsilon >0$ is arbitrary, $F^{(1)}(t,x,v) \ge 0$ over $[0,\mathcal{T}]$.\\
	Suppose that \eqref{uniformbound} holds for $m=0,1,2, \cdots, n-1$. We consider \eqref{uniformbound} for $m =n$. By Lemma \ref{Dampexistence}, for $n$, there exists the solution operator $S_{G^{(n+1)}}(t,s)$ to the system \eqref{iterativesys1} and \eqref{iterativebdry1}. Using the Duhamel principle, the system \eqref{iterativesys1} implies that
	\begin{align*} 
		h^{(n+1)}(t,x,v) &= S_{G^{(n+1)}}(t,0)h_0 + \int_0^t S_{G^{(n+1)}}(t,s)\left[K_wh^{(n)}(s)+w\Gamma^+ (f^{(n)},f^{(n)} )(s)\right]ds.
	\end{align*}
	From \eqref{uniformbound}, Lemma \ref{KESTIMATE}, Corollary \ref{LpGamma+est}, and Lemma \ref{Dampexistence}, we obtain
	\begin{align*}
		\|h^{(n+1)}(t)\|_{L^p_vL^\infty_x} &\le C_1 \mathcal{T}^{\frac{5}{4}}\|h_0\|_{L^p_vL^\infty_x}+C_1 \mathcal{T}^{\frac{5}{4}}\int_0^t \left[\|K_wh^{(n)}(s)\|_{L^p_vL^\infty_x} + \|w\Gamma^+ (f^{(n)},f^{(n)} )(s)\|_{L^p_vL^\infty_x}\right] ds\\
		& \le C_1 \mathcal{T}^{\frac{5}{4}}\|h_0\|_{L^p_vL^\infty_x} + C_pC_1 \mathcal{T}^{\frac{5}{4}}\int_0^t \left[\|h^{(n)}(s)\|_{L^p_vL^\infty_x}+\|h^{(n)}(s)\|_{L^p_vL^\infty_x}^2\right]ds\\
		& \le C_1 \mathcal{T}^{\frac{5}{4}}\|h_0\|_{L^p_vL^\infty_x}\left(1+C_{\mathcal{T},1}t\left(1+\|h_0\|_{L^p_vL^\infty_x}\right)  \right)
	\end{align*}
	for some constant $C_{\mathcal{T},1} \ge 1$ depending only on $\mathcal{T}$. Taking $\hat{t}_1:= \left(C_{\mathcal{T},1}\left(1+\|h_0\|_{L^p_vL^\infty_x}\right)\right)^{-1}$, it holds that
	\begin{align*}
		\|h^{(n+1)}(t)\|_{L^p_vL^\infty_x} \le 2C_1\mathcal{T}^{\frac{5}{4}}\|h_0\|_{L^p_vL^\infty_x}
	\end{align*}
	for $0\le t \le \hat{t}_1$. We conclude the uniform bound for solutions to the system \eqref{iterativesys2}.\\
	We now show the positivity for $F^{(n+1)}$. By a variant of Lemma \ref{RepresentationforDiffuse}, we have
	\begin{align*}
		F^{(n+1)}(t,x,v) &= \mathbf{1}_{\{t_1 \le 0\}} \left(I^{(n)}(t,0)F_0(x-vt,v)+\int_0^t I^{(n)}(t,s) Q^+(F^{(n)},F^{(n)})(s,x-v(t-s),v) ds\right)\\
		& \quad +\mathbf{1}_{\{t_1 > 0\}} c_\mu \mu(v) I^{(n)}(t,t_1)\Bigg\{\sum_{l=1}^{k-1} \int_{\prod_{j=1}^{k-1}\mathcal{V}_j}\mathbf{1}_{\{t_{l+1}\le 0 <t_l\}} F_0(x_l-v_lt_l,v_l)d\Sigma_l^{(n)}(0)\\
		&\qquad +\sum_{l=1}^{k-1} \int_{0}^{t_l}\int_{\prod_{j=1}^{k-1}\mathcal{V}_j} \mathbf{1}_{\{t_{l+1}\le 0 <t_l\}}Q^+(F^{(n)},F^{(n)})(\tau, x_l-v_l(t_l-\tau),v_l)d\Sigma_l^{(n)}(\tau)d\tau\\
		&\qquad +\sum_{l=1}^{k-1} \int_{t_{l+1}}^{t_l}\int_{\prod_{j=1}^{k-1}\mathcal{V}_j} \mathbf{1}_{\{t_{l+1}>0\}}Q^+(F^{(n)},F^{(n)})(\tau, x_l-v_l(t_l-\tau),v_l) d\Sigma_l^{(n)}(\tau)d\tau\\
		&\qquad +\int_{\prod_{j=1}^{k-1}\mathcal{V}_j} \mathbf{1}_{\{t_k>0\}}F^{(n+1)}(t_k,x_k,v_{k-1})d\Sigma_{k-1}^{(n)}(t_k)\Bigg\}.
	\end{align*}
	Here, by the induction assumption, it holds that
	\begin{align*}
		\int_{\mathcal{V}_{k-1}}\left|F^{(n+1)}(t_k,x_k,v_{k-1})\right|\{n(x_k)\cdot v_{k-1}\}dv_{k-1} 
		& \le C_p \mathcal{T}^{\frac{5}{4}} \left(1+\|h_0\|_{L^p_vL^\infty_x}\right),
	\end{align*}
	which implies that for $0\le t \le \hat{t}_1$,
	\begin{align*}
		F^{(n+1)}(t,x,v) \ge -C_p \mathcal{T}^{\frac{5}{4}}\left(1+\|h_0\|_{L^p_vL^\infty_x}\right) \epsilon,
	\end{align*}
	where we have used \eqref{Lsmallineq}. Since $\epsilon >0$ is arbitrary, $F^{(m+1)}(t,x,v) \ge 0$ over $[0,\hat{t}_1]$. By the induction, \eqref{uniformbound} holds for all $m=0,1,2, \cdots$.\\
	\newline
	$\mathbf{(Existence)}$\\
	We denote $\hat{w}(v) = \frac{w(v)}{\sqrt{1+|v|^2}}$ and
	\begin{align*}
		\hat{h}^{(m)}:= \hat{w}(v)f^{(m)}(t,x,v) =(1+|v|^2)^{-1/2}h^{(m)}(t,x,v).
	\end{align*}
	We set
	\begin{align*}
		g^{(m)}(t,x,v) := K_{\hat{w}}\hat{h}^{(m)}(t,x,v)+\hat{w}\Gamma^+ (f^{(m)},f^{(m)} )(t,x,v)
	\end{align*}
	for $m=0,1,2,\cdots$.\\
	\newline
	$\mathbf{Case}$ $t_1(t,x,v) \le 0 : $
	By Duhamel principle, for each $m$, we have
	\begin{align*}
		\hat{h}^{(m+1)}(t,x,v) &= I^{(m)}(t,0) \hat{h}_0(x-vt,v) + \int_0^t I^{(m)}(t,s)g^{(m)}(s,x-v(t-s),v)ds.	
	\end{align*}
	By subtracting $\hat{h}^{(m+2)} - \hat{h}^{(m+1)}$, we obtain
	\begin{equation} \label{local8}
	\begin{aligned}
		&\left|\left(\hat{h}^{(m+2)} - \hat{h}^{(m+1)} \right)(t,x,v)\right|\\
		& \le |I^{(m+1)}(t,0)- I^{(m)}(t,0) ||\hat{h}_0(x-vt,v)|\\
		& \quad +\int_0^t |I^{(m+1)}(t,s)- I^{(m)}(t,s) | \left| g^{(m+1)}(s,x-v(t-s),v)\right|ds\\
		& \quad + \int_0^t I^{(m)}(t,s)\Bigl(|K_{\hat{w}}(\hat{h}^{(m+1)}-\hat{h}^{(m)})(s,x-v(t-s),v)|\\
		& \qquad \qquad +\left|\left(\hat{w}\Gamma^+ (f^{(m+1)},f^{(m+1)} )-\hat{w}\Gamma^+ (f^{(m)},f^{(m)}) \right)(s,x-v(t-s),v)\right|\Bigr)ds\\
		& =:J_1^{(m)} +J_2^{(m)}+J_3^{(m)}.
	\end{aligned}
\end{equation}
	First, let us estimate $J_1^{(m)}$. By a direct calculation, we have
	\begin{equation} \label{local2}
	\begin{aligned}
		|I^{(m+1)}(t,s)- I^{(m)}(t,s) | &\le \left|\int_s^t \left[R(f^{(m+1)})-R(f^{(m)})\right](\tau,x-v(t-\tau),v) d\tau \right| \\
		& \le C \nu(v) \int_s^t  \int_{\R^3} \nu(u)\mu^{1/2}(u) \left|\left(f^{(m+1)}-f^{(m)}\right)(\tau,x-v(t-\tau),u)
	\right|dud\tau \\
		&\le C_p\nu(v) \int_s^t\left\|\left(f^{(m+1)}-f^{(m)}\right)(\tau)\right\|_{L^p_vL^\infty_x}d\tau.
	\end{aligned}
\end{equation}
	Then from \eqref{local2},
	\begin{align} \label{local3}
		\|J_1^{(m)}\|_{L^p_vL^\infty_x} &\le C_p \|h_0\|_{L^p_vL^\infty_x}\int_0^t\left\|\left(f^{(m+1)}-f^{(m)}\right)(\tau)\right\|_{L^p_vL^\infty_x}d\tau\nonumber\\
		& \le C_p t\|h_0\|_{L^p_vL^\infty_x}\sup_{0\le \tau \le t}\left\|\left(f^{(m+1)}-f^{(m)}\right)(\tau)\right\|_{L^p_vL^\infty_x}.
	\end{align}
	Next, let us estimate $J_2^{(m)}$. By Lemma \ref{KESTIMATE}, Corollary \ref{LpGamma+est} and \eqref{local2}, we obtain 
	\begin{align} \label{local4}
		\|J_2^{(m)}\|_{L^p_vL^\infty_x} &\le C_p \int_0^t \left(\left\|K_{w}h^{(m+1)}(s) \right\|_{L^p_vL^\infty_x}+\left\| w\Gamma^+ (f^{(m+1)},f^{(m+1)} )(s)\right\|_{L^p_vL^\infty_x}\right)\nonumber\\
		& \quad \times \int_s^t\left\|\left(f^{(m+1)}-f^{(m)}\right)(\tau)\right\|_{L^p_vL^\infty_x}d\tau ds \nonumber\\
		& \le C_p \int_0^t \left( \|h^{(m+1)}(s)\|_{L^p_vL^\infty_x}+ \|h^{(m+1)}(s)\|_{L^p_vL^\infty_x}^2 \right) \nonumber\\
		& \quad \times \int_s^t\left\|\left(f^{(m+1)}-f^{(m)}\right)(\tau)\right\|_{L^p_vL^\infty_x}d\tau ds\nonumber\\
		& \le C_p t^2 \sup_{0\le s \le t}\left[ \|h^{(m+1)}(s)\|_{L^p_vL^\infty_x}+ \|h^{(m+1)}(s)\|_{L^p_vL^\infty_x}^2 \right] \nonumber\\
		& \quad \times \sup_{0\le \tau \le t}\left\|\left(f^{(m+1)}-f^{(m)}\right)(\tau)\right\|_{L^p_vL^\infty_x}.
	\end{align}
	Combining \eqref{local3} and \eqref{local4}, and using \eqref{uniformbound}, we get for $t \le \hat{t}_1$,
	\begin{align} \label{local5}
		\|J_1^{(m)}+J_2^{(m)}\|_{L^p_vL^\infty_x} &\le C_pt\left(\|h_0\|_{L^p_vL^\infty_x}+t \sup_{0\le s \le t}\left[ \|h^{(m+1)}(s)\|_{L^p_vL^\infty_x}+ \|h^{(m+1)}(s)\|_{L^p_vL^\infty_x}^2 \right] \right)\nonumber\\
		& \quad \times \sup_{0\le \tau \le t}\left\|\left(f^{(m+1)}-f^{(m)}\right)(\tau)\right\|_{L^p_vL^\infty_x}\nonumber\\
		& \le C_{p,\mathcal{T}} t\|h_0\|_{L^p_vL^\infty_x}\sup_{0\le \tau \le t}\left\|\left(f^{(m+1)}-f^{(m)}\right)(\tau)\right\|_{L^p_vL^\infty_x}.
	\end{align}
	Note that
	\begin{align*}
		&\hat{w}\Gamma^+ (f^{(m+1)},f^{(m+1)} )-\hat{w}\Gamma^+ (f^{(m)},f^{(m)} ) \\
		&= \hat{w}\Gamma^+ (f^{(m+1)}-f^{(m)},f^{(m+1)} )+\hat{w}\Gamma^+ (f^{(m)},f^{(m+1)}-f^{(m)} ).
	\end{align*}
	For $J_3^{(m)}$, we use Lemma \ref{KESTIMATE}, Corollary \ref{LpGamma+est} and the positivity for $F^{(m)}$ to obtain
	\begin{align} \label{local6}
		\|J_3^{(m)}\|_{L^p_vL^\infty_x} \le C_p\int_0^t \left(1+\|\hat{h}^{(m)}(s)\|_{L^p_vL^\infty_x} + \|\hat{h}^{(m+1)}(s)\|_{L^p_vL^\infty_x}\right)\|(\hat{h}^{(m+1)}-\hat{h}^{(m)})(s)\|_{L^p_vL^\infty_x}ds.
	\end{align}
	For $t\le \hat{t}_1$, we apply the uniform boundedness \eqref{uniformbound} to \eqref{local6} :
	\begin{align} \label{local7}
		\|J_3^{(m)}\|_{L^p_vL^\infty_x} \le C_{p,\mathcal{T}} t (1+\|h_0\|_{L^p_vL^\infty_x})\sup_{0\le s\le t} \|(\hat{h}^{(m+1)}-\hat{h}^{(m)})(s)\|_{L^p_vL^\infty_x}.
	\end{align}
	Plugging \eqref{local5} and \eqref{local7} into \eqref{local8}, we derive
	\begin{align} \label{local9}
		\sup_{0\le s  \le t}\|(\hat{h}^{(m+2)}-\hat{h}^{(m+1)})(s)\|_{L^p_vL^\infty_x} \le C_{p,\mathcal{T}} t (1+\|h_0\|_{L^p_vL^\infty_x})\sup_{0\le s\le t} \|(\hat{h}^{(m+1)}-\hat{h}^{(m)})(s)\|_{L^p_vL^\infty_x}
	\end{align}
	for all $0\le t \le \hat{t}_1$.\\
	\newline
	$\mathbf{Case}$ $t_1(t,x,v) > 0 : $ By a variant of Lemma \ref{RepresentationforDiffuse}, it holds that
	\begin{align*}
		\hat{h}^{(m+1)}(t,x,v) &= \int_{t_1}^tI^{(m)}(t,s)g^{(m)}(\tau, x-v(t-s),v)ds\\
		&\quad +\frac{I^{(m)}(t,t_1)}{\tilde{w}(v)\sqrt{1+|v|^2}}\sum_{l=1}^{k-1} \int_{\prod_{j=1}^{k-1}\mathcal{V}_j}\mathbf{1}_{\{t_{l+1}\le 0 <t_l\}} \hat{h}_0(x_l-v_lt_l,v_l)d\hat{\Sigma}_l^{(m)}(0)\\
		&\quad +\frac{I^{(m)}(t,t_1)}{\tilde{w}(v)\sqrt{1+|v|^2}} \sum_{l=1}^{k-1} \int_{s}^{t_l}\int_{\prod_{j=1}^{k-1}\mathcal{V}_j} \mathbf{1}_{\{t_{l+1}\le 0 <t_l\}}g^{(m)}(\tau, x_l-v_l(t_l-\tau),v_l)d\hat{\Sigma}_l^{(m)}(\tau)d\tau\\
		&\quad +\frac{I^{(m)}(t,t_1)}{\tilde{w}(v)\sqrt{1+|v|^2}} \sum_{l=1}^{k-1} \int_{t_{l+1}}^{t_l}\int_{\prod_{j=1}^{k-1}\mathcal{V}_j} \mathbf{1}_{\{t_{l+1}>0\}}g^{(m)}(\tau, x_l-v_l(t_l-\tau),v_l) d\hat{\Sigma}_l^{(m)}(\tau)d\tau\\
		&\quad +\frac{I^{(m)}(t,t_1)}{\tilde{w}(v)\sqrt{1+|v|^2}} \int_{\prod_{j=1}^{k-1}\mathcal{V}_j} \mathbf{1}_{\{t_k>0\}}\hat{h}^{(m+1)}(t_k,x_k,v_{k-1})d\hat{\Sigma}_{k-1}^{(m)}(t_k)\\
		&=:J_4^{(m)}+J_5^{(m)}+J_6^{(m)}+J_7^{(m)}+J_8^{(m)},
	\end{align*}
	where 
	\begin{align*}
		d\hat{\Sigma}_{l}^{(m)}(s) = \left\{\prod_{j=l+1}^{k-1}d\sigma_j \right\}\left\{I^{(m)}(t_l,s) \tilde{w}(v_l)\sqrt{1+|v_l|^2}d\sigma_l \right\}\prod_{j=1}^{l-1} \left\{ I^{(m)}(t_j,t_{j+1})d\sigma_j \right\}.
	\end{align*}
	By subtracting $\hat{h}^{(m+2)} - \hat{h}^{(m+1)}$ and taking the norm in $L^p_vL^\infty_x$, we obtain
	\begin{align} \label{local17}
		\|(\hat{h}^{(m+2)} - \hat{h}^{(m+1)})(t)\|_{L^p_vL^\infty_x} \le \sum_{i=4}^8 \|(J_i^{(m+1)}-J_i^{(m)})(t)\|_{L^p_vL^\infty_x}.
	\end{align}
	For $\|J_4^{(m+1)}-J_4^{(m)}\|_{L^p_vL^\infty_x}$, by similar arguments in \eqref{local4} and \eqref{local6}, we have
	\begin{align} \label{local10}
		\|J_4^{(m+1)}-J_4^{(m)}\|_{L^p_vL^\infty_x} \le  C_{p,\mathcal{T}} t (1+\|h_0\|_{L^p_vL^\infty_x})\sup_{0\le s\le t} \|(\hat{h}^{(m+1)}-\hat{h}^{(m)})(s)\|_{L^p_vL^\infty_x}.
	\end{align}
	Before estimating $\|J_i^{(m+1)}-J_i^{(m)}\|_{L^p_vL^\infty_x}$ for $i = 5,6,7,8$, we need to deal with the following difference
	\begin{align*}
		&\Bigg| \sum_{l=1}^{k-1}\int_{\prod_{j=1}^{k-1}\mathcal{V}_j} \mathbf{1}_{\{t_{l+1}\le 0 <t_l\}}g^{(m+1)}(\tau, x_l-v_l(t_l-\tau),v_l)d\hat{\Sigma}_l^{(m+1)}(\tau)\\
		& \quad -\sum_{l=1}^{k-1}\int_{\prod_{j=1}^{k-1}\mathcal{V}_j} \mathbf{1}_{\{t_{l+1}\le 0 <t_l\}}g^{(m)}(\tau, x_l-v_l(t_l-\tau),v_l)d\hat{\Sigma}_l^{(m)}(\tau)\Bigg|\\
		&\le \Bigg| \sum_{l=1}^{k-1}\int_{\prod_{j=1}^{k-1}\mathcal{V}_j} \mathbf{1}_{\{t_{l+1}\le 0 <t_l\}}g^{(m+1)}(\tau, x_l-v_l(t_l-\tau),v_l)\left\{d\hat{\Sigma}_l^{(m+1)}(\tau)-d\hat{\Sigma}_l^{(m)}(\tau) \right\}\Bigg|\\
		& \quad + \left|  \sum_{l=1}^{k-1}\int_{\prod_{j=1}^{k-1}\mathcal{V}_j} \mathbf{1}_{\{t_{l+1}\le 0 <t_l\}}(g^{(m+1)}-g^{(m)})(\tau, x_l-v_l(t_l-\tau),v_l)d\hat{\Sigma}_l^{(m)}(\tau)\right|\\
		&=:L_1 + L_2. 
	\end{align*}
	Let us estimate the term $L_1$:
	\begin{align*}
		L_1 &\le \sum_{l=1}^{k-1}\int_{\prod_{j=1}^{l}\mathcal{V}_j} \left|g^{(m+1)}(\tau, x_l-v_l(t_l-\tau),v_l)\right| \tilde{w}(v_l)\sqrt{1+|v_l|^2}\\
		& \quad \times \Big|I^{(m+1)}(t_l,s)I^{(m+1)}(t_{l-1},t_l) \cdots I^{(m+1)}(t_1,t_2) - I^{(m)}(t_l,s)I^{(m)}(t_{l-1},t_l) \cdots I^{(m)}(t_1,t_2) \Big| d\sigma_l \cdots d\sigma_1.
	\end{align*}
	Note that
	\begin{align} \label{local11}
		&\left|I^{(m+1)}(t_l,s)I^{(m+1)}(t_{l-1},t_l) \cdots I^{(m+1)}(t_1,t_2) - I^{(m)}(t_l,s)I^{(m)}(t_{l-1},t_l) \cdots I^{(m)}(t_1,t_2) \right|\nonumber\\
		& \le \left|I^{(m+1)}(t_l,s)-I^{(m)}(t_l,s) \right|+\left|I^{(m+1)}(t_{l-1},t_l)-I^{(m)}(t_{l-1},t_l) \right|+\cdots+\left|I^{(m+1)}(t_1,t_2)-I^{(m)}(t_1,t_2) \right|\nonumber\\
		& \le C_p\nu(v_l) \int_{s}^{t_l}\left\|\left(f^{(m+1)}-f^{(m)}\right)(\tau)\right\|_{L^p_vL^\infty_x}d\tau+C_p\nu(v_{l-1}) \int^{t_{l-1}}_{t_l}\left\|\left(f^{(m+1)}-f^{(m)}\right)(\tau)\right\|_{L^p_vL^\infty_x}d\tau\nonumber\\
		& \quad +\cdots + C_p\nu(v_1) \int_{t_2}^{t_1}\left\|\left(f^{(m+1)}-f^{(m)}\right)(\tau)\right\|_{L^p_vL^\infty_x}d\tau\nonumber\\
		&\le C_p\left[\nu(v_l)+\nu(v_{l-1})+ \cdots +\nu(v_1)\right] \int_0^t \left\|\left(f^{(m+1)}-f^{(m)}\right)(\tau)\right\|_{L^p_vL^\infty_x}d\tau.
	\end{align}
	By \eqref{local11}, we obtain
	\begin{align} \label{local12}
		L_1 &\le C_p\sum_{l=1}^{k-1}\int_{\prod_{j=1}^{l}\mathcal{V}_j} \left|g^{(m+1)}(\tau, x_l-v_l(t_l-\tau),v_l)\right| \tilde{w}(v_l)\sqrt{1+|v_l|^2}\nonumber\\
		& \quad \times \left[\nu(v_l)+\nu(v_{l-1})+ \cdots +\nu(v_1)\right] \int_0^t \left\|\left(f^{(m+1)}-f^{(m)}\right)(\tau')\right\|_{L^p_vL^\infty_x}d\tau' d\sigma_l \cdots d\sigma_1\nonumber\\
		& \le C_p k\|g^{(m+1)}(\tau)\|_{L^p_vL^\infty_x}\int_0^t \left\|\left(f^{(m+1)}-f^{(m)}\right)(\tau')\right\|_{L^p_vL^\infty_x}d\tau'\nonumber\\
		& \le C_p \mathcal{T}^{\frac{5}{4}}\|g^{(m+1)}(\tau)\|_{L^p_vL^\infty_x}\int_0^t \left\|\left(f^{(m+1)}-f^{(m)}\right)(\tau')\right\|_{L^p_vL^\infty_x}d\tau',
	\end{align}
	where we have used the inequalities
	\begin{align*}
		&\int_{\mathcal{V}_l}\left|g^{(m+1)}(\tau, x_l-v_l(t_l-\tau),v_l)\right|\nu(v_l)\tilde{w}(v_l)\sqrt{1+|v_l|^2} d\sigma_l \le C_p\|g^{(m+1)}(\tau)\|_{L^p_vL^\infty_x},\\
		& \int_{\mathcal{V}_j} \nu(v_j)d\sigma_j \le C \quad  \text{for\ }  j=1,2,\cdots,l-1,
	\end{align*}
	and we have taken $k= \tilde{C}_1 \mathcal{T}^{\frac{5}{4}}$.\\
	Let us estimate the term $L_2$:
	\begin{align} \label{local13}
		L_2 &\le \sum_{l=1}^{k-1}\int_{\prod_{j=1}^{k-1}\mathcal{V}_j} \mathbf{1}_{\{t_{l+1}\le 0 <t_l\}}\left|(g^{(m+1)}-g^{(m)})(\tau, x_l-v_l(t_l-\tau),v_l)\right|d\hat{\Sigma}_l^{(m)}(\tau)\nonumber\\
		& \le C_p\|(g^{(m+1)}-g^{(m)})(\tau)\|_{L^p_vL^\infty_x}.
	\end{align}
	For $\|J_5^{(m+1)}-J_5^{(m)}\|_{L^p_vL^\infty_x}$, by similar arguments in \eqref{local12} and \eqref{local13}, we get
	\begin{align*}
		\left|J_5^{(m+1)}-J_5^{(m)}\right| &\le \frac{\left|I^{(m+1)}(t,t_1)-I^{(m)}(t,t_1)\right|}{\tilde{w}(v)\sqrt{1+|v|^2}}\sum_{l=1}^{k-1} \int_{\prod_{j=1}^{k-1}\mathcal{V}_j}\mathbf{1}_{\{t_{l+1}\le 0 <t_l\}} |\hat{h}_0(x_l-v_lt_l,v_l)|d\hat{\Sigma}_l^{(m+1)}(0)\\
		& \quad +\frac{\left|I^{(m)}(t,t_1)\right|}{\tilde{w}(v)\sqrt{1+|v|^2}}\sum_{l=1}^{k-1} \int_{\prod_{j=1}^{k-1}\mathcal{V}_j}\mathbf{1}_{\{t_{l+1}\le 0 <t_l\}} |\hat{h}_0(x_l-v_lt_l,v_l)|\left\{d\hat{\Sigma}_l^{(m+1)}(0)-d\hat{\Sigma}_l^{(m)}(0)	\right\}\\
		& \le \frac{C_p\nu(v)}{\tilde{w}(v)\sqrt{1+|v|^2}}\|\hat{h}_0\|_{L^p_vL^\infty_x}\int_0^t\left\|\left(f^{(m+1)}-f^{(m)}\right)(\tau)\right\|_{L^p_vL^\infty_x}d\tau\\
		& \quad + \frac{C_{p,\mathcal{T}}}{\tilde{w}(v)\sqrt{1+|v|^2}}\|\hat{h}_0\|_{L^p_vL^\infty_x}\int_0^t \left\|\left(f^{(m+1)}-f^{(m)}\right)(\tau')\right\|_{L^p_vL^\infty_x}d\tau'\\
		& \le t\frac{C_{p,\mathcal{T}}}{\tilde{w}(v)}\|h_0\|_{L^p_vL^\infty_x}\sup_{0\le s\le t} \|(\hat{h}^{(m+1)}-\hat{h}^{(m)})(s)\|_{L^p_vL^\infty_x},
	\end{align*}
	and it follows that
	\begin{align} \label{local14}
		\|J_5^{(m+1)}-J_5^{(m)}\|_{L^p_vL^\infty_x} \le tC_{p,\mathcal{T}}\|h_0\|_{L^p_vL^\infty_x}\sup_{0\le s\le t} \|(\hat{h}^{(m+1)}-\hat{h}^{(m)})(s)\|_{L^p_vL^\infty_x}.
	\end{align}
	For $\|J_6^{(m+1)}-J_6^{(m)}\|_{L^p_vL^\infty_x}$ and $\|J_7^{(m+1)}-J_7^{(m)}\|_{L^p_vL^\infty_x}$, by similar arguments in \eqref{local4}, \eqref{local6}, \eqref{local12}, \eqref{local13}, we have
	\begin{align} \label{local15}
		\|J_6^{(m+1)}-J_6^{(m)}\|_{L^p_vL^\infty_x}+\|J_7^{(m+1)}-J_7^{(m)}\|_{L^p_vL^\infty_x} \le tC_{p,\mathcal{T}} (1+\|h_0\|_{L^p_vL^\infty_x})\sup_{0\le s\le t} \|(\hat{h}^{(m+1)}-\hat{h}^{(m)})(s)\|_{L^p_vL^\infty_x}
	\end{align}
	for $t \le \hat{t}_1$.\\
	For $\|J_8^{(m+1)}-J_8^{(m)}\|_{L^p_vL^\infty_x}$, by Lemma \ref{Lsmall} and similar arguments in \eqref{local12}, \eqref{local13}, we can derive
	\begin{align*}
		\left|J_8^{(m+1)}-J_8^{(m)}\right| &\le \frac{I^{(m+1)}(t,t_1)}{\tilde{w}(v)\sqrt{1+|v|^2}} \int_{\prod_{j=1}^{k-1}\mathcal{V}_j} \mathbf{1}_{\{t_k>0\}}\left|(\hat{h}^{(m+2)}-\hat{h}^{(m+1)})(t_k,x_k,v_{k-1})\right|d\hat{\Sigma}_{k-1}^{(m+1)}(t_k)\\
		& \quad +\frac{\left|I^{(m+1)}(t,t_1)-I^{(m)}(t,t_1)\right|}{\tilde{w}(v)\sqrt{1+|v|^2}} \int_{\prod_{j=1}^{k-1}\mathcal{V}_j} \mathbf{1}_{\{t_k>0\}}\left|\hat{h}^{(m+1)}(t_k,x_k,v_{k-1})\right|d\hat{\Sigma}_{k-1}^{(m+1)}(t_k)\\
		& \quad + \frac{I^{(m)}(t,t_1)}{\tilde{w}(v)\sqrt{1+|v|^2}} \left|\int_{\prod_{j=1}^{k-1}\mathcal{V}_j} \mathbf{1}_{\{t_k>0\}}\hat{h}^{(m+1)}(t_k,x_k,v_{k-1})\left\{d\hat{\Sigma}_{k-1}^{(m+1)}(t_k)-d\hat{\Sigma}_{k-1}^{(m)}(t_k)\right\}\right|\\
		& \le \frac{C_p}{\tilde{w}(v)} \left(\frac{1}{2}\right)^{\tilde{C}_2 \mathcal{T}^{\frac{5}{4}}} \sup_{0\le s\le t} \|(\hat{h}^{(m+2)}-\hat{h}^{(m+1)})(s)\|_{L^p_vL^\infty_x}\\
		& \quad + t\frac{C_{p,\mathcal{T}}}{\tilde{w}(v)}\|h_0\|_{L^p_vL^\infty_x}\sup_{0\le s\le t} \|(\hat{h}^{(m+1)}-\hat{h}^{(m)})(s)\|_{L^p_vL^\infty_x},
	\end{align*}
	which implies that
	\begin{align} \label{local16}
		\|J_8^{(m+1)}-J_8^{(m)}\|_{L^p_vL^\infty_x} &\le C_p\left(\frac{1}{2}\right)^{\tilde{C}_2 \mathcal{T}^{\frac{5}{4}}} \sup_{0\le s\le t} \|(\hat{h}^{(m+2)}-\hat{h}^{(m+1)})(s)\|_{L^p_vL^\infty_x}\nonumber\\
		& \quad +tC_{p,\mathcal{T}}\|h_0\|_{L^p_vL^\infty_x}\sup_{0\le s\le t} \|(\hat{h}^{(m+1)}-\hat{h}^{(m)})(s)\|_{L^p_vL^\infty_x}
	\end{align}
	Plugging \eqref{local10}, \eqref{local14}, \eqref{local15}, \eqref{local16} into \eqref{local17} and then combining \eqref{local9}, we obtain
	\begin{align*}
		\sup_{0\le s  \le t}\|(\hat{h}^{(m+2)}-\hat{h}^{(m+1)})(s)\|_{L^p_vL^\infty_x} &\le tC_{p,\mathcal{T}}  (1+\|h_0\|_{L^p_vL^\infty_x})\sup_{0\le s\le t} \|(\hat{h}^{(m+1)}-\hat{h}^{(m)})(s)\|_{L^p_vL^\infty_x}\\
		& \quad + C_p\left(\frac{1}{2}\right)^{\tilde{C}_2 \mathcal{T}^{\frac{5}{4}}} \sup_{0\le s\le t} \|(\hat{h}^{(m+2)}-\hat{h}^{(m+1)})(s)\|_{L^p_vL^\infty_x}.
	\end{align*}
	Taking $\mathcal{T}>0$ sufficiently large such that $C_p\left(\frac{1}{2}\right)^{\tilde{C}_2 \mathcal{T}^{\frac{5}{4}}} \le \frac{1}{2}$, we deduce
	\begin{align} \label{local20}
		\sup_{0\le s  \le t}\|(\hat{h}^{(m+2)}-\hat{h}^{(m+1)})(s)\|_{L^p_vL^\infty_x}\le tC_{p,\mathcal{T},2}  (1+\|h_0\|_{L^p_vL^\infty_x})\sup_{0\le s\le t} \|(\hat{h}^{(m+1)}-\hat{h}^{(m)})(s)\|_{L^p_vL^\infty_x},
	\end{align}
	where $C_{p,\mathcal{T},2}$ is a positive constant depending on $\mathcal{T}$. We set $\hat{C}_\mathcal{T}:=\max\{C_{\mathcal{T},1},2C_{p,\mathcal{T},2}\}$ and choose
	\begin{align*}
		\hat{t}_0 = \frac{1}{\hat{C}_\mathcal{T}\left(1+\|h_0\|_{L^p_vL^\infty_x}\right)}.
	\end{align*}
	Then it holds that $\hat{t}_0 \le \hat{t}_1$, and it follows from \eqref{local20} that $0\le t \le \hat{t}_0$,
	\begin{align*}
		\sup_{0\le s  \le t}\|(\hat{h}^{(m+2)}-\hat{h}^{(m+1)})(s)\|_{L^p_vL^\infty_x} \le \frac{1}{2} \sup_{0\le s\le t} \|(\hat{h}^{(m+1)}-\hat{h}^{(m)})(s)\|_{L^p_vL^\infty_x}.
	\end{align*}
	Hence $\hat{h}^{(m)}$ is a Cauchy sequence, and $f^{(m)}$ is also a Cauchy sequence. Therefore, there exists a limit $f(t,x,v)$, which is a solution to the equation \eqref{FPBE} with the diffuse reflection boundary condition \eqref{PDRBC}.\\
	\newline	
	$\mathbf{(Uniqueness)}$\\
	Assume that there exists another solution $\tilde{f}(t,x,v)$ to the Boltzmann equation \eqref{FPBE} with the boundary condtion \eqref{PDRBC} and the same initial data as for $f$, and the bound $\sup_{0\le t \le \hat{t}_0}\|w\tilde{f}(t)\|_{L^p_vL^\infty_x}<\infty$.
	By the similar argument to the existence, we obtain
	\begin{align*}
		&\|\hat{w}(f-\tilde{f})(t)\|_{L^p_vL^\infty_x}\\
		& \le C_p\left(1+\sup_{0\le t \le \hat{t}_0}\|w\tilde{f}(t)\|_{L^p_vL^\infty_x}+\sup_{0\le t \le \hat{t}_0}\|wf(t)\|_{L^p_vL^\infty_x} \right)\int_0^t \|\hat{w}(f-\tilde{f})(s)\|_{L^p_vL^\infty_x}ds.
	\end{align*}
	By the Gr\"onwall inequality, the uniqueness is guaranteed.
\end{proof}

\bigskip
\subsection{The a priori assumption and the related lemmas} \label{apriorisamll}
From now on, we first fix $p>0$ which satisfies the condition \eqref{pcondition} and then fix $\beta>\max\Big\{\frac{3p-6}{2p}+2\gamma, \beta(p,\gamma) \Big\}$, where $\beta(p,\gamma)$ is given by
	\begin{align} \label{betacondition1}
	\beta(p,\gamma) := \begin{cases}
 \frac{3p-5}{p}, & \text{if} \quad 0\le \gamma \le \frac{3}{4},p>4\ \text{or}\ \frac{3}{4}<\gamma \le 1, p>\frac{1}{1-\gamma},  \\
 \frac{7p-12}{2p}, & \text{if}\quad \frac{3}{4}<\gamma \le 1, \frac{5}{2-\gamma}< p \le \frac{1}{1-\gamma}.
\end{cases}
\end{align}
The reason for condition on $\beta$ will become clear later. Next, we fix $\mathcal{T}>0$ such that all the results in the previous subsection in Section \ref{smallperturbationproblem} hold and also the following inequality
\begin{align} \label{assump1}
	(C_1 \mathcal{T}^{\frac{5}{4}})^{\frac{1}{\mathcal{T}}} \le e^{\frac{\nu_0}{2}}
\end{align}
holds. Let $f(t,x,v)$ with 
\begin{align} \label{positivityassump}
	F(t,x,v) = \mu(v)+ \mu^{1/2}(v)f(t,x,v) \ge 0
\end{align}
be the solution to the equation \eqref{FPBE} with the boundary condition \eqref{PDRBC} and initial data $f_0(x,v)$ over the time interval $[0,T]$ for $0<T\le \infty$. We set $h(t,x,v)=w(v)f(t,x,v)$. In this section, we impose the a priori assumption : 
\begin{align} \label{Aprioriassump1}
	\sup_{0\le t < \infty } \|h(t)\|_{L^p_vL^\infty_x} \le \bar{\eta} \ll 1,
\end{align}
where $\bar{\eta}$ depends on the initial amplitude $\eta_0>0$ with $\|h_0\|_{L^p_vL^\infty_x} \le \eta_0$, but does not depend on the solution $h$. It will be determined in subsection \ref{endofsmall}.\\
\indent We need to derive the time decay to $\exp\left\{-\int_s^tR(f)(\tau,X_{\mathbf{cl}}(\tau),V_{\mathbf{cl}}(\tau))d\tau\right\}$ in order to extend the local-in-time solutions to global-in-time solutions. Under the a priori assumption and the smallness condition for the initial relative entropy, we can estimate the lower bound for $R(f)$. From this, we may obtain the time decay mentioned previously.

\begin{Lem} \label{Rfest1}
Under the a priori assumption \eqref{Aprioriassump1}, there exists a small constant $\eta_0>0$, with $\|h_0\|_{L^p_vL^\infty_x} \le \eta_0$, such that for any $T_0$, there is a generally small positive constant $\epsilon_1 = \epsilon_1(\bar{\eta},T_0)>0$, depending only on $\bar{\eta}$ and $T_0$, such that if $\mathcal{E}(F_0) \le \epsilon_1$, then we have
	\begin{align*}
		R(f)(t,x,v) \ge \frac{1}{2}\nu(v)
	\end{align*}
	for all $(t,x,v) \in [0,T_0] \times \Omega \times \R^3$.
\end{Lem}
\begin{remark}
	In fact, for this lemma, it suffices that $\beta> \beta(p,\gamma)$, where $\beta(p,\gamma)$ is defined in \eqref{betacondition1}.
\end{remark}
\begin{proof}
 Let $0\le t \le T_0$. We set $h(t,x,v) =w(v)f(t,x,v)$. We recall
	\begin{align*}
		R(f)(t,x,v) &= \int_{\R^3} \int_{\S^2} B(v-u,\omega) [\mu(u) + \mu^{\frac{1}{2}}(u)f(t,x,u)]d\omega du\\
		& = \nu(v)+\int_{\R^3} \int_{\S^2} B(v-u,\omega)\mu^{\frac{1}{2}}(u)f(t,x,u) d\omega du\\
		& \ge \nu(v) -\tilde{C}_3 \nu(v) \int_{\mathbb{R}^3} e^{-\frac{|u|^2}{8}}|f(t,x,u)|du,
	\end{align*}
	where $\tilde{C}_3$ is a constant depending only on the collision kernel $b$. \\
	If it holds that
	\begin{align} \label{aa61}
		\int_{\mathbb{R}^3} e^{-\frac{|u|^2}{8}}|h(t,x,u)|du\le  \frac{1}{2\tilde{C}_3} \quad \text{for all } 0\le t < T_0,\ x\in \Omega,
	\end{align}
	then we can complete the proof of this Lemma. Thus it suffices to show \eqref{aa61}.
	Then by Duhamel principle,
	\begin{align*}
		h(t,x,v) = S_{G_\nu}(t,0)h_0 + \int_0^t S_{G_\nu}(t,s) \left[K_wh(s)+w\Gamma(f,f)(s)\right]ds.
	\end{align*}
	This yields that 
	\begin{align}\label{Rf10}
		&\int_{\mathbb{R}^3}e^{-\frac{|v|^2}{8}}|h(t,x,v)|dv\nonumber\\
		& \le \int_{\mathbb{R}^3} e^{-\frac{|v|^2}{8}} |(S_{G_\nu}(t,0)h_0)(t,x,v)|dv\nonumber\\
		& \quad + \int_{t-\epsilon}^{t}  \int_{\mathbb{R}^3}e^{-\frac{|v|^2}{8}}\left[\left| \left(S_{G_\nu}(t,s)K_wh(s)\right)(t,x,v)\right|+|(S_{G_\nu}(t,s)w\Gamma(f,f)(s))(t,x,v)|\right]dvds\nonumber\\
		& \quad + \int_0^{t-\epsilon}  \int_{\mathbb{R}^3}e^{-\frac{|v|^2}{8}}\left| \left(S_{G_\nu}(t,s)K_wh(s)\right)(t,x,v)\right|dvds\nonumber\\ 
		& \quad + \int_0^{t-\epsilon} \int_{\mathbb{R}^3}e^{-\frac{|v|^2}{8}}|(S_{G_\nu}(t,s)w\Gamma(f,f)(s))(t,x,v)|dvds\nonumber\\
		& =: I_1+I_2+I_3+I_4.
	\end{align}
	First of all, we easily compute that
	\begin{align} \label{Rf11}
		I_1 \le C \|S_{G_\nu}(t,0)h_0\|_{L^p_vL^\infty_x} \le C_p e^{-\frac{\nu_0}{2} t} \|h_0\|_{L^p_vL^\infty_x}.
	\end{align}
	For $I_2$, from Lemma \ref{KESTIMATE} and Corollary \ref{LpGamma+est}, we get
	\begin{align} \label{Rf12}
		I_2 \le \epsilon C_p \sup_{0\le s \le t}\left[\|h(s)\|_{L^p_vL^\infty_x}+\|h(s)\|_{L^p_vL^\infty_x}^2 \right]. 
	\end{align}
	Let us estimate $I_3$ and $I_4$. For $|v|\ge N$, it follows from Lemma \ref{DampLpDecayorgin}, Lemma \ref{KESTIMATE}, and Corollary \ref{LpGamma+est} that
	\begin{align} \label{Rf50}
		&\int_0^{t-\epsilon}  \int_{|v|\ge N}e^{-\frac{|v|^2}{8}}\left[\left| \left(S_{G_\nu}(t,s)K_wh(s)\right)(t,x,v)\right|+|(S_{G_\nu}(t,s)w\Gamma(f,f)(s))(t,x,v)|\right]dvds\nonumber\\
		& \le \frac{C_{p,\gamma}}{N}\sup_{0\le s \le t}\left[\|h(s)\|_{L^p_vL^\infty_x}+\|h(s)\|_{L^p_vL^\infty_x}^2 \right].
	\end{align}
	Then the only remaining term to estimate is
	\begin{align*}
		\int_0^{t-\epsilon}  \int_{|v|\le N}e^{-\frac{|v|^2}{8}}\left[\left| \left(S_{G_\nu}(t,s)K_wh(s)\right)(t,x,v)\right|+|(S_{G_\nu}(t,s)w\Gamma(f,f)(s))(t,x,v)|\right]dvds
	\end{align*}
	By Lemma \ref{RepresentationforDiffuse1}, we deduce 
	\begin{align*}
		&(S_{G_\nu}(t,s)K_wh(s))(t,x,v)\\
		&=\mathbf{1}_{\{t_1 \le s\}} e^{-\nu(v)(t-s)} K_wh(s,x-v(t-s),v)\\
		& \quad + \frac{e^{-\nu(v)(t-t_1)}}{\tilde{w}(v)} \sum_{l=1}^{k-1} \int_{\prod_{j=1}^{k-1} \mathcal{V}_j} \mathbf{1}_{\{t_{l+1} \le s <t_l\}}K_wh(s,x_l-v_l(t_l-s),v_l)d\Sigma_l(s)\\
		& \quad + \frac{e^{-\nu(v)(t-t_1)}}{\tilde{w}(v)}\int_{\prod_{j=1}^{k-1} \mathcal{V}_j}\mathbf{1}_{\{t_k>s\}} (S_{G_\nu}(t,s)K_wh(s))(t_k,x_k,v_{k-1})d\Sigma_{k-1}(t_k)\\
		& =: \tilde{I}_{31}+\tilde{I}_{32}+\tilde{I}_{33},
	\end{align*}
	and
	\begin{align*}
		&(S_{G_\nu}(t,s)w\Gamma(f,f)(s))(t,x,v)\\
		&=\mathbf{1}_{\{t_1 \le s\}} e^{-\nu(v)(t-s)} w\Gamma(f,f)(s,x-v(t-s),v)\\
		& \quad + \frac{e^{-\nu(v)(t-t_1)}}{\tilde{w}(v)} \sum_{l=1}^{k-1} \int_{\prod_{j=1}^{k-1} \mathcal{V}_j} \mathbf{1}_{\{t_{l+1} \le s <t_l\}}w\Gamma(f,f)(s,x_l-v_l(t_l-s),v_l)d\Sigma_l(s)\\
		& \quad + \frac{e^{-\nu(v)(t-t_1)}}{\tilde{w}(v)}\int_{\prod_{j=1}^{k-1} \mathcal{V}_j}\mathbf{1}_{\{t_k>s\}} (S_{G_\nu}(t,s)w\Gamma(f,f)(s))(t_k,x_k,v_{k-1})d\Sigma_{k-1}(t_k)\\
		& =: \tilde{I}_{41}+\tilde{I}_{42}+\tilde{I}_{43}.
	\end{align*}
	We denote
	\begin{align*}
		I_{ij}:=\int_0^{t-\epsilon}  \int_{|v|\le N}e^{-\frac{|v|^2}{8}}|\tilde{I}_{ij}|dvds
	\end{align*}
	for all $i=3,4$ and $j=1,2,3$. Note that from Lemma \ref{DampLpDecayorgin}, Lemma \ref{KESTIMATE}, Lemma \ref{LpGamma+est},
	\begin{align} \label{Rf13}
		&\|(S_{G_\nu}(t,s)K_wh(s))(t_k,x_k,v_{k-1})\|_{L^p_vL^\infty_x} \le C_p e^{-\frac{\nu_0}{2}(t_k-s)} \|h(s)\|_{L^p_vL^\infty_x},\\ \label{Rf141}
		&\|(S_{G_\nu}(t,s)w\Gamma(f,f)(s))(t_k,x_k,v_{k-1})\|_{L^p_vL^\infty_x} \le C_p e^{-\frac{\nu_0}{2}(t_k-s)} \|h(s)\|_{L^p_vL^\infty_x}.
	\end{align}
	For $I_{33}$ and $I_{43}$, since $t-s\ge \epsilon>0$, by Lemma \ref{Lsmall}, it follows from \eqref{Rf13} and \eqref{Rf141} that, for $k=k(\epsilon,T_0)+1$ such that \eqref{Lsmallineq} holds, 
	\begin{align} \label{Rf51}
		I_{33}+I_{43} &\le C_{p}\sup_{0\le s\le t}\left[\|h(s)\|_{L^p_vL^\infty_x}+\|h(s)\|_{L^p_vL^\infty_x}^2 \right] \int_0^{t-\epsilon}e^{-\frac{\nu_0}{2}(t-s)}\int_{|v|\le N} e^{-\frac{|v|^2}{8}}\nonumber\\
		& \quad \times \left(\int_{\prod_{j=1}^{k-2}\mathcal{V}_j} \mathbf{1}_{\{t_{k-1}>s\}} \prod_{j=1}^{k-2}d\sigma_j\right)dvds\nonumber\\
		&\le \epsilon C_p \sup_{0\le s\le t}\left[\|h(s)\|_{L^p_vL^\infty_x}+\|h(s)\|_{L^p_vL^\infty_x}^2 \right],
	\end{align}
	where we have used the H\"older inequality and 
	\begin{align*}
		\left(\int_{\R^3} \left|\mu^{1/2}(v_{k-1}) \frac{|v_{k-1}|}{w(v_{k-1})} \right|^{p'} dv_{k-1}\right)^{1/p'} \le C_p.
	\end{align*}
	For $I_{31}$, it is bounded by
	\begin{align*}
		\int_{t_1}^{t-\epsilon} \int_{|v|\le N}e^{-\frac{|v|^2}{8}}  e^{-\nu(v)(t-s)} \int_{\R^3}|k_w(v,v')| |h(s,X(s),v')|dvds,
	\end{align*}
	where $X(s) : = x-v(t-s)$.
	we divide this term into two cases.\\
	\newline
	$\mathbf{Case\ 1}:$ $|v'|\ge 2N.$\\
	The following holds:
	\begin{align} \label{Rf14}
		|k_w(v,v')| \le e^{-\frac{\epsilon N^2}{8}}|k_w(v,v')|e^{\frac{\epsilon|v-v'|^2}{8}}
	\end{align}
	By Lemma \ref{KESTIMATE}, we have
	\begin{align} \label{Rf15}
		\sup_{|v|\le N}\int_{|v'|\ge 2N} |k_w(v,v')|^{p'}e^{\frac{\epsilon p'|v-v'|^2}{8}}dv' \le \frac{C_p}{(1+|v|)^{p'(1-\gamma)+1}}.
	\end{align}
	Then, from \eqref{Rf14} and \eqref{Rf15}, we deduce that
	\begin{align} \label{Rf16}
		I_{31}^1 
		& \le C_p e^{-\frac{\epsilon N^2}{8}}\int_{t_1}^{t-\epsilon}\int_{|v|\le N}e^{-\frac{|v|^2}{8}}e^{-\nu_0(t-s)}\|h(s)\|_{L^p_vL^\infty_x}dvds\nonumber\\
		& \le C_p e^{-\frac{\epsilon N^2}{8}} \int_0^{t-\epsilon} e^{-\nu_0(t-s)}\|h(s)\|_{L^p_vL^\infty_x}ds\nonumber\\
		& \le C_p e^{-\frac{\epsilon N^2}{8}} \sup_{0\le s \le t} \|h(s)\|_{L^p_vL^\infty_x}.
	\end{align}
	\newline
	$\mathbf{Case\ 2}:$ $|v'|\le 2N.$\\
	We use the Cauchy-Schwartz inequality and Lemma \ref{KESTIMATE} to obtain
	\begin{align*}
		I_{31}^2 &\le \int_{t_1}^{t-\epsilon}\int_{|v|\le N}e^{-\frac{|v|^2}{8}}e^{-\nu(v)(t-s)}\left(\int_{|v'|\le 2N} |k_w(v,v')|^2 dv'\right)^{1/2}\left(\int_{|v'|\le 2N}|h(s,X(s),v')|^2dv'\right)^{1/2}\\
		& \quad \times dvds\\
		& \le C \int_{t_1}^{t-\epsilon}e^{-\nu_0(t-s)} \int_{|v|\le N}e^{-\frac{|v|^2}{8}}\left(\int_{|v'|\le 2N}|h(s,X(s),v')|^2dv'\right)^{1/2}dvds\\
		& \le C\int_{t_1}^{t-\epsilon}e^{-\nu_0(t-s)} \left(\int_{|v|\le N} \int_{|v'|\le 2N}|h(s,X(s),v')|^2dv'dv\right)^{1/2}ds.
	\end{align*}
	We make a change of variables $v \mapsto y:=X(s) = x-v(t-s)$ with $\left|\frac{dy}{dv} \right| =(t-s)^3$:
	\begin{align*}
		I_{31}^2 &\le  C_\epsilon\int_0^{t-\epsilon}e^{-\nu_0(t-s)} \left(\int_{\Omega} \int_{|v'|\le 2N}|h(s,y,v')|^2dv'dy\right)^{1/2}ds.
	\end{align*}
	Here, we use Lemma \ref{L1L2cont} and the interpolation inequality to compute
\begin{equation} \label{RTIRE}
\begin{aligned}
	&\int_{\Omega}\int_{|v'|\le 2N}|h(s,y,v')|^2dv'dy\\
	& \le C_N \int_{\Omega}\int_{|v'|\le 2N}|f(s,y,v')|^2 \mathbf{1}_{\{|f|\le \sqrt{\mu}\}} dv'dy +C_N\int_{\Omega}\int_{|v'|\le 2N}|h(s,y,v')|^2 \mathbf{1}_{\{|f|> \sqrt{\mu}\}}dv'dy\\
	&\le C_N \mathcal{E}(F_0) + C_N\left(\int_{\Omega}\int_{|v'|\le 2N}|h(s,y,v')| \mathbf{1}_{\{|f|> \sqrt{\mu}\}}dv'dy\right)^{\kappa} \left(\int_{\Omega}\int_{|v'|\le 2N}|h(s,y,v')|^p \mathbf{1}_{\{|f|> \sqrt{\mu}\}}dv'dy\right)^{\frac{1-\kappa}{p}}\\
	&\le C_N \mathcal{E}(F_0) + C_{p,N}\|h(s)\|_{L^p_vL^\infty_x}^{1-\kappa}\left(\int_{\Omega}\int_{|v'|\le 2N}\mu^{1/2}(v')|f(s,y,v')| \mathbf{1}_{\{|f|> \sqrt{\mu}\}}dv'dy\right)^{\kappa}\\
	&\le C_N \mathcal{E}(F_0) + C_{p,N}\|h(s)\|_{L^p_vL^\infty_x}^{1-\kappa}\mathcal{E}(F_0)^\kappa,
\end{aligned}
\end{equation}
where
\begin{align} \label{kappaconst}
	\kappa = \frac{1/2 - 1/p}{1-1/p} = \frac{p-2}{2p-2}.
\end{align}
Hence we have
\begin{align} \label{Rf17}
	I_{31}^2 &\le  C_{p,\epsilon,N} \int_0^{t-\epsilon}e^{-\nu_0(t-s)} \left[\mathcal{E}(F_0)^{1/2} + \|h(s)\|_{L^p_vL^\infty_x}^{\frac{1-\kappa}{2}}\mathcal{E}(F_0)^{\kappa/2} \right] ds\nonumber\\
	& \le  C_{p,\epsilon,N} \left[\mathcal{E}(F_0)^{1/2} + \sup_{0 \le s \le t}\|h(s)\|_{L^p_vL^\infty_x}^{\frac{1-\kappa}{2}}\mathcal{E}(F_0)^{\kappa/2} \right]\nonumber\\
	& \le \frac{C_p}{N} \sup_{0 \le s \le t}\|h(s)\|_{L^p_vL^\infty_x}^{1-\kappa} + C_{p,\epsilon,N} \left[\mathcal{E}(F_0)^{1/2}+ \mathcal{E}(F_0)^\kappa \right],
\end{align}
where we have used the Young's inequality.\\
Gathering \eqref{Rf16} and \eqref{Rf17}, we get
\begin{align} \label{Rf18}
	I_{31} \le \frac{C_{p,\epsilon}}{N}\sup_{0\le s \le t}\left[\|h(s)\|_{L^p_vL^\infty_x}+\|h(s)\|_{L^p_vL^\infty_x}^{1-\kappa}\right] + C_{p,\epsilon,N} \left[\mathcal{E}(F_0)^{1/2}+ \mathcal{E}(F_0)^\kappa \right].
\end{align}
For $I_{32}$, it is bounded by
\begin{align*}
	\int_0^{t-\epsilon}  \int_{|v|\le N}e^{-\frac{|v|^2}{8}}\frac{e^{-\nu(v)(t-t_1)}}{\tilde{w}(v)} \sum_{l=1}^{k-1} \int_{\prod_{j=1}^{k-1} \mathcal{V}_j} \mathbf{1}_{\{t_{l+1} \le s <t_l\}}K_wh(s,x_l-v_l(t_l-s),v_l)d\Sigma_l(s)dvds.
\end{align*}
Fix $l$. Then the above term is bounded by
\begin{align*}
	 I_{32l}:=&C\int_{t_{l+1}}^{t_l}e^{-\nu_0(t_l-s)} \int_{|v|\le N}e^{-\frac{|v|^2}{8}}\int_{\prod_{j=1}^{l-1}  \mathcal{V}_j}\int_{\mathcal{V}_l} \int_{\R^3}|k_w(v_l,v')||h(s,x_l-v_l(t_l-s),v')|dv'\\
	 &\quad \times  \mu^{1/4}(v_l)dv_l\prod_{j=1}^{l-1} d\sigma_j dvds.
\end{align*}
we divide this term into four cases.\\
\newline
$\mathbf{Case\ 1}:$ $|v_l|\ge N.$\\
We use the H\"older inequality and Lemma \ref{KESTIMATE} to estimate
\begin{align} \label{Rf19}
	I_{32l}^1 &\le C_p \int_{t_{l+1}}^{t_l}e^{-\nu_0(t_l-s)} \|h(s)\|_{L^p_vL^\infty_x} \int_{|v|\le N}e^{-\frac{|v|^2}{8}}\int_{\prod_{j=1}^{l-1}  \mathcal{V}_j}\int_{|v_l|\ge N}\mu^{1/4}(v_l)dv_l\prod_{j=1}^{l-1} d\sigma_j dvds\nonumber\\
	& \le \frac{C_p}{N}\sup_{0\le s \le t}\|h(s)\|_{L^p_vL^\infty_x}.
\end{align}
\newline
$\mathbf{Case\ 2}:$ $|v_l|\le N, |v'|\ge 2N.$\\
The following is valid:
\begin{align} \label{Rf20}
	|k_w(v_l,v')| ^{p'}\le e^{-\frac{\epsilon}{8}N^2}|k_w(v_l,v')|^{p'}e^{\frac{\epsilon}{8}|v_l-v'|^2},
\end{align}
where $\epsilon$ is sufficiently small.
From Lemma \ref{KESTIMATE}, it holds that
\begin{align} \label{Rf21}
	 \int_{|v'|\ge 2N}|k_w(v_l,v')|^{p'}e^{\frac{\epsilon}{8}|v_l-v'|^2}dv'' \le C_{p} \left(\frac{1}{1+|v_l|}\right)^{p'(1-\gamma)+1}.
\end{align}
Then we use \eqref{Rf20} and \eqref{Rf21} to obtain
\begin{align} \label{Rf22}
	I_{32l}^2&\le C_p e^{-\frac{\epsilon}{8}N^2}\int_{t_{l+1}}^{t_l}e^{-\nu_0(t_l-s)} \|h(s)\|_{L^p_vL^\infty_x} \int_{|v|\le N}e^{-\frac{|v|^2}{8}}\int_{\prod_{j=1}^{l-1}  \mathcal{V}_j}\int_{|v_l|\le N}\mu^{1/4}(v_l)dv_l\prod_{j=1}^{l-1} d\sigma_j dvds\nonumber\\
	& \le C_p e^{-\frac{\epsilon}{8}N^2}\sup_{0\le s \le t}\|h(s)\|_{L^p_vL^\infty_x}.
\end{align}
\newline
$\mathbf{Case\ 3}:$ $t_l-s\le \frac{1}{N}.$\\
We can easily compute
\begin{align} \label{Rf23}
	I_{32l}^3 &\le C_p\int_{t_l-\frac{1}{N}}^{t_l}\|h(s)\|_{L^p_vL^\infty_x} \int_{|v|\le N}e^{-\frac{|v|^2}{8}}\int_{\prod_{j=1}^{l-1}  \mathcal{V}_j}\int_{\mathcal{V}_l} \mu^{1/4}(v_l)dv_l\prod_{j=1}^{l-1} d\sigma_j dvds\nonumber\\
	& \le \frac{C_p}{N}\sup_{0\le s \le t}\|h(s)\|_{L^p_vL^\infty_x},
\end{align}
where we have used Lemma \ref{KESTIMATE}.\\
\newline
$\mathbf{Case\ 4}:$ $t_l-s\ge \frac{1}{N}$ and $|v_l|\le N, |v'|\le 2N$.\\
We use the Cauchy-Schwartz inequality and Lemma \ref{KESTIMATE} to estimate
\begin{align*}
	I_{32l}^4 &\le C\int_{t_{l+1}}^{t_l-\frac{1}{N}}e^{-\nu_0(t_l-s)} \int_{|v|\le N}e^{-\frac{|v|^2}{8}}\int_{\prod_{j=1}^{l-1}  \mathcal{V}_j}\left(\int_{|v_l|\le N} \int_{|v'|\le 2N}|k_w(v_l,v')|^2\mu^{1/2}(v_l)dv'dv_l\right)^{1/2}\\
	 &\quad \times \left(\int_{|v_l|\le N} \int_{|v'|\le 2N}|h(s,x_l-v_l(t_l-s),v')|^2dv' dv_l\right)^{1/2}\prod_{j=1}^{l-1} d\sigma_j dvds\\
	 &\le C_p\int_{t_{l+1}}^{t_l-\frac{1}{N}}e^{-\nu_0(t_l-s)} \int_{|v|\le N}e^{-\frac{|v|^2}{8}}\int_{\prod_{j=1}^{l-1}  \mathcal{V}_j}\left(\int_{|v_l|\le N} \int_{|v'|\le 2N}|h(s,x_l-v_l(t_l-s),v')|^2dv' dv_l\right)^{1/2}\\
	 & \quad \times \prod_{j=1}^{l-1} d\sigma_j dvds.
\end{align*}
We make a change of variables $v_l \mapsto y=x_l-v_l(t_l-s)$ with $\left|\frac{dy}{dv_l}\right| = (t_l-s)^3$:
\begin{align*}
	I_{32l}^4 &\le C_{p,N} \int_{t_{l+1}}^{t_l-\frac{1}{N}}e^{-\nu_0(t_l-s)} \int_{|v|\le N}e^{-\frac{|v|^2}{8}}\int_{\prod_{j=1}^{l-1}  \mathcal{V}_j}\left(\int_{\Omega} \int_{|v'|\le 2N}|h(s,y,v')|^2dv' dy\right)^{1/2}\prod_{j=1}^{l-1} d\sigma_j dvds,
\end{align*}
and It follows from \eqref{RTIRE} that
\begin{align} \label{Rf24}
	I_{32l}^4 &\le C_{p,N} \int_{t_{l+1}}^{t_l}e^{-\nu_0(t_l-s)} \int_{|v|\le N}e^{-\frac{|v|^2}{8}}\int_{\prod_{j=1}^{l-1}  \mathcal{V}_j}\left[\mathcal{E}(F_0)^{1/2} + \|h(s)\|_{L^p_vL^\infty_x}^{\frac{1-\kappa}{2}}\mathcal{E}(F_0)^{\kappa/2} \right]\prod_{j=1}^{l-1} d\sigma_j dvds\nonumber\\
	& \le \frac{C_p}{N}\sup_{0\le s\le t}\|h(s)\|_{L^p_vL^\infty_x}^{1-\kappa}+ C_{p,N}\left[\mathcal{E}(F_0)^{1/2}+ \mathcal{E}(F_0)^{\kappa}\right],
\end{align}
where we have used the Young's inequality.
Combining \eqref{Rf19}, \eqref{Rf22}, \eqref{Rf23}, \eqref{Rf24} and summing over $1\le l \le k(\epsilon,T_0)-1$, we obtain
\begin{align} \label{Rf25}
	I_{32}\le \frac{C_{p,\epsilon,T_0}}{N}\sup_{0\le s \le t}\left[\|h(s)\|_{L^p_vL^\infty_x}+\|h(s)\|_{L^p_vL^\infty_x}^{1-\kappa}\right] + C_{p,\epsilon,N,T_0} \left[\mathcal{E}(F_0)^{1/2}+ \mathcal{E}(F_0)^\kappa \right].
\end{align}
For $I_{41}$, we first split it into two terms:
\begin{align*}
	&\int_{t_1}^{t-\epsilon} \int_{|v|\le N}e^{-\frac{|v|^2}{8}}e^{-\nu(v)(t-s)}|w\Gamma(f,f)(s,X(s),v)|dvds \\
	& \le \int_{t_1}^{t-\epsilon} \int_{|v|\le N}e^{-\frac{|v|^2}{8}}e^{-\nu(v)(t-s)}|w\Gamma^+(f,f)(s,X(s),v)|dvds\\
	& \quad +\int_{t_1}^{t-\epsilon} \int_{|v|\le N}e^{-\frac{|v|^2}{8}}e^{-\nu(v)(t-s)}|w\Gamma^-(f,f)(s,X(s),v)|dvds\\
	& =: I_{411}+I_{412},
\end{align*}
where $X(s)=x-v(t-s)$.
Firstly, let us estimate $I_{411}$. From the estimate \eqref{Gamma+estimatecombine}, the term $I_{411}$ is bounded by	
\begin{align*}
	C_{p,\gamma}\int_{t_1}^{t-\epsilon}e^{-\nu_0(t-s)}\|h(s)\|_{L^p_vL^\infty_x} \int_{|v|\le N}e^{-\frac{|v|^2}{8}} \left(\int_{\R^3}(1+|\eta|)^{e(p,\gamma)}|f(s,X(s),\eta)|^2 d\eta \right)^{1/2}dvds.
\end{align*}
 We divide this term into two cases.\\
\newline
$\mathbf{Case\ 1}:$ $|\eta|\ge N$.\\
Using the H\"older inequality and thanks to the condition \eqref{betacondition1}, we deduce
\begin{align} \label{Rf26}
	I_{411}^1 &\le \frac{C_{\beta,p,\gamma}}{N^{1/2}}\int_{t_1}^{t-\epsilon}e^{-\nu_0(t-s)}\|h(s)\|_{L^p_vL^\infty_x} \int_{|v|\le N}e^{-\frac{|v|^2}{8}} \left(\int_{|\eta|\ge N}(1+|\eta|)^{e(p,\gamma)+1-2\beta}|h(s,X(s),\eta)|^2 d\eta \right)^{1/2}dvds	 \nonumber\\
	& \le \frac{C_{p,\gamma}}{N^{1/2}}\int_0^{t-\epsilon}e^{-\nu_0(t-s)}\|h(s)\|_{L^p_vL^\infty_x}^2ds \int_{|v|\le N}e^{-\frac{|v|^2}{8}}dv \nonumber\\
	& \le \frac{C_{p,\gamma}}{N^{1/2}}\sup_{0\le s \le t}\|h(s)\|_{L^p_vL^\infty_x}^2.
\end{align}
\newline
$\mathbf{Case\ 2}:$ $|\eta|\le N$.\\
We have
\begin{align*}
	I_{411}^2 &\le C_{p,\gamma}\int_{t_1}^{t-\epsilon}e^{-\nu_0(t-s)}\|h(s)\|_{L^p_vL^\infty_x} \int_{|v|\le N}e^{-\frac{|v|^2}{8}} \left(\int_{|\eta|\le N}|h(s,X(s),\eta)|^2 d\eta \right)^{1/2}dvds\\
	& \le C_{p,\gamma}\int_{t_1}^{t-\epsilon}e^{-\nu_0(t-s)}\|h(s)\|_{L^p_vL^\infty_x} \left(\int_{|v|\le N} \int_{|\eta|\le N}|h(s,X(s),\eta)|^2 d\eta dv\right)^{1/2}ds.
\end{align*}
We make a change of variables $v \mapsto y:=X(s) = x-v(t-s)$ with $\left|\frac{dy}{dv} \right| =(t-sx)^3$:
	\begin{align*}
		I_{411}^2 &\le C_{p,\gamma,\epsilon}\int_{t_1}^{t-\epsilon}e^{-\nu_0(t-s)}\|h(s)\|_{L^p_vL^\infty_x} \left(\int_{\Omega} \int_{|\eta|\le N}|h(s,y,\eta)|^2d\eta dy\right)^{1/2}ds.
	\end{align*}
From \eqref{RTIRE}, we get
\begin{align} \label{Rf27}
	I_{411}^2 &\le C_{p,\gamma,\epsilon}\int_0^{t-\epsilon}e^{-\nu_0(t-s)}\|h(s)\|_{L^p_vL^\infty_x}\left[\mathcal{E}(F_0)^{1/2} + \|h(s)\|_{L^p_vL^\infty_x}^{\frac{1-\kappa}{2}}\mathcal{E}(F_0)^{\kappa/2} \right] ds \nonumber\\
	& \le \frac{C_p}{N} \sup_{0 \le s \le t}\left[\|h(s)\|_{L^p_vL^\infty_x}^2+\|h(s)\|_{L^p_vL^\infty_x}^{3-\kappa}\right] + C_{p,\gamma,\epsilon,N} \left[\mathcal{E}(F_0)+ \mathcal{E}(F_0)^\kappa \right],
\end{align}
where we have used the Young's inequality.\\
Gathering \eqref{Rf26} and \eqref{Rf27}, we obtain
\begin{align} \label{Rf28}
	I_{411} \le \frac{C_{p,\gamma}}{N^{1/2}}\sup_{0 \le s \le t}\left[\|h(s)\|_{L^p_vL^\infty_x}^2+\|h(s)\|_{L^p_vL^\infty_x}^{3-\kappa}\right]+ C_{p,\gamma,\epsilon,N} \left[\mathcal{E}(F_0)+ \mathcal{E}(F_0)^\kappa \right].
\end{align}
Next, let us estimate $I_{412}$. Then the term $I_{412}$ is bounded by	
\begin{align*}
	C_{p,\gamma}\int_{t_1}^{t-\epsilon}e^{-\nu_0(t-s)} \int_{|v| \le N}e^{-\frac{|v|^2}{8}} \nu(v)|h(s,X(s),v)|\left(\int_{\R^3}\mu^{1/4}(u)|f(s,X(s),u)| du \right)dvds.
\end{align*}
We divide this term into two cases.\\
\newline
$\mathbf{Case\ 1}:$  $|u| \ge N$.\\
We can easily compute 
\begin{align} \label{Rf29}
	I_{412}^1 &\le \frac{C_{p,\gamma}}{N}\int_{t_1}^{t-\epsilon}e^{-\nu_0(t-s)} \|h(s)\|_{L^p_vL^\infty_x} \int_{|v|\le N} e^{-\frac{|v|^2}{16}}|h(s,X(s),v)|dvds \nonumber\\	
	& \le \frac{C_{p,\gamma}}{N}\int_{t_1}^{t-\epsilon}e^{-\nu_0(t-s)} \|h(s)\|_{L^p_vL^\infty_x}^2ds \nonumber\\
	& \le \frac{C_{p,\gamma}}{N} \sup_{0\le s \le t}\|h(s)\|_{L^p_vL^\infty_x}^2.
\end{align}
\newline
$\mathbf{Case\ 2}:$  $|u| \le N$.\\
We use the H\"older inequality to derive
\begin{align*}
	I_{412}^2 &\le C_{p,\gamma}\int_{t_1}^{t-\epsilon}e^{-\nu_0(t-s)} \int_{|v|\le N}e^{-\frac{|v|^2}{16}} |h(s,X(s),v)|\left(\int_{|u|\le N}|h(s,X(s),u)|^2 du \right)^{1/2}dvds \nonumber\\
	& \le C_{p,\gamma}\int_{t_1}^{t-\epsilon}e^{-\nu_0(t-s)} \|h(s)\|_{L^p_vL^\infty_x} \left(\int_{|v|\le N}\int_{|u|\le N}|h(s,X(s),u)|^2 du dv\right)^{1/2}ds.
\end{align*}
We make a change of variables $v \mapsto y:=X(s) = x-v(t-s)$ with $\left|\frac{dy}{dv} \right| =(t-s)^3$:
	\begin{align*}
		I_{412}^2 &\le C_{p,\gamma,\epsilon}\int_0^{t-\epsilon}e^{-\nu_0(t-s)}\|h(s)\|_{L^p_vL^\infty_x} \left(\int_{\Omega} \int_{|u|\le N}|h(s,y,u)|^2du dy\right)^{1/2}ds.
	\end{align*}
From \eqref{RTIRE}, we get
\begin{align} \label{Rf30}
	I_{412}^2 &\le C_{p,\gamma,\epsilon}\int_0^{t-\epsilon}e^{-\nu_0(t-s)}\|h(s)\|_{L^p_vL^\infty_x}\left[\mathcal{E}(F_0)^{1/2} + \|h(s)\|_{L^p_vL^\infty_x}^{\frac{1-\kappa}{2}}\mathcal{E}(F_0)^{\kappa/2} \right] ds \nonumber\\
	& \le \frac{C_p}{N} \sup_{0 \le s \le t}\left[\|h(s)\|_{L^p_vL^\infty_x}^2+\|h(s)\|_{L^p_vL^\infty_x}^{3-\kappa}\right] + C_{p,\gamma,\epsilon,N} \left[\mathcal{E}(F_0)+ \mathcal{E}(F_0)^\kappa \right],
\end{align}
where we have used the Young's inequality.\\
Gathering \eqref{Rf29} and \eqref{Rf30}, we obtain
\begin{align} \label{Rf31}
	I_{412} \le \frac{C_{p,\gamma}}{N}\sup_{0 \le s \le t}\left[\|h(s)\|_{L^p_vL^\infty_x}^2+\|h(s)\|_{L^p_vL^\infty_x}^{3-\kappa}\right]+ C_{p,\gamma,\epsilon,N} \left[\mathcal{E}(F_0)+ \mathcal{E}(F_0)^\kappa \right].
\end{align}
For $I_{42}$, we split it into two terms:
\begin{align*}
	 &\int_0^{t-\epsilon} \int_{|v|\le N} e^{-\frac{|v|^2}{8}} \frac{e^{-\nu(v)(t-t_1)}}{\tilde{w}(v)} \sum_{l=1}^{k-1} \int_{\prod_{j=1}^{k-1} \mathcal{V}_j} \mathbf{1}_{\{t_{l+1} \le s <t_l\}}|w\Gamma(f,f)(s,x_l-v_l(t_l-s),v_l)|d\Sigma_l(s)dvds\\
	 &\le \int_0^{t-\epsilon} \int_{|v|\le N} e^{-\frac{|v|^2}{8}} e^{-\nu(v)(t-t_1)}\sum_{l=1}^{k-1} \int_{\prod_{j=1}^{k-1} \mathcal{V}_j} \mathbf{1}_{\{t_{l+1} \le s <t_l\}}|w\Gamma^+(f,f)(s,x_l-v_l(t_l-s),v_l)|d\Sigma_l(s)dvds\\
	 & \quad +\int_0^{t-\epsilon} \int_{|v|\le N} e^{-\frac{|v|^2}{8}} e^{-\nu(v)(t-t_1)} \sum_{l=1}^{k-1} \int_{\prod_{j=1}^{k-1} \mathcal{V}_j} \mathbf{1}_{\{t_{l+1} \le s <t_l\}}|w\Gamma^-(f,f)(s,x_l-v_l(t_l-s),v_l)|d\Sigma_l(s)dvds\\
	 &=:I_{421}+I_{422}.
\end{align*}
Fix $l$. Then the $l$-th term for $I_{421}$ is bounded by
\begin{align*}
	I_{421l}&:=\int_{t_{l+1}}^{t_l} e^{-\nu_0(t_l-s)}\int_{|v|\le N} e^{-\frac{|v|^2}{8}} \int_{\prod_{j=1}^{l-1} \mathcal{V}_j} \int_{\mathcal{V}_l} |w\Gamma^+(f,f)(s,x_l-v_l(t_l-s),v_l)|\mu^{1/4}(v_l)dv_l\\
	& \quad \times \prod_{j=1}^{l-1} d\sigma_j dvds\\
	& \le C_{p,\gamma}\int_{t_{l+1}}^{t_l} e^{-\nu_0(t_l-s)} \|h(s)\|_{L^p_vL^\infty_x}\int_{|v|\le N} e^{-\frac{|v|^2}{8}}\int_{\prod_{j=1}^{l-1} \mathcal{V}_j}  \\
	& \quad  \times \int_{\mathcal{V}_l} \left(\int_{\R^3}(1+|\eta|)^{e(p,\gamma)}|f(s,x_l-v_l(t_l-s),\eta)|^2 d\eta \right)^{1/2}\mu^{1/4}(v_l)dv_l\prod_{j=1}^{l-1} d\sigma_j dvds,
\end{align*}
where we have used \eqref{Gamma+estimatecombine}. We divide it into four cases.\\
\newline
$\mathbf{Case\ 1}:$ $|v_l|\ge N.$\\
We use the H\"older inequality to estimate
\begin{align} \label{Rf32}
	I_{421l}^1 &\le C_{p,\gamma}\int_{t_{l+1}}^{t_l} e^{-\nu_0(t_l-s)} \|h(s)\|_{L^p_vL^\infty_x}^2\int_{|v|\le N} e^{-\frac{|v|^2}{8}}\int_{\prod_{j=1}^{l-1} \mathcal{V}_j} \int_{|v_l|\ge N}\mu^{1/4}(v_l)dv_l\prod_{j=1}^{l-1} d\sigma_j dvds\nonumber\\
	& \le \frac{C_{p,\gamma}}{N}\int_{t_{l+1}}^{t_l}e^{-\nu_0(t_l-s)} \|h(s)\|_{L^p_vL^\infty_x}^2ds\nonumber\\
	& \le \frac{C_{p,\gamma}}{N} \sup_{0\le s \le t}\|h(s)\|_{L^p_vL^\infty_x}^2.
\end{align}
\newline
$\mathbf{Case\ 2}:$ $|v_l|\le N, |\eta|\ge N.$\\
We use the H\"older inequality to obtain
\begin{align} \label{Rf33}
	I_{421l}^2 &\le \frac{C_{p,\gamma}}{N^{1/2}}\int_{t_{l+1}}^{t_l} e^{-\nu_0(t_l-s)} \|h(s)\|_{L^p_vL^\infty_x}^2\int_{|v|\le N} e^{-\frac{|v|^2}{8}}\int_{\prod_{j=1}^{l-1} \mathcal{V}_j} \int_{|v_l|\le N}\mu^{1/4}(v_l)dv_l\prod_{j=1}^{l-1} d\sigma_j dvds\nonumber\\
	& \le \frac{C_{p,\gamma}}{N^{1/2}}\sup_{0\le s \le t}\|h(s)\|_{L^p_vL^\infty_x}^2.
\end{align}
\newline
$\mathbf{Case\ 3}:$ $t_l-s\le \frac{1}{N}.$\\
We can easily compute
\begin{align} \label{Rf34}
	I_{421l}^3 &\le  C_{p,\gamma}\int_{t_l-\frac{1}{N}}^{t_l}  \|h(s)\|_{L^p_vL^\infty_x}^2\int_{|v|\le N} e^{-\frac{|v|^2}{8}}\int_{\prod_{j=1}^{l-1} \mathcal{V}_j} \int_{\mathcal{V}_l}\mu^{1/4}(v_l)dv_l\prod_{j=1}^{l-1} d\sigma_j dvds\nonumber\\
	& \le C_{p,\gamma}\int_{t_l-\frac{1}{N}}^{t_l}  \|h(s)\|_{L^p_vL^\infty_x}^2ds\nonumber\\
	& \le \frac{C_{p,\gamma}}{N}\sup_{0\le s \le t}\|h(s)\|_{L^p_vL^\infty_x}^2.
\end{align}
\newline
$\mathbf{Case\ 4}:$ $t_l-s\ge \frac{1}{N}$ and $|v_l|\le N, |\eta|\le N$.\\
We use the Cauchy-Schwartz inequality to estimate
\begin{align*}
	I_{421l}^4 &\le C_{p,\gamma}\int_{t_{l+1}}^{t_l-\frac{1}{N}} e^{-\nu_0(t_l-s)} \|h(s)\|_{L^p_vL^\infty_x}\int_{|v|\le N} e^{-\frac{|v|^2}{8}}\int_{\prod_{j=1}^{l-1} \mathcal{V}_j} \\
	& \quad  \times \left(\int_{|v_l|\le N} \int_{|\eta|\le N}|h(s,x_l-v_l(t_l-s),\eta)|^2 d\eta dv_l\right)^{1/2}\prod_{j=1}^{l-1} d\sigma_j dvds.
\end{align*}
We make a change of variables $v_l \mapsto y=x_l-v_l(t_l-s)$ with $\left|\frac{dy}{dv_l}\right| = (t_l-s)^3$:
\begin{align*}
	I_{421l}^4 &\le C_{p,\gamma,N}\int_{t_{l+1}}^{t_l-\frac{1}{N}} e^{-\nu_0(t_l-s)} \|h(s)\|_{L^p_vL^\infty_x}\int_{|v|\le N} e^{-\frac{|v|^2}{8}}\int_{\prod_{j=1}^{l-1} \mathcal{V}_j} \left(\int_{\Omega} \int_{|\eta|\le N}|h(s,y,\eta)|^2 d\eta dy\right)^{1/2}\\
	&\quad \times \prod_{j=1}^{l-1} d\sigma_j dvds.
\end{align*}
It follows from \eqref{RTIRE} that
\begin{align} \label{Rf35}
	I_{421l}^4 &\le C_{p,\gamma,N}\int_{t_{l+1}}^{t_l} e^{-\nu_0(t_l-s)} \|h(s)\|_{L^p_vL^\infty_x}\int_{|v|\le N} e^{-\frac{|v|^2}{8}}\int_{\prod_{j=1}^{l-1} \mathcal{V}_j} \left[\mathcal{E}(F_0)^{1/2} + \|h(s)\|_{L^p_vL^\infty_x}^{\frac{1-\kappa}{2}}\mathcal{E}(F_0)^{\kappa/2}\right]\nonumber\\
	&\quad \times \prod_{j=1}^{l-1} d\sigma_j dvds\nonumber\\
	& \le C_{p,\gamma,N}\int_{t_{l+1}}^{t_l} e^{-\nu_0(t_l-s)} \|h(s)\|_{L^p_vL^\infty_x}\left[\mathcal{E}(F_0)^{1/2} + \|h(s)\|_{L^p_vL^\infty_x}^{\frac{1-\kappa}{2}}\mathcal{E}(F_0)^{\kappa/2}\right]ds\nonumber\\
	& \le \frac{C_p}{N}\sup_{0\le s\le t}\left[\|h(s)\|_{L^p_vL^\infty_x}^2+\|h(s)\|_{L^p_vL^\infty_x}^{3-\kappa}\right]+ C_{p,\gamma,N}\left[\mathcal{E}(F_0)+ \mathcal{E}(F_0)^{\kappa}\right],
\end{align}
where we have used the Young's inequality.\\
Now, the $l$-th term for $I_{422}$ is bounded by
\begin{align*}
	I_{422l}&:=\int_{t_{l+1}}^{t_l} e^{-\nu_0(t_l-s)}\int_{|v|\le N} e^{-\frac{|v|^2}{8}} \int_{\prod_{j=1}^{l-1} \mathcal{V}_j} \int_{\mathcal{V}_l} |w\Gamma^-(f,f)(s,x_l-v_l(t_l-s),v_l)|\mu^{1/4}(v_l)dv_l\\
	& \quad \times \prod_{j=1}^{l-1} d\sigma_j dvds\\
	& \le C_{p}\int_{t_{l+1}}^{t_l} e^{-\nu_0(t_l-s)} \int_{|v|\le N} e^{-\frac{|v|^2}{8}}\int_{\prod_{j=1}^{l-1} \mathcal{V}_j} \int_{\mathcal{V}_l} \nu(v_l)\mu^{1/4}(v_l)|h(s,x_l-v_l(t_l-s),v_l)|dv_l \\
	& \quad  \times \left(\int_{\R^3}\mu^{1/4}(u)|f(s,x_l-v_l(t_l-s),u)| du \right)\prod_{j=1}^{l-1} d\sigma_j dvds.
\end{align*}
 We divide it into four cases.\\
\newline
$\mathbf{Case\ 1}:$ $|v_l|\ge N.$\\
We use the H\"older inequality to estimate
\begin{align} \label{Rf36}
	I_{422l}^1 &\le C_{p}\int_{t_{l+1}}^{t_l} e^{-\nu_0(t_l-s)} \|h(s)\|_{L^p_vL^\infty_x}\int_{|v|\le N} e^{-\frac{|v|^2}{8}}\int_{\prod_{j=1}^{l-1} \mathcal{V}_j}\nonumber \\
	& \quad  \times  \int_{|v_l| \ge N} \nu(v_l)\mu^{1/4}(v_l)|h(s,x_l-v_l(t_l-s),v_l)|dv_l\prod_{j=1}^{l-1} d\sigma_j dvds \nonumber\\
	& \le \frac{C_p}{N} \int_{t_{l+1}}^{t_l} e^{-\nu_0(t_l-s)} \|h(s)\|_{L^p_vL^\infty_x}^2ds\nonumber\\
	& \le \frac{C_p}{N} \sup_{0\le s \le t}\|h(s)\|_{L^p_vL^\infty_x}^2.
\end{align}
\newline
$\mathbf{Case\ 2}:$ $|v_l|\le N, |u|\ge N.$\\
We use the H\"older inequality to obtain
\begin{align} \label{Rf37}
	I_{422l}^2 & \le \frac{C_p}{N}\int_{t_{l+1}}^{t_l} e^{-\nu_0(t_l-s)} \|h(s)\|_{L^p_vL^\infty_x}\int_{|v|\le N} e^{-\frac{|v|^2}{8}}\int_{\prod_{j=1}^{l-1} \mathcal{V}_j}\nonumber \\
	& \quad  \times  \int_{|v_l| \le N} \nu(v_l)\mu^{1/4}(v_l)|h(s,x_l-v_l(t_l-s),v_l)|dv_l\prod_{j=1}^{l-1} d\sigma_j dvds \nonumber\\
	& \le  \frac{C_p}{N} \sup_{0\le s \le t}\|h(s)\|_{L^p_vL^\infty_x}^2.
\end{align}
\newline
$\mathbf{Case\ 3}:$ $t_l-s\le \frac{1}{N}.$\\
We can easily compute
\begin{align} \label{Rf38}
	I_{422l}^3 & \le C_{p}\int_{t_l-\frac{1}{N}}^{t_l}  \|h(s)\|_{L^p_vL^\infty_x}^2ds \le \frac{C_{p}}{N}\sup_{0\le s \le t}\|h(s)\|_{L^p_vL^\infty_x}^2.
\end{align}
\newline
$\mathbf{Case\ 4}:$ $t_l-s\ge \frac{1}{N}$ and $|v_l|\le N, |\eta|\le N$.\\
By a similar way to derive \eqref{Rf35}, we deduce
\begin{align} \label{Rf39}
	I_{422l}^3 \le \frac{C_p}{N}\sup_{0\le s\le t}\left[\|h(s)\|_{L^p_vL^\infty_x}^2+\|h(s)\|_{L^p_vL^\infty_x}^{3-\kappa}\right]+ C_{p, N}\left[\mathcal{E}(F_0)+ \mathcal{E}(F_0)^{\kappa}\right].
\end{align}
Gathering \eqref{Rf32}, \eqref{Rf33}, \eqref{Rf34}, \eqref{Rf35}, \eqref{Rf36}, \eqref{Rf37}, \eqref{Rf38} ,\eqref{Rf39}, and summing over $1\le l \le k(\epsilon,T_0)-1$, we obtain
\begin{align} \label{Rf40}
	I_{42} \le \frac{C_{p,\epsilon,T_0}}{N^{1/2}}\sup_{0\le s \le t}\left[ \|h(s)\|_{L^p_vL^\infty_x}^2+\|h(s)\|_{L^p_vL^\infty_x}^{3-\kappa}\right]+C_{p,N,\epsilon,T_0}\left[\mathcal{E}(F_0)+ \mathcal{E}(F_0)^{\kappa}\right]
\end{align}
Inserting \eqref{Rf11}, \eqref{Rf12}, \eqref{Rf50}, \eqref{Rf51}, \eqref{Rf18}, \eqref{Rf25}, \eqref{Rf28}, \eqref{Rf31}, \eqref{Rf40} into \eqref{Rf10}, it follows that
\begin{align*}
	\int_{\mathbb{R}^3}e^{-\frac{|v|^2}{8}}|h(t,x,v)|dv &\le e^{-\frac{\nu_0}{2} t} \|h_0\|_{L^p_vL^\infty_x}\\
	& \quad  + C_p \left(\epsilon+\frac{C_{p,\gamma,\epsilon,T_0}}{N^{1/2}} \right)  \sup_{0\le s \le t}\Bigl[\|h(s)\|_{L^p_vL^\infty_x}+\|h(s)\|_{L^p_vL^\infty_x}^2+\|h(s)\|_{L^p_vL^\infty_x}^{1-\kappa}\\
	& \qquad +\|h(s)\|_{L^p_vL^\infty_x}^{3-\kappa} \Bigr]\\
	& \quad +C_{p,\gamma,\epsilon,N,T_0} \left[\mathcal{E}(F_0)^{1/2}+\mathcal{E}(F_0)+ \mathcal{E}(F_0)^\kappa \right]\\
	& \le  \|h_0\|_{L^p_vL^\infty_x}+C_p \left(\epsilon+\frac{C_{p,\gamma,\epsilon,T_0}}{N^{1/2}} \right) \left[\bar{\eta} + \bar{\eta}^2 + \bar{\eta}^{1-\kappa}+ \bar{\eta}^{3-\kappa} \right]\\
	& \quad +C_{p,\gamma,\epsilon,N,T_0} \left[\mathcal{E}(F_0)^{1/2}+\mathcal{E}(F_0)+ \mathcal{E}(F_0)^\kappa \right].
\end{align*}
Now we first choose suitably small $\epsilon>0$, then take $N>0$ large enough, and last choose sufficiently small $\epsilon_1$ with $\mathcal{E}(F_0) \le \epsilon_1$ so that
\begin{align*}
	C_p \left(\epsilon+\frac{C_{p,\gamma,\epsilon,T_0}}{N^{1/2}} \right) \left[\bar{\eta} + \bar{\eta}^2 + \bar{\eta}^{1-\kappa}+ \bar{\eta}^{3-\kappa} \right]
	 +C_{p,\gamma,\epsilon,N,T_0} \left[\mathcal{E}(F_0)^{1/2}+\mathcal{E}(F_0)+ \mathcal{E}(F_0)^\kappa \right] \le \frac{1}{4\tilde{C}_3}.
\end{align*}
We also choose small $\eta_0>0$ such that
\begin{align} \label{etacond10}
	\|h_0\|_{L^p_vL^\infty_x} \le \eta_0 \le \frac{1}{4\tilde{C}_3}.
\end{align}
Therefore, we conclude \eqref{aa61}, and we completes the proof of this lemma.
\end{proof}

\bigskip

\subsection{$L^2$ estimate}
In this subsection, our aim is to prove the $L^2$ estimate to the Boltzmann equation \eqref{FPBER}. To derive the exponential time-decay in $L^p_vL^\infty_x$ for the equation \eqref{WFPBE}, we need to consider the exponential $L^2$ decay. We define the $L^2_{v}$ projection $P$ of $f$ corresponding to operator $L$ as
\begin{align} \label{LPro}
	Pf(t,x,v) = a(t,x) \mu^{1/2}(v) + b(t,x)\cdot v\mu^{1/2}(v) + c(t,x)\frac{|v|^2-3}{\sqrt{6}}\mu^{1/2}(v),
\end{align}
where
\begin{align*}
	& a(t,x) = \int_{\mathbb{R}^3}f(t,x,v) \mu^{1/2}(v)dv,\\
	& b(t,x) = \int_{\mathbb{R}^3}vf(t,x,v) \mu^{1/2}(v)dv,\\
	& c(t,x) = \int_{\mathbb{R}^3}\frac{|v|^2-3}{\sqrt{6}}f(t,x,v) \mu^{1/2}(v)dv.
\end{align*}
It is well-known the operator $L$ satisfies the $L^2$ coercivity $(Lf,f)_{L_v^2} \ge C_L\|(I-P)f\|_{L_v^2}^2$ for all $f$ in $L_v^2$, where $C_L$ is a generic constant. We also define the $L^2_v$ projection $P_\gamma$ of $f$ on the boundary $\gamma$ as
\begin{align} \label{Pgamm}
	P_\gamma f = c_\mu \mu^{1/2}(v)\int_{n(x)\cdot v' >0} f(t,x,v')\mu^{1/2}(v')\{n(x)\cdot v'\}dv'.
\end{align}
\newline
\indent The following lemma states the $L^2_{x,v}$ bound for $Pf$ by $(I-P)f$ and the effects of the boundary. The lemma gives the key estimate to derive the exponential decay in $L^2_{x,v}$. The proof of this lemma is left in Section \ref{Appendix}.
\begin{Lem} \label{coer}
	Assume that $f\in L^2_{x,v}$ is a solution to the equation \eqref{FPBER} satisfying the mass conservation \eqref{massconserv}. There exists a function $G(t)$ such that for all $0\le s \le t$, $G(s) \lesssim \|f(s)\|_{L^2_{x,v}}^2$ and
	\begin{align*}
		\int_s^t \|Pf(\tau)\|_{L^2_{x,v}}^2d\tau &\lesssim G(t)-G(s) + \int_s^t\left[\|(I-P)f(\tau)\|_{L^2_{x,v}}^2+\left\|\left(I-P_\gamma\right)(f)(\tau)\right\|^2_{L_{\gamma_+}^2}\right]d\tau\\
		& \quad +\int_s^t \|\Gamma(f,f)(\tau)\|_{L^2_{x,v}}^2d\tau.
	\end{align*}
\end{Lem} 

\bigskip

Before proving Lemma \ref{nonlinear L^2 decay}, we need to estimate the time integration to the nonlinear term $\Gamma(f,f). $ Let $p>2$ and $\beta > \frac{3p-6}{2p}+2\gamma$. We claim that
\begin{align*}
	\int_s^t \|e^{\lambda_2 \tau}\Gamma(f,f)(\tau)\|_{L^2_{x,v}}^2d\tau \lesssim  \sup_{s \le \tau \le t}\|wf(\tau)\|_{L^p_vL^\infty_x}^2 \int_{s}^t \|e^{\lambda_2 \tau}f(\tau)\|_{L^2_{x,v}}^2 d\tau.
\end{align*}
For simplicity, we set $s=0$.
First, let us estimate
\begin{align*}
	\int_0^t \|e^{\lambda_2 \tau} \Gamma^-(f,f)(\tau)\|_{L^2_{x,v}}^2d\tau.
\end{align*}
By the definition of $\Gamma^-$, the above term is bounded by
\begin{align} \label{Gamma-L2}
	\int^t_0 \|e^{\lambda_2 \tau}\Gamma^-(f,f)(\tau)\|_{L^2_{x,v}}^2d\tau 
	&\lesssim \int_0^t e^{2\lambda_2 \tau} \int_{\Omega \times \R^3}|\nu(v) f(v)|^2 \left|\int_{\R^3\times \S^2} b(\cos\theta)\nu(u)\mu^{1/2}(u)|f(u)|d\omega du\right|^2\nonumber\\
	& \quad \times dxdv d\tau\nonumber\\
	& \lesssim \int_0^t e^{2\lambda_2 \tau} \int_{\Omega \times \R^3}|\nu(v) f(v)|^2 \left|\int_{\R^3}\nu(u)\mu^{1/2}(u)|f(u)| du\right|^2dxdvd\tau\nonumber\\
	& \lesssim \int_0^t e^{2\lambda_2 \tau} \|wf(\tau)\|_{L^p_v L^\infty_x}^2 \|f(\tau)\|_{L^2_{x,v}}^2d\tau\nonumber\\
	& \lesssim \sup_{0\le \tau \le t}\|wf(\tau)\|_{L^p_v L^\infty_x}^2\int_0^t \|e^{\lambda_2 \tau}f(\tau)\|_{L^2_{x,v}}^2 d\tau,
\end{align}
where $\beta>\frac{3p-6}{2p}+\gamma$.\\

Before we estimate the term $\Gamma^+$, we give some facts about the pre-post velocities. Let 
\begin{align*}
	\mathbf{k} = \begin{cases}
\frac{v-u}{|v-u|}, & \text{if }v\not=u  \\
(1,0,0), & \text{if }v=u.
\end{cases}
\end{align*}
The spherical coordinate gives
\begin{align*}
	\omega = \cos \varphi \mathbf{k}+\sin \varphi \rho,
\end{align*}
where $\varphi \in [0, \pi/2]$ and $\rho \in \S^1(\mathbf{k})=\{\rho \in \S^2 : \rho \cdot \mathbf{k}=0\}$.
Then we have
\begin{align*}
	u' = \sin^2 \frac{\varphi}{2} v+ \cos^2\frac{\varphi}{2}u-\frac{|v-u|}{2}\sin \varphi\rho,
\end{align*}
and
\begin{align*}
	|u'-v| = |v-u| \cos \frac{\varphi}{2}.
\end{align*}
Then it follows from $\varphi \in (0,\pi/2]$ that
\begin{align} \label{somefact1}
	|u|^2 \le (|u-v|+|v|)^2 \le \left(\frac{|u'-v|}{\cos\frac{\pi}{4}}+|v|\right)^2 \le (\sqrt{2}|u'|+(\sqrt{2}+1)|v|)^2 \lesssim |u'|^2+|v|^2.
\end{align}
Now, let us estimate the following time integration for the term $\Gamma^+$ :
\begin{align*}
	\int_0^t \|e^{\lambda_2 \tau}\Gamma^+(f,f)(\tau)\|_{L^2_{x,v}}^2 d\tau.
\end{align*}
By the definition for $\Gamma^+$, we have
\begin{align*}
	|\Gamma^+(f,f)(v)| &\le \int_{\R^3} \int_{\S^2}b(\cos \theta)|v-u|^\gamma \mu^{1/2}|f(v')||f(u')|d\omega du\\
	& \le \left(\int_{\R^3} \int_{\S^2}b(\cos \theta)|f(u')|^2d\omega du\right)^{1/2}\left(\int_{\R^3} \int_{\S^2}b(\cos \theta)|v-u|^{2\gamma}\mu(u) |f(u')|^2d\omega du\right)^{1/2}
\end{align*}
By making the change of variables $u \mapsto u'$ and using \eqref{Jacobianuu'}, we have
\begin{align*}
	\left(\int_{\R^3}\int_{\S^2}b(\cos \theta) \left|f(u')\right|^2 d\omega du\right)^{1/2} \le  \|f\|_{L^2_v},
\end{align*}
which implies that
\begin{align*}
	|\Gamma^+(f,f)(v)|  \lesssim \|f\|_{L^2_v}\left(\int_{\R^3} \int_{\S^2}b(\cos \theta)|v-u|^{2\gamma}\mu(u) |f(v')|^2d\omega du\right)^{1/2}.
\end{align*}
Taking the $L^2_{v}$ norm to the above inequality, we obtain
\begin{align} \label{ee1} 
	\|\Gamma^+(f,f)\|_{L^2_v} &\lesssim \|f\|_{L^2_v} \left(\int_{\R^3} \int_{\R^3} \int_{\S^2}b(\cos \theta)|v-u|^{2\gamma}\mu(u) |f(v')|^2d\omega dudv\right)^{1/2}\nonumber\\
	& \lesssim \|f\|_{L^2_v} \left(\int_{\R^3} \int_{\R^3} \int_{\S^2}b(\cos \theta)|v-u|^{2\gamma}\mu(u') |f(v)|^2d\omega dudv\right)^{1/2}\nonumber\\
	& \lesssim \|f\|_{L^2_v} \left(\int_{\R^3} \langle v \rangle^{2\gamma}|f(v)|^2 \int_{\R^3} \int_{\S^2}b(\cos \theta)\langle u \rangle^{2\gamma}\mu(u') d\omega dudv\right)^{1/2}
\end{align}
where we have used a change of variables.
Here, we can derive
\begin{align*}
	\int_{\R^3} \int_{\S^2}b(\cos \theta)\langle u \rangle^{2\gamma}\mu(u') d\omega du & \le  C_\gamma \int_{\R^3} \int_{\S^2}b(\cos \theta)\left(\langle u' \rangle^{2\gamma}+\langle v \rangle^{2\gamma} \right) \mu(u') d\omega du\\
	& \le C_\gamma \langle v \rangle^{2\gamma}\int_{\R^3} \int_{\S^2}b(\cos \theta)\langle u' \rangle^{2\gamma}  \mu(u') d\omega du,
\end{align*}
where we have used \eqref{somefact1}, and by making the change of variables $u \mapsto u'$ and using \eqref{Jacobianuu'}, it yields
\begin{align*}
	\int_{\R^3} \int_{\S^2}b(\cos \theta)\langle u \rangle^{2\gamma}\mu(u') d\omega du \le  C_\gamma \langle v \rangle^{2\gamma} .
\end{align*}
Thus \eqref{ee1} becomes
\begin{align} \label{betacondL2}
	\|\Gamma^+(f,f)\|_{L^2_v} &\le C_{\gamma} \|f\|_{L^2_v} \left(\int_{\R^3} \langle v \rangle^{4\gamma}|f(v)|^2dv\right)^{1/2}\nonumber\\
	& \le  C_{\gamma} \|f\|_{L^2_v} \left(\int_{\R^3} \left|\frac{\langle v \rangle^{4\gamma}}{w(v)^2} \right|^{\frac{p}{p-2}}dv\right)^{\frac{p-2}{2p}}\|wf\|_{L^p_vL^\infty_x}\\
	& \le  C_{p,\gamma,\beta} \|f\|_{L^2_v}\|wf\|_{L^p_vL^\infty_x} \nonumber
\end{align}
provided that $\beta>\frac{3p-6}{2p}+2\gamma$.
Hence we can derive
\begin{align} \label{Gamma+L2}
	\int^t_0 \|e^{\lambda_2 \tau}\Gamma^+(f,f)(\tau)\|_{L^2_{x,v}}^2d\tau 
	& \le C_{p, \gamma,\beta}\int_0^t e^{2\lambda_2 \tau} \|f(\tau)\|_{L^2_{x,v}}^2 \|wf(\tau)\|^2_{L^p_vL^\infty_x}dxdv d\tau\nonumber\\
	& \le  C_{p, \gamma ,\beta}\sup_{0\le \tau \le t}\|wf(\tau)\|_{L^p_v L^\infty_x}^2\int_0^t \|e^{\lambda_2 \tau}f(\tau)\|_{L^2_{x,v}}^2 d\tau.
\end{align}

\bigskip

\begin{Lem} \label{nonlinear L^2 decay}
	Let $p>2$ and $\beta>\frac{3p-6}{2p}+2\gamma$.
Suppose that $f$ solves the equation \eqref{FPBE} satisfying the mass conservation \eqref{massconserv}. Then there are $\lambda_2>0$ and $C_2>0$ with $\lambda_2 \le \nu_0/16$ such that for $0\le s \le t$,
	\begin{align*}
		\|e^{\lambda_2 t}f(t)\|_{L^2_{x,v}}^2+ \int_{s}^t \|e^{\lambda_2 \tau}f(\tau)\|_{L^2_{x,v}}^2 d\tau &\le C_2\|e^{\lambda_2 s}f(s)\|_{L^2_{x,v}}^2\\
		& \quad  + C_2\sup_{s \le \tau \le t}\|wf(\tau)\|_{L^p_vL^\infty_x}^2 \int_{s}^t \|e^{\lambda_2 \tau}f(\tau)\|_{L^2_{x,v}}^2 d\tau.
	\end{align*}	
\end{Lem}
\begin{remark}
	$C_2$ depends only on $p$, $\beta$, and $\gamma$.
\end{remark}
\begin{proof}
	For simplicity, we set $s=0$. Multiplying $e^{2\lambda_2 t} f$ by \eqref{FPBER}, we can derive the $L^2$ estimate:
	\begin{align*}
		&\|e^{\lambda_2 t} f(t)\|_{L^2_{x,v}}^2-\|f(0)\|_{L^2_{x,v}}^2+\int_0^t \|e^{\lambda_2 \tau}(I-P_\gamma)f(\tau)\|_{L^2_{\gamma_+}}^2d\tau + 2\int_0^t \int_{\Omega}\int_{\R^3} e^{2\lambda_2 \tau}f(\tau) Lf(\tau)dxdv d\tau\\
		& = 2 \int_0^t \int_{\Omega}\int_{\R^3} e^{2\lambda_2 \tau}f(\tau) \Gamma(f,f)(\tau)dxdv d\tau + \lambda_2 \int_0^t \|e^{\lambda_2 \tau} f(\tau)\|_{L^2_{x,v}}^2d\tau.
	\end{align*}
		From the coercivity $\left(L f,f\right)_{L_v^2} \ge C_L \left\|\left(I-P\right)f\right\|_{L^2_v}^2$, we deduce
	\begin{align} \label{L2dd1}
		&\|e^{\lambda_2 t} f(t)\|_{L^2_{x,v}}^2-\|f(0)\|_{L^2_{x,v}}^2+\int_0^t \|e^{\lambda_2 \tau}(I-P_\gamma)f(\tau)\|_{L^2_{\gamma_+}}^2d\tau\nonumber\\
		 &\le  -2C_L \int_0^t\|e^{\lambda_2 \tau}\left(I-P\right)f(\tau)\|_{L^2_{x,v}}^2d\tau+ 2  \int_0^t \int_{\Omega}\int_{\R^3} e^{2\lambda_2 \tau}f(\tau) \Gamma(f,f)(\tau)dxdv d\tau \nonumber\\
		& \quad + \lambda_2 \int_0^t \|e^{\lambda_2 \tau} f(\tau)\|_{L^2_{x,v}}^2d\tau\nonumber\\
		&\le  -2C_L \int_0^t\|e^{\lambda_2 \tau}\left(I-P\right)f(\tau)\|_{L^2_{x,v}}^2d\tau+  \epsilon \int_0^t \|e^{\lambda_2 \tau}f(\tau)\|_{L^2_{x,v}}^2d\tau  \nonumber\\
		& \quad +C(\epsilon)\int_0^t  \|e^{\lambda_2 \tau}\Gamma(f,f)(\tau)\|_	{L^2_{x,v}}^2 d\tau+ \lambda_2 \int_0^t \|e^{\lambda_2 \tau} f(\tau)\|_{L^2_{x,v}}^2d\tau,
	\end{align}
	where we have used the Young's inequality and $\epsilon$ can be small enough.
	Applying Lemma \ref{coer} to $e^{\lambda_2t}f$, we have
	\begin{align}\label{L2dd2}
		\int_0^t\|e^{\lambda_2 \tau}Pf(\tau)\|^2_{L^2_{x,v}}d\tau &\le C\|e^{\lambda_2 t}f(t)\|^2_{L^2_{x,v}} + C\|f(0)\|^2_{L^2_{x,v}}\nonumber\\
		& \quad +C\int_0^t\left[\|e^{\lambda_2 \tau}\left(I-P\right)f(\tau)\|^2_{L^2_{x,v}}+\|e^{\lambda_2 \tau}(I-P_\gamma)f(\tau)\|_{L^2_{\gamma_+}}^2\right]d\tau\nonumber\\
		&\quad +C\lambda_2^2 \int_0^t\|e^{\lambda_2 \tau}f(\tau)\|^2_{L^2_{x,v}}ds+C\int_0^t \|e^{\lambda_2 \tau}\Gamma(f,f)(\tau)\|_{L^2_{x,v}}^2d\tau,
	\end{align}
	where $C$ is a constant.
	For $\delta >0$, \eqref{L2dd1} + $\delta \times$\eqref{L2dd2} yields
	\begin{align*}
		&(1-\delta C) \|e^{\lambda_2 t}f(t)\|^2_{L^2_{x,v}}+ (2C_L-\delta C) \int_0^t\|e^{\lambda_2 \tau}\left(I-P\right)f(\tau)\|_{L^2_{x,v}}^2d\tau+ \delta \int_0^t\|e^{\lambda_2 \tau}Pf(\tau)\|^2_{L^2_{x,v}}d\tau\\ 
		&+ (1-\delta C)\int_0^t\|e^{\lambda_2 \tau}(I-P_\gamma)f(\tau)\|_{L^2_{\gamma_+}}^2d\tau\\ 
		&\le (1+\delta C) \|f(0)\|^2_{L^2_{x,v}} + (\epsilon +\lambda_2 + \delta C\lambda_2^2)\int_0^t \|e^{\lambda_2 \tau} f(\tau)\|_{L^2_{x,v}}^2d\tau\\
		& \quad + (\delta C+C(\epsilon)) \int_0^t \|e^{\lambda_2 \tau}\Gamma(f,f)(\tau)\|_{L^2_{x,v}}^2d\tau.
	\end{align*}
	Firstly, choosing $\delta >0$ sufficiently small  such that $2C_L-\delta C>\delta$, we get
	\begin{align*}
		&(1-\delta C) \|e^{\lambda_2 t}f(t)\|^2_{L^2_{x,v}}+ \delta \int_0^t\|e^{\lambda_2 \tau}f(\tau)\|^2_{L^2_{x,v}}d\tau+ (1-\delta C)\int_0^t\|e^{\lambda_2 \tau}(I-P_\gamma)f(\tau)\|_{L^2_{\gamma_+}}^2d\tau\\ 
		&\le (1+\delta C) \|f(0)\|^2_{L^2_{x,v}} + (\epsilon +\lambda_2 + \delta C\lambda_2^2)\int_0^t \|e^{\lambda_2 \tau} f(\tau)\|_{L^2_{x,v}}^2d\tau\\
		& \quad + (\delta C+C(\epsilon)) \int_0^t \|e^{\lambda_2 \tau}\Gamma(f,f)(\tau)\|_{L^2_{x,v}}^2d\tau.
	\end{align*}
	Next, choosing $\lambda_2>0$ and $\epsilon>0$ small enough such that 
	\begin{align*} 
		\delta/2 >\epsilon+ \lambda_2 + \delta C \lambda_2^2,
	\end{align*}
	we obtain
	\begin{align*}
		\|e^{\lambda_2 t}f(t)\|^2_{L^2_{x,v}}+\int_0^t\|e^{\lambda_2 \tau}f(\tau)\|^2_{L^2_{x,v}}d\tau \lesssim  \|f(0)\|^2_{L^2_{x,v}}+\sup_{0\le \tau \le t}\|wf(\tau)\|_{L^p_v L^\infty_x}^2\int_0^t \|e^{\lambda_2 \tau}f(\tau)\|_{L^2_{x,v}}^2 d\tau,
	\end{align*}
	where we have used the estimates \eqref{Gamma-L2} and \eqref{Gamma+L2}.
\end{proof}

\bigskip

\begin{remark}
	If $\sup_{s \le \tau \le t}\|wf(\tau)\|_{L^p_vL^\infty_x}^2$ is sufficiently small due to the a priori assumption \eqref{Aprioriassump1}, from Lemma \ref{nonlinear L^2 decay}, we obtain for $0 \le s \le t$,
\begin{align*}
	\|e^{\lambda_2 t}f(t)\|_{L^2_{x,v}}^2+ \int_{s}^t \|e^{\lambda_2 \tau}f(\tau)\|_{L^2_{x,v}(\nu)}^2 d\tau\lesssim \|e^{\lambda_2 s}f(s)\|_{L^2_{x,v}}^2,
\end{align*}
and thus we can derive the $L^2$ decay for $f$:
\begin{align*}
	\|e^{\lambda_2 t}f(t)\|_{L^2_{x,v}}\lesssim \|f(0)\|_{L^2_{x,v}}
\end{align*}
for all $t \ge 0$.
\end{remark}

\bigskip

\subsection{Nonlinear $L^p_vL^\infty_x$ estimate} Setting $h=wf$, the equation \eqref{FPBE} can be written as
\begin{align}\label{WFPBE2}
	\partial_t h + v \cdot \nabla_x h + R(f) h=K_wh + w\Gamma^+\left(f ,f \right).
\end{align}
Let $S_{G_f}(t,s)$ be the solution operator to the equation \eqref{WFPBE2} with the diffuse reflection boundary condition \eqref{WPDRBC}. 
Recall that we fix $p$, $\beta$, $\mathcal{T}$ in subsection \ref{apriorisamll}. To show the global-in-time existence to the equation \eqref{FPBER}, we need to deal with the nonlinear $L^p_vL^\infty_x$ estimate for the equation \eqref{WFPBE2}. This subsection is devoted to deriving the nonlinear $L^p_vL^\infty_x$ estimate.
\begin{Lem} \label{LALinftyestimate1} 
	 Let $(t,x,v)\in (0,T_0]\times \Omega \times \R^3$. Assume $\mathcal{E}(F_0) \le \epsilon_1 = \epsilon_1(\bar{\eta},T_0)$, where $\epsilon_1$ is determined in Lemma \ref{Rfest1}. Under the a priori assumption \eqref{Aprioriassump1}, it holds that
\begin{align*}
	|h(t,x,v)| &\le \left|(S_{G_f}(t,0)h_0)(t,x,v) \right|+\int_{t-\epsilon}^{t}\left|(S_{G_f}(t,s)K_wh(s))(t,x,v) \right|ds\\
	& \quad +\int_{t-\epsilon}^t\left|(S_{G_f}(t,s)w\Gamma^+(f,f)(s))(t,x,v) \right|ds\\
	& \quad + \int_0^{t-\epsilon} \mathbf{1}_{\{t_1\le s\}}e^{-\frac{\nu_0}{2}(t-s)} \left|(K_wh)(s,x-v(t-s),v) \right|ds\\
	& \quad + \int_0^{t-\epsilon} \mathbf{1}_{\{t_1\le s\}}e^{-\frac{\nu_0}{2}(t-s)} \left|w\Gamma^+(f,f)(s,x-v(t-s),v) \right|ds\\
	& \quad +\frac{1}{\tilde{w}(v)}\left(C_p\epsilon + \frac{C_{p,\epsilon,T_0}}{N^{1/2}} \right)\sup_{0\le s\le t}\left[\|h(s)\|_{L^p_vL^\infty_x}+\|h(s)\|_{L^p_vL^\infty_x}^2+\|h(s)\|_{L^p_vL^\infty_x}^{\frac{p}{2p-2}}+\|h(s)\|_{L^p_vL^\infty_x}^{\frac{5p-4}{2p-2}}\right]\\
	& \quad +\frac{C_{p,N,\Omega,T_0}}{\tilde{w}(v)}\left[\mathcal{E}(F_0)^{1/2}+\mathcal{E}(F_0)+\mathcal{E}(F_0)^{\frac{p-2}{2p-2}}\right],
\end{align*}	
where $\epsilon>0$ can be arbitrary small and $N>0$ can be arbitrary large.
\end{Lem}
\begin{proof}
	Let $(t,x,v) \in (0, T_0] \times \Omega \times \R^3$. We apply the Duhamel principle to the weighted Boltzmann equation \eqref{WFPBE2}:
	\begin{equation} \label{LA10}
	\begin{aligned}
			h(t,x,v) &= S_{G_f}(t,0) h_0+ \int_0^t \left[S_{G_f}(t,s)K_wh(s)+S_{G_f}(t,s)w\Gamma^+(f,f)(s)\right] ds\\
		& = S_{G_f}(t,0) h_0+ \int_{t-\epsilon}^t \left[S_{G_f}(t,s)K_wh(s)+S_{G_f}(t,s)w\Gamma^+(f,f)(s)\right] ds\\
		& \quad +  \int_0^{t-\epsilon} S_{G_f}(t,s)K_wh(s)ds+\int_0^{t-\epsilon}S_{G_f}(t,s)w\Gamma^+(f,f)(s) ds,
	\end{aligned}
	\end{equation}
	where $\epsilon>0$ is to be fixed sufficiently small later.\\
	Now, let us estimate the last term in \eqref{LA10}:\
	\begin{align*}
		 \int_0^{t-\epsilon} S_{G_f}(t,s)K_wh(s)ds+\int_0^{t-\epsilon}S_{G_f}(t,s)w\Gamma^+(f,f)(s) ds.
	\end{align*}
	First of all, let us estimate
	\begin{align*}
		\int_0^{t-\epsilon} S_{G_f}(t,s)K_wh(s)ds.
	\end{align*}
	From Lemma \ref{RepresentationforDiffuse}, we have
	\begin{equation} \label{LA18}
	\begin{aligned}
		&\left(S_{G_f}(t,s)K_wh(s)\right)(t,x,v)\\
		& =\mathbf{1}_{\{t_1\le s\}} I^f(t,s)(K_wh)(s,x-v(t-s),v)\\
		&\quad +\frac{I^f(t,t_1)}{\tilde{w}(v)} \sum_{l=1}^{k-1} \int_{\prod_{j=1}^{k-1}\mathcal{V}_j}\mathbf{1}_{\{t_{l+1}\le s <t_l\}} (K_wh)(s,x_l-v_l(t_l-s),v_l)d\Sigma_l^f(s)\\
		&\quad +\frac{I^f(t,t_1)}{\tilde{w}(v)} \int_{\prod_{j=1}^{k-1}\mathcal{V}_j} \mathbf{1}_{\{t_k>s\}}(S_{G_f}(t,s)K_wh(s))(t_k,x_k,v_{k-1})d\Sigma_{k-1}^f(t_k)\\
		& =: I_{1}+I_{2}+I_{3}.
\end{aligned}
\end{equation}
Note that from the a priori assumption \eqref{Aprioriassump1} and Lemma \ref{Rfest1} , we have
\begin{align*}
		R(f)(t,x,v) \ge \frac{1}{2}\nu(v)
	\end{align*}
	for all $(t,x,v) \in [0,T_0] \times \Omega \times \R^3$. Thus it holds that
	\begin{align} \label{Ibound}
		I^f(t,s) \le e^{-\frac{\nu_0}{2}(t-s)}
	\end{align}
	for all $0\le s \le t \le T_0$.\\
	For $I_3$, from \eqref{DampLpdecay3} and Lemma \ref{KESTIMATE}, we get
	\begin{align} \label{LA11}
		\|(S_{G_f}(t,s)K_wh(s))(t_k,x_k,v_{k-1})\|_{L^p_vL^\infty_x} \le C_{\beta,p} e^{-\frac{\nu_0}{4}(t_k-s)} \|h(s)\|_{L^p_vL^\infty_x}.
	\end{align}
	Since $t-s\ge \epsilon>0$, by Lemma \ref{Lsmall}, it follows from \eqref{LA11} that, for $k=k(\epsilon,T_0)+1$ such that \eqref{Lsmallineq} holds, 
	\begin{align*}
		|I_3| &\le C_{\beta,p}\frac{e^{-\frac{\nu_0}{4}(t-s)}}{\tilde{w}(v)} \|h(s)\|_{L^p_vL^\infty_x} \left(\int_{\prod_{j=1}^{k-2}\mathcal{V}_j} \mathbf{1}_{\{t_{k-1}>s\}} \prod_{j=1}^{k-2}d\sigma_j\right)\\
		&\le \epsilon C_{\beta,p}\frac{e^{-\frac{\nu_0}{4}(t-s)}}{\tilde{w}(v)} \|h(s)\|_{L^p_vL^\infty_x},
	\end{align*}
	where we have used the H\"older inequality and 
	\begin{align*}
		\left(\int_{\R^3} \left|\mu^{1/2}(v_{k-1}) \frac{|v_{k-1}|}{w(v_{k-1})} \right|^{p'} dv_{k-1}\right)^{1/p'} \le C_p.
	\end{align*}
	Thus we can deduce that
	\begin{align} \label{LA16}
		\int_0^{t-\epsilon} |I_3| ds \le \frac{\epsilon C_{\beta,p}}{\tilde{w}(v)} \sup_{0\le s \le t}\|h(s)\|_{L^p_vL^\infty_x}.
	\end{align}
	For $I_2$, fix $l$ and we denote each term $l=1,2,\cdots, k-1$ in the summation by $I_{2l}$. Then from \eqref{Ibound}, we have
	\begin{align*}
		\int_0^{t-\epsilon} |I_{2l}| ds \le \int_{t_{l+1}}^{t_l} \frac{e^{-\frac{\nu_0}{2}(t-s)}}{\tilde{w}(v)}  \int_{\prod_{j=1}^{l-1}\mathcal{V}_j} \int_{\mathcal{V}_{l}} \tilde{w}(v_l)\int_{\R^3}|k_w(v_l,v')||h(s,x_l-v_l(t_l-s),v')|dv'd\sigma_l\prod_{j=1}^{l-1}d\sigma_j ds.
	\end{align*}
	We divide this term into four cases.\\
	\newline
	$\mathbf{Case\ 1:}$ $|v_l| \ge N$.\\
	We use the H\"older inequality and Lemma \ref{KESTIMATE} to estimate
\begin{align} \label{LA12}
	\int_0^{t-\epsilon}|I_{2l}^1| ds 
	&\le \int_{t_{l+1}}^{t_l} \frac{e^{-\frac{\nu_0}{2}(t-s)}}{\tilde{w}(v)}\int_{\prod_{j=1}^{l-1}\mathcal{V}_j} \int_{|v_l|\ge N} \tilde{w}(v_l)\int_{\R^3}|k_w(v_l,v')||h(s,x_l-v_l(t_l-s),v')|dv'd\sigma_l\prod_{j=1}^{l-1}d\sigma_j ds\nonumber\\
	&\le C_p \int_{t_{l+1}}^{t_l} \frac{e^{-\frac{\nu_0}{2}(t-s)}}{\tilde{w}(v)}\|h(s)\|_{L^p_vL^\infty_x}\int_{\prod_{j=1}^{l-1}\mathcal{V}_j} \int_{|v_l|\ge N} \tilde{w}(v_l) d\sigma_l\prod_{j=1}^{l-1}d\sigma_j ds\nonumber\\
	&\le \frac{C_p}{N}\int_{0}^{t} \frac{e^{-\frac{\nu_0}{2}(t-s)}}{\tilde{w}(v)}\|h(s)\|_{L^p_vL^\infty_x}ds\nonumber \\
	&\le \frac{C_p}{N}\frac{1}{{\tilde{w}(v)}}\sup_{0\le s \le t} \|h(s)\|_{L^p_vL^\infty_x}.
\end{align}
	\newline
	$\mathbf{Case\ 2:}$ $|v_l| \le N$, $|v'| \ge 2N$.\\
	The following is valid:
\begin{align} \label{LA41}
	|k_w(v_l,v')| ^{p'}\le e^{-\frac{\epsilon}{8}N^2}|k_w(v_l,v')|^{p'}e^{\frac{\epsilon}{8}|v_l-v'|^2},
\end{align}
where $\epsilon$ is sufficiently small.
From Lemma \ref{KESTIMATE}, it holds that
\begin{align} \label{LA42}
	 \int_{|v'|\ge 2N}|k_w(v_l,v')|^{p'}e^{\frac{\epsilon}{8}|v_l-v'|^2}dv'' \le C_{\beta,p} \left(\frac{1}{1+|v_l|}\right)^{p'+1},
\end{align}
where $p'$ is the exponent conjugate of $p$.
Then we use the H\"older inequality to obtain
\begin{align} \label{LA13}
	\int_0^{t-\epsilon}|I_{2l}^2| ds 
	&\le \int_{t_{l+1}}^{t_l} \frac{e^{-\frac{\nu_0}{2}(t-s)}}{\tilde{w}(v)}\int_{\prod_{j=1}^{l-1}\mathcal{V}_j} \int_{|v_l|\le N} \tilde{w}(v_l)\int_{|v'| \ge 2N}|k_w(v_l,v')||h(s,x_l-v_l(t_l-s),v')|dv'\nonumber\\
	& \quad \times d\sigma_l\prod_{j=1}^{l-1}d\sigma_j ds\nonumber\\
	&\le C_pe^{-\frac{\epsilon N^2}{8}} \int_{t_{l+1}}^{t_l} \frac{e^{-\frac{\nu_0}{2}(t-s)}}{\tilde{w}(v)}\|h(s)\|_{L^p_vL^\infty_x}\int_{\prod_{j=1}^{l-1}\mathcal{V}_j} \int_{|v_l|\le N} \tilde{w}(v_l) d\sigma_l\prod_{j=1}^{l-1}d\sigma_j ds\nonumber\\
	&\le C_pe^{-\frac{\epsilon N^2}{8}}\int_{0}^{t} \frac{e^{-\frac{\nu_0}{2}(t-s)}}{\tilde{w}(v)}\|h(s)\|_{L^p_vL^\infty_x}ds \nonumber\\
	&\le C_pe^{-\frac{\epsilon N^2}{8}}\frac{1}{{\tilde{w}(v)}}\sup_{0\le s \le t} \|h(s)\|_{L^p_vL^\infty_x}.
\end{align}
\newline
$\mathbf{Case\ 3:}$ $t_l-s \le \frac{1}{N}$.\\
We can easily compute
\begin{align} \label{LA14}
	\int_0^{t-\epsilon}|I_{2l}^3| ds 
	&\le \int_{t_{l}-\frac{1}{N}}^{t_l} \frac{e^{-\frac{\nu_0}{2}(t-s)}}{\tilde{w}(v)}\int_{\prod_{j=1}^{l-1}\mathcal{V}_j} \int_{\mathcal{V}_l} \tilde{w}(v_l)\int_{\R^3}|k_w(v_l,v')||h(s,x_l-v_l(t_l-s),v')|dv'\nonumber\\
	&\quad  \times d\sigma_l\prod_{j=1}^{l-1}d\sigma_j ds\nonumber\\
	&\le C_p \int_{t_l-\frac{1}{N}}^{t_l} \frac{1}{\tilde{w}(v)}\|h(s)\|_{L^p_vL^\infty_x}ds\nonumber\\
	&\le \frac{C_p}{N}\frac{1}{{\tilde{w}(v)}}\sup_{0\le s \le t} \|h(s)\|_{L^p_vL^\infty_x}.
\end{align}
	\newline
	$\mathbf{Case\ 4:}$ $t_l-s \ge \frac{1}{N} \text{ and }|v_l|\le N$, $|v'| \le 2N$.\\
We use the Cauchy-Schwartz inequality and Lemma \ref{KESTIMATE} to estimate
\begin{align*}
	\int_0^{t-\epsilon}|I_{2l}^4| ds &\le \int_{t_{l+1}}^{t_l-\frac{1}{N}} \frac{e^{-\frac{\nu_0}{2}(t-s)}}{\tilde{w}(v)}\int_{\prod_{j=1}^{l-1}\mathcal{V}_j} \int_{|v_l| \le N} \tilde{w}(v_l)\int_{|v'|\le 2N}|k_w(v_l,v')||h(s,x_l-v_l(t_l-s),v')|dv'\\
	& \quad \times d\sigma_l\prod_{j=1}^{l-1}d\sigma_j ds\\
	& \le \int_{t_{l+1}}^{t_l-\frac{1}{N}} \frac{e^{-\frac{\nu_0}{2}(t-s)}}{\tilde{w}(v)}\int_{\prod_{j=1}^{l-1}\mathcal{V}_j} \int_{|v_l| \le N} \tilde{w}(v_l)\left(\int_{|v'|\le 2N}|h(s,x_l-v_l(t_l-s),v')|^2dv'\right)^{1/2}\nonumber\\
	& \quad \times  d\sigma_l\prod_{j=1}^{l-1}d\sigma_j ds\\
	& \le C_N\int_{t_{l+1}}^{t_l-\frac{1}{N}} \frac{e^{-\frac{\nu_0}{2}(t-s)}}{\tilde{w}(v)}\int_{\prod_{j=1}^{l-1}\mathcal{V}_j} \left(\int_{|v_l| \le N} \int_{|v'|\le 2N}|h(s,x_l-v_l(t_l-s),v')|^2dv'dv_l\right)^{1/2}\nonumber \\
	& \quad \times  \prod_{j=1}^{l-1}d\sigma_j ds
\end{align*}
where we have used
\begin{align*}
	\int_{\R^3}|k_w(v_l,v')|^2dv' \le C_\beta. 
\end{align*}
We make a change of variables $v_l \mapsto y=x_l-v_l(t_l-s)$ with $\left|\frac{dy}{dv_l}\right| = (t_l-s)^3$:
\begin{align*}
	\int_0^{t-\epsilon}|I_{2l}^4| ds \le C_N\int_{t_{l+1}}^{t_l-\frac{1}{N}} \frac{e^{-\frac{\nu_0}{2}(t-s)}}{\tilde{w}(v)}\int_{\prod_{j=1}^{l-1}\mathcal{V}_j} \left(\int_{\Omega}\int_{|v'|\le 2N}|h(s,y,v')|^2dv'dy\right)^{1/2}\prod_{j=1}^{l-1}d\sigma_j ds  .   
\end{align*}
It follows from \eqref{RTIRE} that
\begin{align} \label{LA15}
	\int_0^{t-\epsilon}|I_{2l}^4| ds &\le C_{p,N,\Omega} \int_0^t \frac{e^{-\nu_0(t-s)}}{\tilde{w}(v)}\left[\mathcal{E}(F_0)^{1/2} + \|h(s)\|_{L^p_vL^\infty_x}^{\frac{1-\kappa}{2}}\mathcal{E}(F_0)^{\kappa/2}\right]ds\nonumber\\
	& \le \frac{C_{p,N,\Omega}}{\tilde{w}(v)}\left[\mathcal{E}(F_0)^{1/2} + \mathcal{E}(F_0)^{\kappa/2}\sup_{0\le s\le t}\|h(s)\|_{L^p_vL^\infty_x}^{\frac{1-\kappa}{2}}\right]\nonumber\\
	& \le \frac{1}{\tilde{w}(v)}\left[ \frac{C_p}{N}\sup_{0\le s\le t}\|h(s)\|_{L^p_vL^\infty_x}^{1-\kappa}+ C_{p,N,\Omega}\left[\mathcal{E}(F_0)^{1/2}+ \mathcal{E}(F_0)^{\kappa}\right]\right],
\end{align}
where we have used the Young's inequality.\\
Gathering \eqref{LA12}, \eqref{LA13}, \eqref{LA14}, \eqref{LA15} and summing over $1\le l \le k(\epsilon,T_0)-1$, we obtain
\begin{equation} \label{LA17}
\begin{aligned}
	\int_0^{t-\epsilon} |I_2| ds &\le  \frac{C_{p,\epsilon,T_0}}{N}\frac{1}{{\tilde{w}(v)}}\sup_{0\le s \le t}\left[ \|h(s)\|_{L^p_vL^\infty_x}+\|h(s)\|_{L^p_vL^\infty_x}^{1-\kappa}\right]\\
	& \quad +\frac{1}{\tilde{w}(v)}C_{p,N,\Omega,\epsilon,T_0}\left[\mathcal{E}(F_0)^{1/2}+ \mathcal{E}(F_0)^{\kappa}\right].
\end{aligned}
\end{equation}
Combining \eqref{LA18}, \eqref{LA16}, \eqref{LA17}, we deduce that
\begin{equation} \label{LA31}
	\begin{aligned}
	\int_0^{t-\epsilon} |S_{G_f}(t,s)K_wh(s)|ds &\le \int_0^{t-\epsilon}\mathbf{1}_{\{t_1\le s\}} e^{-\frac{\nu_0}{2}(t-s)}\left|(K_wh)(s,x-v(t-s),v)\right|ds \\
	& \quad + \frac{C_p}{\tilde{w}(v)}\left(\epsilon+\frac{C_{p,\epsilon,T_0}}{N} \right)\sup_{0\le s \le t}\left[ \|h(s)\|_{L^p_vL^\infty_x}+\|h(s)\|_{L^p_vL^\infty_x}^{1-\kappa}\right]\\
	& \quad +\frac{C_{p,N,\Omega,\epsilon,T_0}}{\tilde{w}(v)}\left[\mathcal{E}(F_0)^{1/2}+ \mathcal{E}(F_0)^{\kappa}\right].
\end{aligned}
\end{equation}
Next, let us estimate
\begin{align*}
	\int_0^{t-\epsilon}S_{G_f}(t,s)w\Gamma^+(f,f)(s) ds.
\end{align*}
From Lemma \ref{RepresentationforDiffuse}, we have
\begin{equation} \label{LA30}
	\begin{aligned}
		&\left(S_{G_f}(t,s)w\Gamma^+(f,f)(s)\right)(t,x,v)\\
		& =\mathbf{1}_{\{t_1\le s\}} I^f(t,s)w\Gamma^+(f,f)(s,x-v(t-s),v)\\
		&\quad +\frac{I^f(t,t_1)}{\tilde{w}(v)} \sum_{l=1}^{k-1} \int_{\prod_{j=1}^{k-1}\mathcal{V}_j}\mathbf{1}_{\{t_{l+1}\le s <t_l\}} w\Gamma^+(f,f)(s,x_l-v_l(t_l-s),v_l)d\Sigma_l^f(s)\\
		&\quad +\frac{I^f(t,t_1)}{\tilde{w}(v)} \int_{\prod_{j=1}^{k-1}\mathcal{V}_j} \mathbf{1}_{\{t_k>s\}}(S_{G_f}(t,s)w\Gamma^+(f,f)(s))(t_k,x_k,v_{k-1})d\Sigma_{k-1}^f(t_k)\\
		& =: I_{4}+I_{5}+I_{6}.
\end{aligned}
\end{equation}
	For $I_6$, from \eqref{DampLpdecay3} and Corollary \ref{LpGamma+est}, we get
	\begin{align} \label{LA19}
		\|(S_{G_f}(t,s)w\Gamma^+(f,f)(s))(t_k,x_k,v_{k-1})\|_{L^p_vL^\infty_x} \le C_{\beta,p} e^{-\frac{\nu_0}{4}(t_k-s)} \|h(s)\|_{L^p_vL^\infty_x}^2.
	\end{align}
	Since $t-s\ge \epsilon>0$, by Lemma \ref{Lsmall}, it follows from \eqref{LA19} that, for $k=k(\epsilon,T_0)+1$ such that \eqref{Lsmallineq} holds, 
	\begin{align*}
		|I_6| &\le C_{\beta,p}\frac{e^{-\frac{\nu_0}{2}(t-s)}}{\tilde{w}(v)} \|h(s)\|_{L^p_vL^\infty_x}^2 \left(\int_{\prod_{j=1}^{k-2}\mathcal{V}_j} \mathbf{1}_{\{t_{k-1}>s\}} \prod_{j=1}^{k-2}d\sigma_j\right)\\
		&\le \epsilon C_{\beta,p}\frac{e^{-\frac{\nu_0}{2}(t-s)}}{\tilde{w}(v)} \|h(s)\|_{L^p_vL^\infty_x}^2,
	\end{align*}
	where we have used the H\"older inequality and 
	\begin{align*}
		\left(\int_{\R^3} \left|\mu^{1/2}(v_{k-1}) \frac{|v_{k-1}|}{w(v_{k-1})} \right|^{p'} dv_{k-1}\right)^{1/p'} \le C_p.
	\end{align*}
	Thus we can deduce that
	\begin{align} \label{LA20}
		\int_0^{t-\epsilon} |I_6| ds \le \frac{\epsilon C_{\beta,p}}{\tilde{w}(v)} \sup_{0\le s \le t}\|h(s)\|_{L^p_vL^\infty_x}^2.
	\end{align}
	For $I_5$, fix $l$ and we denote each term $l=1,2,\cdots, k-1$ in the summation by $I_{5l}$. Then we have
	\begin{align*}
		\int_0^{t-\epsilon} |I_{5l}| ds &\le \int_{t_{l+1}}^{t_l} \frac{e^{-\frac{\nu_0}{2}(t-s)}}{\tilde{w}(v)}  \int_{\prod_{j=1}^{l-1}\mathcal{V}_j} \int_{\mathcal{V}_{l}}  \left|w\Gamma^+(f,f)(s,x_l-v_l(t_l-s),v_l)\right| \tilde{w}(v_l)d\sigma_l\prod_{j=1}^{l-1}d\sigma_j ds\
	\end{align*}
	First, let us estimate $I_{5l}$. By \eqref{Gamma+estimatecombine}, the term $I_{5l}$ is bounded by
	\begin{align*}
		&C_{p,\gamma}\int_{t_{l+1}}^{t_l} \frac{e^{-\frac{\nu_0}{2}(t-s)}}{\tilde{w}(v)} \|h(s)\|_{L^p_vL^\infty_x} \int_{\prod_{j=1}^{l-1}\mathcal{V}_j} \int_{\mathcal{V}_{l}}  \left(\int_{\R^3}(1+|\eta|)^{e(p,\gamma)}|f(s,x_l-v_l(t_l-s),\eta)|^2 d\eta \right)^{1/2} \tilde{w}(v_l)d\sigma_l\\
		&\quad  \times \prod_{j=1}^{l-1}d\sigma_j ds.
	\end{align*}
	We divide this term into four cases.\\
	\newline
	$\mathbf{Case\ 1:}$ $|v_l| \ge N$.\\
	We use the H\"older inequality to estimate
\begin{align} \label{LA21}
	|I_{5l}^1|
	&\le C_{p,\gamma}\int_{t_{l+1}}^{t_l} \frac{e^{-\frac{\nu_0}{2}(t-s)}}{\tilde{w}(v)} \|h(s)\|_{L^p_vL^\infty_x} \int_{\prod_{j=1}^{l-1}\mathcal{V}_j} \int_{|v_l| \ge N}  \left(\int_{\R^3}(1+|\eta|)^{e(p,\gamma)}|f(s,x_l-v_l(t_l-s),\eta)|^2 d\eta \right)^{1/2} \nonumber \\
	& \quad \times \tilde{w}(v_l)d\sigma_l\prod_{j=1}^{l-1}d\sigma_j ds\nonumber\\
	&\le C_{p,\gamma}\int_{t_{l+1}}^{t_l} \frac{e^{-\frac{\nu_0}{2}(t-s)}}{\tilde{w}(v)} \|h(s)\|_{L^p_vL^\infty_x}^2 \int_{\prod_{j=1}^{l-1}\mathcal{V}_j} \int_{|v_l| \ge N}  \left(\int_{\R^3}(1+|\eta|)^{\left(e(p,\gamma)-2\beta\right)\cdot\frac{p}{p-2}} d\eta \right)^{\frac{p-2}{2p}}\tilde{w}(v_l)d\sigma_l    \nonumber\\
	& \quad \times \prod_{j=1}^{l-1}d\sigma_j ds\nonumber\\
	&\le \frac{C_{p,\gamma}}{N}\int_{0}^{t} \frac{e^{-\frac{\nu_0}{2}(t-s)}}{\tilde{w}(v)}\|h(s)\|_{L^p_vL^\infty_x}^2ds\nonumber \\
	&\le \frac{C_{p,\gamma}}{N}\frac{1}{{\tilde{w}(v)}}\sup_{0\le s \le t} \|h(s)\|_{L^p_vL^\infty_x}^2,
\end{align}
where we have used
\begin{align*}
	\int_{\R^3}(1+|\eta|)^{\left(e(p,\gamma)-2\beta\right)\cdot\frac{p}{p-2}} d\eta <C_{\beta, p}
\end{align*} 
since $\beta$ satisfies $\left(e(p,\gamma)-2\beta\right)\frac{p}{p-2}<-3$.\\
	\newline
	$\mathbf{Case\ 2:}$ $|v_l| \le N$, $|\eta| \ge N$.\\
We use the H\"older inequality to obtain
\begin{align} \label{LA22}
	|I_{5l}^2| 
	&\le  \frac{C_{p,\gamma}}{N^{1/2}}\int_{t_{l+1}}^{t_l} \frac{e^{-\frac{\nu_0}{2}(t-s)}}{\tilde{w}(v)} \|h(s)\|_{L^p_vL^\infty_x}^2 \int_{\prod_{j=1}^{l-1}\mathcal{V}_j} \int_{|v_l| \le N}  \left(\int_{|\eta|\ge N}(1+|\eta|)^{\left(e(p,\gamma)+1-2\beta\right)\cdot\frac{p}{p-2}} d\eta \right)^{\frac{p-2}{2p}}    \nonumber\\
	& \quad \times \tilde{w}(v_l)d\sigma_l\prod_{j=1}^{l-1}d\sigma_j ds\nonumber\\
	&\le \frac{C_{p,\gamma}}{N^{1/2}} \int_{0}^{t} \frac{e^{-\frac{\nu_0}{2}(t-s)}}{\tilde{w}(v)}\|h(s)\|_{L^p_vL^\infty_x}^2ds\nonumber\\
	&\le\frac{C_{p,\gamma}}{N^{1/2}}\frac{1}{{\tilde{w}(v)}}\sup_{0\le s \le t} \|h(s)\|_{L^p_vL^\infty_x}^2,
\end{align}
where we have used
\begin{align*}
	\int_{\R^3}(1+|\eta|)^{\left(e(p,\gamma)+1-2\beta\right)\cdot\frac{p}{p-2}} d\eta <C_{\beta, p}
\end{align*} 
since $\beta$ satisfies $\left(e(p,\gamma)+1-2\beta\right)\frac{p}{p-2}<-3$.\\
\newline
$\mathbf{Case\ 3:}$ $t_l-s \le \frac{1}{N}$.\\
We can easily compute
\begin{align} \label{LA23}
	|I_{5l}^3| 
	&\le C_{p,\gamma}\int_{t_{l}-\frac{1}{N}}^{t_l} \frac{e^{-\frac{\nu_0}{2}(t-s)}}{\tilde{w}(v)} \|h(s)\|_{L^p_vL^\infty_x} \int_{\prod_{j=1}^{l-1}\mathcal{V}_j} \int_{\mathcal{V}_l}  \left(\int_{\R^3}(1+|\eta|)^{e(p,\gamma)}|f(s,x_l-v_l(t_l-s),\eta)|^2 d\eta \right)^{1/2} \nonumber \\
	& \quad \times \tilde{w}(v_l)d\sigma_l\prod_{j=1}^{l-1}d\sigma_j ds\nonumber\\
	&\le C_{p,\gamma} \int_{t_l-\frac{1}{N}}^{t_l} \frac{1}{\tilde{w}(v)}\|h(s)\|_{L^p_vL^\infty_x}^2ds\nonumber\\
	&\le \frac{C_{p,\gamma}}{N}\frac{1}{{\tilde{w}(v)}}\sup_{0\le s \le t} \|h(s)\|_{L^p_vL^\infty_x}^2,
\end{align}
where we have used the H\"older inequality.\\
	\newline
	$\mathbf{Case\ 4:}$ $t_l-s \ge \frac{1}{N} \text{ and }|v_l|\le N$, $|\eta| \le N$.\\
We use the Cauchy-Schwartz inequality to estimate
\begin{align*}
	|I_{5l}^4| 
	&\le C_{p,\gamma,N} \int_{t_{l+1}}^{t_l-\frac{1}{N}} \frac{e^{-\frac{\nu_0}{2}(t-s)}}{\tilde{w}(v)} \|h(s)\|_{L^p_vL^\infty_x} \int_{\prod_{j=1}^{l-1}\mathcal{V}_j} \left(\int_{|v_l| \le N}  \int_{|\eta|\le N}|h(s,x_l-v_l(t_l-s),\eta)|^2 d\eta dv_l \right)^{1/2}\nonumber \\
	& \quad \times \prod_{j=1}^{l-1}d\sigma_j ds.
\end{align*}
We make a change of variables $v_l \mapsto y=x_l-v_l(t_l-s)$ with $\left|\frac{dy}{dv_l}\right| = (t_l-s)^3$:
\begin{align*}
	|I_{5l}^4| &\le C_{p,\gamma,N}\int_{t_{l+1}}^{t_l-\frac{1}{N}} \frac{e^{-\frac{\nu_0}{2}(t-s)}}{\tilde{w}(v)} \|h(s)\|_{L^p_vL^\infty_x} \int_{\prod_{j=1}^{l-1}\mathcal{V}_j} \left(\int_{\Omega}  \int_{|\eta|\le N}|h(s,y,\eta)|^2 d\eta dy \right)^{1/2} \prod_{j=1}^{l-1}d\sigma_j ds .   
\end{align*}
It follows from \eqref{RTIRE} that
\begin{align} \label{LA24}
	|I_{5l}^4| &\le C_{p,\gamma,N,\Omega} \int_0^t \frac{e^{-\frac{\nu_0}{2}(t-s)}}{\tilde{w}(v)}\|h(s)\|_{L^p_vL^\infty_x}\left[\mathcal{E}(F_0)^{1/2} + \|h(s)\|_{L^p_vL^\infty_x}^{\frac{1-\kappa}{2}}\mathcal{E}(F_0)^{\kappa/2}\right]ds\nonumber\\
	& \le \frac{C_{p,\gamma,N,\Omega}}{\tilde{w}(v)}\left[\mathcal{E}(F_0)^{1/2}\sup_{0\le s\le t}\|h(s)\|_{L^p_vL^\infty_x} + \mathcal{E}(F_0)^{\kappa/2}\sup_{0\le s\le t}\|h(s)\|_{L^p_vL^\infty_x}^{\frac{3-\kappa}{2}}\right]\nonumber\\
	& \le \frac{1}{\tilde{w}(v)}\left[ \frac{C_p}{N}\sup_{0\le s\le t}\left[\|h(s)\|_{L^p_vL^\infty_x}^2+\|h(s)\|_{L^p_vL^\infty_x}^{3-\kappa}\right]+ C_{p,\gamma,N,\Omega}\left[\mathcal{E}(F_0)+ \mathcal{E}(F_0)^{\kappa}\right]\right],
\end{align}
where we have used the Young's inequality.\\
Combining \eqref{LA21}, \eqref{LA22}, \eqref{LA23}, \eqref{LA24}, and summing over $1\le l \le k-1$, we obtain
\begin{equation} \label{LA29}
\begin{aligned}
	\int_0^{t-\epsilon} |I_{5}| ds &\le  \frac{C_{p,\epsilon,T_0}}{N^{1/2}}\frac{1}{{\tilde{w}(v)}}\sup_{0\le s \le t}\left[ \|h(s)\|_{L^p_vL^\infty_x}^2+\|h(s)\|_{L^p_vL^\infty_x}^{3-\kappa}\right]\\
	& \quad +\frac{1}{\tilde{w}(v)}C_{p,\gamma,N,\Omega,\epsilon,T_0}\left[\mathcal{E}(F_0)+ \mathcal{E}(F_0)^{\kappa}\right].
\end{aligned}
\end{equation}
Gathering \eqref{LA30}, \eqref{LA20}, \eqref{LA29}, we deduce that
\begin{equation} \label{LA32}
	\begin{aligned}
\int_0^{t-\epsilon} |S_{G_f}(t,s)w\Gamma^+(f,f)(s)|ds &\le \int_0^{t-\epsilon}\mathbf{1}_{\{t_1\le s\}} e^{-\frac{\nu_0}{2} (t-s)}\left|w\Gamma^+(f,f)(s,x-v(t-s),v)\right|ds \\
	& \quad + \frac{C_{p,\gamma}}{\tilde{w}(v)}\left(\epsilon+\frac{C_{p,\gamma,\epsilon,T_0}}{N^{1/2}} \right)\sup_{0\le s \le t}\left[ \|h(s)\|_{L^p_vL^\infty_x}^2+\|h(s)\|_{L^p_vL^\infty_x}^{3-\kappa}\right]\\
	& \quad +\frac{C_{p,,\gamma,N,\Omega,\epsilon,T_0}}{\tilde{w}(v)}\left[\mathcal{E}(F_0)^{1/2}+\mathcal{E}(F_0)+ \mathcal{E}(F_0)^{\kappa}\right].
\end{aligned}
\end{equation}
Therefore, we conclude the result of this lemma by inserting \eqref{LA31} and \eqref{LA32} into \eqref{LA10}.
\end{proof}

\bigskip

\begin{Lem} \label{LALinftyestimate2}
	Assume $\mathcal{E}(F_0) \le \epsilon_1 = \epsilon_1(\bar{\eta},T_0)$, where $\epsilon_1$ is determined in Lemma \ref{Rfest1}. Under the a priori assumption \eqref{Aprioriassump1}, there exists a generic constant $C_3 \ge 1$, depending on $p$, $\beta$, and $\gamma$, such that
\begin{align*}
	\|h(t)\|_{L^p_vL^\infty_x} &\le C_3\|h_0\|_{L^p_vL^\infty_x}\left(1+\int_0^t \|h(s)\|_{L^p_vL^\infty_x}ds \right)e^{-\frac{\nu_0}{4}t}\\
	& \quad + C_3\bigg\{\left(\epsilon + \frac{C_{p,,\gamma,\epsilon,T_0}}{N^{1/2}} \right) \sup_{0\le s\le t}\biggl[\|h(s)\|_{L^p_vL^\infty_x}+\|h(s)\|_{L^p_vL^\infty_x}^2+\|h(s)\|_{L^p_vL^\infty_x}^3+\|h(s)\|_{L^p_vL^\infty_x}^4\nonumber\\
	& \qquad+\|h(s)\|_{L^p_vL^\infty_x}^{\frac{p}{2p-2}}+\|h(s)\|_{L^p_vL^\infty_x}^{\frac{3p-2}{2p-2}}+\|h(s)\|_{L^p_vL^\infty_x}^{\frac{5p-4}{2p-2}}+\|h(s)\|_{L^p_vL^\infty_x}^{\frac{7p-6}{2p-2}}+\|h(s)\|_{L^p_vL^\infty_x}^{\frac{9p-8}{2p-2}}\biggr]\\
	& \qquad +C_{p,\gamma,N,\epsilon,\Omega,T_0}\left[\mathcal{E}(F_0)+\mathcal{E}(F_0)^2+\mathcal{E}(F_0)^{\frac{p-2}{2p-2}}+\mathcal{E}(F_0)^{\frac{p-2}{p-1}}\right]\bigg\}
\end{align*}
for all $0 \le t \le T_0$, where $\epsilon>0$ can be arbitrary small and $N>0$ can be arbitrary large.
\end{Lem}
\begin{proof}
	In order to derive our result in this lemma, we need to focus on estimating the following terms in Lemma \ref{LALinftyestimate1}:
	\begin{align*}
		&\int_0^{t-\epsilon} \mathbf{1}_{\{t_1\le s\}}e^{-\frac{\nu_0}{2}(t-s)}\left|(K_wh)(s,x-v(t-s),v) \right|ds\\
		&+ \int_0^{t-\epsilon} \mathbf{1}_{\{t_1\le s\}}e^{-\frac{\nu_0}{2}(t-s)} \left|w\Gamma^+(f,f)(s,x-v(t-s),v) \right|ds\\
		&=: J_1+J_2.
	\end{align*}
	Let $(t,x,v) \in (0, T_0] \times \Omega \times \R^3$. Denote $X(s):=x-v(t-s)$.\\
	First of all, let us estimate $J_1$.
	\begin{align*}
		J_1 \le \int_0^{t-\epsilon} \mathbf{1}_{\{t_1\le s\}}e^{-\frac{\nu_0}{2}(t-s)} \int_{\R^3}|k_w(v,v')||h(s,X(s),v')|dv'ds
	\end{align*}
	and we apply the estimate in Lemma \ref{LALinftyestimate1} to bound the term $|h(s,X(s),v')|$ in the above integral
	\begin{align*}
		J_1 &\le \int_{t_1}^{t-\epsilon} e^{-\frac{\nu_0}{2}(t-s)} \int_{\R^3}|k_w(v,v')||(S_{G_f}(s,0)h_0)(s,X(s),v')|dv'ds\\
		& \quad + \int_{t_1}^{t-\epsilon} e^{-\frac{\nu_0}{2}(t-s)} \int_{\R^3}|k_w(v,v')|\bigg\{\int_{s-\epsilon}^s|(S_{G_f}(s,s')K_wh(s'))(s,X(s),v')|ds'\\
		& \qquad +\int_{s-\epsilon}^s|(S_{G_f}(s,s')w\Gamma^+(f,f)(s'))(s,X(s),v')|ds'\bigg\}dv'ds \\
		& \quad + C_{p,\gamma} \int_0^{t-\epsilon} e^{-\frac{\nu_0}{2}(t-s)} \int_{\R^3}|k_w(v,v')|\frac{1}{\tilde{w}(v')}\Biggl\{\left(\epsilon + \frac{C_{p,\gamma,\epsilon,T_0}}{N^{1/2}} \right)\sup_{0\le s\le t}\biggl[\|h(s)\|_{L^p_vL^\infty_x}+\|h(s)\|_{L^p_vL^\infty_x}^2\\
		& \qquad +\|h(s)\|_{L^p_vL^\infty_x}^{\frac{p}{2p-2}}+\|h(s)\|_{L^p_vL^\infty_x}^{\frac{5p-4}{2p-2}}\biggr]+C_{p,\gamma,N,\Omega,T_0}\left[\mathcal{E}(F_0)^{1/2}+\mathcal{E}(F_0)+\mathcal{E}(F_0)^{\frac{p-2}{2p-2}}\right]\Biggr\} dv'ds\\
		&\quad +\int_0^{t-\epsilon}\int_0^{s-\epsilon} \mathbf{1}_{\{t_1\le s\}}\mathbf{1}_{\{t_1'\le s'\}}e^{-\frac{\nu_0}{2}(t-s')}\int_{\R^3}\int_{\R^3}|k_w(v,v')||k_w(v',v'')||h(s',X(s)-v'(s-s'),v'')|\\
		&\qquad \times dv''dv'ds'ds\\
		& \quad +\int_0^{t-\epsilon}\int_0^{s-\epsilon} \mathbf{1}_{\{t_1\le s\}}\mathbf{1}_{\{t_1'\le s'\}}e^{-\frac{\nu_0}{2}(t-s')}\int_{\R^3}|k_w(v,v')||w\Gamma^+(f,f)(s',X(s)-v'(s-s'),v')|dv'ds'ds\\
		& =:J_{11}+J_{12}+J_{13}+J_{14}+J_{15},
	\end{align*}
	where we have denoted $t_1':=t_1(s,X(s),v')$.\\
	For $J_{11}$, we use Lemma \ref{KESTIMATE} and \eqref{DampLpdecay3} to obtain
	\begin{align*} 
		J_{11} &\le C_{p,\gamma} \int_0^{t-\epsilon} e^{-\frac{\nu_0}{2}(t-s)} e^{-\frac{\nu_0}{4}s} \left(\frac{1}{1+|v|}\right)^{\frac{p'(1-\gamma)+1}{p'}}\|h_0\|_{L^p_vL^\infty_x}ds\nonumber\\
		&\le C_{p,\gamma} \left(\frac{1}{1+|v|}\right)^{\frac{p'(1-\gamma)+1}{p'}}e^{-\frac{\nu_0}{4}t}\|h_0\|_{L^p_vL^\infty_x},
	\end{align*}  
	which implies
	\begin{align}\label{LA33}
		\|J_{11}\|_{L^p_vL^\infty_x} \le C_{p,\gamma} e^{-\frac{\nu_0}{4}t}\|h_0\|_{L^p_vL^\infty_x}.
	\end{align}
	For $J_{12}$, we use Lemma \ref{KESTIMATE} and \eqref{LA11}, \eqref{LA19} to obtain
	\begin{align*}
		J_{12} &\le C_{p,\gamma}\int_0^{t-\epsilon} e^{-\frac{\nu_0}{2}(t-s)}\int_{s-\epsilon}^s \left(\frac{1}{1+|v|}\right)^{\frac{p'(1-\gamma)+1}{p'}} e^{-\frac{\nu_0}{4}(s-s')}\left(\|h(s')\|_{L^p_vL^\infty_x}+\|h(s')\|_{L^p_vL^\infty_x}^2 \right)ds'ds\\
		& \le \epsilon C_{p,\gamma} \left(\frac{1}{1+|v|}\right)^{\frac{p'(1-\gamma)+1}{p'}}\sup_{0\le s'\le t}\left[\|h(s')\|_{L^p_vL^\infty_x}+\|h(s')\|_{L^p_vL^\infty_x}^2 \right],
	\end{align*}
	which implies
	\begin{align} \label{LA34}
		\|J_{12}\|_{L^p_vL^\infty_x} \le \epsilon C_{p,\gamma} \sup_{0\le s'\le t}\left[\|h(s')\|_{L^p_vL^\infty_x}+\|h(s')\|_{L^p_vL^\infty_x}^2 \right].
	\end{align}
	For $J_{13}$, we use Lemma \ref{KESTIMATE} to get
	\begin{align*}
		J_{13} &\le \frac{C_{p,\gamma}}{(1+|v|)^{2-\gamma}} \Biggl\{\left(\epsilon + \frac{C_{p,\gamma,\epsilon,T_0}}{N^{1/2}} \right)\sup_{0\le s\le t}\biggl[\|h(s)\|_{L^p_vL^\infty_x}+\|h(s)\|_{L^p_vL^\infty_x}^2+\|h(s)\|_{L^p_vL^\infty_x}^{\frac{p}{2p-2}}+\|h(s)\|_{L^p_vL^\infty_x}^{\frac{5p-4}{2p-2}}\biggr]\\
		& \quad +C_{p,\gamma,N,\Omega,T_0}\left[\mathcal{E}(F_0)^{1/2}+\mathcal{E}(F_0)+\mathcal{E}(F_0)^{\frac{p-2}{2p-2}}\right]\Biggr\},
	\end{align*}
	which implies
	\begin{align} \label{LA35}
		\|J_{13}\|_{L^p_vL^\infty_x} &\le C_{p,\gamma}\left(\epsilon + \frac{C_{p,\gamma,\epsilon,T_0}}{N^{1/2}} \right)\sup_{0\le s\le t}\biggl[\|h(s)\|_{L^p_vL^\infty_x}+\|h(s)\|_{L^p_vL^\infty_x}^2+\|h(s)\|_{L^p_vL^\infty_x}^{\frac{p}{2p-2}}+\|h(s)\|_{L^p_vL^\infty_x}^{\frac{5p-4}{2p-2}}\biggr]\nonumber\\
		& \quad +C_{p,\gamma,N,\Omega,T_0}\left[\mathcal{E}(F_0)^{1/2}+\mathcal{E}(F_0)+\mathcal{E}(F_0)^{\frac{p-2}{2p-2}}\right].
	\end{align}
	For $J_{14}$, we divide this term into three cases.\\
	\newline
	$\mathbf{Case\ 1:}$ $|v|\ge N$.\\
	We use Lemma \ref{KESTIMATE} to estimate
\begin{align*}
	|J_{14}^1| &\le  \int_{t_1}^t\int_{t_1'}^s e^{-\frac{\nu_0}{2}(t-s')} \mathbf{1}_{\{|v|\ge N\}}  \int_{\R^3} \int_{\R^3} |k_w(v,v')||k_w(v',v'')| |h(s',X(s)-v'(s-s'),v'')|dv''dv'ds'ds\\
	& \le C_p \sup_{0\le s' \le t}\|h(s')\|_{L^p_v L^\infty_x} \mathbf{1}_{\{|v|\ge N\}}\int_0^t\int_0^s e^{-\frac{\nu_0}{2}(t-s')}   \int_{v' \in \R^3} |k_w(v,v')|\left(\frac{1}{1+|v'|}\right)^{\frac{p'(1-\gamma)+1}{p'}}dv'ds'ds\\
	& \le\frac{C_p}{(1+|v|)^{\frac{p'(1-\gamma)+1}{p'}+2-\gamma}}\mathbf{1}_{\{|v|\ge N\}}\sup_{0\le s' \le t}\|h(s')\|_{L^p_v L^\infty_x}.
\end{align*}
Thus we obtain
\begin{align} \label{LA36}
		\|J_{14}^1\|_{L^p_vL^\infty_x}\le \frac{C_p}{N}\sup_{0\le s'\le t}\|h(s')\|_{L^p_vL^\infty_x}.
\end{align}
\newline
$\mathbf{Case\ 2:}$ $|v|\le N, |v'|\ge 2N$ or $|v'|\le 2N, |v''|\ge 3N$.\\
Either one of the following are valid:
\begin{equation} \label{class0}
	\begin{aligned}
	&|k_w(v,v')| \le e^{-\frac{\epsilon}{8}N^2}|k_w(v,v')|e^{\frac{\epsilon}{8}|v-v'|^2}\\
	\text{or} \quad &|k_w(v',v'')| ^{p'}\le e^{-\frac{\epsilon}{8}N^2}|k_w(v',v'')|^{p'}e^{\frac{\epsilon}{8}|v'-v''|^2},
\end{aligned}
\end{equation}
where $\epsilon$ is sufficiently small.
From Lemma \ref{KESTIMATE}, either one of the following are valid:
\begin{equation} \label{class1}
	\begin{aligned}
	&\int_{|v'| \ge 2N}|k_w(v,v')|e^{\frac{\epsilon}{8}|v-v'|^2}\left(\frac{1}{1+|v'|}\right)^{\frac{p'(1-\gamma)+1}{p'}}dv' \lesssim \left(\frac{1}{1+|v|}\right)^{\frac{p'(1-\gamma)+1}{p'}+2-\gamma}\\
	\text{or} \quad & \int_{|v''|\ge 3N}|k_w(v',v'')|^{p'}e^{\frac{\epsilon}{8}|v'-v''|^2}dv'' \lesssim \left(\frac{1}{1+|v'|}\right)^{p'(1-\gamma)+1}.
\end{aligned}
\end{equation}
From \eqref{class0} and \eqref{class1}, we obtain
\begin{align*}
	J_{14}^2 &\le \int_{t_1}^t\int_{t_1'}^s e^{-\frac{\nu_0}{2}(t-s')}\mathbf{1}_{\{|v|\le N\}}   \int_{|v'| \ge 2N}\int_{v'' \in \R^3} |k_w(v,v')||k_w(v',v'')| |h(s',X(s)-v'(s-s'),v'')|\\
	& \qquad \times dv''dv'ds'ds\\
	& \quad + \int_{t_1}^t\int_{t_1'}^s e^{-\frac{\nu_0}{2}(t-s')}  \int_{|v'|\le 2N}\int_{|v''| \ge 3N} |k_w(v,v')||k_w(v',v'')| |h(s',X(s)-v'(s-s'),v'')|\\
	& \qquad \times dv''dv'ds'ds\\
	& \le \frac{C_p}{(1+|v|)^{\frac{p'(1-\gamma)+1}{p'}+2-\gamma}} e^{-\frac{\epsilon }{8}N^2}\sup_{0\le s'\le t}\|h(s')\|_{L^p_vL^\infty_x}.
\end{align*}
Thus we obtain
	\begin{align} \label{LA37}
		\|J_{14}^2\|_{L^p_vL^\infty_x}\le C_pe^{-\frac{\epsilon }{8}N^2} \sup_{0\le s'\le t}\|h(s')\|_{L^p_vL^\infty_x}.
	\end{align}
	\newline
	$\mathbf{Case\ 3:}$ $|v|\le N, |v'|\le 2N, |v''|\le 3N$.\\
	Using the Cauchy-Schwartz inequality and Lemma \ref{KESTIMATE}, we deduce
	\begin{align*}
		J_{14}^3 &\le \int_{t_1}^{t-\epsilon}\int_{t_1'}^{s-\epsilon} e^{-\frac{\nu_0}{2}(t-s')}\mathbf{1}_{\{|v|\le N\}}\left(\int_{|v'|\le 2N}\int_{|v''|\le 3N}|k_w(v,v')|^2|k_w(v',v'')|^2dv''dv'\right)^{1/2} \\
		& \quad \times \left(\int_{|v'|\le 2N}\int_{|v''|\le 3N}|h(s',X(s)-v'(s-s'),v'')|^2dv''dv'  \right)^{1/2}ds'ds\\
		& \le \frac{C_p}{(1+|v|)^{6-4\gamma}} \int_{t_1}^{t-\epsilon}\int_{t_1'}^{s-\epsilon} e^{-\frac{\nu_0}{2}(t-s')}\mathbf{1}_{\{|v|\le N\}}\\
		& \quad \times\left(\int_{|v'|\le 2N}\int_{|v''|\le 3N}|h(s',X(s)-v'(s-s'),v'')|^2dv''dv'  \right)^{1/2} ds'ds.
	\end{align*}
	We make a change of variables $v' \mapsto y=X(s)-v'(s-s')$ with $\left|\frac{dy}{dv'}\right| = (s-s')^3$:
	\begin{align*}
		J_{14}^3  \le  \frac{C_{p,N,\epsilon}}{(1+|v|)^{6-4\gamma}}\int_{t_1}^{t-\epsilon}\int_{t_1'}^{s-\epsilon} e^{-\frac{\nu_0}{2}(t-s')}\left(\int_{\Omega}\int_{|v''|\le 3N}|h(s',y,v'')|^2dv''dy  \right)^{1/2}ds'ds.
	\end{align*}
	It follows from \eqref{RTIRE} that 
	\begin{align*}
		J_{14}^3 \le \frac{1}{(1+|v|)^{6-4\gamma}}\left[ \frac{C_{p,\epsilon}}{N}\sup_{0\le s\le t}\|h(s)\|_{L^p_vL^\infty_x}^{1-\kappa}+ C_{p,N,\epsilon,\Omega}\left[\mathcal{E}(F_0)^{1/2}+ \mathcal{E}(F_0)^{\kappa}\right]\right],
	\end{align*}
	which implies that
	\begin{align} \label{LA38}
		\|J_{14}^3\|_{L^p_vL^\infty_x} \le \frac{C_{p,\epsilon}}{N}\sup_{0\le s\le t}\|h(s)\|_{L^p_vL^\infty_x}^{1-\kappa}+ C_{p,N,\epsilon,\Omega}\left[\mathcal{E}(F_0)^{1/2}+ \mathcal{E}(F_0)^{\kappa}\right].
	\end{align}
	Gathering \eqref{LA36}, \eqref{LA37}, \eqref{LA38}, we derive
	\begin{align} \label{LA39}
		\|J_{14}\|_{L^p_vL^\infty_x} \le \frac{C_{p,\epsilon}}{N}\sup_{0\le s\le t}\left[\|h(s)\|_{L^p_vL^\infty_x}+\|h(s)\|_{L^p_vL^\infty_x}^{1-\kappa}\right]+ C_{p,N,\epsilon,\Omega}\left(\mathcal{E}(F_0)^{1/2}+ \mathcal{E}(F_0)^{\kappa}\right).
	\end{align}
	\newline
	For $J_{15}$, by \eqref{Gamma+estimatecombine}, it is bounded by
	\begin{align*}
		&C_{p,\gamma}\int_{t_1}^{t-\epsilon}\int_{t_1'}^{s-\epsilon} e^{-\frac{\nu_0}{2}(t-s')}\|h(s')\|_{L^p_vL^\infty_x}\int_{\R^3}|k_w(v,v')|\frac{1}{(1+|v'|)^{\mathcal{M}}}\\
		& \quad \times \left(\int_{\R^3}(1+|\eta|)^{e(p,\gamma)}|f(s',X(s)-v'(s-s'),\eta)|^2 d\eta \right)^{1/2}dv'ds'ds,
	\end{align*}
	where $\mathcal{M}:=\min\left\{\frac{1}{p'},\frac{2-p'\gamma}{p'}\right\}$.
	We divide this term into four cases.\\
	\newline
	$\mathbf{Case\ 1:}$ $|v|\ge N$.\\
	We use Lemma \ref{KESTIMATE} and a similar argument in \eqref{LA21} to obtain
	\begin{align*}
		J_{15}^1 &\le C_{p,\gamma} \int_{t_1}^{t-\epsilon}\int_{t_1'}^{s-\epsilon}e^{-\frac{\nu_0}{2}(t-s')}\mathbf{1}_{\{|v|\ge N\}}\|h(s')\|_{L^p_vL^\infty_x}^2\int_{\R^3}|k_w(v,v')|\frac{1}{(1+|v'|)^{\mathcal{M}}}dv'ds'ds\\
		& \le C_{p,\gamma}\int_{0}^{t}\int_{0}^{s}e^{-\frac{\nu_0}{2}(t-s')}\mathbf{1}_{\{|v|\ge N\}}\frac{\|h(s')\|_{L^p_vL^\infty_x}^2}{(1+|v|)^{2-\gamma+\mathcal{M}}}ds'ds\\
		& \le \frac{C_{p,\gamma}}{(1+|v|)^{2-\gamma+\mathcal{M}}}\mathbf{1}_{\{|v|\ge N\}}\sup_{0\le s' \le t}\|h(s')\|_{L^p_vL^\infty_x}^2,
	\end{align*}
	and it follows that
	\begin{align}  \label{LA40}
		\|J_{15}^1\|_{L^p_vL^\infty_x} \le \frac{C_{p,\gamma}}{N^{2-\gamma}}\sup_{0\le s' \le t}\|h(s')\|_{L^p_vL^\infty_x}^2.	
	\end{align}
	\newline
	$\mathbf{Case\ 2:}$ $|v|\le N, |v'|\ge 2N$.\\
	Using Lemma \ref{KESTIMATE} and similar arguments in \eqref{LA41}, \eqref{LA42}, \eqref{LA13}, we deduce that
	\begin{align*} 
		J_{15}^2 & \le C_{p,\gamma} \int_{t_1}^{t-\epsilon}\int_{t_1'}^{s-\epsilon} e^{-\frac{\nu_0}{2}(t-s')}\mathbf{1}_{\{|v|\le N\}}\|h(s')\|_{L^p_vL^\infty_x}^2 \int_{|v'|\ge 2N}|k_w(v,v')|\frac{1}{(1+|v'|)^{\mathcal{M}}}dv'ds'ds\\
		& \le  C_{p,\gamma} e^{-\frac{\epsilon N^2}{8}}\frac{1}{(1+|v|)^{2-\gamma+\mathcal{M}}}\sup_{0\le s' \le t}\|h(s')\|_{L^p_vL^\infty_x}^2,
	\end{align*}	
	and it follows that
	\begin{align} \label{LA43}
		\|J_{15}^2\|_{L^p_vL^\infty_x} \le C_{p,\gamma}e^{-\frac{\epsilon N^2}{8}}\sup_{0\le s' \le t}\|h(s')\|_{L^p_vL^\infty_x}^2.
	\end{align}
	\newline
	$\mathbf{Case\ 3:}$ $|v|\le N, |v'|\le 2N, |\eta| \ge N$.\\
	We apply Lemma \ref{KESTIMATE} and a similar argument in \eqref{LA22} to obtain
	\begin{align*}
		J_{15}^3 &\le \frac{C_{p,\gamma}}{N^{1/2}}\int_{t_1}^{t-\epsilon}\int_{t_1'}^{s-\epsilon} e^{-\frac{\nu_0}{2}(t-s')}\mathbf{1}_{\{|v|\le N\}}\|h(s')\|_{L^p_vL^\infty_x}^2\int_{|v'|\le 2N}|k_w(v,v')|\frac{1}{(1+|v'|)^{\mathcal{M}}}dv'ds'ds\\
		& \le \frac{C_{p,\gamma}}{N^{1/2}}\frac{1}{(1+|v|)^{2-\gamma+\mathcal{M}}}\sup_{0\le s' \le t}\|h(s')\|_{L^p_vL^\infty_x}^2,
	\end{align*}
	which implies that
	\begin{align} \label{LA44}
		\|J_{15}^3\|_{L^p_vL^\infty_x} \le \frac{C_{p,\gamma}}{N^{1/2}}\sup_{0\le s' \le t}\|h(s')\|_{L^p_vL^\infty_x}^2.
	\end{align}
	\newline
	$\mathbf{Case\ 4:}$ $|v|\le N, |v'|\le 2N, |\eta| \le N$.\\
	Using the Cauchy-Schwartz inequality and Lemma \ref{KESTIMATE}, we get
	\begin{align*}
		J_{15}^4 &\le C_{p,\gamma}\int_{t_1}^{t-\epsilon}\int_{t_1'}^{s-\epsilon} e^{-\frac{\nu_0}{2}(t-s')}\|h(s')\|_{L^p_vL^\infty_x}\left(\int_{|v'|\le 2N}|k_w(v,v')|^2dv'\right)^{1/2}\\
		& \quad \times \left(\int_{|v'|\le 2N}\int_{|\eta|\le N}|h(s',X(s)-v'(s-s'),\eta)|^2 d\eta dv' \right)^{1/2}ds'ds\\
		& \le \frac{C_{p,\gamma}}{(1+|v|)^{\frac{3-2\gamma}{2}}} \int_{t_1}^{t-\epsilon}\int_{t_1'}^{s-\epsilon} e^{-\frac{\nu_0}{2}(t-s')}\|h(s')\|_{L^p_vL^\infty_x}\\
		& \quad \times\left(\int_{|v'|\le 2N}\int_{|\eta|\le N}|h(s',X(s)-v'(s-s'),\eta)|^2 d\eta dv' \right)^{1/2} ds'ds
	\end{align*}
	We make a change of variables $v' \mapsto y=X(s)-v'(s-s')$ with $\left|\frac{dy}{dv'}\right| = (s-s')^3$:
	\begin{align*}
		J_{15}^4  \le  \frac{C_{p,\gamma,\epsilon}}{(1+|v|)^{\frac{3-2\gamma}{2}}}\int_{t_1}^{t-\epsilon}\int_{t_1'}^{s-\epsilon} e^{-\frac{\nu_0}{2}(t-s')}\|h(s')\|_{L^p_vL^\infty_x}\left(\int_{\Omega}\int_{|\eta|\le N}|h(s',y,\eta)|^2d\eta dy  \right)^{1/2}ds'ds.
	\end{align*}
	It follows from \eqref{RTIRE} that 
	\begin{align*}
		J_{15}^4 \le \frac{1}{(1+|v|)^{\frac{3-2\gamma}{2}}}\left[ \frac{C_{p,\gamma,\epsilon}}{N}\sup_{0\le s\le t}\left[\|h(s)\|_{L^p_vL^\infty_x}^2+\|h(s)\|_{L^p_vL^\infty_x}^{3-\kappa}\right]+ C_{p,\gamma,N,\epsilon,\Omega}\left[\mathcal{E}(F_0)+ \mathcal{E}(F_0)^{\kappa}\right]\right],
	\end{align*}
	which implies that
	\begin{align} \label{LA45}
		\|J_{15}^4\|_{L^p_vL^\infty_x} \le \frac{C_{p,\gamma,\epsilon}}{N}\sup_{0\le s\le t}\left[\|h(s)\|_{L^p_vL^\infty_x}^2+\|h(s)\|_{L^p_vL^\infty_x}^{3-\kappa}\right]+ C_{p,\gamma,N,\epsilon ,\Omega}\left[\mathcal{E}(F_0)+ \mathcal{E}(F_0)^{\kappa}\right].
	\end{align}
	Gathering \eqref{LA40}, \eqref{LA43}, \eqref{LA44}, \eqref{LA45}, we derive
	\begin{align} \label{LA46}
		\|J_{15}\|_{L^p_vL^\infty_x} \le \frac{C_{p,\gamma,\epsilon}}{N^{1/2}}\sup_{0\le s\le t}\left[\|h(s)\|_{L^p_vL^\infty_x}^2+\|h(s)\|_{L^p_vL^\infty_x}^{3-\kappa}\right]+ C_{p,\gamma,N,\epsilon ,\Omega}\left[\mathcal{E}(F_0)+ \mathcal{E}(F_0)^{\kappa}\right].
	\end{align}
	\newline
	Inserting \eqref{LA33}, \eqref{LA34}, \eqref{LA35}, \eqref{LA39}, \eqref{LA46} into $J_1$, we derive
	\begin{align} \label{LA52}
		\|J_1\|_{L^p_vL^\infty_x} &\le C_p e^{-\frac{\nu_0}{4}t} \|h_0\|_{L^p_vL^\infty_x}\nonumber\\
		& \quad  + C_{p,\gamma}\left(\epsilon + \frac{C_{p,\gamma,\epsilon,T_0}}{N^{1/2}} \right)\sup_{0\le s\le t}\biggl[\|h(s)\|_{L^p_vL^\infty_x}+\|h(s)\|_{L^p_vL^\infty_x}^2+\|h(s)\|_{L^p_vL^\infty_x}^{\frac{p}{2p-2}}+\|h(s)\|_{L^p_vL^\infty_x}^{\frac{5p-4}{2p-2}}\biggr]\nonumber\\
		& \quad +C_{p,\gamma,N,\epsilon,\Omega,T_0}\left[\mathcal{E}(F_0)^{1/2}+\mathcal{E}(F_0)+\mathcal{E}(F_0)^{\frac{p-2}{2p-2}}\right]
	\end{align}
	Next, let us estimate $J_2$. From \eqref{Gamma+estimatecombine}, we have
	\begin{align} \label{LA53}
		J_{2} &\le C_{p,\gamma}\int_{t_1}^{t-\epsilon} e^{-\frac{\nu_0}{2}(t-s)} \frac{\|h(s)\|_{L^p_vL^\infty_x}}{(1+|v|)^{\mathcal{M}}} \left(\left(\int_{|\eta|\ge N}+\int_{|\eta|\le N}\right)(1+|\eta|)^{e(p,\gamma)}|f(s,X(s),\eta)|^2 d\eta \right)^{1/2} ds\nonumber\\
		& \le  \frac{C_{p,\gamma}}{N}\frac{1}{(1+|v|)^{\mathcal{M}}}\sup_{0\le s\le t}\|h(s)\|_{L^p_vL^\infty_x}^2\nonumber\\
		& \quad + C_{p,\gamma}\int_{t_1}^{t-\epsilon} e^{-\frac{\nu_0}{2}(t-s)} \frac{\|h(s)\|_{L^p_vL^\infty_x}}{(1+|v|)^{\mathcal{M}}} \left(\int_{|\eta|\le N}(1+|\eta|)^{e(p,\gamma)-2\beta}|h(s,X(s),\eta)|^2 d\eta \right)^{1/2} ds,
	\end{align}
	where $\mathcal{M}:=\min\left\{\frac{1}{p'},\frac{2-p'\gamma}{p'}\right\}$.
	We apply the estimate in Lemma \ref{LALinftyestimate1} to bound the term $|h(s,X(s),\eta)|$ in the last term in the above
	\begin{align*}
		&\left(\int_{|\eta|\le N}(1+|\eta|)^{e(p,\gamma)-2\beta}|h(s,X(s),\eta)|^2 d\eta \right)^{1/2} \\
		& \le \left(\int_{|\eta|\le N}(1+|\eta|)^{e(p,\gamma)-2\beta}|(S_{G_f}(s,0)h_0)(s,X(s),\eta)|^2 d\eta \right)^{1/2}\\
		& \quad + \left(\int_{|\eta|\le N}(1+|\eta|)^{e(p,\gamma)-2\beta}\left|\int_{s-\epsilon}^s\left|S_{G_f}(s,s')K_wh(s')\right|+\left|S_{G_f}(s,s')w\Gamma^+(f,f)(s')\right|ds'\right|^2 d\eta \right)^{1/2}\\
		& \quad + C_{p,\gamma}\left(\int_{|\eta|\le N}(1+|\eta|)^{e(p,\gamma)-2\beta}\frac{1}{\tilde{w}(\eta)^2} d\eta \right)^{1/2}\Biggl\{\left(\epsilon + \frac{C_{p,\epsilon,T_0}}{N^{1/2}} \right)\sup_{0\le s\le t}\biggl[\|h(s)\|_{L^p_vL^\infty_x}+\|h(s)\|_{L^p_vL^\infty_x}^2\\
		& \qquad +\|h(s)\|_{L^p_vL^\infty_x}^{\frac{p}{2p-2}}+\|h(s)\|_{L^p_vL^\infty_x}^{\frac{5p-4}{2p-2}}\biggr]+C_{p,N,\Omega,T_0}\left[\mathcal{E}(F_0)^{1/2}+\mathcal{E}(F_0)+\mathcal{E}(F_0)^{\frac{p-2}{2p-2}}\right]\Biggr\}\\
		& \quad +\left(\int_{|\eta|\le N}(1+|\eta|)^{e(p,\gamma)-2\beta}\left|\int_{0}^{s-\epsilon}\mathbf{1}_{\{t_1' \le s'\}}e^{-\frac{\nu_0}{2}(s-s')}(K_wh)(s',X(s)-\eta(s-s'),\eta)ds'\right|^2 d\eta \right)^{1/2}\\
		& \quad +\left(\int_{|\eta|\le N}(1+|\eta|)^{e(p,\gamma)-2\beta}\left|\int_{0}^{s-\epsilon}\mathbf{1}_{\{t_1' \le s'\}}e^{-\frac{\nu_0}{2}(s-s')}w\Gamma^+(f,f)(s',X(s)-\eta(s-s'),\eta)ds'\right|^2 d\eta \right)^{1/2}\\
		&=: \tilde{L}_1+\tilde{L}_2+\tilde{L}_3+\tilde{L}_4+\tilde{L}_5.
	\end{align*}
	We define
	\begin{align*}
		L_{i} = \int_{t_1}^{t-\epsilon}e^{-\frac{\nu_0}{2}(t-s)}\frac{\|h(s)\|_{L^p_vL^\infty_x}}{(1+|v|)^{\mathcal{M}}}\tilde{L}_ids.
	\end{align*}
	For $L_1$, we use \eqref{DampLpdecay3} to obtain
	\begin{align*}
		L_1 \le C_{p,\gamma} \int_{t_1}^{t-\epsilon} e^{-\frac{\nu_0}{2}(t-s)}e^{-\frac{\nu_0}{4}s}\frac{\|h(s)\|_{L^p_vL^\infty_x}}{(1+|v|)^{\mathcal{M}}}\|h_0\|_{L^p_vL^\infty_x}ds,
	\end{align*}
	which implies that
	\begin{align} \label{LA54}
		\|L_1\|_{L^p_vL^\infty_x} \le C_{p,\gamma} e^{-\frac{\nu_0}{4}t}\|h_0\|_{L^p_vL^\infty_x} \int_0^t\|h(s)\|_{L^p_vL^\infty_x}ds.
	\end{align}
	For $L_2$, we use the Minkowski integral inequality and \eqref{LA11}, \eqref{LA19} to get
	\begin{align*}
		L_2 &\le \epsilon C_{p,\gamma}\int_{t_1}^{t-\epsilon} e^{-\frac{\nu_0}{2}(t-s)}\frac{\|h(s)\|_{L^p_vL^\infty_x}}{(1+|v|)^{\mathcal{M}}}\sup_{0\le s' \le t}\left[\|h(s')\|_{L^p_vL^\infty_x}+\|h(s')\|_{L^p_vL^\infty_x}^2 \right]ds\\
		& \le \frac{\epsilon C_{p,\gamma}}{(1+|v|)^{\mathcal{M}}}\sup_{0\le s \le t}\left[\|h(s)\|_{L^p_vL^\infty_x}^2+\|h(s)\|_{L^p_vL^\infty_x}^3\right],
	\end{align*} 
	and it follows that
	\begin{align} \label{LA55}
		\|L_2\|_{L^p_vL^\infty_x} \le \epsilon C_{p,\gamma}	\sup_{0\le s \le t}\left[\|h(s)\|_{L^p_vL^\infty_x}^2+\|h(s)\|_{L^p_vL^\infty_x}^3\right].
	\end{align}
	For $L_3$, we easily compute
	\begin{align} \label{LA56}
		\|L_3\|_{L^p_vL^\infty_x} &\le C_{p,\gamma}\Biggl\{\left(\epsilon + \frac{C_{p,\gamma,\epsilon,T_0}}{N^{1/2}} \right)\sup_{0\le s\le t}\biggl[\|h(s)\|_{L^p_vL^\infty_x}^2+\|h(s)\|_{L^p_vL^\infty_x}^3+\|h(s)\|_{L^p_vL^\infty_x}^{\frac{3p-2}{2p-2}}+\|h(s)\|_{L^p_vL^\infty_x}^{\frac{7p-6}{2p-2}}\biggr]\nonumber\\
		& \qquad +C_{p,\gamma,N,\Omega,T_0}\left[\mathcal{E}(F_0)+\mathcal{E}(F_0)^2+\mathcal{E}(F_0)^{\frac{p-2}{p-1}}\right]\Biggr\}.
	\end{align}
	For $L_4$, we first need to estimate $\tilde{L}_4$. We use the H\"older inequality and Lemma \ref{KESTIMATE} to obtain
	\begin{align*}
		\tilde{L}_4 &\le C_{p,\gamma}\left(\int_{|\eta|\le N}(1+|\eta|)^{e(p,\gamma)-2\beta}\int_{t_1'}^{s-\epsilon}e^{-\frac{\nu_0}{2}(s-s')}\left|(K_wh)(s',X(s)-\eta(s-s'),\eta)\right|^2ds' d\eta \right)^{\frac{1}{2}}\\
		& \le C_{p,\gamma} e^{-\frac{\epsilon N^2}{8}} \sup_{0\le s \le t}\|h(s)\|_{L^p_vL^\infty_x}\\
		& \quad + C_{p,\gamma,N}\left(\int_{t_1'}^{s-\epsilon}e^{-\frac{\nu_0}{2}(s-s')}\int_{|\eta|\le N}\int_{|v'|\le 2N}|h(s',X(s)-\eta(s-s'),v')|^2dv' d\eta ds'  \right)^{1/2}\\
		& \le C_{p,\gamma} e^{-\frac{\epsilon N^2}{8}} \sup_{0\le s \le t}\|h(s)\|_{L^p_vL^\infty_x}+ C_{p,\gamma,N} \mathcal{E}(F_0)^{1/2} + C_{p,\gamma,N,\Omega}\sup_{0\le s \le t}\|h(s)\|_{L^p_vL^\infty_x}^{\frac{1-\kappa}{2}}\mathcal{E}(F_0)^{\kappa/2},
	\end{align*}
	where we have used \eqref{RTIRE} in the last term.
	It follow that
	\begin{align} \label{LA57}
		\|L_4\|_{L^p_vL^\infty_x} &\le C_{p,\gamma} e^{-\frac{\epsilon N^2}{8}} \sup_{0\le s \le t}\|h(s)\|_{L^p_vL^\infty_x}^2+ C_{p,\gamma,N} \mathcal{E}(F_0)^{1/2}\sup_{0\le s \le t}\|h(s)\|_{L^p_vL^\infty_x}\nonumber\\
		& \quad  + C_{p,\gamma,N,\Omega}\sup_{0\le s \le t}\|h(s)\|_{L^p_vL^\infty_x}^{\frac{3-\kappa}{2}}\mathcal{E}(F_0)^{\kappa/2} \nonumber\\
		& \le \frac{C_{p,\gamma,\epsilon}}{N} \sup_{0\le s \le t}\left[\|h(s)\|_{L^p_vL^\infty_x}^2 +\|h(s)\|_{L^p_vL^\infty_x}^{3-\kappa}\right]+C_{p,\gamma,N,\Omega}\left[ \mathcal{E}(F_0)+\mathcal{E}(F_0)^{\kappa}\right]
	\end{align}
	For $L_5$, we first need to estimate $\tilde{L}_5$. We use the H\"older inequality and similar arguments in \eqref{LA21}, \eqref{LA22}, \eqref{LA24} to deduce
	\begin{align*}
		\tilde{L}_5 &\le C_{p,\gamma}\left(\int_{t_1'}^{s-\epsilon}e^{-\frac{\nu_0}{2}(s-s')}\int_{|\eta|\le N}(1+|\eta|)^{e(p,\gamma)-2\beta}\left|w\Gamma^+(f,f)(s',X(s)-\eta(s-s'),\eta)\right|^2d\eta ds'  \right)^{1/2}\\
		& \le \frac{C_{p,\gamma}}{N}\sup_{0\le s \le t}\|h(s)\|_{L^p_vL^\infty_x}^2+C_{p,\gamma,N,\Omega}\left[\mathcal{E}(F_0)^{1/2}\sup_{0\le s\le t}\|h(s)\|_{L^p_vL^\infty_x} + \mathcal{E}(F_0)^{\kappa/2}\sup_{0\le s\le t}\|h(s)\|_{L^p_vL^\infty_x}^{\frac{3-\kappa}{2}}\right],
	\end{align*}
	and it follows that
	\begin{align} \label{LA58}
		\|L_5\|_{L^p_vL^\infty_x} \le \frac{C_{p,\gamma}}{N} \sup_{0\le s \le t}\left[\|h(s)\|_{L^p_vL^\infty_x}^3+\|h(s)\|_{L^p_vL^\infty_x}^4 +\|h(s)\|_{L^p_vL^\infty_x}^{5-\kappa}\right]+C_{p,\gamma,N,\Omega}\left[ \mathcal{E}(F_0)+\mathcal{E}(F_0)^{\kappa}\right].
	\end{align}
	Inserting \eqref{LA54}, \eqref{LA55}, \eqref{LA56}, \eqref{LA57}, \eqref{LA58} into \eqref{LA53}, we derive
	\begin{align} \label{LA59}
		\|J_{2}\|_{L^p_vL^\infty_x} &\le C_p e^{-\frac{\nu_0}{4}t}\|h_0\|_{L^p_vL^\infty_x} \int_0^t\|h(s)\|_{L^p_vL^\infty_x}ds\nonumber\\
		& \quad + C_{p,\gamma}\Biggl\{\left(\epsilon + \frac{C_{p,\gamma,\epsilon,T_0}}{N^{1/2}} \right)\sup_{0\le s\le t}\biggl[\|h(s)\|_{L^p_vL^\infty_x}^2+\|h(s)\|_{L^p_vL^\infty_x}^3+\|h(s)\|_{L^p_vL^\infty_x}^4\nonumber\\
		& \qquad+\|h(s)\|_{L^p_vL^\infty_x}^{\frac{3p-2}{2p-2}}+\|h(s)\|_{L^p_vL^\infty_x}^{\frac{5p-4}{2p-2}}+\|h(s)\|_{L^p_vL^\infty_x}^{\frac{7p-6}{2p-2}}+\|h(s)\|_{L^p_vL^\infty_x}^{\frac{9p-8}{2p-2}}\biggr]\nonumber\\
		& \qquad +C_{p,\gamma,N,\epsilon,\Omega,T_0}\left[\mathcal{E}(F_0)+\mathcal{E}(F_0)^2+\mathcal{E}(F_0)^{\frac{p-2}{2p-2}}+\mathcal{E}(F_0)^{\frac{p-2}{p-1}}\right]\Biggr\}.
	\end{align}
	From the Minkowski integral inequality, we can easily compute that
	\begin{align} \label{LA62}
		 \left\|(S_{G_f}(t,0)h_0)(t,x,v) \right\|_{L^p_vL^\infty_x} \le C_p e^{-\frac{\nu_0}{4}t}\|h_0\|_{L^p_vL^\infty_x}
	\end{align}
	and
	\begin{align} \label{LA63}
		&\left\|\int_{t-\epsilon}^{t}\left|(S_{G_f}(t,s)K_wh(s))(t,x,v) \right|+\left|(S_{G_f}(t,s)w\Gamma^+(f,f)(s))(t,x,v) \right|ds\right\|_{L^p_vL^\infty_x}\nonumber\\
		& \le \epsilon C_p \sup_{0\le s \le t}\left[\|h(s)\|_{L^p_vL^\infty_x}+\|h(s)\|_{L^p_vL^\infty_x}^2\right],
	\end{align}
	where we have used the estimate \eqref{DampLpdecay3}.
	Plugging \eqref{LA52}, \eqref{LA59}, \eqref{LA62}, \eqref{LA63} into the estimate in Lemma \ref{LALinftyestimate1}, we completes the proof of this lemma.
\end{proof}

\bigskip

\subsection{Small perturbation result with smallness for relative entropy} \label{endofsmall}

 We now demonstrate the small perturbation problem with the smallness condition for relative entropy.  

\begin{proof} [\textbf{Proof of Theorem \ref{Mainresult1}}]
	Set $C_4 :=\max\{C_1 \mathcal{T}^{\frac{5}{4}},C_3 \}$ and let us fix $\eta_0>0$ to satisfy the conditions \eqref{etacond10}, \eqref{etacond1}, \eqref{etacond2}, and \eqref{etacond3}, which will be specified later. Take
		\begin{align} \label{LA64}
			\bar{\eta} := 4 C_4^2 \eta_0.
		\end{align}
		Using the a priori assumption \eqref{Aprioriassump1} and Lemma \ref{LALinftyestimate2} and taking $T_0=\infty$ in Lemma \ref{LALinftyestimate2}, we get
		\begin{align} \label{LA66}
			\|h(t)\|_{_{L^p_vL^\infty_x}} \le C_4 \eta_0 \left(1+\int_0^t \|h(s)\|_{L^p_vL^\infty_x}ds\right)e^{-\frac{\nu_0}{4}t}+E \quad \text{for all } t\ge 0,
		\end{align}
		where
		\begin{align*}
			E:&= C_4\bigg\{\left(\epsilon + \frac{C_{p,\gamma,\epsilon}}{N^{1/2}} \right) \biggl[\bar{\eta}+\bar{\eta}^2+\bar{\eta}^3+\bar{\eta}^4+\bar{\eta}^{\frac{p}{2p-2}}+\bar{\eta}^{\frac{3p-2}{2p-2}}+\bar{\eta}^{\frac{5p-4}{2p-2}}+\bar{\eta}^{\frac{7p-6}{2p-2}}+\bar{\eta}^{\frac{9p-8}{2p-2}}\biggr]\\
			& \qquad +C_{p,\gamma,N,\Omega}\left[\mathcal{E}(F_0)+\mathcal{E}(F_0)^2+\mathcal{E}(F_0)^{\frac{p-2}{2p-2}}+\mathcal{E}(F_0)^{\frac{p-2}{p-1}}\right]\bigg\}.
		\end{align*}
		We define
		\begin{align*}
			G(t) := 1+\int_0^t \|h(s)\|_{L^p_vL^\infty_x}ds.
		\end{align*}
		Then the inequality \eqref{LA66} becomes
		\begin{align} \label{LA67}
			G'(t) \le C_4 \eta_0 e^{-\frac{\nu_0}{4}t} G(t) + E.
		\end{align}
		By Gr\"{o}nwall's inequality, \eqref{LA67} implies that
		\begin{equation} \label{LA68}
		\begin{aligned}
			G(t) &\le (1+Et) \exp\left\{\frac{4}{\nu_0}C_4\eta_0  \left(1-e^{-\frac{\nu_0}{4}t}\right)\right\} \le (1+Et) \exp\left\{\frac{4}{\nu_0}C_4\eta_0 \right\}
		\end{aligned}
		\end{equation}
		for all $t\ge 0$.\\
		Inserting \eqref{LA68} into \eqref{LA66} and since 
		\begin{align} \label{etacond1}
			\exp\left\{\frac{4}{\nu_0}C_4\eta_0 \right\}<2,
		\end{align}
		we deduce for $t\ge 0$
		\begin{align*}
			\|h(t)\|_{L^p_vL^\infty_x} &\le C_4 \eta_0  \exp\left\{\frac{4}{\nu_0}C_4\eta_0 \right\}(1+Et)e^{-\frac{\nu_0}{4}t}+E\\
			& \le \frac{1}{4C_4}\bar{\eta}(1+Et)e^{-\frac{\nu_0}{4}t} +E\\
			& \le \frac{1}{4C_4}\bar{\eta}\left(1+\frac{8}{\nu_0} E\right)e^{-\frac{\nu_0}{8}t} +E.
		\end{align*}
		We first choose $\epsilon >0$ small enough, then choose $N>0$ sufficiently large, and assume $\mathcal{E}(F_0) \le \epsilon_0(\eta_0)$ with small enough $\epsilon_0(\eta_0)<\epsilon_1$, which is determined in Lemma \ref{Rfest1},  so that
		\begin{align*}
			E \le \min\left\{\frac{\nu_0}{32  } , \frac{\eta_0}{8}\right\},
		\end{align*}
		and it follows that
		\begin{align} \label{LA69}
			\|h(t)\|_{L^p_vL^\infty_x} \le \frac{5}{16C_4}\bar{\eta}e^{-\frac{\nu_0}{8}t} + \frac{\eta_0}{8} \le \frac{1}{2C_4} \bar{\eta}
		\end{align}
		for all $ t \ge 0$.\\
		Hence we have closed the a priori assumption over $t \in [0,\infty)$ if $\mathcal{E}(F_0) \le \epsilon_0$.\\
		We claim that a solution to the Boltzmann equation \eqref{Boltzmanneq} extends into time interval $[0,\infty)$. From Lemma \ref{Localexistence}, there exists the Boltzmann solution $F(t) \ge 0$ to \eqref{Boltzmanneq} on $[0,\hat{t}_0]$ such that
		\begin{align} \label{LA70}
			\sup_{0\le t \le  \hat{t}_0} \|h(t)\|_{L^p_vL^\infty_x} \le 2C_4 \|h_0\|_{L^p_vL^\infty_x} \le \frac{1}{2C_4} \bar{\eta}.
		\end{align}
		We define $t_*:=\left(\hat{C}_{\mathcal{T}}\left[1+(2C_4)^{-1}\bar{\eta}\right]\right)^{-1}>0$, where $\hat{C}_{\mathcal{T}}$ is a constant in Theorem \ref{Localexistence}. Taking $t=\hat{t}_0$ as the initial time, it follows from \eqref{LA70} and Theorem \ref{Localexistence} that we can extend the Boltzmann equation solution $F(t) \ge 0$ into time interval $[0, \hat{t}_0+t_*]$ satisfying
		\begin{align*}
			\sup_{\hat{t}_0 \le t \le \hat{t}_0+t_*} \|h(t)\|_{L^p_vL^\infty_x} \le 2C_4 \|h(\hat{t}_0)\|_{L^p_vL^\infty_x} \le \bar{\eta}.
		\end{align*}
		Thus we have
		\begin{align} \label{LA71}
			\sup_{0 \le t \le \hat{t}_0+t_*} \|h(t)\|_{L^p_vL^\infty_x} \le \bar{\eta}.
		\end{align}
		Note that \eqref{LA71} means $h(t)$ satisfies the a priori assumption \eqref{Aprioriassump1} over $[0,\hat{t}_0+t_*]$.
		From \eqref{LA69}, we can obtain
		\begin{align*}
			\sup_{0 \le t \le \hat{t}_0+t_*} \|h(t)\|_{L^p_vL^\infty_x} \le\frac{1}{2C_4}\bar{\eta}.
		\end{align*}
		Repeating the same process infinitely many times, we can derive that there exists the Boltzmann equation solution $F(t)\ge 0$ on the time interval $[0,\infty)$ such that
		\begin{align*}
			\sup_{0\le t < \infty} \|h(t)\|_{L^p_vL^\infty_x} \le \frac{1}{2C_4} \bar{\eta}.
		\end{align*}
		\newline
		$\mathbf{(Exponential\ decay)}$\\
		We first show that there exist $\tilde{T}_0>0$ and $C_{p,\Omega,\tilde{T}_0}>0$ such that 
		\begin{equation} \label{finitetimeest}
	\begin{aligned}
			\|h(\tilde{T}_0)\|_{L^p_vL^\infty_x} \le e^{-\frac{\nu_0}{16} \tilde{T}_0}\|h_0\|_{L^p_vL^\infty_x} +C_{p,\Omega,\tilde{T}_0} \int_0^{\tilde{T}_0}\|f(s')\|_{L^2_{x,v}}ds'.
		\end{aligned}
\end{equation}
		Let $0\le t \le \tilde{T}_0$. Note that from the a priori assumption \eqref{Aprioriassump1} and Lemma \ref{Rfest1}, we have
\begin{align*}
		R(f)(t,x,v) \ge \frac{1}{2}\nu(v)
	\end{align*}
	for all $(t,x,v) \in [0,\tilde{T}_0] \times \Omega \times \R^3$. Thus it holds that
	\begin{align} \label{Ibound2}
		I^f(t,s) \le e^{-\frac{\nu_0}{2}(t-s)}
	\end{align}
	for all $0\le s \le t \le \tilde{T}_0$. We apply the Duhamel principle to the equation \eqref{WFPBE2} to obtain
		\begin{align*}
			h(t,x,v) &= (S_{G_f}(t,0) h_0)(t,x,v)\\
			& \quad  +\int_0^t (S_{G_f}(t,s)w\Gamma^+(f,f)(s))(t,x,v)ds\\
			& \quad + \int_0^t (S_{G_f}(t,s)K_wh(s))(t,x,v)ds.
		\end{align*}
		From Lemma \ref{Rfest1}, we derive
	\begin{align*}
			|h(t,x,v)| &\le \left|(S_{G_f}(t,0) h_0)(t,x,v)\right|\\
			& \quad  +\int_0^t \left|(S_{G_f}(t,s)w\Gamma^+(f,f)(s))(t,x,v)\right|ds\\
			& \quad + \int_0^t \left|(S_{G_f}(t,s)K_wh(s))(t,x,v)\right|ds\\
			&=:I_1+I_2+I_3.
		\end{align*}
		For $I_1$, by \eqref{DampLpdecay3}, we can easily get
	\begin{align} \label{df1}
		\|I_1\|_{L^p_vL^\infty_x} \le \left\|S_{G_f}(t,0)h_0\right\|_{L^p_vL^\infty_x} \le C_{p,\gamma,\mathcal{T}} e^{-\frac{\nu_0}{4}t}\|h_0\|_{L^p_vL^\infty_x}.
	\end{align}
		For $I_2$, thanks to \eqref{DampLpdecay3} and Corollary \ref{LpGamma+est}, we obtain
		\begin{align} \label{df2}
			\|I_2\|_{L^p_vL^\infty_x} &\le \int_0^t \left\|(S_{G_f}(t,s)w\Gamma^+(f,f))(s)\right\|_{L^p_vL^\infty_x}ds\nonumber\\
			& \le C_{p,\gamma,\mathcal{T}} \int_0^t e^{-\frac{\nu_0}{4}(t-s)} \|h(s)\|^2_{L^p_vL^\infty_x}ds\nonumber\\
			& \le C_{p,\gamma,\mathcal{T}}\sup_{0\le s \le t}\left[e^{-\frac{\nu_0}{8}(t-s)}\|h(s)\|^2_{L^p_vL^\infty_x} \right],
		\end{align}
		For $I_3$, we use the Duhamel principle once again to $h(s)$:
	\begin{align*}
			I_3 &\le \int_0^t \left|(S_{G_f}(t,s)K_wS_{G_f}(s,0)h_0)(t,x,v)\right|ds\\
			& \quad + \int_0^t \int_0^s\left|(S_{G_f}(t,s)K_wS_{G_f}(s,s')w\Gamma^+(f,f)(s'))(t,x,v)\right|ds'ds\\
			& \quad + \int_0^t \int_0^s\left|(S_{G_f}(t,s)K_wS_{G_f}(s,s')K_wh(s'))(t,x,v)\right|ds'ds\\
			&=:I_{31}+I_{32}+I_{33}.
		\end{align*}
		For $I_{31}$, it follows from \eqref{DampLpdecay3} and Lemma \ref{KESTIMATE} that
		\begin{align} \label{df4}
			\|I_{31}\|_{L^p_vL^\infty_x} &\le  \int_0^t \left\|S_{G_f}(t,s)K_wS_{G_f}(s,0)h_0\right\|_{L^p_vL^\infty_x}ds	 \nonumber\\
			& \le C_{p,\gamma,\mathcal{T}} t e^{-\frac{\nu_0}{4}t}\|h_0\|_{L^p_vL^\infty_x}.	
		\end{align}
		For $I_{32}$, we use \eqref{DampLpdecay3}, Lemma \ref{KESTIMATE}, and Corollary \ref{LpGamma+est} to obtain
		\begin{align} \label{df5}
			\|I_{32}\|_{L^p_vL^\infty_x} &\le\int_0^t \int_0^s\left\|S_{G_f}(t,s)K_wS_{G_f}(s,s')w\Gamma^+(f,f)(s')\right\|_{L^p_vL^\infty_x}ds'ds\nonumber\\
			& \le  C_{p,\gamma,\mathcal{T}} \int_0^t \int_0^s e^{-\frac{\nu_0}{4}(t-s')}\|h(s')\|_{L^p_vL^\infty_x}ds'ds\nonumber\\
			& \le C_{p,\gamma,\mathcal{T}}\sup_{0\le s' \le t}\left[e^{-\frac{\nu_0}{8}(t-s')}\|h(s')\|^2_{L^p_vL^\infty_x} \right].
		\end{align}
		It remains to estimate $I_{33}$. We can divide $I_{33}$ into three integrations as follows:
	\begin{align*} 
		I_{33}
		&\le\int_{t-\epsilon}^t\int_0^s \left|S_{G_f}(t,s)K_wS_{G_f}(s,s')K_wh(s')\right|ds'ds\nonumber\\
		& \quad + \int_0^{t-\epsilon}\int_{s-\epsilon}^s  \left|S_{G_f}(t,s)K_wS_{G_f}(s,s')K_wh(s')\right|ds'ds\nonumber\\
		& \quad + \int_0^{t-\epsilon}\int_0^{s-\epsilon} \left|S_{G_f}(t,s)K_wS_{G_f}(s,s')K_wh(s')\right|ds'ds\nonumber\\
		&=: J_1+J_2+J_3.
	\end{align*}
	For $J_1$, we can compute 
	\begin{align} \label{df8}
		\|J_1\|_{L^p_vL^\infty_x} 
		& \le C_{p,\gamma,\mathcal{T}} \int_{t-\epsilon}^t\int_0^se^{-\frac{\nu_0}{4}(t-s')}\|h(s')\|_{L^p_vL^\infty_x}ds'ds\nonumber\\
		& \le C_{p,\gamma,\mathcal{T}} \sup_{0 \le s' \le t }\left[e^{-\frac{\nu_0}{8}(t-s')}\|h(s')\|_{L^p_vL^\infty_x}\right] \int_{t-\epsilon}^{t}\int_0^s e^{-\frac{\nu_0}{8}(t-s')}ds'ds\nonumber\\
		& \le \epsilon C_{p,\gamma,\mathcal{T}} \sup_{0 \le s' \le t }\left[e^{-\frac{\nu_0}{8}(t-s')}\|h(s')\|_{L^p_vL^\infty_x}\right],
	\end{align}
	where we have used the Minkowski's integral inequality, the estimate \eqref{DampLpdecay3}, and Lemma \ref{KESTIMATE}.
	Here, we can calculate
	\begin{align*}
		\int_{t-\epsilon}^{t}\int_0^s e^{-\frac{\nu_0}{8}(t-s')}ds'ds \lesssim \int_{t-\epsilon}^{t} e^{-\frac{\nu_0}{8}(t-s)}ds \lesssim\epsilon.
	\end{align*}
	For $J_2$, by a similar way to estimate $J_1$, we can deduce
	\begin{align} \label{df9}
		\|J_2\|_{L^p_vL^\infty_x} & \le C_p \sup_{0 \le s' \le t }\left[e^{-\frac{\nu_0}{8}(t-s')}\|h(s')\|_{L^p_vL^\infty_x}\right] \int_{0}^{t}\int_{s-\epsilon}^s e^{-\frac{\nu_0}{8}(t-s')}ds'ds\nonumber\\
		& \le \epsilon C_p \sup_{0 \le s' \le t }\left[e^{-\frac{\nu_0}{8}(t-s')}\|h(s')\|_{L^p_vL^\infty_x}\right].
	\end{align}
	We now consider the term $J_3$. Fix $(t,x,v)$ so that $(x,v) \notin \gamma_0$. From Lemma \ref{RepresentationforDiffuse}, we have
	\begin{align*}
		&S_{G_f}(t,s)K_wS_{G_f}(s,s')K_wh(s')\\
		& = I^f(t,s)\mathbf{1}_{\{t_{1}\le s \}}\left\{K_wS_{G_f}(s,s')K_wh(s')\right\}(s,x-v(t-s),v)\\
		&\quad +  \frac{I^f(t,t_1)}{\tilde{w}(v)} \sum_{l=1}^{k-1}\int_{\prod_{j=1}^{k-1}\mathcal{V}_j} \mathbf{1}_{\{t_{l+1}\le s < t_l\}} \left\{K_wS_{G_f}(s,s')K_wh(s')\right\}(s,x_l-v_l(t_l-s),v_l) d\Sigma_l^f(s)\\
		&\quad + \frac{I^f(t,t_1)}{\tilde{w}(v)}\int_{\prod_{j=1}^{k-1}\mathcal{V}_j} \mathbf{1}_{\{t_k>s\}} \left\{S_{G_f}(t,s)K_wS_{G_f}(s,s')K_wh(s')\right\}(t_k,x_k,v_{k-1}) d\Sigma_{k-1}^f(t_k)\\
		&=: \tilde{M}_1+\tilde{M}_2+\tilde{M}_3,
	\end{align*}
	where
	\begin{align*}
		d\Sigma_{l}^f(s) = \left\{\prod_{j=l+1}^{k-1}d\sigma_j \right\}\left\{I^f(t_l,s) \tilde{w}(v_l)d\sigma_l \right\}\prod_{j=1}^{l-1} \left\{ I^f(t_j,t_{j+1})d\sigma_j \right\},
	\end{align*}
	and the exponential factor in $d \Sigma_l^f(s)$ is bounded by $e^{-\frac{\nu_0}{2}(t_1-s)}$.\\
	Then the term $J_3$ can be written as
	\begin{align*} 
		J_3 &\le  \int_0^{t-\epsilon}\int_0^{s-\epsilon} |\tilde{M}_1| ds'ds+\int_0^{t-\epsilon}\int_0^{s-\epsilon} |\tilde{M}_2| ds'ds+\int_0^{t-\epsilon}\int_0^{s-\epsilon} |\tilde{M}_3| ds'ds\nonumber\\
		&=:M_1+M_2+M_3.
	\end{align*}
	For $M_3$, by \eqref{DampLpdecay3} and Lemma \ref{KESTIMATE}, we compute
	\begin{equation} \label{ppp0}
	\begin{aligned}
		&\left\|\left\{S_{G_f}(t,s)K_wS_{G_f}(s,s')K_wh(s')\right\}(t_k)\right\|_{L^p_vL^\infty_x}\le C_{p,\gamma,\mathcal{T}} e^{-\frac{\nu_0}{4} (t_k-s')}\left\|h(s')\right\|_{L^p_vL^\infty_x}.
	\end{aligned}
	\end{equation}
	Since $t-s \ge \epsilon >0$, from Lemma \ref{Lsmall}, it follows from \eqref{Ibound2} and \eqref{ppp0} that, for $k=k(\epsilon,\tilde{T}_0)+1$ such that \eqref{Lsmallineq} holds,
	\begin{align*}
		|\tilde{M}_3| &\le C_{p,\gamma,\mathcal{T}}\frac{e^{-\frac{\nu_0}{2}(t-t_k)}}{\tilde{w}(v)}e^{-\frac{\nu_0}{4} (t_k-s')}\left\|h(s')\right\|_{L^p_vL^\infty_x}\int_{\prod_{j=1}^{k-2}\mathcal{V}_j}  \mathbf{1}_{\{t_{k}>s\}}\left(\int_{\R^3}\left|\frac{|v_{k-1}|}{w(v_{k-1})}\mu^{1/2}(v_{k-1}) \right|^{p'} dv_{k-1}\right)^{1/p'} \\
		& \quad \times  \prod_{j=1}^{k-2}d\sigma_j\\
		& \le C_{p,\gamma,\mathcal{T}} \frac{e^{-\frac{\nu_0}{4} (t-s')}}{\tilde{w}(v)}\left\|h(s')\right\|_{L^p_vL^\infty_x}\left(\int_{\prod_{j=1}^{k-2}\mathcal{V}_j} \mathbf{1}_{\{t_{k-1}>s\}} \prod_{j=1}^{k-2}d\sigma_j\right)\\
		& \le \epsilon C_{p,\gamma,\mathcal{T}} \frac{e^{-\frac{\nu_0}{4} (t-s')}}{\tilde{w}(v)}\left\|h(s')\right\|_{L^p_vL^\infty_x},
	\end{align*}
	which implies that
	\begin{align*}
		M_3 &\le \frac{\epsilon C_{p,\gamma,\mathcal{T}}}{\tilde{w}(v)} \int_0^{t}\int_0^{s}e^{-\frac{\nu_0}{4} (t-s')}\left\|h(s')\right\|_{L^p_vL^\infty_x}ds'ds\\
		&\le \frac{\epsilon C_{p,\gamma,\mathcal{T}}}{\tilde{w}(v)} \sup_{0 \le s' \le t }\left[e^{-\frac{\nu_0}{8}(t-s')}\|h(s')\|_{L^p_vL^\infty_x}\right].
	\end{align*}
	This yields
	\begin{align} \label{df10}
		\left\|M_3\right\|_{L^p_vL^\infty_x}\le \epsilon C_{p,\gamma,\mathcal{T}}\sup_{0 \le s' \le t  }\left[e^{-\frac{\nu_0}{8}(t-s')}\|h(s')\|_{L^p_vL^\infty_x}\right].
	\end{align}
	We now consider the terms $M_1$ and $M_2$. We denote $X_{\mathbf{cl}}(s) = x_l-v_l(t_l-s)$, $t_{l'}' = t_{l'}(s,X_{\mathbf{cl}}(s),v')$, $x_{l'}' = X_{\mathbf{cl}}(t_{l'}',X_{\mathbf{cl}}(s),v')$, and $v_{l'}' \in \mathcal{V}'_{l'}$ for $1\le l' \le k-1$. We derive from Lemma \ref{RepresentationforDiffuse} the formula for $K_wS_{G_f}(s,s')K_wh(s')$:
	\begin{align*}
		&\left\{K_wS_{G_f}(s,s')K_wh(s')\right\}(s,X_{\mathbf{cl}}(s),v_l)\\
		& = \int_{\mathbb{R}^3} k_w(v_l,v')\left\{S_{G_f}(s,s')K_wh(s')\right\}(s,X_{\mathbf{cl}}(s),v')dv'\\
		& = \int_{\mathbb{R}^3} k_w(v_l,v')I^f(s,s') \mathbf{1}_{\{t_1' \le s'\}}\left\{K_wh(s')\right\}(s',X_{\mathbf{cl}}(s)-v'(s-s'),v')dv'\\
		& \quad + \int_{\mathbb{R}^3} k_w(v_l,v')\frac{I^f(s,t_1')}{\tilde{w}(v')} \sum_{l'=1}^{k-1}\int_{\prod_{j=1}^{k-1}\mathcal{V}_j'} \mathbf{1}_{\{t_{l'+1}'\le s' < t_{l'}'\}}\left\{\int_{\mathbb{R}^3}k_w(v_{l'}',v'')h(s',x_{l'}'-v_{l'}'(t_{l'}'-s'),v'')dv''\right\}\\
		& \qquad \times d\Sigma_{l'}'^f(s')dv'\\
		& \quad +  \int_{\mathbb{R}^3}k_w(v_l,v')\frac{I^f(s,t_1')}{\tilde{w}(v')} \int_{\prod_{j=1}^{k-1}\mathcal{V}_j'}\mathbf{1}_{\{t_k' > s'\}}\left\{S_{G_f}(s,s')K_wh(s')\right\}(t_k',x_k',v_{k-1}')d\Sigma_{k-1}'^f(t_k')dv'\\
		&=: L_1(v_l)+L_2(v_l)+L_3(v_l),
	\end{align*}
	where 
	\begin{align*}
		d\Sigma_{l'}'^f(s) = \left\{\prod_{j=l'+1}^{k-1}d\sigma'_j \right\}\left\{I^f(t'_{l'},s') \tilde{w}(v'_{l'})d\sigma'_{l'} \right\}\prod_{j=1}^{l'-1} \left\{ I^f(t'_j,t'_{j+1})d\sigma'_j \right\}	
	\end{align*}
	with  $d\sigma'_j :=c_\mu \mu(v'_j) \{n(x'_j) \cdot v'_j\}dv'_j$.\\
	For $L_3$, from \eqref{DampLpdecay3} and Lemma \ref{KESTIMATE}, we have
	\begin{align} \label{ppp10}
		\left\|\left\{S_{G_f}(s,s')K_wh(s')\right\}(t_k')\right\|_{L^p_vL^\infty_x} \le C_{p,\gamma,\mathcal{T}}e^{-\frac{\nu_0}{4}(t_k'-s')} \|h(s')\|_{L^p_vL^\infty_x}.
	\end{align}
	From Lemma \ref{Lsmall}, it follows from \eqref{Ibound2} and \eqref{ppp10} that, for $k=k(\epsilon,\tilde{T}_0)+1$ such that \eqref{Lsmallineq} holds,
	\begin{align} \label{L3contribution}
		|L_3(v_l)| &\le  C_{p,\gamma,\mathcal{T}} \int_{\mathbb{R}^3}|k_w(v_l,v')|  \frac{e^{-\frac{\nu_0}{2}(s-t_k')}}{\tilde{w}(v')}e^{-\frac{\nu_0}{4}(t_k'-s')} \|h(s')\|_{L^p_vL^\infty_x}\nonumber\\
		&\quad \times \int_{\prod_{j=1}^{k-2}\mathcal{V}_j'}\mathbf{1}_{\{t_k' > s'\}} \left(\int_{\R^3}\left|\frac{|v_{k-1}'|}{w(v_{k-1}')}\mu^{1/2}(v_{k-1}') \right|^{p'} dv_{k-1}'\right)^{1/p'} \prod_{j=1}^{k-2}d\sigma_j' dv'\nonumber\\
		&\le C_{p,\gamma,\mathcal{T}} e^{-\frac{\nu_0}{4}(s-s')} \|h(s')\|_{L^p_vL^\infty_x}\int_{\mathbb{R}^3}|k_w(v_l,v')|\left(\int_{\prod_{j=1}^{k-2}\mathcal{V}_j'} \mathbf{1}_{\{t_{k-1}'>s'\}} \prod_{j=1}^{k-2}d\sigma_j'\right) dv'\nonumber\\
		&\le \epsilon C_{p,\gamma,\mathcal{T}} e^{-\frac{\nu_0}{4}(s-s')} \|h(s')\|_{L^p_vL^\infty_x}\int_{\mathbb{R}^3}|k_w(v_l,v')| dv'\nonumber\\
		&\le \epsilon \frac{C_{p,\gamma,\mathcal{T}}}{(1+|v_l|)^{2-\gamma}} e^{-\frac{\nu_0}{4}(s-s')} \|h(s')\|_{L^p_vL^\infty_x}.
	\end{align}
	Using the estimate \eqref{Ibound2} and \eqref{L3contribution}, the contribution of $M_1$ and $L_3$ is 
	\begin{align*}
		&\int_0^{t-\epsilon}\int_0^{s-\epsilon} I^f(t,s)\mathbf{1}_{\{t_{1}\le s \}} |L_3(v)|ds'ds\\
		&\le \epsilon \frac{C_{p,\gamma,\mathcal{T}}}{(1+|v|)^{2-\gamma}} \sup_{0 \le s' \le t  }\left[e^{-\frac{\nu_0}{8}(t-s')}\|h(s')\|_{L^p_vL^\infty_x}\right] \int_0^{t}\int_0^{s}  e^{-\frac{\nu_0}{8} (t-s')} ds'ds\\
		&\le \epsilon \frac{C_{p,\gamma,\mathcal{T}}}{(1+|v|)^{2-\gamma}} \sup_{0 \le s' \le t  }\left[e^{-\frac{\nu_0}{8}(t-s')}\|h(s')\|_{L^p_vL^\infty_x}\right],
	\end{align*}
	which implies that
	\begin{align} \label{df11}
		\left\| \int_0^{t-\epsilon}\int_0^{s-\epsilon} I^f(t,s)\mathbf{1}_{\{t_{1}\le s \}} L_3(v)ds'ds\right\|_{L^p_vL^\infty_x} \le \epsilon C_{p,\gamma,\mathcal{T}} \sup_{0 \le s' \le t  }\left[e^{-\frac{\nu_0}{8}(t-s')}\|h(s')\|_{L^p_vL^\infty_x}\right].  
	\end{align}
	Using the estimate \eqref{Ibound2} and \eqref{L3contribution}, the contribution of $M_2$ and $L_3$ also is
	\begin{align*}
		&\int_0^{t-\epsilon}\int_0^{s-\epsilon}\frac{I^f(t,t_1)}{\tilde{w}(v)} \sum_{l=1}^{k-1}\int_{\prod_{j=1}^{k-1}\mathcal{V}_j} \mathbf{1}_{\{t_{l+1}\le s < t_l\}} |L_3(v_l)|d\Sigma_l^f(s)ds'ds\\
		& \le \epsilon C_{p,\gamma,\mathcal{T}}  \sum_{l=1}^{k-1}\int_{t_{l+1}}^{t_l}\int_0^{s} \frac{e^{-\frac{\nu_0}{2}(t-t_1)}}{\tilde{w}(v)} \int_{\prod_{j=1}^{l}\mathcal{V}_j} e^{-\frac{\nu_0}{4} (t_1-s')}\|h(s')\|_{L^p_vL^\infty_x}\tilde{w}(x_l,v_l)\prod_{j=1}^{l}d\sigma_jds'ds\\
		& \le \epsilon \frac{C_{p,\gamma,\mathcal{T}}}{\tilde{w}(v)} \sum_{l=1}^{k-1}\int_{t_{l+1}}^{t_l}\int_0^{s} e^{-\frac{\nu_0}{4} (t-s')}\|h(s')\|_{L^p_vL^\infty_x}ds'ds\\
		& \le \epsilon \frac{C_{p,\gamma,\mathcal{T}}}{\tilde{w}(v)}\sup_{0 \le s' \le t  }\left[e^{-\frac{\nu_0}{8}(t-s')}\|h(s')\|_{L^p_vL^\infty_x}\right]\int_0^{t}\int_0^{s}  e^{-\frac{\nu_0}{8} (t-s')} ds'ds\\ 
		& \le \epsilon \frac{C_{p,\gamma,\mathcal{T}}}{\tilde{w}(v)}\sup_{0 \le s' \le t  }\left[e^{-\frac{\nu_0}{8}(t-s')}\|h(s')\|_{L^p_vL^\infty_x}\right],
	\end{align*}
	which implies that
	\begin{align} \label{df12}
		& \left\| \int_0^{t-\epsilon}\int_0^{s-\epsilon}\frac{I^f(t,t_1)}{\tilde{w}(v)} \sum_{l=1}^{k-1}\int_{\prod_{j=1}^{k-1}\mathcal{V}_j} \mathbf{1}_{\{t_{l+1}\le s < t_l\}} L_3(v_l)d\Sigma_l^f(s)ds'ds \right\|_{L^p_vL^\infty_x}\nonumber\\
		& \le \epsilon C_{p,\gamma,\mathcal{T}} \sup_{0 \le s' \le t  }\left[e^{-\frac{\nu_0}{8}(t-s')}\|h(s')\|_{L^p_vL^\infty_x}\right]. 
	\end{align}
	The remaining terms are
	\begin{equation} \label{RMC}
	\begin{aligned}
		&\int_0^{t-\epsilon}\int_0^{s-\epsilon} I^f(t,s) \mathbf{1}_{\{t_{1}\le s \}}\int_{\mathbb{R}^3}k_w(v,v') I^f(s,s') \mathbf{1}_{\{t_1' \le s'\}} \left\{K_wh(s')\right\}(s',x-v(t-s)-v'(s-s'),v')\\
		&\quad \times dv'ds'ds\\
		&+\int_0^{t-\epsilon}\int_0^{s-\epsilon}  I^f(t,t_1)  \mathbf{1}_{\{t_{1}\le s \}}\int_{\mathbb{R}^3} k_w(v,v')\frac{I^f(s,t_1')}{\tilde{w}(v')}\\ 
		&\quad \times \sum_{l'=1}^{k-1}\int_{\prod_{j=1}^{k-1}\mathcal{V}_j'} \mathbf{1}_{\{t_{l'+1}'\le s' < t_{l'}'\}} \left\{\int_{\mathbb{R}^3}k_w(v_{l'}',v'')h(s',x_{l'}'-v_{l'}'(t_{l'}'-s'),v'')dv''\right\}d\Sigma_{l'}'^f(s')dv'ds'ds\\
		&+\int_0^{t-\epsilon}\int_0^{s-\epsilon} \frac{I^f(t,t_1)}{\tilde{w}(v)} \sum_{l=1}^{k-1}\int_{\prod_{j=1}^{k-1}\mathcal{V}_j} \mathbf{1}_{\{t_{l+1}\le s < t_l\}} \int_{\mathbb{R}^3}  k_w(v_l,v')I^f(s,s') \mathbf{1}_{\{t_1' \le s'\}}\\
		&\quad \times \left\{K_wh(s')\right\}(s',x_l-v_l(t_l-s)-v'(s-s'),v')dv'd\Sigma_l^f(s)ds'ds\\
		&+\int_0^{t-\epsilon}\int_0^{s-\epsilon} \frac{I^f(t,t_1)}{\tilde{w}(v)} \sum_{l=1}^{k-1}\int_{\prod_{j=1}^{k-1}\mathcal{V}_j} \mathbf{1}_{\{t_{l+1}\le s < t_l\}} \int_{\mathbb{R}^3} k_w(v_l,v')\frac{I^f(s,t_1')}{\tilde{w}(v')}\\
		&\quad \times \sum_{l'=1}^{k-1}\int_{\prod_{j=1}^{k-1}\mathcal{V}_j'} \mathbf{1}_{\{t_{l'+1}'\le s' < t_{l'}'\}}\biggl\{\int_{\mathbb{R}^3}k_w(v_{l'}',v'')h(s',x_{l'}'-v_{l'}'(t_{l'}'-s'),v'')dv''\biggr\}d\Sigma_{l'}'^f(s')dv'd\Sigma_l^f(s)ds'ds\\
		&=: R_1+R_2+R_3+R_4.
	\end{aligned}
	\end{equation}
	\newline
	We denote $X'(s')=x-v(t-s)-v'(s-s')$.\\
	For $R_1$ in \eqref{RMC} :
	\begin{align*}
		&\int_{t_1}^{t-\epsilon}\int_{t_1'}^{s-\epsilon}   I^f(t,s) \int_{\mathbb{R}^3} k_w(v,v')   I^f(s,s')\int_{\mathbb{R}^3} k_w(v',v'')h(s',X'(s'),v'')dv''dv'ds'ds \nonumber \\ 
		&\le \int_{t_1}^{t-\epsilon}\int_{t_1'}^{s-\epsilon}e^{-\frac{\nu_0}{2}(t-s')}\int_{\mathbb{R}^3} \int_{\mathbb{R}^3}|k_w(v,v')| |k_w(v',v'')||h(s',X'(s'),v'')|dv''dv'ds'ds,
	\end{align*}
	where we have used the estimate \eqref{Ibound2}. By a similar way in \eqref{LA36} and \eqref{LA37}, we obtain
	\begin{align} \label{class11}
	R_1 &\le \frac{C_{p,\epsilon}}{N} \sup_{0 \le s' \le t  }\left[e^{-\frac{\nu_0}{4}(t-s')}\|h(s')\|_{L^p_vL^\infty_x}\right]\nonumber\\
	& \quad  +\int_{t_1}^{t-\epsilon}\int_{t_1'}^{s-\epsilon} e^{-\frac{\nu_0}{2}(t-s')}  \mathbf{1}_{\{|v|\le N\}} \int_{|v'|\le 2N} \int_{|v''|\le 3N} |k_w(v,v')||k_w(v',v'')| |h(s',X'(s'),v'')|dv''dv'\nonumber\\
	& \qquad \times ds'ds\nonumber\\
	& =: R_{11}+R_{12}.
\end{align}
Since $|k_w(v,v')|^{q}$ has possible integrable singularity of $\frac{1}{|v-v'|^{q}}$, where $1\le q<3$, we can choose $k_N(v,v')$ smooth with compact support such that
\begin{align} \label{df13}
	\sup_{|v|\le 3N} \int_{|v'|\le 3N}\left|k_N(v,v')- k_w(v,v')\right|^{q}dv' \le \frac{1}{N^{q+1}}. 
\end{align}
We split
\begin{align*}
	k_w(v,v')k_w(v',v'') &= \left\{k_w(v,v')- k_N(v,v')\right\}k_w(v',v'')\\
	&\quad +\left\{k_w(v',v'')- k_N(v',v'')\right\}k_N(v,v')\\
	& \quad +k_N(v,v')k_N(v',v'').
\end{align*}
We use Lemma \ref{KESTIMATE} and \eqref{df13} to estimate
\begin{align*}
	&\int_{t_1}^{t-\epsilon}\int_{t_1'}^{s-\epsilon} e^{-\frac{\nu_0}{2}(t-s')} \mathbf{1}_{\{|v|\le N\}} \int_{|v'|\le 2N} \int_{|v''|\le 3N} |k_w(v,v')-k_N(v,v')||k_w(v',v'')| |h(s',X'(s'),v'')|\\
	& \quad \times dv''dv'ds'ds\\
	& \le  C_p \mathbf{1}_{\{|v|\le N\}}\sup_{0 \le s' \le t  }\left[e^{-\frac{\nu_0}{4}(t-s')}\|h(s')\|_{L^p_vL^\infty_x}\right]\int_0^t\int_0^s e^{-\frac{\nu_0}{4}(t-s')}ds'ds \\
	&\quad \times \left\{\sup_{|v|\le N}\int_{|v'| \le 2N}|k_w(v,v')-k_N(v,v')|dv'\right\}\\
	& \le \frac{C_p}{N^{2}} \mathbf{1}_{\{|v|\le N\}}\sup_{0 \le s' \le t  }\left[e^{-\frac{\nu_0}{4}(t-s')}\|h(s')\|_{L^p_vL^\infty_x}\right].
\end{align*}
Similarly, we derive
\begin{align*}
	&\int_{t_1}^{t-\epsilon}\int_{t_1'}^{s-\epsilon} e^{-\frac{\nu_0}{2}(t-s')} \mathbf{1}_{\{|v|\le N\}} \int_{|v'|\le 2N} \int_{|v''|\le 3N} |k_N(v,v')||k_w(v',v'')-k_N(v',v'')| |h(s',X'(s'),v'')|\\
	& \quad \times dv''dv'ds'ds\\
	& \le \frac{C_p}{N^{\frac{p'+1}{p'}}} \mathbf{1}_{\{|v|\le N\}}\sup_{0 \le s' \le t  }\left[e^{-\frac{\nu_0}{4}(t-s')}\|h(s')\|_{L^p_vL^\infty_x}\right]\int_0^t\int_0^s e^{-\frac{\nu_0}{4}(t-s')}ds'ds  \left\{\int_{|v'|\le 2N}|k_N(v,v')|dv'\right\}\\
	& \le \frac{C_p}{N^{\frac{p'+1}{p'}}} \mathbf{1}_{\{|v|\le N\}}\sup_{0 \le s' \le t  }\left[e^{-\frac{\nu_0}{4}(t-s')}\|h(s')\|_{L^p_vL^\infty_x}\right].
\end{align*}
This implies
\begin{align*}
	R_{12} &\le \frac{C_p}{N^{\frac{p'+1}{p'}}} \mathbf{1}_{\{|v|\le N\}}\sup_{0 \le s' \le t  }\left[e^{-\frac{\nu_0}{4}(t-s')}\|h(s')\|_{L^p_vL^\infty_x}\right]\\
	& \quad +C_{p,N} \mathbf{1}_{\{|v|\le N\}}\int_{t_1}^{t-\epsilon}\int_{t_1'}^{s-\epsilon} e^{-\frac{\nu_0}{2}(t-s')}\int_{|v'|\le 2N} \int_{|v''|\le 3N} |h(s',X'(s'),v'')|dv''dv'ds'ds\\
	& =:A_1+A_2.
\end{align*}
We can easily get
\begin{align} \label{class13}
	\|A_1\|_{L^p_vL^\infty_x} &\le \frac{C_p}{N^{\frac{p'+1}{p'}}}\sup_{0 \le s' \le t  }\left[e^{-\frac{\nu_0}{4}(t-s')}\|h(s')\|_{L^p_vL^\infty_x}\right] \left(\int_{|v|\le N} dv\right)^{1/p}\nonumber\\
	& \le C_p\left(\frac{1}{N}\right)^{\frac{p'+1}{p'}-\frac{3}{p}}\sup_{0 \le s' \le t  }\left[e^{-\frac{\nu_0}{4}(t-s')}\|h(s')\|_{L^p_vL^\infty_x}\right].
\end{align}
We make a change of variables $v' \mapsto y=X'(s')=x-v(t-s)-v'(s-s')$ with $\left|\frac{dy}{dv'}\right| = (s-s')^3$:
\begin{align*}
	A_2 &\le C_{p,N,\epsilon} \mathbf{1}_{\{|v|\le N\}}\int_{t_1}^{t-\epsilon}\int_{t_1'}^{s-\epsilon} e^{-\frac{\nu_0}{2}(t-s')} \int_{\Omega} \int_{|v''|\le 3N} |h(s',y,v'')|dv''dyds'ds\\
	& \le C_{p,N,\epsilon,\Omega} \mathbf{1}_{\{|v|\le N\}}\int_0^t\int_0^{s} e^{-\frac{\nu_0}{2}(t-s')}\|f(s')\|_{L^2_{x,v}}ds'ds\\ 
	& \le C_{p,N,\epsilon,\Omega} \mathbf{1}_{\{|v|\le N\}} \int_0^{t} \|f(s')\|_{L^2_{x,v}}ds'     
\end{align*}
which implies that
\begin{align} \label{class14}
	\|A_2\|_{L^p_vL^\infty_x} \le C_{p,N,\epsilon,\Omega} \int_0^{t} \|f(s')\|_{L^2_{x,v}}ds'  
\end{align}
Inserting \eqref{class13} and \eqref{class14} into \eqref{class11}, we can derive
	\begin{align} \label{df14}
		\|R_1\|_{L^p_vL^\infty_x} \le \frac{C_{p,\epsilon}}{N}\sup_{0 \le s' \le t  }\left[e^{-\frac{\nu_0}{4}(t-s')}\|h(s')\|_{L^p_vL^\infty_x}\right]+C_{p,N,\epsilon,\Omega}\int_0^{t}\|f(s')\|_{L^2_{x,v}}ds'.
	\end{align}
	\newline
	For $R_2$ in \eqref{RMC} :
	\begin{align*} 
		&\int_{t_1}^{t-\epsilon} I^f(t,s) \int_{\mathbb{R}^3} k_w(v,v')\frac{I^f(s,t_1')}{\tilde{w}(v')}\nonumber\\
		&\quad \times \sum_{l'=1}^{k-1}\int_{t_{l'+1}'}^{t_{l'}'}\int_{\prod_{j=1}^{k-1}\mathcal{V}_j'}\left\{\int_{\mathbb{R}^3}k_w(v_{l'}',v'')h(s',x_{l'}'-v_{l'}'(t_{l'}'-s'),v'')dv''\right\}d\Sigma_{l'}'^f(s')ds'dv'ds \nonumber \\
		& \le \int_{t_1}^{t-\epsilon} e^{-\frac{\nu_0}{2}(t-t_1')}\int_{\mathbb{R}^3} |k_w(v,v')|\frac{1}{\tilde{w}(v')}\sum_{l'=1}^{k-1}\int_{t_{l'+1}'}^{t_{l'}'}e^{-\frac{\nu_0}{2}(t_1'-s')}\int_{\prod_{j=1}^{k-1}\mathcal{V}_j'} \nonumber \\ 
		&\quad \times  \Biggl\{\int_{\mathbb{R}^3}|k_w(v_{l'}',v'')| |h(s',x_{l'}'-v_{l'}'(t_{l'}'-s'),v'')|dv''\Biggr\}\left\{\prod_{j=l'+1}^{k-1}d\sigma_j'\right\}\left\{\tilde{w}(v'_{l'})d\sigma_{l'}'\right\}\nonumber\\
		& \quad \times \left\{\prod_{j=1}^{l'-1}d\sigma_j'\right\}ds'dv'ds,
	\end{align*}
	where we have used the estimate \eqref{Ibound2}. Fix $l'$. We divide this term into four cases.\\
	\newline
	$\mathbf{Case\ 1:}$ $|v| \ge N$ or $|v_{l'}'| \ge N$.\\
	From Lemma \ref{KESTIMATE}, we have
	\begin{align*}
		\int_{|v_{l'}'|\ge N} \left(\int_{\mathbb{R}^3} |k_w(v_{l'}',v'')||h(s',x_{l'}'-v_{l'}'(t_{l'}'-s'),v'')|dv''\right)\tilde{w}(v_{l'}')d\sigma_{l'}' \le \frac{C_p}{N} \|h(s')\|_{L^p_vL^\infty_x}.
	\end{align*}
	Then $R_2$ in the case $|v_{l'}'| \ge N$ is bounded by
	\begin{align*} 
		&\frac{C_p}{N} \int_0^{t-\epsilon} \int_{t_{l'+1}'}^{t_{l'}'} e^{-\frac{\nu_0}{2}(t-s')}\|h(s')\|_{L^p_vL^\infty_x} \int_{\mathbb{R}^3}|k_w(v,v')|dv'ds'ds\\
		& \le\frac{C_p}{N}\left(\frac{1}{1+|v|}\right)^{2-\gamma} \int_0^t\int_0^s e^{-\frac{\nu_0}{2}(t-s')}\|h(s')\|_{L^p_vL^\infty_x} ds'ds\\
		& \le \frac{C_p}{N}\left(\frac{1}{1+|v|}\right)^{2-\gamma} \sup_{0 \le s' \le t  }\left[e^{-\frac{\nu_0}{4}(t-s')}\|h(s')\|_{L^p_vL^\infty_x}\right],
	\end{align*}
	which implies the $L^p_vL^\infty_x$ norm of $R_2$ in this case is bounded by
	\begin{align} \label{class21}
		\frac{C_p}{N} \sup_{0 \le s' \le t  }\left[e^{-\frac{\nu_0}{4}(t-s')}\|h(s')\|_{L^p_vL^\infty_x}\right].
	\end{align}
	By Lemma \ref{KESTIMATE}, $R_2$ in the case $|v| \ge N$ is bounded by 
	\begin{align*}
		& C_p \int_0^{t-\epsilon}\int_{t_{l'+1}'}^{t_{l'}'}e^{-\frac{\nu_0}{2}(t-s')}\mathbf{1}_{\{|v|\ge N\}}\|h(s')\|_{L^p_vL^\infty_x}\int_{\mathbb{R}^3}|k_w(v,v')|\frac{1}{\tilde{w}(v')} dv'ds'ds\\
		& \le \frac{C_p}{N}\left(\frac{1}{1+|v|}\right)^2 \int_0^{t}\int_{0}^{s} e^{-\frac{\nu_0}{2}(t-s')} \|h(s')\|_{L^p_vL^\infty_x} ds'ds\\
		&\le \frac{C_p}{N}\left(\frac{1}{1+|v|}\right)^2 \sup_{0 \le s' \le t  }\left[e^{-\frac{\nu_0}{4}(t-s')}\|h(s')\|_{L^p_vL^\infty_x}\right],
	\end{align*}
	which implies the $L^p_vL^\infty_x$ norm of $R_2$ in this case is bounded by
	\begin{align} \label{class22}
		\frac{C_p}{N} \sup_{0 \le s' \le t  }\left[e^{-\frac{\nu_0}{4}(t-s')}\|h(s')\|_{L^p_vL^\infty_x}\right].
	\end{align}
	\newline
	$\mathbf{Case\ 2:}$ $|v| \le N$, $|v'| \ge 2N$ or $|v_{l'}'| \le N$, $|v''| \ge 2N$.\\
	Note that either $|v-v'| \ge N$ or $|v_{l'}'-v''| \ge N$. Then we use \eqref{class0} and \eqref{class1} to bound $R_2$ in the case $|v| \le N, |v'|\ge 2N$ by
	\begin{align*} 
		&C_p \int_0^{t-\epsilon}\int_{t_{l'+1}'}^{t_{l'}'}e^{-\frac{\nu_0}{2}(t-s')}\|h(s')\|_{L^p_vL^\infty_x}\mathbf{1}_{\{|v|\le N\}} \int_{|v'| \ge 2N}|k_w(v,v')| dv'ds'ds\nonumber\\
		&\le C_p e^{-\frac{\epsilon N^2}{8}} \left(\frac{1}{1+|v|}\right)^2\int_0^t\int_0^s e^{-\frac{\nu_0}{2}(t-s')}\|h(s')\|_{L^p_vL^\infty_x}ds'ds\nonumber\\
		&\le  C_p e^{-\frac{\epsilon N^2}{8}} \left(\frac{1}{1+|v|}\right)^2\sup_{0 \le s' \le t  }\left[e^{-\frac{\nu_0}{4}(t-s')}\|h(s')\|_{L^p_vL^\infty_x}\right],
	\end{align*}
	which implies the $L^p_vL^\infty_x$ norm of $R_2$ in this case is bounded by
	\begin{align} \label{class23}
		 C_p e^{-\frac{\epsilon N^2}{8}} \sup_{0 \le s' \le t  }\left[e^{-\frac{\nu_0}{4}(t-s')}\|h(s')\|_{L^p_vL^\infty_x}\right].
	\end{align}
	Similarly, we use \eqref{class0} and \eqref{class1} to bound $R_2$ in the case $|v_{l'}'| \le N, |v''| \ge 2N$ by 
	\begin{align*}
		&C_p \int_0^{t-\epsilon}\int_{t_{l'+1}'}^{t_{l'}'}e^{-\frac{\nu_0}{2}(t-s')}\|h(s')\|_{L^p_vL^\infty_x}\int_{\mathbb{R}^3}|k_w(v,v')| \int_{\prod_{j=1}^{k-1}\mathcal{V}_j' ,\ |v_{l'}'|\le N}\Biggl\{\int_{|v''| \ge 2N}|k_w(v_{l'}',v'')|^{p'}dv''\Biggr\}^{1/p'}\nonumber\\
		&\quad  \times  \left\{\prod_{j=l'+1}^{k-1}d\sigma_j'\right\}\left\{\tilde{w}(v_{l'}')d\sigma_{l'}'\right\}\left\{\prod_{j=1}^{l'-1}d\sigma_j'\right\}dv'ds'ds\nonumber\\
		&\le C_p e^{-\frac{\epsilon N^2}{8}} \int_0^t\int_0^se^{-\frac{\nu_0}{2}(t-s')}\|h(s')\|_{L^p_vL^\infty_x} \left\{\int_{\mathbb{R}^3}|k_w(v,v')| dv'\right\}ds'ds \nonumber\\
		&\le C_p e^{-\frac{\epsilon N^2}{8}} \left(\frac{1}{1+|v|}\right)^{2-\gamma} \sup_{0 \le s' \le t  }\left[e^{-\frac{\nu_0}{4}(t-s')}\|h(s')\|_{L^p_vL^\infty_x}\right],
	\end{align*}
	which implies the $L^p_vL^\infty_x$ norm of $R_2$ in this case is bounded by
	\begin{align} \label{class24}
		 C_p e^{-\frac{\epsilon N^2}{8}} \sup_{0 \le s' \le t  }\left[e^{-\frac{\nu_0}{4}(t-s')}\|h(s')\|_{L^p_vL^\infty_x}\right].
	\end{align}
	\newline
	$\mathbf{Case\ 3:}$ $t_{l'}'-s'\le \frac{1}{N}$.\\
	We easily obtain
	\begin{align}\label{class25}
		\|R_2^3\|_{L^p_vL^\infty_x} \le \frac{C_p}{N}\sup_{0 \le s' \le t  }\left[e^{-\frac{\nu_0}{4}(t-s')}\|h(s')\|_{L^p_vL^\infty_x}\right].
	\end{align}
	\newline
	$\mathbf{Case\ 4:}$ $t_{l'}'-s'\ge \frac{1}{N}$ and $|v|\le N$, $|v'| \le 2N$, $|v_{l'}'| \le N$, $|v''| \le 2N$.\\
	We split
\begin{align*}
	k_w(v,v')k_w(v_{l'}',v'') &= \left\{k_w(v,v')- k_N(v,v')\right\}k_w(v_{l'}',v'')\\
	&\quad +\left\{k_w(v_{l'}',v'')- k_N(v_{l'}',v'')\right\}k_N(v,v')\\
	& \quad +k_N(v,v')k_N(v_{l'}',v'').
\end{align*}
We use Lemma \ref{KESTIMATE} and \eqref{df13} to estimate
\begin{align*}
	&\int_{t_1}^{t-\epsilon}\int_{t_{l'+1}'}^{t_{l'}'-\frac{1}{N}} e^{-\frac{\nu_0}{2}(t-s')}\mathbf{1}_{\{|v|\le N\}}\int_{|v'|\le 2N} |k_w(v,v')-k_N(v,v')|\frac{1}{\tilde{w}(v')}  \int_{\prod_{j=1}^{k-1}\mathcal{V}_j', |v_{l'}'| \le N}\nonumber \\
	&\quad \times \Biggl\{\int_{|v''|\le2N}|k_w(v_{l'}',v'')| |h(s',x_{l'}'-v_{l'}'(t_{l'}'-s'),v'')|dv''\Biggr\}\left\{\prod_{j=l'+1}^{k-1}d\sigma_j'\right\}\left\{\tilde{w}(v'_{l'})d\sigma_{l'}'\right\}\\
	& \quad \times \left\{\prod_{j=1}^{l'-1}d\sigma_j'\right\}ds'dv'ds\\
	&\le C_p\int_{0}^{t}\int_{0}^{s} e^{-\frac{\nu_0}{2}(t-s')}\|h(s')\|_{L^p_vL^\infty_x}\mathbf{1}_{\{|v|\le N\}}\int_{|v'|\le 2N} |k_w(v,v')-k_N(v,v')|\frac{1}{\tilde{w}(v')}dv'ds'ds\\
	& \le \frac{C_p}{N^{2}} \mathbf{1}_{\{|v|\le N\}}\sup_{0 \le s' \le t  }\left[e^{-\frac{\nu_0}{4}(t-s')}\|h(s')\|_{L^p_vL^\infty_x}\right].
\end{align*}
Similarly, we derive
\begin{align*}
	&\int_{t_1}^{t-\epsilon}\int_{t_{l'+1}'}^{t_{l'}'-\frac{1}{N}} e^{-\frac{\nu_0}{2}(t-s')}\mathbf{1}_{\{|v|\le N\}}\int_{|v'|\le 2N} |k_N(v,v')|\frac{1}{\tilde{w}(v')}  \int_{\prod_{j=1}^{k-1}\mathcal{V}_j', |v_{l'}'| \le N}\nonumber \\
	&\quad \times  \Biggl\{\int_{|v''|\le2N}|k_w(v_{l'}',v'')-k_N(v_{l'}',v'')||h(s',x_{l'}'-v_{l'}'(t_{l'}'-s'),v'')|dv''\Biggr\}\left\{\prod_{j=l'+1}^{k-1}d\sigma_j'\right\}\\
	& \quad \times \left\{\tilde{w}(v'_{l'})d\sigma_{l'}'\right\}\left\{\prod_{j=1}^{l'-1}d\sigma_j'\right\}ds'dv'ds\\
	&\le C_p\int_0^{t}\int_0^{s} e^{-\frac{\nu_0}{2}(t-s')} \|h(s')\|_{L^p_v L^\infty_x}\mathbf{1}_{\{|v|\le N\}} \int_{|v'|\le 2N} |k_N(v,v')|\frac{1}{\tilde{w}(v')}\int_{\prod_{j=1}^{k-1}\mathcal{V}_j', |v_{l'}'| \le N}\\
	& \quad \times  \Biggl\{\sup_{|v_{l'}'|\le N}\int_{|v''|\le2N}|k_w(v_{l'}',v'')-k_N(v_{l'}',v'')|^{p'}dv''\Biggr\}^{1/p'}\left\{\prod_{j=l'+1}^{k-1}d\sigma_j'\right\}\left\{\tilde{w}(v'_{l'})d\sigma_{l'}'\right\}\\
	& \quad \times \left\{\prod_{j=1}^{l'-1}d\sigma_j'\right\}ds'dv'ds\\
	& \le \frac{C_p}{N^{\frac{p'+1}{p'}}}\int_0^{t}\int_0^{s} e^{-\frac{\nu_0}{2}(t-s')} \|h(s')\|_{L^p_v L^\infty_x}\mathbf{1}_{\{|v|\le N\}} \int_{|v'|\le 2N} |k_N(v,v')|\frac{1}{\tilde{w}(v')}dv'ds'ds\\
	& \le \frac{C_p}{N^{\frac{p'+1}{p'}}} \mathbf{1}_{\{|v|\le N\}}\sup_{0 \le s' \le t  }\left[e^{-\frac{\nu_0}{4}(t-s')}\|h(s')\|_{L^p_vL^\infty_x}\right].
\end{align*}
This implies
\begin{align*}
	R_2^4 &\le \frac{C_p}{N^{\frac{p'+1}{p'}}} \mathbf{1}_{\{|v|\le N\}}\sup_{0 \le s' \le t  }\left[e^{-\frac{\nu_0}{4}(t-s')}\|h(s')\|_{L^p_vL^\infty_x}\right]\\
	& \quad +C_{p,N} \mathbf{1}_{\{|v|\le N\}}\int_{t_1}^{t-\epsilon}\int_{t_{l'+1}'}^{t_{l'}'-\frac{1}{N}} e^{-\frac{\nu_0}{2}(t-s')}\int_{|v'|\le 2N}\frac{1}{\tilde{w}(v')}\int_{\prod_{j=1}^{k-1}\mathcal{V}_j', |v_{l'}'|\le N}\\
		&\quad \times  \Biggl\{\int_{|v''|\le 2N} |h(s',x_{l'}'-v_{l'}'(t_{l'}'-s'),v'')|dv''\Biggr\}\left\{\prod_{j=l'+1}^{k-1}d\sigma_j'\right\}\left\{\tilde{w}(v'_{l'})d\sigma_{l'}'\right\}\left\{\prod_{j=1}^{l'-1}d\sigma_j'\right\}ds'dv'ds\\
	& =:A_3+A_4.
\end{align*}
We can easily get
\begin{align} \label{class26}
	\|A_3\|_{L^p_vL^\infty_x} &\le \frac{C_p}{N^{\frac{p'+1}{p'}}}\sup_{0 \le s' \le t  }\left[e^{-\frac{\nu_0}{4}(t-s')}\|h(s')\|_{L^p_vL^\infty_x}\right] \left(\int_{|v|\le N} dv\right)^{1/p}\nonumber\\
	& \le \frac{C_p}{N}\sup_{0 \le s' \le t  }\left[e^{-\frac{\nu_0}{4}(t-s')}\|h(s')\|_{L^p_vL^\infty_x}\right].
\end{align}
since $p>4$.
We make a change of variables $v_{l'}' \mapsto y=x_{l'}'-v_{l'}'(t_{l'}'-s')$ with $\left|\frac{dy}{dv_{l'}'}\right| = (t_{l'}'-s')^3$:
\begin{align*}
	A_4 
	& \le C_{p,N} \mathbf{1}_{\{|v|\le N\}}\int_{t_1}^{t-\epsilon}\int_{t_{l'+1}'}^{t_{l'}'-\frac{1}{N}} e^{-\frac{\nu_0}{2}(t-s')}\int_{|v'|\le 2N}\frac{1}{\tilde{w}(v')}\int_{\prod_{j=1}^{l'-1}\mathcal{V}_j'}\Biggl\{\int_{\Omega}\int_{|v''|\le 2N} |h(s',y,v'')|dv''dy\Biggr\}\\
		&\quad \times  \left\{\prod_{j=1}^{l'-1}d\sigma_j'\right\}ds'dv'ds\\
	& \le C_{p,N,\Omega} \mathbf{1}_{\{|v|\le N\}}\int_0^t\int_0^{s} e^{-\frac{\nu_0}{2}(t-s')}\|f(s')\|_{L^2_{x,v}}ds'ds\\ 
	& \le C_{p,N,\Omega} \mathbf{1}_{\{|v|\le N\}} \int_0^{t} \|f(s')\|_{L^2_{x,v}}ds',     
\end{align*}
which implies that
\begin{align} \label{class27}
	\|A_4\|_{L^p_vL^\infty_x} \le C_{p,N,\Omega} \int_0^{t} \|f(s')\|_{L^2_{x,v}}ds'  
\end{align}
Combining \eqref{class21}, \eqref{class22}, \eqref{class23}, \eqref{class24}, \eqref{class25}, \eqref{class26}, \eqref{class27} and summing over $1\le l' \le k(\epsilon,\tilde{T}_0)-1$, we can derive
	\begin{equation} \label{df15}
	\begin{aligned}
		\|R_2\|_{L^p_vL^\infty_x} \le \frac{C_{p,\epsilon,\tilde{T}_0}}{N}\sup_{0 \le s' \le t  }\left[e^{-\frac{\nu_0}{4}(t-s')}\|h(s')\|_{L^p_vL^\infty_x}\right]+C_{p,N,\Omega,\tilde{T}_0}\int_0^{t}\|f(s')\|_{L^2_{x,v}}ds'.
	\end{aligned}
	\end{equation}
	\newline
	We denote $X_l'(s')=x_l-v_l(t_l-s)-v'(s-s')$.\\
	For $R_3$ in \eqref{RMC} : 
	\begin{align*} 
		& \sum_{l=1}^{k-1} \int_{t_{l+1}}^{t_l}\int_{t_1'}^{s-\epsilon} \frac{I^f(t,t_1)}{\tilde{w}(v)}\int_{\prod_{j=1}^{k-1}\mathcal{V}_j}  \int_{\mathbb{R}^3} k_w(v_l,v') I^f(s,s')\int_{\mathbb{R}^3}k_w(v',v'')h(s',X_l'(s'),v'') dv''dv'\nonumber\\
		& \quad \times d\Sigma_l^f(s)ds'ds \nonumber \\
		&\le \sum_{l=1}^{k-1} \int_{t_{l+1}}^{t_l}\int_{t_1'}^{s-\epsilon} \frac{e^{-\frac{\nu_0}{2}(t-s')}}{\tilde{w}(v)}\int_{\prod_{j=1}^{k-1}\mathcal{V}_j}  \int_{\mathbb{R}^3}|k_w(v_l,v')| \int_{\mathbb{R}^3}|k_w(v',v'')||h(s',X_l'(s'),v'')| dv''dv'\nonumber \\  
		&\quad \times\left\{\prod_{j=l+1}^{k-1}d\sigma_j\right\}\left\{\tilde{w}(v_{l})d\sigma_{l}\right\}\left\{\prod_{j=1}^{l-1}d\sigma_j\right\}ds'ds,
	\end{align*}
	where we have used the estimate \eqref{Ibound2}. Fix $l$. We divide this term into three cases.\\
	$\mathbf{Case\ 1:}$ $|v_l| \ge N$.\\
	We use Lemma \ref{KESTIMATE} to estimate
\begin{align*}
	|R_3^1| 
	&\le \int_{t_{l+1}}^{t_l}\int_{t_1'}^{s-\epsilon} \frac{e^{-\frac{\nu_0}{2}(t-s')}}{\tilde{w}(v)}\int_{\prod_{j=1}^{k-1}\mathcal{V}_j,|v_l|\ge N}  \int_{\mathbb{R}^3}|k_w(v_l,v')| \int_{\mathbb{R}^3}|k_w(v',v'')||h(s',X_l'(s'),v'')| dv''dv'\nonumber \\  
		&\quad \times\left\{\prod_{j=l+1}^{k-1}d\sigma_j\right\}\left\{\tilde{w}(v_{l})d\sigma_{l}\right\}\left\{\prod_{j=1}^{l-1}d\sigma_j\right\}ds'ds\\
	&\le \frac{C_p}{N}\int_0^t\int_0^s \frac{e^{-\frac{\nu_0}{2}(t-s')}}{\tilde{w}(v)} \|h(s')\|_{L^p_v L^\infty_x}ds'ds\\
	& \le \frac{C_p}{N}\frac{1}{\tilde{w}(v)} \sup_{0 \le s' \le t  }\left[e^{-\frac{\nu_0}{4}(t-s')}\|h(s')\|_{L^p_vL^\infty_x}\right].
\end{align*}
Thus we obtain
\begin{align} \label{class31}
	\|R_3^1\|_{L^p_v  L^\infty_x} \le \frac{C_p}{N} \sup_{0 \le s' \le t  }\left[e^{-\frac{\nu_0}{4}(t-s')}\|h(s')\|_{L^p_vL^\infty_x}\right].
\end{align}
	\newline
	$\mathbf{Case\ 2:}$ $|v_l| \le N$, $|v'| \ge 2N$ or $|v'| \le 2N$, $|v''| \ge 3N$.\\
	We apply \eqref{class0} and \eqref{class1} to deduce
\begin{align*}
	R_3^2 &\le \int_0^t\int_0^s \frac{e^{-\frac{\nu_0}{2}(t-s')}}{\tilde{w}(v)}\int_{\prod_{j=1}^{k-1}\mathcal{V}_j,|v_l|\le N}  \int_{|v'|\ge 2N}|k_w(v_l,v')| \int_{\mathbb{R}^3}|k_w(v',v'')||h(s',X_l'(s'),v'')| dv''dv'\nonumber \\  
		&\qquad \times\left\{\prod_{j=l+1}^{k-1}d\sigma_j\right\}\left\{\tilde{w}(v_{l})d\sigma_{l}\right\}\left\{\prod_{j=1}^{l-1}d\sigma_j\right\}ds'ds\\
	& \quad + \int_0^t\int_0^s \frac{e^{-\frac{\nu_0}{2}(t-s')}}{\tilde{w}(v)}\int_{\prod_{j=1}^{k-1}\mathcal{V}_j}  \int_{|v'|\le 2N}|k_w(v_l,v')| \int_{|v''|\ge 3N}|k_w(v',v'')||h(s',X_l'(s'),v'')| dv''dv'\nonumber \\  
	&\qquad \times\left\{\prod_{j=l+1}^{k-1}d\sigma_j\right\}\left\{\tilde{w}(v_{l})d\sigma_{l}\right\}\left\{\prod_{j=1}^{l-1}d\sigma_j\right\}ds'ds\\	& \le \frac{C_p}{\tilde{w}(v)} e^{-\frac{\epsilon N^2}{8}}\sup_{0 \le s' \le t  }\left[e^{-\frac{\nu_0}{4}(t-s')}\|h(s')\|_{L^p_vL^\infty_x}\right].
\end{align*}
Thus we obtain
\begin{align} \label{class32}
	\|R_3^2\|_{L^p_vL^\infty_x} \le C_pe^{-\frac{\epsilon N^2}{8}}\sup_{0 \le s' \le t  }\left[e^{-\frac{\nu_0}{4}(t-s')}\|h(s')\|_{L^p_vL^\infty_x}\right].
\end{align}
	\newline
	$\mathbf{Case\ 3:}$ $|v_l|\le N$, $|v'| \le 2N$, $|v''| \le 3N$.\\
	We split
\begin{align*}
	k_w(v_l,v')k_w(v',v'') &= \left\{k_w(v_l,v')- k_N(v_l,v')\right\}k_w(v',v'')\\
	&\quad +\left\{k_w(v',v'')- k_N(v',v'')\right\}k_N(v_l,v')\\
	& \quad +k_N(v_l,v')k_N(v',v'').
\end{align*}
We use Lemma \ref{KESTIMATE} and \eqref{df13} to estimate
\begin{align*}
	&\int_{t_{l+1}}^{t_l}\int_{t_1'}^{s-\epsilon} \frac{e^{-\frac{\nu_0}{2}(t-s')}}{\tilde{w}(v)}\int_{\prod_{j=1}^{k-1}\mathcal{V}_j,|v_l|\le N}  \int_{|v'|\le 2N}\int_{|v''|\le 3N}|k_w(v_l,v')-k_N(v_l,v')| |k_w(v',v'')| \nonumber \\  
		&\quad \times |h(s',X_l'(s'),v'')|dv''dv'\left\{\prod_{j=l+1}^{k-1}d\sigma_j\right\}\left\{\tilde{w}(v_{l})d\sigma_{l}\right\}\left\{\prod_{j=1}^{l-1}d\sigma_j\right\}ds'ds\\
	& \le \frac{C_p}{N^{2}}\frac{1}{\tilde{w}(v)} \sup_{0 \le s' \le t  }\left[e^{-\frac{\nu_0}{4}(t-s')}\|h(s')\|_{L^p_vL^\infty_x}\right].
\end{align*}
Similarly, we derive
\begin{align*}
	&\int_{t_{l+1}}^{t_l}\int_{t_1'}^{s-\epsilon} \frac{e^{-\frac{\nu_0}{2}(t-s')}}{\tilde{w}(v)}\int_{\prod_{j=1}^{k-1}\mathcal{V}_j,|v_l|\le N}  \int_{|v'|\le 2N}\int_{|v''|\le 3N}|k_N(v_l,v')| |k_w(v',v'')-k_N(v',v'')| \nonumber \\  
		&\quad \times |h(s',X_l'(s'),v'')|dv''dv'\left\{\prod_{j=l+1}^{k-1}d\sigma_j\right\}\left\{\tilde{w}(v_{l})d\sigma_{l}\right\}\left\{\prod_{j=1}^{l-1}d\sigma_j\right\}ds'ds\\
	& \le \frac{C_p}{N^{\frac{p'+1}{p'}}} \frac{1}{\tilde{w}(v)}\sup_{0 \le s' \le t  }\left[e^{-\frac{\nu_0}{4}(t-s')}\|h(s')\|_{L^p_vL^\infty_x}\right].
\end{align*}
This implies
\begin{align*}
	R_3^3 &\le \frac{C_p}{N} \frac{1}{\tilde{w}(v)}\sup_{0 \le s' \le t  }\left[e^{-\frac{\nu_0}{4}(t-s')}\|h(s')\|_{L^p_vL^\infty_x}\right]\\
	& \quad +C_{p,N} \int_{t_{l+1}}^{t_l}\int_{t_1'}^{s-\epsilon} \frac{e^{-\frac{\nu_0}{2}(t-s')}}{\tilde{w}(v)}\int_{\prod_{j=1}^{k-1}\mathcal{V}_j,|v_l|\le N}  \int_{|v'|\le 2N}\int_{|v''|\le 3N}|h(s',X_l'(s'),v'')| dv''dv'\nonumber \\  
		&\qquad \times \left\{\prod_{j=l+1}^{k-1}d\sigma_j\right\}\left\{\tilde{w}(v_{l})d\sigma_{l}\right\}\left\{\prod_{j=1}^{l-1}d\sigma_j\right\}ds'ds\\
	& =:A_5+A_6.
\end{align*}
We can easily get
\begin{align} \label{class33}
	\|A_5\|_{L^p_vL^\infty_x} &\le \frac{C_p}{N}\sup_{0 \le s' \le t  }\left[e^{-\frac{\nu_0}{4}(t-s')}\|h(s')\|_{L^p_vL^\infty_x}\right].
\end{align}
We make a change of variables $v' \mapsto y=X_{l'}'(s')=x_l-v_l(t_l-s)-v'(s-s')$ with $\left|\frac{dy}{dv'}\right| = (s-s')^3$:
\begin{align*}
	A_6 &\le C_{p,N,\epsilon} \int_{t_{l+1}}^{t_l}\int_{t_1'}^{s-\epsilon} \frac{e^{-\frac{\nu_0}{2}(t-s')}}{\tilde{w}(v)}\int_{\prod_{j=1}^{k-1}\mathcal{V}_j,|v_l|\le N}  \int_{\Omega}\int_{|v''|\le 3N}|h(s',y,v'')| dv''dy\nonumber \\  
		&\quad \times \left\{\prod_{j=l+1}^{k-1}d\sigma_j\right\}\left\{\tilde{w}(v_{l})d\sigma_{l}\right\}\left\{\prod_{j=1}^{l-1}d\sigma_j\right\}ds'ds \\
	& \le \frac{C_{p,N,\epsilon,\Omega}}{\tilde{w}(v)} \int_0^t\int_0^{s}e^{-\frac{\nu_0}{2}(t-s')}\|f(s')\|_{L^2_{x,v}}ds'ds\\ 
	& \le \frac{C_{p,N,\epsilon,\Omega}}{\tilde{w}(v)}  \int_0^{t} \|f(s')\|_{L^2_{x,v}}ds',    
\end{align*}
which implies that
\begin{align} \label{class34}
	\|A_6\|_{L^p_vL^\infty_x} \le C_{p,N,\epsilon,\Omega} \int_0^{t} \|f(s')\|_{L^2_{x,v}}ds'  
\end{align}
Combining the bounds \eqref{class31}, \eqref{class32}, \eqref{class33}, \eqref{class34} and summing over $1\le l' \le k(\epsilon,\tilde{T}_0)-1$, we derive
	\begin{equation} \label{df16}
	\begin{aligned}
		\|R_3\|_{L^p_vL^\infty_x} \le \frac{C_{p,\epsilon,\tilde{T}_0}}{N}\sup_{0 \le s' \le t  }\left[e^{-\frac{\nu_0}{4}(t-s')}\|h(s')\|_{L^p_vL^\infty_x}\right]+C_{p,N,\epsilon,\Omega,\tilde{T}_0}\int_0^{t}\|f(s')\|_{L^2_{x,v}}ds'.
	\end{aligned}
	\end{equation}
	\newline
	For $R_4$ in \eqref{RMC} :
	\begin{align} 
		&\sum_{l=1}^{k-1}\int_{t_{l+1}}^{t_l}\sum_{l'=1}^{k-1}\int_{t_{l'+1}'}^{t_{l'}'}\frac{I^f(t,t_1)}{\tilde{w}(v)}\int_{\prod_{j=1}^{k-1}\mathcal{V}_j} \int_{\mathbb{R}^3} k_w(v_l,v') \frac{I^f(s,t_1')}{\tilde{w}(v')}  \nonumber \\
		&\quad \times \int_{\prod_{j=1}^{k-1}\mathcal{V}_j'}\left\{\int_{\mathbb{R}^3}k_w(v_{l'}',v'')h(s',x_{l'}'-v_{l'}'(t_{l'}'-s'),v'')dv''\right\} d\Sigma_{l'}'^f(s')dv'd\Sigma_l^f(s)ds'ds \nonumber \\
		&\le \sum_{l=1}^{k-1}\int_{t_{l+1}}^{t_l}\sum_{l'=1}^{k-1}\int_{t_{l'+1}'}^{t_{l'}'}\frac{e^{-\frac{\nu_0}{2}(t-s')}}{\tilde{w}(v)}\int_{\prod_{j=1}^{k-1}\mathcal{V}_j} \int_{\mathbb{R}^3} |k_w(v_l,v')|\frac{1}{\tilde{w}(v')}  \nonumber \\
		& \quad \times  \int_{\prod_{j=1}^{k-1}\mathcal{V}_j'}\left\{\int_{\mathbb{R}^3}|k_w(v_{l'}',v'')||h(s',x_{l'}'-v_{l'}'(t_{l'}'-s'),v'')|dv''\right\} \left\{\prod_{j=l'+1}^{k-1}d\sigma_j'\right\}\left\{\tilde{w}(v_{l'}')d\sigma_{l'}'\right\} \nonumber\\ \label{FFTMC}
		& \quad \times \left\{\prod_{j=1}^{l'-1}d\sigma_j'\right\}dv'\left\{\prod_{j=l+1}^{k-1}d\sigma_j\right\}\left\{\tilde{w}(v_{l})d\sigma_{l}\right\}\left\{\prod_{j=1}^{l-1}d\sigma_j\right\}ds'ds,
	\end{align}
	where we have used the estimate \eqref{Ibound2}. Fix $l,l'$. We divide this term into four cases.\\
	\newline
	$\mathbf{Case\ 1:}$ $|v_l| \ge N$ or $|v_{l'}'| \ge N$.\\
	Using similar approaches as in \eqref{class21} and \eqref{class31}, we obtain
	\begin{align} \label{class42}
		\|R_4^1\|_{L^p_vL^\infty_x}\le \frac{C_p}{N} \sup_{0 \le s' \le t  }\left[e^{-\frac{\nu_0}{4}(t-s')}\|h(s')\|_{L^p_vL^\infty_x}\right].
	\end{align}
	\newline
	$\mathbf{Case\ 2:}$ $|v_l| \le N$, $|v'| \ge 2N$ or $|v_{l'}'| \le N$, $|v''| \ge 2N$.\\
	It follows from arguments similar to those used to derive \eqref{class24} and \eqref{class32} that
	\begin{align} \label{class44}
		\|R_4^2\|_{L^p_vL^\infty_x} \le  C_p e^{-\frac{\epsilon N^2}{8}} \sup_{0 \le s' \le t  }\left[e^{-\frac{\nu_0}{4}(t-s')}\|h(s')\|_{L^p_vL^\infty_x}\right].
	\end{align}
	\newline
	$\mathbf{Case\ 3:}$ $t_{l'}'-s' \le \frac{1}{N}$.\\
	We can easily compute
	\begin{align}\label{class45}
		\|R_4^3\|_{L^p_vL^\infty_x} \le \frac{C_p}{N}\sup_{0 \le s' \le t  }\left[e^{-\frac{\nu_0}{4}(t-s')}\|h(s')\|_{L^p_vL^\infty_x}\right].
	\end{align}
	\newline
	$\mathbf{Case\ 4:}$ $t_{l'}'-s' \ge \frac{1}{N}$ and $|v_l|\le N$, $|v'| \le 2N$, $|v_{l'}'| \le N$, $|v''| \le 2N$.\\
	From similar ways to obtain \eqref{class26}, \eqref{class27}, \eqref{class33}, \eqref{class34}, we have 
\begin{align} \label{class46}
	\|R_4^4\|_{L^p_vL^\infty_x} &\le \frac{C_p}{N}\sup_{0 \le s' \le t  }\left[e^{-\frac{\nu_0}{4}(t-s')}\|h(s')\|_{L^p_vL^\infty_x}\right]+C_{p,N,\Omega} \int_0^{t} \|f(s')\|_{L^2_{x,v}}ds'.
\end{align}
Combining \eqref{class42}, \eqref{class44}, \eqref{class45}, \eqref{class46} and summing over $1\le l \le k(\epsilon,\tilde{T}_0)-1$, $1\le l' \le k(\epsilon,\tilde{T}_0)-1$, it holds that
	\begin{equation} \label{df17}
	\begin{aligned}
		\|R_4\|_{L^p_vL^\infty_x} \le \frac{C_{p,\epsilon,\tilde{T}_0}}{N}\sup_{0 \le s' \le t  }\left[e^{-\frac{\nu_0}{4}(t-s')}\|h(s')\|_{L^p_vL^\infty_x}\right]+C_{p,N,\Omega,\epsilon,\tilde{T}_0}\int_0^{t}\|f(s')\|_{L^2_{x,v}}ds'.
	\end{aligned}
	\end{equation}
	\newline
	Gathering \eqref{df1}, \eqref{df2}, \eqref{df4}, \eqref{df5}, \eqref{df8}, \eqref{df9}, \eqref{df10}, \eqref{df11}, \eqref{df12}, \eqref{df14}, \eqref{df15}, \eqref{df16}, and \eqref{df17}, we deduce for $0\le t \le \tilde{T}_0$,
	\begin{align} \label{df18}
		\|h(t)\|_{L^p_vL^\infty_x} 
		& \le  C_{p,\gamma,\mathcal{T}}(1+t) e^{-\frac{\nu_0}{4}t}\|h_0\|_{L^p_vL^\infty_x}+C_{p,\gamma,\mathcal{T}}\sup_{0 \le s' \le t  }\left[e^{-\frac{\nu_0}{8}(t-s')}\|h(s')\|_{L^p_vL^\infty_x}^2\right]\nonumber\\
		& \quad + \left( \epsilon C_{p,\gamma,\mathcal{T}} +\frac{C_{p,\epsilon,\tilde{T}_0}}{N}\right) \sup_{0 \le s' \le t  }\left[e^{-\frac{\nu_0}{8}(t-s')}\|h(s')\|_{L^p_vL^\infty_x}\right]\nonumber\\
		& \quad +C_{p,N,\epsilon,\Omega,\tilde{T}_0}\int_0^{t}\|f(s')\|_{L^2_{x,v}}ds'.
	\end{align}
	By the a priori assumption \eqref{Aprioriassump1}, the estimate \eqref{df18} becomes
	\begin{align*}
		\|h(t)\|_{L^p_vL^\infty_x} 
		& \le  C_{p,\gamma,\mathcal{T}}e^{-\frac{\nu_0}{8}t}\|h_0\|_{L^p_vL^\infty_x}\\
		& \quad + \left( \epsilon C_{p,\gamma,\mathcal{T}} +\frac{C_{p,\epsilon,\tilde{T}_0}}{N}+C_5 \bar{\eta} \right) \sup_{0 \le s' \le t  }\left[e^{-\frac{\nu_0}{8}(t-s')}\|h(s')\|_{L^p_vL^\infty_x}\right]\\
		& \quad +C_{p,N,\epsilon,\Omega,\tilde{T}_0}\int_0^{t}\|f(s')\|_{L^2_{x,v}}ds',
	\end{align*}
	where $C_5$ is a positive constant depending only on $p$, $\gamma$, and $\mathcal{T}$. Note that $C_5$ is independent of $\tilde{T}_0$. This implies that
	\begin{align*}
		\sup_{0\le t \le \tilde{T}_0}\left[e^{\frac{\nu_0}{8}t}\|h(t)\|_{L^p_v L^\infty_x}\right]
		& \le  C_{p,\gamma,\mathcal{T}}\|h_0\|_{L^p_vL^\infty_x}\\
		& \quad + \left( \epsilon C_{p,\gamma,\mathcal{T}}+\frac{C_{p,\epsilon,\tilde{T}_0}}{N}+C_5 \bar{\eta}\right) \sup_{0\le s'\le \tilde{T}_0}\left[e^{\frac{\nu_0}{8}s'}\|h(s')\|_{L^p_v L^\infty_x}\right]\\
		& \quad +C_{p,N,\epsilon,\Omega,\tilde{T}_0}\int_0^{\tilde{T}_0}\|f(s')\|_{L^2_{x,v}}ds'.
	\end{align*}
	Choosing first $\epsilon>0$ small enough and then choosing $N$ sufficiently large such that $\epsilon C_{p,\gamma,\mathcal{T}}+\frac{C_{p,\epsilon,\tilde{T}_0}}{N}< \frac{1}{4}$, and since 
	\begin{align} \label{etacond2}
		C_5 \bar{\eta} = 4C_5C_4^2 \eta_0 <1/4,
	\end{align}
	 we obtain
	\begin{align*}
		\sup_{0\le t \le \tilde{T}_0}\left[e^{\frac{\nu_0}{8}t}\|h(t)\|_{L^p_v L^\infty_x}\right]
		& \le  C_{p,\gamma,\mathcal{T}}\|h_0\|_{L^p_vL^\infty_x} +C_{p,\Omega,\tilde{T}_0}\int_0^{\tilde{T}_0}\|f(s')\|_{L^2_{x,v}}ds'.
	\end{align*}
	This yields that
	\begin{align*}
		\|h(\tilde{T}_0)\|_{L^p_vL^\infty_x} \le C_{p,\gamma,\mathcal{T}}e^{-\frac{\nu_0}{8}\tilde{T}_0}\|h_0\|_{L^p_vL^\infty_x} +C_{p,\Omega,\tilde{T}_0} \int_0^{\tilde{T}_0}\|f(s')\|_{L^2_{x,v}}ds'.
	\end{align*}
	Choosing large $\tilde{T}_0 >0$ so that
	\begin{align*} 
		C_{p,\gamma,\mathcal{T}}e^{-\frac{\nu_0}{8}\tilde{T}_0} \le e^{-\frac{\nu_0}{16} \tilde{T}_0},
	\end{align*}
	we conclude the estimate \eqref{finitetimeest}. Note that $\tilde{T}_0$ depends only $p$, $\gamma$, $\beta$, and $\mathcal{T}$, where $\mathcal{T}$ is a fixed constant satisfying \eqref{assump1}. \\
	Using Lemma \ref{nonlinear L^2 decay}, thanks to the condition 
	\begin{align} \label{etacond3}
		C_2 \bar{\eta}^2= 16C_2C_4^4 \eta_0^2 < 1/2,
	\end{align}
	where $C_2$ is defined in Lemma \ref{nonlinear L^2 decay}, we can derive the $L^2$ decay as follows:
	\begin{align} \label{L2decayorigin}
		\|f(t)\|_{L^2_{x,v}} \lesssim e^{-\lambda_2 t}\|f_0\|_{L^2_{x,v}}
	\end{align}
	for all $t\ge 0$. Note that $\lambda_2 \le \nu_0/16$. For any $n\ge 1$, we apply the estimate \eqref{finitetimeest} repeatedly to $h(l\tilde{T}_0)$ for $l = n-1,n-2, \cdots,0$:
	\begin{align*}
		\|h(n\tilde{T}_0)\|_{L^p_vL^\infty_x} &\le e^{-\frac{\nu_0}{16} \tilde{T}_0}\|h((n-1)\tilde{T}_0)\|_{L^p_vL^\infty_x}+C_{p,\Omega,\tilde{T}_0} \int_0^{\tilde{T}_0} \|f((n-1)\tilde{T}_0+s)\|_{L^2_{x,v}}ds\\
		& = e^{-\frac{\nu_0}{16} \tilde{T}_0}\|h((n-1)\tilde{T}_0)\|_{L^p_vL^\infty_x}+C_{p,\Omega,\tilde{T}_0} \int_{(n-1)\tilde{T}_0}^{n\tilde{T}_0} \|f(s)\|_{L^2_{x,v}}ds\\
		& \le e^{-2\frac{\nu_0}{16} \tilde{T}_0}\|h((n-2)\tilde{T}_0)\|_{L^p_vL^\infty_x}+e^{-\frac{\nu_0}{16} \tilde{T}_0}C_{p,\Omega,\tilde{T}_0}\int_{(n-2)\tilde{T}_0}^{(n-1)\tilde{T}_0}\|f(s)\|_{L^2_{x,v}}ds\\
		& \quad +C_{p,\Omega,\tilde{T}_0} \int_{(n-1)\tilde{T}_0}^{n\tilde{T}_0} \|f(s)\|_{L^2_{x,v}}ds\\
		& \le e^{-n\frac{\nu_0}{16}\tilde{T}_0}\|h_0\|_{L^p_vL^\infty_x} + C_{p,\Omega,\tilde{T}_0}\sum_{l=0}^{n-1}e^{-l\frac{\nu_0}{16} \tilde{T}_0} \int_{(n-l-1)\tilde{T}_0}^{(n-l)\tilde{T}_0} \|f(s)\|_{L^2_{x,v}}ds.
	\end{align*}
	From the assumption \eqref{L2decayorigin}, in the interval $(n-l-1)\tilde{T}_0 \le s \le (n-l)\tilde{T}_0$, we have $\|f(s)\|_{L^2_{x,v}}\lesssim e^{-\lambda_2s}\|f_0\|_{L^2_{x,v}}\lesssim e^{-\lambda_2(n-l-1)\tilde{T}_0}\|f_0\|_{L^2_{x,v}}$. Then we obtain
	\begin{align*}
		\|h(n\tilde{T}_0)\|_{L^p_vL^\infty_x} &\le e^{-n\lambda_2\tilde{T}_0}\|h_0\|_{L^p_vL^\infty_x}+ C_{p,\Omega,\tilde{T}_0}\sum_{l=0}^{n-1}e^{-l\lambda_2 \tilde{T}_0} \int_{(n-l-1)\tilde{T}_0}^{(n-l)\tilde{T}_0} e^{-\lambda_2(n-l-1)\tilde{T}_0}\|f_0\|_{L^2_{x,v}}ds\\
		& \le e^{-n\lambda_2\tilde{T}_0}\|h_0\|_{L^p_vL^\infty_x}+C_{p,\Omega,\tilde{T}_0}e^{\lambda_2 \tilde{T}_0} n\tilde{T}_0e^{-n\lambda_2 \tilde{T}_0}\|f_0\|_{L^2_{x,v}},
	\end{align*}
	where we have used the fact $\lambda_2 \le \nu_0/16$.
	Since $\beta > \frac{3(p-2)}{2p}$, we use the H\"older inequality to derive 
	\begin{align*}
		\|f_0\|_{L^2_{x,v}} \le C_{\Omega} \|h_0\|_{L^p_vL^\infty_x} \|w^{-2}\|_{L^{\frac{p}{p-2}}_v} \le C_{\Omega,p,\beta} \|h_0\|_{L^p_vL^\infty_x},
	\end{align*}
	and $\|h(n\tilde{T}_0)\|_{L^p_vL^\infty_x}$ is bounded by
	\begin{align*}
		&e^{-n\lambda_2\tilde{T}_0}\|h_0\|_{L^p_vL^\infty_x}+C_{p,\Omega,\tilde{T}_0}e^{\lambda_2 \tilde{T}_0} n\tilde{T}_0e^{-n\lambda_2 \tilde{T}_0}\|h_0\|_{L^p_vL^\infty_x}\\
		& \le e^{-n\lambda_2\tilde{T}_0}\|h_0\|_{L^p_vL^\infty_x}+C_{p,\Omega,\tilde{T}_0,\lambda_2}e^{-\frac{n\lambda_2 \tilde{T}_0}{2}}\|h_0\|_{L^p_vL^\infty_x}\\
		& \le C_{p,\Omega,\tilde{T}_0,\lambda_2} e^{-\frac{n\lambda_2 \tilde{T}_0}{2}}\|h_0\|_{L^p_vL^\infty_x}      ,
	\end{align*}
	where we have used $n\tilde{T}_0 e^{-n\lambda_2\tilde{T}_0} \le e^{-\frac{n\lambda_2\tilde{T}_0}{2}}$. For $t \ge 0$, there exists $n \ge 1$ such that $n\tilde{T}_0 \le t \le (n+1)\tilde{T}_0$, and it follows that
	\begin{align*}
		\|h(t)\|_{L^p_vL^\infty_x} &\le C\|h(n\tilde{T}_0)\|_{L^p_vL^\infty_x} \le C_{p,\Omega,\tilde{T}_0,\lambda_2} e^{-\frac{n\lambda_2 \tilde{T}_0}{2}}\|h_0\|_{L^p_vL^\infty_x}\\ 
		& \le  C_{p,\Omega,\tilde{T}_0,\lambda_2} e^{-\frac{\lambda_2 }{2}t}\|h_0\|_{L^p_vL^\infty_x},
	\end{align*}
	where $e^{-\frac{n\lambda_2 \tilde{T}_0}{2}} \le e^{\frac{\lambda_2\tilde{T}_0}{2}}e^{-\frac{\lambda_2 }{2}t}$. Therefore, we conclude the exponential time decay.
\end{proof}

\bigskip

\section{Large amplitude problem} \label{Largeamplitudeproblem}
Our main goal in this section is to establish the global-in-time existence to the Boltzmann equation \eqref{FPBER} with initial data of large amplitude in $L^p_vL^\infty_x$. 
\subsection{The a priori assumption in the large amplitude problem}
Recall that we fix $p$, $\beta$, $\mathcal{T}$ in subsection \ref{apriorisamll}. 
Let $f(t,x,v)$ with 
\begin{align*} 
	F(t,x,v) = \mu(v)+ \mu^{1/2}(v)f(t,x,v) \ge 0
\end{align*}
be the solution to the equation \eqref{FPBE} with initial data $f_0(x,v)$ over the time interval $[0,T)$ for $0<T\le \infty$. We set $h(t,x,v)=w(v)f(t,x,v)$. Throughout this section, we make the a priori assumption : 
\begin{align} \label{Aprioriassump2}
	\sup_{0\le t < T } \|h(t)\|_{L^p_vL^\infty_x} \le \bar{M},
\end{align}
where $\bar{M}\ge 1$ is a large constant depending only on the initial large amplitude $M_0$, not on the solution $h$ and $T$. $\bar{M}$ and $T$ will be determined later. See \eqref{LA1} and \eqref{LA2}.\\
\indent As in the small data problem, we need to estimate a lower bound for the term $R(f)$. Unlike Lemma \ref{Rfest1}, it is slightly different because initial data has a large amplitude. The difference will be controlled by the exponential time decay for the term associated with the initial data.\begin{Lem} \label{Rfest2}
	Under the a priori assumption \eqref{Aprioriassump2}, there exists a generic constant $C_8 \ge 1$ such that for given $T_0 > \tilde{t}$ with 
	\begin{align} \label{smalltimeslot}
		\tilde{t}:= \frac{2}{\nu_0} \log (C_8 M_0)>0,
	\end{align}
   there is a generally small positive constant $\epsilon_2 = \epsilon_2(\bar{M},T_0)>0$, depending only on $\bar{M}$ and $T_0$, such that if $\mathcal{E}(F_0) \le \epsilon_2$, then we have
	\begin{align*}
		R(f)(t,x,v) \ge \frac{1}{2}\nu(v)
	\end{align*}
	for all $(t,x,v) \in [\tilde{t},T_0) \times \Omega \times \R^3$.
\end{Lem}
\begin{proof}
Let $0\le t \le T_0$. We set $h(t,x,v) =w(v)f(t,x,v)$. We recall
	\begin{align*}
		R(f)(t,x,v) &= \int_{\mathbb{R}^3 \times \mathbb{S}^2} B(v-u,\omega) [\mu(u) + \mu^{\frac{1}{2}}(u)f(t,x,u)]d\omega du\\
		& = \nu(v)+\int_{\mathbb{R}^3 \times \mathbb{S}^2} B(v-u,\omega)\mu^{\frac{1}{2}}(u)f(t,x,u) d\omega du\\
		& \ge \nu(v) -\tilde{C}_3 \nu(v) \int_{\mathbb{R}^3} e^{-\frac{|u|^2}{8}}|f(t,x,u)|du,
	\end{align*}
	where $\tilde{C}_3$ is a constant depending only on the collision kernel $b$. \\
	If it holds that
	\begin{align} \label{RRf1}
		\int_{\mathbb{R}^3} e^{-\frac{|u|^2}{8}}|h(t,x,u)|du\le  \frac{1}{2\tilde{C}_3} \quad \text{for all } \tilde{t} \le t \le T_0,\ x\in \Omega,
	\end{align}
	then we can complete the proof of this Lemma. Thus it suffices to show \eqref{aa61}.
	Then by Duhamel principle,
	\begin{align*}
		h(t,x,v) = S_{G_\nu}(t,0)h_0 + \int_0^t S_{G_\nu}(t,s) \left[K_wh(s)+w\Gamma(f,f)(s)\right]ds.
	\end{align*}
	This yields that
	\begin{equation} \label{RRf2}
	\begin{aligned}
		&\int_{\mathbb{R}^3}e^{-\frac{|v|^2}{8}}|h(t,x,v)|dv\\
		& \le \int_{\mathbb{R}^3} e^{-\frac{|v|^2}{8}} |(S_{G_\nu}(t,0)h_0)(t,x,v)|dv\\
		& \quad + \int_{t-\epsilon}^{t}  \int_{\mathbb{R}^3}e^{-\frac{|v|^2}{8}}\left[\left| \left(S_{G_\nu}(t,s)K_wh(s)\right)(t,x,v)\right|+|(S_{G_\nu}(t,s)w\Gamma(f,f)(s))(t,x,v)|\right]dvds\\
		& \quad + \int_0^{t-\epsilon}  \int_{\mathbb{R}^3}e^{-\frac{|v|^2}{8}}\left| \left(S_{G_\nu}(t,s)K_wh(s)\right)(t,x,v)\right|dvds\\ 
		& \quad + \int_0^{t-\epsilon} \int_{\mathbb{R}^3}e^{-\frac{|v|^2}{8}}|(S_{G_\nu}(t,s)w\Gamma(f,f)(s))(t,x,v)|dvds.
	\end{aligned}
	\end{equation}
	By the same argument in the proof of Lemma \ref{Rfest1} and the a priori assumption \eqref{Aprioriassump2}, we can estimate \eqref{RRf2}:
\begin{align*}
	\int_{\mathbb{R}^3}e^{-\frac{|v|^2}{8}}|h(t,x,v)|dv &\le e^{-\frac{\nu_0}{2} t} \|h_0\|_{L^p_vL^\infty_x}\\
	& \quad  + C_p \left(\epsilon+\frac{C_{p,\gamma,\epsilon,T_0}}{N^{1/2}} \right)  \sup_{0\le s \le t}\Bigl[\|h(s)\|_{L^p_vL^\infty_x}+\|h(s)\|_{L^p_vL^\infty_x}^2+\|h(s)\|_{L^p_vL^\infty_x}^{1-\kappa}\\
	& \qquad +\|h(s)\|_{L^p_vL^\infty_x}^{3-\kappa} \Bigr]\\
	& \quad +C_{p,\gamma,\epsilon,N,T_0} \left[\mathcal{E}(F_0)^{1/2}+\mathcal{E}(F_0)+ \mathcal{E}(F_0)^\kappa \right]\\
	& \le  e^{-\frac{\nu_0}{2} t}\|h_0\|_{L^p_vL^\infty_x}+C_p \left(\epsilon+\frac{C_{p,\gamma,\epsilon,T_0}}{N^{1/2}} \right) \left[\bar{M} + \bar{M}^2 + \bar{M}^{1-\kappa}+ \bar{M}^{3-\kappa} \right]\\
	& \quad +C_{p,\gamma,\epsilon,N,T_0} \left[\mathcal{E}(F_0)^{1/2}+\mathcal{E}(F_0)+ \mathcal{E}(F_0)^\kappa \right],
\end{align*}
where $\kappa$ is defined in \eqref{kappaconst}. Taking
\begin{align*}
	\tilde{t}:= \frac{2}{\nu_0} \log(4\tilde{C}_3M_0)=\frac{2}{\nu_0} \log(C_8M_0),
\end{align*}
it holds that for $t \ge\tilde{t}$
\begin{align*}
	 e^{-\frac{\nu_0}{2} t}\|h_0\|_{L^p_vL^\infty_x} \le  e^{-\frac{\nu_0}{2} t}M_0 \le \frac{1}{4\tilde{C}_3}.
\end{align*}
We first choose suitably small $\epsilon>0$, then take $N>0$ large enough, and last choose sufficiently small $\epsilon_1$ with $\mathcal{E}(F_0) \le \epsilon_2$ so that
\begin{align*}
	C_p \left(\epsilon+\frac{C_{p,\gamma,\epsilon,T_0}}{N^{1/2}} \right) \left[\bar{M} + \bar{M}^2 + \bar{M}^{1-\kappa}+ \bar{M}^{3-\kappa} \right]
	 +C_{p,\gamma,\epsilon,N,T_0} \left[\mathcal{E}(F_0)^{1/2}+\mathcal{E}(F_0)+ \mathcal{E}(F_0)^\kappa \right] \le \frac{1}{4\tilde{C}_3}.
\end{align*}
Therefore, we conclude \eqref{RRf1}, and we completes the proof of this lemma.
\end{proof}

\bigskip
Recall that $S_{G_f}(t,s)$ is the solution operator to the equation \eqref{WFPBE2} with the diffuse reflection boundary condition \eqref{WPDRBC}. To deal with the nonlinear $L^p_vL^\infty_x$ estimate in the large amplitude problem, we need to derive the exponential time decay for the operator $S_{G_f}(t,s)$. By using the estimate \eqref{assump1}, we can show the below lemma.
\bigskip
\begin{Lem} \label{L61}
	Assume the a priori assumption \eqref{Aprioriassump2}. Then there exists a constant $C_{p,\gamma,\mathcal{T}}>0$ so that if $\mathcal{E}(F_0) \le \epsilon_2$, where $\epsilon_2 = \epsilon_2(\bar{M},T_0)$ is determined in Lemma \ref{Rfest2},  
	\begin{align*}
		\|S_{G_f}(t,s)h(s)\|_{L^p_vL^\infty_x} \le C_{p,\gamma,\mathcal{T}} e^{\frac{3}{4}\nu_0\tilde{t}}e^{-\frac{1}{4}\nu_0(t-s)} \|h(s)\|_{L^p_vL^\infty_x}
	\end{align*}
	for all $0\le s \le t \le T_0$.
\end{Lem}
\begin{proof}
	Suppose that $\mathcal{E}(F_0) \le \epsilon_2(\bar{M},T_0)$. By Lemma \ref{Rfest2}, we have
	\begin{equation} \label{LL1}
	\begin{aligned}
		R(f)(t,x,v) \ge 
		\begin{cases}
			0 & \text{if }t\in [0,\tilde{t}], \\
			\frac{1}{2}\nu(v) \quad & \text{if }t\in (\tilde{t},T_0],
		\end{cases}
	\end{aligned}
	\end{equation}
	for all $(x,v)\in \Omega\times \mathbb{R}^3$.\\
	\newline
	$\mathbf{Case\ of \ t \in [0,\tilde{t}]:}$ We know that $t/ \mathcal{T}\in [m,m+1)$ for some $m \in \{0,1,\cdots\ , [\tilde{t}/\mathcal{T}]\}$, where $\tilde{t}$ is determined in Lemma \ref{Rfest2}. Thus it follows from Lemma \ref{Dampexistence} and \eqref{assump1} that
	\begin{align*}
		\|S_{G_f}(t,0)h_0\|_{L^p_vL^\infty_x} &\le C_1\mathcal{T}^{\frac{5}{4}}\|S_{G_f}(m\mathcal{T},0)h_0\|_{L^p_vL^\infty_x} \le C_1 \mathcal{T}^{\frac{5}{4}}\left(C_1\mathcal{T}^{\frac{5}{4}}\right)^m \|h_0\|_{L^p_vL^\infty_x} \\
		& \le C_1 \mathcal{T}^{\frac{5}{4}}\left(C_1\mathcal{T}^{\frac{5}{4}}\right)^{\tilde{t}/\mathcal{T}}\|h_0\|_{L^p_vL^\infty_x}\le C_{p,\gamma,\mathcal{T}} e^{\frac{\nu_0}{2}\tilde{t}}\|h_0\|_{L^p_vL^\infty_x},
	\end{align*}
	and we derive
	\begin{align*}
		\|S_{G_f}(t,0)h_0\|_{L^p_vL^\infty_x} \le C_{p,\gamma,\mathcal{T}} e^{\frac{3}{4}\nu_0\tilde{t}}e^{-\frac{1}{4}\nu_0t}\|h_0\|_{L^p_vL^\infty_x}.
	\end{align*}
	\newline
	$\mathbf{Case\ of \ t \in (\tilde{t},T_0]:}$ We note that $S_{G_f}(t,0)h_0 = S_{G_f}(t,\tilde{t})S_{G_f}(\tilde{t},0)h_0$. From \eqref{LL1} and Lemma \ref{DampLpdecay}, we obtain
	\begin{align*}
		\|S_{G_f}(t,0)h_0\|_{L^p_vL^\infty_x} &\le C_{p,\gamma,\mathcal{T}} e^{-\frac{\nu_0}{4}(t-\tilde{t})}\|S_{G_f}(\tilde{t},0)h_0\|_{L^p_vL^\infty_x}.
	\end{align*}
	We apply the previous case to get
	\begin{align*}
		\|S_{G_f}(t,0)h_0\|_{L^p_vL^\infty_x} \le C_{p,\gamma,\mathcal{T}} e^{-\frac{\nu_0}{4}(t-\tilde{t})} e^{\frac{1}{2}\nu_0\tilde{t}}\|h_0\|_{L^p_vL^\infty_x}\le C_{p,\gamma,\mathcal{T}} e^{\frac{3}{4}\nu_0\tilde{t}}e^{-\frac{1}{4}\nu_0t}\|h_0\|_{L^p_vL^\infty_x}.
	\end{align*}
	\newline
	Gathering two cases, we conclude that
	\begin{align*}
		\|S_{G_f}(t,0)h_0\|_{L^p_vL^\infty_x} \le C_{p,\gamma,\mathcal{T}} e^{\frac{3}{4}\nu_0\tilde{t}}e^{-\frac{1}{4}\nu_0t}\|h_0\|_{L^p_vL^\infty_x}.
	\end{align*}
	for all $0 \le t \le T_0$.
\end{proof}
\bigskip
We now demonstrate the nonlinear $L^p_vL^\infty_x$ estimate in the large amplitude problem using the above lemma.

\bigskip

\begin{Lem} \label{LALinftyestimate3} 
	Assume that $\mathcal{E}(F_0)\le \epsilon_2=\epsilon_2(\bar{M},T_0)$. Under the a priori assumption \eqref{Aprioriassump2}, there exists a generic constant $C_7 \ge 1$, depending on $p$, $\gamma$, $\beta$, such that
\begin{align*}
	\|h(t)\|_{L^p_vL^\infty_x} &\le C_7e^{2\nu_0 \tilde{t}}\|h_0\|_{L^p_vL^\infty_x}\left(1+\int_0^t \|h(s)\|_{L^p_vL^\infty_x}ds \right)e^{-\frac{\nu_0}{8}t}\\
	& \quad + C_7e^{2\nu_0 \tilde{t}}\bigg\{\left(\epsilon + \frac{C_{p,,\gamma,\epsilon,T_0}}{N^{1/2}} \right) \sup_{0\le s\le t}\biggl[\|h(s)\|_{L^p_vL^\infty_x}+\|h(s)\|_{L^p_vL^\infty_x}^2+\|h(s)\|_{L^p_vL^\infty_x}^3\nonumber\\
	& \qquad+\|h(s)\|_{L^p_vL^\infty_x}^4+\|h(s)\|_{L^p_vL^\infty_x}^{\frac{p}{2p-2}}+\|h(s)\|_{L^p_vL^\infty_x}^{\frac{3p-2}{2p-2}}+\|h(s)\|_{L^p_vL^\infty_x}^{\frac{5p-4}{2p-2}}+\|h(s)\|_{L^p_vL^\infty_x}^{\frac{7p-6}{2p-2}}\\
	& \qquad+\|h(s)\|_{L^p_vL^\infty_x}^{\frac{9p-8}{2p-2}}\biggr] +C_{p,\gamma,N,\epsilon,\Omega,T_0}\left[\mathcal{E}(F_0)+\mathcal{E}(F_0)^2+\mathcal{E}(F_0)^{\frac{p-2}{2p-2}}+\mathcal{E}(F_0)^{\frac{p-2}{p-1}}\right]\bigg\}
\end{align*}
for all $0 \le t \le T_0$, where $\epsilon>0$ can be arbitrary small and $N>0$ can be arbitrary large.
\end{Lem}
\begin{proof}
	Let $(t,x,v) \in (0,T_0]\times \Omega \times \R^3$. We apply the Duhamel principle to the weighted Boltzmann equation \eqref{WFPBE}:
	\begin{align*}
			h(t,x,v) &= S_{G_f}(t,0) h_0+ \int_0^t \left[S_{G_f}(t,s)K_wh(s)+S_{G_f}(t,s)w\Gamma^+(f,f)(s)\right] ds\\
		& = S_{G_f}(t,0) h_0+ \int_{t-\epsilon}^t \left[S_{G_f}(t,s)K_wh(s)+S_{G_f}(t,s)w\Gamma^+(f,f)(s)\right] ds\\
		& \quad +  \int_0^{t-\epsilon} S_{G_f}(t,s)K_wh(s)ds+\int_0^{t-\epsilon}S_{G_f}(t,s)w\Gamma^+(f,f)(s) ds,
	\end{align*}
	where $\epsilon>0$ is to be fixed sufficiently small later. Note that from the a priori assumption \eqref{Aprioriassump2} and Lemma \ref{Rfest2} , we have
\begin{align*}
		R(f)(t,x,v) \ge \frac{1}{2}\nu(v)
	\end{align*}
	for all $(t,x,v) \in [\tilde{t},T_0] \times \Omega \times \R^3$. Thus it holds that
	\begin{align} \label{IboundLA}
		I^f(t,s) \le e^{\frac{\nu_0}{4}\tilde{t}}e^{-\frac{\nu_0}{4}(t-s)}
	\end{align}
	for all $0\le s \le t \le T_0$.\\
By similar approaches in Lemma \ref{LALinftyestimate1} and using Lemma \ref{L61} instead of Lemma \ref{DampLpdecay}, we deduce
	\begin{align} \label{llp1}
	|h(t,x,v)| &\le \left|(S_{G_f}(t,0)h_0)(t,x,v) \right|+\int_{t-\epsilon}^{t}\left|(S_{G_f}(t,s)K_wh(s))(t,x,v) \right|ds\nonumber\\
	& \quad +\int_{t-\epsilon}^t\left|(S_{G_f}(t,s)w\Gamma^+(f,f)(s))(t,x,v) \right|ds\nonumber\\
	& \quad + e^{\frac{\nu_0}{4}\tilde{t}}\int_0^{t-\epsilon} \mathbf{1}_{\{t_1\le s\}}e^{-\frac{\nu_0}{4}(t-s)} \left|(K_wh)(s,x-v(t-s),v) \right|ds\nonumber\\
	& \quad + e^{\frac{\nu_0}{4}\tilde{t}}\int_0^{t-\epsilon} \mathbf{1}_{\{t_1\le s\}}e^{-\frac{\nu_0}{2}(t-s)} \left|w\Gamma^+(f,f)(s,x-v(t-s),v) \right|ds\nonumber\\
	& \quad +\frac{e^{\nu_0\tilde{t}}}{\tilde{w}(v)}\left(C_p\epsilon + \frac{C_{p,\epsilon,T_0}}{N^{1/2}} \right)\sup_{0\le s\le t}\left[\|h(s)\|_{L^p_vL^\infty_x}+\|h(s)\|_{L^p_vL^\infty_x}^2+\|h(s)\|_{L^p_vL^\infty_x}^{\frac{p}{2p-2}}+\|h(s)\|_{L^p_vL^\infty_x}^{\frac{5p-4}{2p-2}}\right]\nonumber\\
	& \quad +e^{\nu_0\tilde{t}}\frac{C_{p,N,\Omega,T_0}}{\tilde{w}(v)}\left[\mathcal{E}(F_0)^{1/2}+\mathcal{E}(F_0)+\mathcal{E}(F_0)^{\frac{p-2}{2p-2}}\right],
\end{align}	
where $\epsilon>0$ can be arbitrary small and $N>0$ can be arbitrary large. Using similar ways to derive the estimate in Lemma \ref{LALinftyestimate2} and applying \eqref{IboundLA}, \eqref{llp1}, and Lemma \ref{L61}, we can conclude the result for this Lemma. 
\end{proof}

\bigskip
\subsection{Proof of the main theorem}
\begin{proof}[\textbf{Proof of Theorem \ref{Mainresult2}}]
	Take $C_6 := \max\left\{C_1 \mathcal{T}^{\frac{5}{4}},C_0,C_7\right\} >1$, and let
		\begin{align} \label{LA1}
			\bar{M} := 4 C_6^2M_0 \exp\left\{2\nu_0 \tilde{t}+\frac{8}{\nu_0}C_6M_0e^{2\nu_0 \tilde{t}} \right\} = 4C_6^2C_8^4M_0^5 \exp\left\{\frac{8}{\nu_0}C_6C_8^4M_0^5 \right\},
		\end{align}
		and
		\begin{align} \label{LA2}
			T_0 := \frac{16}{\nu_0}\left(\log\bar{M}+|\log\eta_0|\right).
		\end{align}
		By the a priori assumption \eqref{Aprioriassump2} and Lemma \ref{LALinftyestimate2}, we get
		\begin{align} \label{LA3}
			\|h(t)\|_{_{L^p_vL^\infty_x}} \le C_6C_8^4 M_0^5 \left(1+\int_0^t \|h(s)\|_{L^p_vL^\infty_x}ds\right)e^{-\frac{\nu_0}{8}t}+\tilde{E} \quad \text{for all } 0\le t \le T_0,
		\end{align}
		where
		\begin{align*}
			\tilde{E}:&= C_6C_8^4 M_0^4 \bigg\{\left(\epsilon + \frac{C_{p,\gamma,\epsilon,T_0}}{N^{1/2}} \right) \biggl[\bar{M}+\bar{M}^2+\bar{M}^3+\bar{M}^4+\bar{M}^{\frac{p}{2p-2}}+\bar{M}^{\frac{3p-2}{2p-2}}+\bar{M}^{\frac{5p-4}{2p-2}}\\
			& \qquad +\bar{M}^{\frac{7p-6}{2p-2}}+\bar{M}^{\frac{9p-8}{2p-2}}\biggr]+C_{p,\gamma,\epsilon,N,T_0,\Omega}\left[\mathcal{E}(F_0)+\mathcal{E}(F_0)^2+\mathcal{E}(F_0)^{\frac{p-2}{2p-2}}+\mathcal{E}(F_0)^{\frac{p-2}{p-1}}\right]\bigg\}.
		\end{align*}
		We define
		\begin{align*}
			G(t) := 1+\int_0^t \|h(s)\|_{L^p_vL^\infty_x}ds.
		\end{align*}
		Then the inequality \eqref{LA3} becomes
		\begin{align} \label{LA4}
			G'(t) \le C_6C_8^4 M_0^5  e^{-\frac{\nu_0}{8}t} G(t) + \tilde{E}.
		\end{align}
		By Gr\"{o}nwall's inequality, \eqref{LA4} implies that
		\begin{equation} \label{LA5}
		\begin{aligned}
			G(t) &\le (1+\tilde{E}t) \exp\left\{\frac{8}{\nu_0}C_6C_8^4 M_0^5  \left(1-e^{-\frac{\nu_0}{8}t}\right)\right\} \le (1+\tilde{E}t) \exp\left\{\frac{8}{\nu_0}C_6C_8^4 M_0^5 \right\}
		\end{aligned}
		\end{equation}
		for all $0\le t \le T_0$.\\
		Inserting \eqref{LA5} into \eqref{LA3}, we deduce for $0\le t \le T_0$
		\begin{align*}
			\|h(t)\|_{L^p_vL^\infty_x} &\le C_6C_8^4 M_0^5  \exp\left\{\frac{8}{\nu_0}C_6C_8^4 M_0^5 \right\}(1+\tilde{E}t)e^{-\frac{\nu_0}{8}t}+\tilde{E}\\
			& \le \frac{1}{4C_6}\bar{M}(1+\tilde{E}t)e^{-\frac{\nu_0}{8}t} +\tilde{E}\\
			& \le \frac{1}{4C_6}\bar{M}\left(1+\frac{16}{\nu_0} \tilde{E}\right)e^{-\frac{\nu_0}{16}t} +\tilde{E}.
		\end{align*}
		We first choose $\epsilon >0$ small enough, then choose $N>0$ sufficiently large, and assume $\mathcal{E}(F_0) \le \epsilon_3(\eta_0, M_0)$ with small enough $\epsilon_3(\eta_0, M_0)$ so that
		\begin{align*}
			\tilde{E} \le \min\left\{\frac{\nu_0}{64} , \frac{\eta_0}{8}\right\},
		\end{align*}
		and it follows that
		\begin{align} \label{LA6}
			\|h(t)\|_{L^p_vL^\infty_x} \le \frac{5}{16C_6}\bar{M}e^{-\frac{\nu_0}{16}t} + \frac{\eta_0}{8} \le \frac{1}{2C_6} \bar{M}
		\end{align}
		for all $0\le t \le T_0$. Hence we have closed the a priori assumption \eqref{Aprioriassump2} over $t \in [0,T_0]$ if $\mathcal{E}(F_0) \le \min\{\epsilon_2,\epsilon_3\}$.\\
		We claim that a solution to the Boltzmann equation \eqref{Boltzmanneq} extends into time interval $[0,T_0]$. From Lemma \ref{Localexistence}, there exists the Boltzmann solution $F(t) \ge 0$ to \eqref{Boltzmanneq} on $[0,\hat{t}_0]$ such that
		\begin{align} \label{LA7}
			\sup_{0\le t \le  \hat{t}_0} \|h(t)\|_{L^p_vL^\infty_x} \le 2 \|h_0\|_{L^p_vL^\infty_x} \le \frac{1}{2C_6} \bar{M}.
		\end{align}
		We define $t_*:=\left(\hat{C}_\mathcal{T}\left[1+(2C_6)^{-1}\bar{M}\right]\right)^{-1}>0$, where $\hat{C}_\mathcal{T}$ is a constant in Theorem \ref{Localexistence}. Taking $t=\hat{t}_0$ as the initial time, it follows from \eqref{LA7} and Theorem \ref{Localexistence} that we can extend the Boltzmann equation solution $F(t) \ge 0$ into time interval $[0, \hat{t}_0+t_*]$ satisfying
		\begin{align*}
			\sup_{\hat{t}_0 \le t \le \hat{t}_0+t_*} \|h(t)\|_{L^p_vL^\infty_x} \le 2 \|h(\hat{t}_0)\|_{L^p_vL^\infty_x} \le \bar{M}.
		\end{align*}
		Thus we have
		\begin{align} \label{LA8}
			\sup_{0 \le t \le \hat{t}_0+t_*} \|h(t)\|_{L^p_vL^\infty_x} \le \bar{M}.
		\end{align}
		Note that \eqref{LA8} means $h(t)$ satisfies the a priori assumption \eqref{Aprioriassump2} over $[0,\hat{t}_0+t_*]$.
		From \eqref{LA6}, we can obtain
		\begin{align*}
			\sup_{0 \le t \le \hat{t}_0+t_*} \|h(t)\|_{L^p_vL^\infty_x} \le\frac{1}{2C_6}\bar{M}.
		\end{align*}
		Repeating the same process for finite times, we can derive that there exists the Boltzmann equation solution $F(t)\ge 0$ on the time interval $[0,T_0]$ such that
		\begin{align*}
			\sup_{0\le t \le T_0} \|h(t)\|_{L^p_vL^\infty_x} \le \frac{1}{2C_6} \bar{M}.
		\end{align*}
		Let us consider the case $[T_0,\infty)$. Assume $\mathcal{E}(F_0) \le \tilde{\epsilon}_0:=\min\{\epsilon_0,\epsilon_2,\epsilon_3\}$. From \eqref{LA6}, we get
		\begin{align*}
			\|h(T_0)\|_{L^p_vL^\infty_x} \le \frac{5}{16C_6}\bar{M}e^{-\frac{\nu_0}{16}T_0} + \frac{\eta_0}{8} \le \frac{5 \eta_0}{16C_6}+\frac{\eta_0}{8} < \frac{\eta_0}{2},
		\end{align*}
		and from Lemma \ref{relativedecrease},
		\begin{align*}
			\mathcal{E}(F(T_0)) \le \mathcal{E}(F_0)<\tilde{\epsilon}_0<\epsilon_0.
		\end{align*}
		Taking $t=T_0$ as the initial time and using Theorem \ref{Mainresult1}, we conclude that there exists the Boltzmann equation solution $F(t)\ge 0$ on $[0,\infty)$. Therefore, we have proven the global existence and uniqueness of the Boltzmann equation \eqref{Boltzmanneq}.\\
		It remains to show the exponential decay of the Boltzmann solution $f(t)$ in the weighted $L^p_vL^\infty_x$ space. Assume $\mathcal{E}(F_0) \le \tilde{\epsilon}_0$. By Theorem \ref{Mainresult1}, for all $t \ge T_0$,
		\begin{align} \label{LA9}
			\|h(t)\|_{L^p_vL^\infty_x} \le C_0 \|h(T_0)\|_{L^p_vL^\infty_x} e^{-\lambda_0(t-T_0)} \le C_0 \eta_0 e^{-\lambda_0(t-T_0)}.
		\end{align}
		 Taking $\tilde{C}_L := 4C_6^3C_8^4$ and $\lambda_L:= \min \left\{\lambda_0,\frac{\nu_0}{16}\right\}$, it follows from \eqref{LA6} and \eqref{LA9} that 
		 \begin{align*}
		 	\|h(t)\|_{L^p_vL^\infty_x} &\le \max\left\{\frac{1}{2},C_0\right\}\bar{M}e^{-\lambda_L t}\\
		 	& \le C_6 \bar{M} e^{-\lambda_L t}\\
		 	& \le 4C_6^3 C_8^4 M_05\exp\left\{ \frac{8}{\nu_0}C_6C_8^4M_0^5 \right\} e^{-\lambda_L t}\\
		 	& \le \tilde{C}_L M_0^5\exp\left\{ \frac{2\tilde{C}_L}{\nu_0}M_0^5 \right\}e^{-\lambda_L t}
		 \end{align*}
		 for all $t \ge 0$.
		 Therefore, we completes the proof of this theorem.
\end{proof}

\bigskip

\section{Appendix} \label{Appendix}
In this section, we present the proof of Lemma \ref{coer}. By choosing suitable test functions in a weak formulation and using the elliptic estimate, we will demonstrate this lemma.

\begin{proof}[\textbf{Proof of Lemma \ref{coer}}]

For simplicity, we set $s=0$. We will choose suitable test functions $\psi=\psi(t,x,v) \in H_{x,v}^1.$\\
	We deduce the weak formulation for the equation \eqref{FPBER}: 
	\begin{align} \label{weakformulation}
		&\int_{\Omega \times \mathbb{R}^3} \psi(t)f(t)dxdv - \int_{\Omega \times \mathbb{R}^3} \psi(0)f(0)dxdv\nonumber\\
		& = \int_0^t \int_{\Omega \times \mathbb{R}^3} f \left(v\cdot \nabla_x \psi \right)dxdvds-\int_0^t\int_{\gamma} \psi f \{n(x) \cdot v\}dS(x)dvds -\int_0^t \int_{\Omega \times \mathbb{R}^3}Lf \psi dxdvds \nonumber\\
		& \quad+\int_0^t \int_{\Omega \times \mathbb{R}^3}\Gamma(f,f)\psi dxdvds+\int_0^t \int_{\Omega \times \mathbb{R}^3} f (\partial_t \psi)dxdvds.
	\end{align}
	Decomposing $f=Pf+\left(I-P\right)f$ and using the fact $L\left(Pf\right) = 0$, the weak form \eqref{weakformulation} becomes
	\begin{align}\label{wwf}
		&\underbrace{-\int_0^t \int_{\Omega \times \mathbb{R}^3}\left(v\cdot \nabla_x \psi \right) Pf dxdvds}_{=:\mathcal{W}_1 }\nonumber\\
		& =\underbrace{\int_{\Omega \times \mathbb{R}^3} \psi(0) f(0) dxdv}_{=:\mathcal{W}_2 } - \underbrace{\int_{\Omega \times \mathbb{R}^3} \psi(t)f(t) dxdv}_{=:\mathcal{W}_3}+\underbrace{\int_0^t \int_{\Omega \times \mathbb{R}^3} \left(v\cdot \nabla_x \psi \right) (I-P)f dxdvds}_{=:\mathcal{W}_4}\nonumber\\
		& \quad -\underbrace{\int_0^t \int_{\Omega \times \mathbb{R}^3} L\left(\left(I-P\right)f\right) \psi dxdvds}_{=:\mathcal{W}_5 }-\underbrace{\int_0^t\int_{\gamma} \psi f \{n(x) \cdot v\}dS(x)dvds}_{=:\mathcal{W}_6}\nonumber\\
		& \quad +\underbrace{\int_0^t \int_{\Omega \times \mathbb{R}^3} f\left(\partial_t\psi\right)dxdvds}_{=:\mathcal{W}_7} +\underbrace{\int_0^t \int_{\Omega \times \mathbb{R}^3} \Gamma(f,f)\psi dxdvds}_{=:\mathcal{W}_8 }.
	\end{align}
	We consider $Pf = a(t,x)\mu^{1/2}(v)+b(t,x)\cdot v\mu^{1/2}(v) + c(t,x)\frac{|v|^2-3}{\sqrt{6}}\mu^{1/2}(v)$, and will derive the estimates for $a(t,x)$, $b(t,x)$, $c(t,x)$ from the weak formulation \eqref{wwf}.\\
	\newline
	$\mathbf{(Estimate\ for \ c(t,x))}$\\
	We can choose the test function
	\begin{align*}
		\psi = \psi_c(t,x,v) = \left(|v|^2-\beta_c\right)\mu^{1/2}(v)v\cdot \nabla_x \phi_c(t,x),
	\end{align*}
	where the function $\phi_c$ satisfies 
	\begin{align*}
		\begin{cases}
		-\Delta_x \phi_c(t,x) = c(t,x)\\
		\phi_c |_{\partial \Omega}=0,
	\end{cases}
	\end{align*}
	and $\beta_c>0$ is chosen such that
	\begin{align} \label{constbetac}
		\int_{\mathbb{R}^3} \left(|v|^2-\beta_c\right)\mu(v)v_i^2dv=0 \quad \text{for} \quad i=1,2,3.
	\end{align}
	We can get $\beta_c = 5$ from the simple computation.\\
	From the standard elliptic estimate, we have
	\begin{align} \label{ellipticest1}
		\|\phi_c(t)\|_{H_x^2} \lesssim \|c(t)\|_{L_x^2} 
	\end{align}
	for all $t\ge 0$.\\
	For $\mathcal{W}_2$ and $\mathcal{W}_3$, from the elliptic estimate \eqref{ellipticest1}, we obtain 
	\begin{align} \label{L2C1}
		\left|\mathcal{W}_3\right| &\le \|f(t)\|_{L^2_{x,v}}\left(\int_{\Omega \times \mathbb{R}^3}\left|\left(|v|^2-\beta_c\right)\mu^{1/2}(v)v\cdot \nabla_x\phi_c(t,x)\right|^2dxdv\right)^{\frac{1}{2}}\nonumber\\
		&\lesssim  \|f(t)\|_{L^2_{x,v}}\|\phi_c(t)\|_{H_x^2}\nonumber \\
		&\lesssim \|f(t)\|_{L^2_{x,v}}\|c(t)\|_{L_x^2}\nonumber\\
		&\lesssim \|f(t)\|_{L^2_{x,v}}^2,
	\end{align}
	where $\int_{\mathbb{R}^3}\left|\left(|v|^2-\beta_c\right)\mu^{\frac{1}{2}}(v)v\right|^2dv$ is finite. Similarly, we can derive
	\begin{align} \label{L2C2}
		|\mathcal{W}_2| \lesssim \|f(0)\|_{L^2_{x,v}}^2
	\end{align}  
	For $\mathcal{W}_4$, from the elliptic estimate \eqref{ellipticest1}, we obtain
	\begin{align} \label{L2C3}
		&\left|\int_{\Omega\times \mathbb{R}^3} \left(v\cdot \nabla_x \psi \right) (I-P)f(t) dxdv\right|\nonumber\\
		& \le \sum_{i,j=1}^3\left|\int_{\Omega\times \mathbb{R}^3}\left(|v|^2-\beta_c\right)\mu^{1/2}(v)v_iv_j\partial_{ij}\phi_c(t,x)(I-P)f(t)dxdv\right|\nonumber\\
		& \le \sum_{i,j=1}^3\left(\int_{\Omega \times \mathbb{R}^3}\left|\left(|v|^2-\beta_c\right)\mu^{1/2}(v)v_iv_j\partial_{ij}\phi_c(t,x)\right|^2dxdv\right)^{1/2}\left\|(I-P)f(t)\right\|_{L^2_{x,v}}\nonumber\\
		& \lesssim \|\phi_c(t)\|_{H_x^2}\left\|(I-P)f(t)\right\|_{L^2_{x,v}}\nonumber\\
		& \lesssim \|c(t)\|_{L^2_x}\left\|(I-P)f(t)\right\|_{L^2_{x,v}},
	\end{align}
	where $\int_{\mathbb{R}^3}\left|\left(|v|^2-\beta_c\right)\mu^{\frac{1}{2}}(v)v_iv_j\right|^2dv$ is finite.\\
	Thus we deduce that
	\begin{align} \label{L2C4}
		|\mathcal{W}_4| \lesssim \int_0^t \|c(s)\|_{L^2_x}\left\|(I-P)f(s)\right\|_{L^2_{x,v}}ds.
	\end{align}
	For $\mathcal{W}_5$, because $L =\nu-K$, it is bounded by
	\begin{align*}
		&\left|\int_{\Omega \times \mathbb{R}^3}\psi\nu(v)(I-P)f(t)dxdv\right| + \left|\int_{\Omega \times \mathbb{R}^3}\psi K\left((I-P)f(t)\right)dxdv\right|\\
		&=: I_1+I_2.
	\end{align*}
	We can easily get
	\begin{align*}
		I_1 \lesssim \|c(t)\|_{L_x^2}\left\|(I-P)f(t)\right\|_{L^2_{x,v}}.
	\end{align*}
	Since $K$ is bounded in $L^2_{x,v}$, we also get
	\begin{align*}
		I_2 &\le \left(\int_{\Omega \times \mathbb{R}^3} \left|\left(|v|^2-\beta_c\right)\mu^{1/2}(v)v\cdot \nabla_x \phi_c(t,x)\right|^2dxdv\right)^{1/2}\left\|K\left((I-P)f(t)\right)\right\|_{L^2_{x,v}}\\
		& \lesssim\|\phi_c(t)\|_{H_x^2}\left\|(I-P)f(t)\right\|_{L^2_{x,v}}\\
		& \lesssim\|c(t)\|_{L^2_x}\left\|(I-P)f(t)\right\|_{L^2_{x,v}},
	\end{align*}
	where we have used the elliptic estimate \eqref{ellipticest1}.
	Gathering $I_1$, $I_2$ and integrating from $0$ to $t$, we obtain
	\begin{align} \label{L2C5}
	|\mathcal{W}_5| \lesssim \int_0^t \|c(s)\|_{L^2_x}\left\|(I-P)f(s)\right\|_{L^2_{x,v}}ds.
	\end{align}
	For $\mathcal{W}_6$, we can decompose it into two terms:
	\begin{align*}
		\mathcal{W}_6 = \int_{\gamma_+}\psi f\{n(x) \cdot v\}dS(x)dv +\int_{\gamma_-}\psi f\{n(x) \cdot v\}dS(x)dv.
	\end{align*}
	Decomposing $f = P_\gamma f + (I-P_\gamma)f$, we get
	\begin{align*}
		\mathcal{W}_6 =\int_{\gamma_+} \psi \left[(I-P_\gamma)f\right]\{n(x) \cdot v\}dS(x)dv + \int_{\gamma} \psi \left(P_\gamma f\right)\{n(x) \cdot v\}dS(x)dv.
	\end{align*}
	Setting $z(t,x) =c_\mu \int_{n(x) \cdot v' >0}f(x,v')\mu^{1/2}(v')\{n(x)\cdot v'\}dv'$, we obtain 
	\begin{align*}
		\int_\gamma \psi (P_\gamma f)\{n(x) \cdot v\}dS(x)dv
		&=\sum_{i=1}^3 \left(\int_{\mathbb{R}^3} \left(|v|^2-\beta_c\right)\mu(v)v_i^2 dv\right)\left(\int_{\partial \Omega}\partial_{x_i} \phi_c(t,x)z(t,x)n_i(x)dS(x)\right)\\
		& =0,
	\end{align*}
	where we have used the oddness for integral and the condition \eqref{constbetac} for $\beta_c$.\\
	Then we can simplify
	\begin{align*}
		|\mathcal{W}_6| 
		&\le \sum_{i,k=1}^3\left(\int_{\gamma_+} \left|\left(|v|^2-\beta_c\right)\mu^{1/2}(v)v_iv_k n_k(x)\partial_{x_i} \phi_c(t,x)\right|^2 dS(x)dv\right)^{1/2}\|(I-P_\gamma)f(t)\|_{L^2_{\gamma_+}}\\
		&\lesssim \left\|\partial_{x_i} \phi_c(t)\right\|_{L^2(\partial\Omega)}\|(I-P_\gamma)f(t)\|_{L^2_{\gamma_+}},
	\end{align*}
	where $\int_{\mathbb{R}^3}\left|\left(|v|^2-\beta_c\right)\mu^{1/2}(v)v_iv_k\right|^2dv$ is finite.\\
	By the trace theorem, the above is bounded by
	\begin{align*}
		\left\|\phi_c(t)\right\|_{H^2_x}\|(I-P_\gamma)f(t)\|_{L^2_{\gamma_+}} \lesssim \|c(t)\|_{L^2_x}\|(I-P_\gamma)f(t)\|_{L^2_{\gamma_+}}.
	\end{align*}
	Thus, we deduce
	\begin{align} \label{L2C16}
		|\mathcal{W}_6|\lesssim \int_0^t \|c(s)\|_{L^2_x}\|(I-P_\gamma)f(s)\|_{L^2_{\gamma_+}}ds.
	\end{align}

	For $\mathcal{W}_8$, we apply the Cauchy-Schwarz inequality to obtain
	\begin{align*}
		\left|\int_{\Omega \times \mathbb{R}^3}\Gamma(f,f)(t) \psi  dxdv\right| &\lesssim \left(\int_{\Omega \times \mathbb{R}^3} \left|\left(|v|^2-\beta_c\right)\mu^{1/2}(v)v\cdot \nabla_x \phi_c(t,x)\right|^2dxdv\right)^{1/2} \| \Gamma(f,f)(t)\|_{L^2_{x,v}} \\
		& \lesssim\|\phi_c(t)\|_{H_x^2}\| \Gamma(f,f)(t)\|_{L^2_{x,v}} \\ 
		& \lesssim \|c(t)\|_{L^2_x}\| \Gamma(f,f)(t)\|_{L^2_{x,v}},
	\end{align*}
	which implies that
	\begin{align} \label{L2C6}
		|\mathcal{W}_8| \lesssim \int_0^t \|c(s)\|_{L^2_x}\| \Gamma(f,f)(s)\|_{L^2_{x,v}}ds.
	\end{align}
	For $\mathcal{W}_1$,	 from the construction \eqref{constbetac} for $\beta_c$ and the oddness of integration in $v$, we deduce
	\begin{align*}
		\mathcal{W}_1
		& = \int_0^t \int_{\Omega \times \mathbb{R}^3}\left(|v|^2-\beta_c\right)\mu(v) \left\{\sum_{i,j=1}^3v_iv_j\partial_{ij}\phi_c(s,x)\right\}\left[a(s,x)+b(s,x)\cdot v+c(s,x)\frac{|v|^2-3}{\sqrt{6}}\right]dxdvds\\
		& = \sum_{i=1}^3 \int_0^t \int_{\Omega \times \mathbb{R}^3} \left(|v|^2-\beta_c\right)\mu(v)v_i^2\partial_{ii}\phi_c(s,x) \left[c(s,x)\frac{|v|^2-3}{\sqrt{6}}\right]dxdvds\\
		& = \sum_{i=1}^3 \int_0^t \left(\int_{\mathbb{R}^3} \left(|v|^2-\beta_c\right)\mu(v)v_i^2\frac{|v|^2-3}{\sqrt{6}}dv\right)\left(\int_{\Omega}\left[\partial_{ii}\phi_c(s,x)\right] c(s,x)dx\right)ds.
	\end{align*}
	Here, for $i=1,2,3$,
	\begin{align*}
		\int_{\mathbb{R}^3} \left(|v|^2-\beta_c\right)\mu(v)v_i^2\frac{|v|^2-3}{\sqrt{6}}dv &= \frac{1}{3}\int_{\mathbb{R}^3} \left(|v|^2-\beta_c\right)|v|^2\frac{|v|^2-3}{\sqrt{6}}\mu(v)dv\\
		&=\frac{1}{3\sqrt{6}}(7\cdot5\cdot3-8\cdot5\cdot3+15\cdot3)=:A>0.
	\end{align*}
	This yields
	\begin{align} \label{L2C7}
		\mathcal{W}_1 &=-A\sum_{i=1}^3 \int_0^t \left(\int_{\Omega}\left[\partial_{ii}\phi_c(s,x)\right] c(s,x)dx\right)ds \nonumber\\
		& = A\int_0^t\|c(s)\|_{L_x^2}^2ds. 
	\end{align}
	For $\mathcal{W}_7$, we decompose $f=Pf+\left(I-P\right)f$ to get
	\begin{align*}
		\mathcal{W}_7
		&=\sum_{i=1}^3\int_0^t\int_{\Omega \times \mathbb{R}^3}\left(|v|^2-\beta_c\right)\mu(v)v_i\left[\partial_t\partial_{x_i}\phi_c\right]\left[a(s,x)+b(s,x)\cdot v+c(s,x)\frac{|v|^2-3}{\sqrt{6}}\right]dxdvds\\
		&\quad +\sum_{i=1}^3\int_0^t\int_{\Omega \times \mathbb{R}^3} \left(|v|^2-\beta_c\right)\mu^{1/2}(v)v_i\left[\partial_t\partial_{x_i}\phi_c\right](I-P)f(s)dxdvds.
	\end{align*}
	From the construction for $\beta_c$ and the oddness of integration in $v$, the above expression becomes
	\begin{align*}
		&\sum_{i=1}^3\int_0^t\int_{\Omega \times \mathbb{R}^3}\left(|v|^2-\beta_c\right)\mu(v)v_i^2\left[\partial_t\partial_{x_i}\phi_c\right]b_i(s,x)dxdvds\\
		&\quad +\sum_{i=1}^3\int_0^t\int_{\Omega \times \mathbb{R}^3} \left(|v|^2-\beta_c\right)\mu^{1/2}(v)v_i\left[\partial_t\partial_{x_i}\phi_c\right](I-P)f(s)dxdvds\\
		& = \sum_{i=1}^3\int_0^t\int_{\Omega \times \mathbb{R}^3} \left(|v|^2-\beta_c\right)\mu^{1/2}(v)v_i\left[\partial_t\partial_{x_i}\phi_c\right](I-P)f(s)dxdvds.
	\end{align*}
	Lately, we will demonstrate the estimate \eqref{L2C15} for $\nabla_x \partial_t \phi_c$ as following:
	\begin{align*}
		\left\|\nabla_x \partial_t \phi_c(t)\right\|_{L^2_x} \lesssim \|b(t)\|_{L_x^2}+\|(I-P)f(t)\|_{L_{x,v}^2}.
	\end{align*}
	By the estimate \eqref{L2C15} for $\nabla_x \partial_t \phi_c$, we get
	\begin{align} \label{L2C8}
		\mathcal{W}_7 = &\int_0^t\int_{\Omega \times \mathbb{R}^3} \left(|v|^2-\beta_c\right)\mu^{1/2}(v)v_i\left[\partial_t\partial_{x_i}\phi_c\right](I-P)f(s)dxdvds \nonumber\\
		& \lesssim\int_0^t \left\|\nabla_x \partial_t \phi_c(s)\right\|_{L^2_x}\|(I-P)f(s)\|_{L_{x,v}^2}ds  \nonumber\\
		& \lesssim\int_0^t \left(\|b(s)\|_{L_x^2}+\|(I-P)f(s)\|_{L_{x,v}^2}\right)\|(I-P)f(s)\|_{L_{x,v}^2}ds  \nonumber\\
		& \le \epsilon \int_0^t\|b(s)\|_{L_x^2}^2ds+C(\epsilon)\int_0^t\|(I-P)f(s)\|_{L_{x,v}^2}^2ds, 
	\end{align}
	where we have used the Young's inequality.\\
	Gathering \eqref{L2C1}, \eqref{L2C2}, \eqref{L2C3}, \eqref{L2C4}, \eqref{L2C5}, \eqref{L2C16}, \eqref{L2C6}, \eqref{L2C7}, \eqref{L2C8}, we obtain
	\begin{align*}
		A\int_0^t\|c(s)\|_{L_x^2}^2ds &\le G_f^c(t)-G_f^c(0) +C\int_0^t \|c(s)\|_{L^2_x}\left\|(I-P)f(s)\right\|_{L^2_{x,v}}ds+\\
		& \quad +\int_0^t \|c(s)\|_{L^2_x}\|(I-P_\gamma)f(s)\|_{L^2_{\gamma_+}}ds+C\int_0^t\|c(s)\|_{L^2_x}\|\Gamma(f,f)(s)\|_{L^2_{x,v}}ds\\
		& \quad + \epsilon \int_0^t\|b(s)\|_{L_x^2}^2ds+C(\epsilon)\int_0^t\|(I-P)f(s)\|_{L_{x,v}^2}^2ds\\
		&\le  G_f^c(t)-G_f^c(0) +\epsilon\int_0^t \|c(s)\|_{L^2_x}^2ds + \epsilon \int_0^t\|b(s)\|_{L_x^2}^2ds\\
		&\quad +C(\epsilon)\int_0^t\left[\|(I-P)f(s)\|_{L_{x,v}^2}^2+\|(I-P_\gamma)f(s)\|_{L^2_{\gamma_+}}^2 \right]ds\\
		& \quad +C(\epsilon)\int_0^t \|\Gamma(f,f)(s)\|_{L^2_{x,v}}^2ds ,
	\end{align*}
	where we have used the Young's inequality and $G_f^c(s) \lesssim \|f(s)\|_{L^2_{x,v}}^2$.\\
	Thus, choosing sufficiently small $\epsilon >0$, we conclude that
	 
	\begin{align} \label{co1}
		\int_0^t\|c(s)\|_{L_x^2}^2ds &\le G_f^c(t)-G_f^c(0) +\epsilon \int_0^t\|b(s)\|_{L_x^2}^2ds \nonumber\\
		&\quad +C(\epsilon)\int_0^t\left[\|(I-P)f(s)\|_{L_{x,v}^2}^2+\|(I-P_\gamma)f(s)\|_{L^2_{\gamma_+}}^2 \right]ds\nonumber\\
		& \quad +C(\epsilon) \int_0^t\|\Gamma(f,f)(s)\|_{L^2_{x,v}}^2ds. 
	\end{align}
	\newline
	$\mathbf{(Estimate\ for \ \nabla_x\partial_t\phi_c)}$\\
	We consider the weak formulation for the equation \eqref{FPBER} over $[t,t+\epsilon]$:
	\begin{align} \label{wf}
		&\underbrace{\int_{\Omega \times \mathbb{R}^3} \psi(x,v)f(t+\epsilon)dxdv}_{=:\bar{\mathcal{W}}_1} - \underbrace{\int_{\Omega \times \mathbb{R}^3} \psi(x,v)f(t)dxdv}_{=:\bar{\mathcal{W}}_2}\nonumber\\
		& = \underbrace{\int_t^{t+\epsilon} \int_{\Omega \times \mathbb{R}^3} f\left(v\cdot \nabla_x \psi \right)dxdvds}_{=:\bar{\mathcal{W}}_3} -\underbrace{\int_t^{t+\epsilon}\int_{\gamma} \psi f \{n(x) \cdot v\}dS(x)dvds}_{=:\bar{\mathcal{W}}_4}- \underbrace{\int_t^{t+\epsilon} \int_{\Omega \times \mathbb{R}^3}Lf \psi dxdvds}_{=:\bar{\mathcal{W}}_5}\nonumber\\
		&\quad +\underbrace{\int_t^{t+\epsilon} \int_{\Omega \times \mathbb{R}^3}\Gamma(f,f)\psi dxdvds}_{=:\bar{\mathcal{W}}_6} +\underbrace{\int_t^{t+\epsilon} \int_{\Omega \times \mathbb{R}^3} f (\partial_t \psi)dxdvds}_{=:\bar{\mathcal{W}}_7}. 
	\end{align}
	We choose the test function
	\begin{align*}
		\psi = \psi(x,v) = \phi(x) \frac{|v|^2-3}{\sqrt{6}}\mu^{\frac{1}{2}}(v),
	\end{align*}
	where $\phi(x)$ depends only on $x$, choosing later.\\
	For $\bar{\mathcal{W}}_1$ and $\bar{\mathcal{W}}_2$, we can easily get
	\begin{align} \label{L2C9}
		\bar{\mathcal{W}}_1 = \int_{\Omega \times \mathbb{R}^3} \phi(x)  \frac{|v|^2-3}{\sqrt{6}}\mu^{1/2}(v)f(t+\epsilon)dxdv = \int_{\Omega}\phi(x)c(t+\epsilon,x)dx.
	\end{align}
	Similarly, we obtain
	\begin{align} \label{L2C10}
		\bar{\mathcal{W}}_2 = \int_{\Omega \times \mathbb{R}^3} \psi(x,v)h(t)dxdv = \int_{\Omega} \phi(x)c(t,x)dx.
	\end{align}
	For $\bar{\mathcal{W}}_3$, we decompose $f=Pf+\left(I-P\right)f$ to get
	\begin{align} \label{L2C11}
		\bar{\mathcal{W}}_3
		&= \int_{t}^{t+\epsilon} \int_{\Omega \times \mathbb{R}^3} \sum_{i=1}^3 v_i \partial_{x_i}\phi(x) \frac{|v|^2-3}{\sqrt{6}}\mu(v) \left[a+b\cdot v +\frac{|v|^2-3}{\sqrt{6}}c\right]dxdvds \nonumber\\
		& \quad +\int_{t}^{t+\epsilon} \int_{\Omega \times \mathbb{R}^3} v\cdot \nabla_x \phi(x) \frac{|v|^2-3}{\sqrt{6}}\mu^{1/2}(v)\left[\left(I-P\right)f(s)\right]dxdvds\nonumber\\
		&= \sum_{i=1}^3	\frac{\sqrt{6}}{3}\int_t^{t+\epsilon}\left(\int_{\Omega}\partial_{x_i} \phi(x) b_i(s,x)dx\right)ds\nonumber\\
		& \quad + \int_{t}^{t+\epsilon} \int_{\Omega \times \mathbb{R}^3} v\cdot \nabla_x \phi(x) \frac{|v|^2-3}{\sqrt{6}}\mu^{1/2}(v)\left[\left(I-P\right)f(s)\right]dxdvds,
	\end{align}
	where we have used the oddness of integration in $v$.\\
		Since $\int_{\Omega} \partial_t c(t,x)=0$ by the energy conservation, for fixed $t\ge 0$, we define $\phi(x) = \Phi_c(x)$ with
	\begin{align*}
			\begin{cases}
			-\Delta_x \Phi_c(x) = \partial_t c(t,x)\\
			\Phi_c|_{\partial \Omega} = 0.
		\end{cases}
	\end{align*}
	Then we have for fixed t,
	\begin{align} \label{L2C17}
		\Phi_c (x) = -\Delta_x ^{-1}\partial_t c(t,x) = \partial_t \phi_c(t,x).
	\end{align}
	For $\bar{\mathcal{W}_4}$, 	from the fact $\Phi_c|_{\partial \Omega} = 0$, 
	\begin{align} \label{L2C18}
		\bar{\mathcal{W}_4}=0.
	\end{align} 
	For $\bar{\mathcal{W}_5}$, it holds that
	\begin{align*}
		\int_{\mathbb{R}^3} Lf \frac{|v|^2-3}{\sqrt{6}}\mu^{1/2}(v)dv=0,
	\end{align*}
	and we obtain
	\begin{align} \label{L2C12}
		\bar{\mathcal{W}}_5 = 0.
	\end{align}
	For $\bar{\mathcal{W}_6}$, it holds that
	\begin{align*}
		\int_{\R^3} \Gamma(f,f) \frac{|v|^2-3}{\sqrt{6}}\mu^{1/2}(v)dv=0,
	\end{align*}
	and it follows that
	\begin{align} \label{L2C13}
		\bar{\mathcal{W}}_6 = 0.
	\end{align}
	For $\bar{\mathcal{W}_7}$, we easily get
	\begin{align} \label{L2C14}
		\bar{\mathcal{W}}_7 = \int_t^{t+\epsilon}\int_{\Omega \times \mathbb{R}^3}(\partial_t \psi)fdxdvds=0
	\end{align}
	since $\psi$ is independent of $t$.\\
	Gathering \eqref{L2C9}, \eqref{L2C10}, \eqref{L2C11}, \eqref{L2C18}, \eqref{L2C12}, \eqref{L2C13}, \eqref{L2C14} and taking the difference quotient in \eqref{wf}, we can derive for all $t \ge 0$,
	\begin{align*}
		\int_{\Omega} \phi(x) \partial_t c(t,x) dx &= \frac{\sqrt{6}}{3}\int_{\Omega}b(t,x)\cdot \nabla_x \phi(x) dx\\
		&\quad + \int_{\Omega \times \mathbb{R}^3} v\cdot \nabla_x \phi(x) \frac{|v|^2-3}{\sqrt{6}}\mu^{1/2}(v)\left[\left(I-P\right)f(t)\right]dxdv
	\end{align*}
	From the above equality and \eqref{L2C17}, we have for all $t \ge 0$,
	\begin{align*}
		\|\nabla_x \partial_t \phi_c(t)\|_{L^2_x}^2 &= \int_\Omega \left|\nabla_x \Phi_c(x)\right|^2 dx = -\int_\Omega \Phi_c(x) \left(\Delta_x \Phi_c(x)\right)dx= \int_\Omega \phi(x) \partial_t c(t,x) dx\\
		& \lesssim \epsilon \left\{\|\nabla_x \Phi_c\|_{L^2_x}^2+\|\Phi_c\|_{L^2_x}^2\right\} + \|b(t)\|_{L^2_x}^2+\| \left(I-P\right)f(t)\|_{L^2_{x,v}}^2,
	\end{align*}
	where we have used the integration by parts and the Young's inequality.\\
	We use the Poincar\'e inequality to obtain
	\begin{align*}
		\|\nabla_x \partial_t \phi_c(t)\|_{L^2_x}^2 \lesssim \epsilon \|\nabla_x \Phi_c\|_{L^2_x}^2 + \|b(t)\|_{L^2_x}^2+\| \left(I-P\right)f(t)\|_{L^2_{x,v}}^2.
	\end{align*}
	Choosing sufficiently small $\epsilon >0$, we have for all $t \ge 0$,
	\begin{align} \label{L2C15}
		\|\nabla_x \partial_t \phi_c(t)\|_{L^2_x} \lesssim \|b(t)\|_{L^2_x}+\| \left(I-P\right)f(t)\|_{L^2_{x,v}}.
	\end{align}
	\newline
	$\mathbf{(Estimate\ for\ b(t,x))}$\\
	In order to derive the estimate for $b(t,x)$, we first consider the estimate for $(\partial_{x_i} \partial_{x_j} \Delta_x^{-1}b_j)b_i$ for $i,j =1,2,3$. Similar to the estimate for $c(t,x)$, we consider the weak formulation \eqref{wwf} for the equation \eqref{FPBER}. Fix $i,j$. We choose the test function
	\begin{align*}
		\psi = \psi_b^{i,j}(t,x,v) = (v_i^2-\beta_b)\mu^{1/2}(v) \partial_{x_j}\phi_b^j(t,x),
	\end{align*}
	where
	\begin{align*}
		\begin{cases}
			-\Delta_x\phi_b^j(t,x) = b_j(t,x)\\
			\phi_b^j|_{\partial \Omega}=0
		\end{cases}
	\end{align*}
	and $\beta_b>0$ is chosen such that for all $i=1,2,3$,
	\begin{align*}
		\int_{\mathbb{R}^3}\left[(v_i)^2-\beta_b\right]\mu(v)dv = \frac{1}{3}\int_{\mathbb{R}^3}\left(|v|^2-3\beta_b\right)\mu(v)dv=0.
	\end{align*}
	Note that for all $i \not= k$,
	\begin{align*}
		&\int_{\mathbb{R}^3}(v_i^2-\beta_b)v_k^2\mu(v)dv = \int_{\mathbb{R}^3}(v_1^2-1)v_2^2\mu(v)dv=0,\\
		&\int_{\mathbb{R}^3}(v_i^2-\beta_b)v_i^2\mu(v)dv = \frac{1}{\sqrt{2\pi}}\int_{\mathbb{R}}(v_1^4-v_1^2)e^{-\frac{v_1^2}{2}}dv_1=2.
	\end{align*}
	From the standard elliptic estimate, we get
	\begin{align} \label{ellipticest2} 
		\|\phi_b^j(t)\|_{H_x^2} \lesssim \|b_j(t)\|_{L_x^2}
	\end{align}
	for all $t\ge 0$.\\
	For $\mathcal{W}_2$ and $\mathcal{W}_3$, we use the H\"{o}lder inequality and the elliptic estimate \eqref{ellipticest2} to obtain
	\begin{align} \label{L2B1}
		|\mathcal{W}_2| \lesssim \|f(0)\|_{L^2_{x,v}}^2, \quad |\mathcal{W}_3| \lesssim \|f(t)\|_{L^2_{x,v}}^2.
	\end{align}
	For $\mathcal{W}_4$, using the elliptic estimate \eqref{ellipticest2}, we obtain
	\begin{align} \label{L2B2}
		|\mathcal{W}_4| \lesssim  \int_0^t \|b(s)\|_{L^2_x} \|(I-P)f(s)\|_{L^2_{x,v}}ds.
	\end{align}
	For $\mathcal{W}_5$, we can easily get
	\begin{align} \label{L2B3}
		|\mathcal{W}_5| \lesssim \int_0^t	\|b(s)\|_{L^2_x} \|(I-P)f(s)\|_{L^2_{x,v}}ds.
	\end{align}
	For $\mathcal{W}_6$, we decompose it into two terms for $P_\gamma f$ and $(I-P_\gamma)f$, and use the trace theorem to obtain
	\begin{align} \label{L2B13}
		|\mathcal{W}_6|\lesssim \int_0^t \|b(s)\|_{L^2_x}\|(I-P_\gamma)f(s)\|_{L^2_{\gamma_+}}ds.
	\end{align}
	For $\mathcal{W}_8$, we use the Cauchy-Schwartz inequality to derive
	\begin{align} \label{L2B4}
		|\mathcal{W}_8| \lesssim \int_0^t \|b(s)\|_{L^2_x}\| \Gamma(f,f)(s)\|_{L^2_{x,v}}ds.
	\end{align}
	For $\mathcal{W}_1$, from the construction for $\beta_b$ and the oddness of integration in $v$, we deduce
	\begin{equation} \label{L2B5}
	\begin{aligned}
		\mathcal{W}_1
		& = \sum_{k=1}^3\int_0^t \int_{\Omega \times \mathbb{R}^3}\left(v_i^2-\beta_b\right)\mu(v)v_k \partial_{kj}\phi_b^j(s,x)\left[a(s,x)+b(s,x)\cdot v+c(s,x)\frac{|v|^2-3}{\sqrt{6}}\right]dxdvds\\
		& = \sum_{k=1}^3 \int_0^t \int_{\Omega \times \mathbb{R}^3}\left(v_i^2-\beta_b\right)\mu(v)v_k^2 \left[\partial_{kj}\phi_b^j(s,x)\right]b_k(s,x)dxdvds\\
		& = 2\int_0^t \left(\int_{\Omega}\left[\partial_{ij}\phi_b^j(s,x)\right] b_i(s,x)dx\right)ds\\
		& = -2\int_0^t \left(\int_{\Omega}\left(\partial_{ij}\Delta_x^{-1} b_j\right)(s,x) b_i(s,x)dx\right)ds
	\end{aligned}
	\end{equation}
	For $\mathcal{W}_7$, we decompose $f=Pf+\left(I-P\right)f$ to get
	\begin{align*}
		&\int_0^t \int_{\Omega\times \mathbb{R}^3} f\left(\partial_t\psi\right)dxdvds\\ 
		&=\int_0^t\int_{\Omega \times \mathbb{R}^3}(v_i^2-\beta_b)\mu(v)\left[\partial_t\partial_{x_j}\phi_b^j\right]\left[a(s,x)+b(s,x)\cdot v+c(s,x)\frac{|v|^2-3}{\sqrt{6}}\right]dxdvds\\
		&\quad +\int_0^t\int_{\Omega \times \mathbb{R}^3} (v_i^2-\beta_b)\mu^{\frac{1}{2}}(v)\left[\partial_t\partial_{x_j}\phi_b^j\right](I-P)f(s)dxdvds.
	\end{align*}
	From the construction for $\beta_b$ and the oddness of integration in $v$, the above expression becomes
	\begin{align*}
		\mathcal{W}_7 = &\frac{2}{\sqrt{6}}\int_0^t\int_{\Omega \times \mathbb{R}^3} \left[\partial_t\partial_{x_j}\phi_b^j(s,x)\right]c(s,x)dxdvds\\
		&+\int_0^t\int_{\Omega \times \mathbb{R}^3}(v_i^2-\beta_b)\mu^{\frac{1}{2}}(v) \left[\partial_t\partial_{x_j}\phi_b^j(s,x)\right] (I-P)f(s)dxdvds
	\end{align*}
	Lately, we will demonstrate the estimate \eqref{L2B30} for $\nabla_x \partial_t \phi_b^j$ as following :
	\begin{align*}
		\left\|\nabla_x \partial_t \phi_b^j(t)\right\|_{L^2_x} \lesssim \|a(t)\|_{L_x^2}+\|c(t)\|_{L_x^2}+\|(I-P)f(t)\|_{L_{x,v}^2}.
	\end{align*}
	By the estimate \eqref{L2B30} for $\nabla_x \partial_t \phi_b^j$, we get
	\begin{align} \label{L2B6}	
		|\mathcal{W}_7|
		& \lesssim \int_0^t \left\|\nabla_x \partial_t \phi_b^j(s)\right\|_{L_x^2} \left(\|(I-P)f(s)\|_{L_{x,v}^2} + \|c(s)\|_{L_x^2}\right)ds \nonumber\\
		& \lesssim \int_0^t \left(\|a(s)\|_{L_x^2}+\|c(s)\|_{L_x^2}+\|(I-P)f(s)\|_{L_{x,v}^2}\right) \left(\|(I-P)f(s)\|_{L_{x,v}^2} + \|c(s)\|_{L_x^2}\right)ds\nonumber\\
		& \lesssim \epsilon \int_0^t \|a(s)\|^2_{L_x^2}ds + C(\epsilon)\int_0^t\|c(s)\|^2_{L_x^2}ds + C(\epsilon)\int_0^t \|(I-P)f(s)\|_{L_{x,v}^2}^2ds,
	\end{align}
	where we have used the Young's inequality.\\
	Gathering \eqref{L2B1}, \eqref{L2B2}, \eqref{L2B3}, \eqref{L2B13}, \eqref{L2B4}, \eqref{L2B5}, \eqref{L2B6}, we obtain
	\begin{align} \label{bo1}
		&\int_0^t\int_{\Omega}\left(\partial_{x_i}\partial_{x_j} \Delta_x^{-1} b_j\right)(s,x)b_i(s,x)dxds\nonumber\\ 
		&\le G_f^{b}(t)-G_f^{b}(0) +C\int_0^t \|b(s)\|_{L^2_x}\left\|(I-P)f(s)\right\|_{L^2_{x,v}}ds+C\int_0^t \|b(s)\|_{L^2_x}\|(I-P_\gamma)f(s)\|_{L^2_{\gamma_+}}ds\nonumber\\
		&\quad  +C\int_0^t\|b(s)\|_{L^2_x}\| \Gamma(f,f)(s)\|_{L^2_{x,v}}ds+ \epsilon \int_0^t\|a(s)\|_{L_x^2}^2ds + C(\epsilon)\int_0^t\|c(s)\|^2_{L_x^2}ds\nonumber\\
		&\quad + C(\epsilon)\int_0^t \|(I-P)f(s)\|_{L_{x,v}^2}^2ds\nonumber\\
		&\le G_f^{b}(t)-G_f^{b}(0) +\epsilon\int_0^t \|b(s)\|_{L^2_x}^2ds + \epsilon \int_0^t\|a(s)\|_{L_x^2}^2ds\nonumber\\
		&\quad +C(\epsilon)\int_0^t\|(I-P)f(s)\|_{L_{x,v}^2}^2ds+C(\epsilon)\int_0^t\|(I-P_\gamma)f(s)\|_{L^2_{\gamma_+}}^2ds \nonumber\\
		& \quad  +C(\epsilon)\int_0^t\| \Gamma(f,f)(s)\|_{L^2_{x,v}}^2ds +C(\epsilon)\int_0^t \|c(s)\|_{L^2_x}^2ds,
	\end{align}
	where we have used the Young's inequality and $G_f^b(s) \lesssim \|f(s)\|_{L^2_{x,v}}^2$.\\
	We now consider the estimate for $(\partial_j\partial_j\Delta_x^{-1}b_i)b_i$ for $i \not= j$. We choose the test function
	\begin{align*}
		\psi = \psi_b^{i,j}(t,x,v) = |v|^2v_iv_j \mu^{\frac{1}{2}}(v) \partial_j\phi_b^i(t,x),
	\end{align*}
	where
	\begin{align*}
		\begin{cases}
			-\Delta_x \phi_b^i(t,x) = b_i(t,x)\\
			\phi_b^i |_{\partial \Omega}=0
		\end{cases}.
	\end{align*}
	For $\mathcal{W}_2$ and $\mathcal{W}_3$, using the elliptic estimate \eqref{ellipticest2}, we obtain
	\begin{align} \label{L2B7}
		|\mathcal{W}_2| \lesssim \|f(0)\|_{L^2_{x,v}}^2, \quad |\mathcal{W}_3| \lesssim \|f(t)\|_{L^2_{x,v}}^2
	\end{align}
	For $\mathcal{W}_4$, using the elliptic estimate \eqref{ellipticest2}, we obtain
	\begin{align} \label{L2B8}
		|\mathcal{W}_4| \lesssim  \int_0^t \|b(s)\|_{L^2_x} \|(I-P)f(s)\|_{L^2_{x,v}}ds.
	\end{align}
	For $\mathcal{W}_5$, we can easily get
	\begin{align} \label{L2B9}
		|\mathcal{W}_5| \lesssim \int_0^t	\|b(s)\|_{L^2_x} \|(I-P)f(s)\|_{L^2_{x,v}}ds.
	\end{align}
	For $\mathcal{W}_6$, we decompose it into two terms for $P_\gamma f$ and $(I-P_\gamma)f$, and use the trace theorem to derive 
	\begin{align} \label{L2B14}
		|\mathcal{W}_6|\lesssim \int_0^t \|b(s)\|_{L^2_x}\|(I-P_\gamma)f(s)\|_{L^2_{\gamma_+}}ds.
	\end{align}
	For $\mathcal{W}_8$, we use the Cauchy-Schwartz inequality to derive
	\begin{align} \label{L2B10}
		|\mathcal{W}_8| \lesssim \int_0^t \|b(s)\|_{L^2_x}\|\Gamma(f,f)(s)\|_{L^2_{x,v}}ds.
	\end{align}
	For $\mathcal{W}_1$, from the oddness of integration in $v$, we deduce that
	\begin{align} \label{L2B11}
		\mathcal{W}_1
		& = \sum_{k=1}^3\int_0^t \int_{\Omega \times \mathbb{R}^3}|v|^2v_iv_jv_k\mu(v) \partial_{kj}\phi_b^i(s,x)\left[a(s,x)+b(s,x)\cdot v+c(s,x)\frac{|v|^2-3}{\sqrt{6}}\right]dxdvds\nonumber\\
		& = \int_0^t \int_{\Omega \times \mathbb{R}^3} |v|^2v_i^2v_j^2\mu(v)\left[\partial_{ij}\phi^i_b(s,x) b_j(s,x) + \partial_{jj} \phi_b^i(s,x)b_i(s,x) \right]dxdvds\nonumber\\
		& = 7\int_0^t \int_{\Omega}\partial_{ij}\phi^i_b(s,x) b_j(s,x) + \partial_{jj} \phi_b^i(s,x)b_i(s,x)dxds\nonumber\\
		& = -7\int_0^t \int_{\Omega}\left(\partial_{ij}\Delta_x^{-1} b_i\right)(s,x) b_j(s,x)+\left(\partial_{jj}\Delta_x^{-1}b_i\right)(s,x)b_i(s,x)dxds.
	\end{align}
	For $\mathcal{W}_7$, we decompose $f=Pf+\left(I-P\right)f$ to get
	\begin{align*}
		&\int_0^t \int_{\Omega \times \mathbb{R}^3} f\left(\partial_t\psi\right)dxdvds\\ 
		&=\int_0^t\int_{\Omega \times \mathbb{R}^3}|v|^2v_iv_j\mu(v)\left[\partial_t\partial_{x_j}\phi_b^i\right]\left[a(s,x)+b(s,x)\cdot v+c(s,x)\frac{|v|^2-3}{\sqrt{6}}\right]dxdvds\\
		&\quad +\int_0^t\int_{\Omega \times \mathbb{R}^3} |v|^2v_iv_j\mu^{\frac{1}{2}}(v)\left[\partial_t\partial_{x_j}\phi_b^i\right](I-P)f(s)dxdvds.
	\end{align*}
	From the oddness of integration in $v$, the above expression becomes
	\begin{align*}
		\mathcal{W}_7 = \int_0^t\int_{\Omega \times \mathbb{R}^3}|v|^2v_iv_j\mu^{\frac{1}{2}}(v) \left[\partial_t\partial_{x_j}\phi_b^i(s,x)\right] (I-P)f(s)dxdvds
	\end{align*}
	Lately, we will demonstrate the estimate \eqref{L2B30} for $\nabla_x \partial_t \phi_b^i$ as following :
	\begin{align*} 
		\left\|\nabla_x \partial_t \phi_b^i(t)\right\|_{L^2_x} \lesssim \|a(t)\|_{L_x^2}+\|c(t)\|_{L_x^2}+\|(I-P)f(t)\|_{L_{x,v}^2}.
	\end{align*}
	By the estimate \eqref{L2B30} for $\nabla_x \partial_t \phi_b^i$, we get
	\begin{align} \label{L2B12}
		\mathcal{W}_7
		& \lesssim \int_0^t \left\|\nabla_x \partial_t \phi_b^i(s)\right\|_{L_x^2}\|(I-P)f(s)\|_{L_{x,v}^2} ds\nonumber\\
		& \lesssim \int_0^t \left(\|a(s)\|_{L_x^2}+\|c(s)\|_{L_x^2}+\|(I-P)f(s)\|_{L_{x,v}^2}\right)\|(I-P)f(s)\|_{L_{x,v}^2}ds\nonumber\\
		& \lesssim \epsilon \int_0^t \|a(s)\|^2_{L_x^2}ds + \epsilon \int_0^t\|c(s)\|^2_{L_x^2}ds + C(\epsilon)\int_0^t \|(I-P)f(s)\|_{L_{x,v}^2}^2ds ,
	\end{align}
	where we have used the Young's inequality.\\
	Gathering \eqref{L2B7}, \eqref{L2B8}, \eqref{L2B9}, \eqref{L2B14},  \eqref{L2B10}, \eqref{L2B11}, \eqref{L2B12}, we obtain for all $i\not= j$,
	\begin{align} \label{bo2}
		&\int_0^t\int_{\Omega}\left(\partial_{ij}\Delta_x^{-1} b_i\right)(s,x) b_j(s,x)+\left(\partial_{jj}\Delta_x^{-1}b_i\right)(s,x)b_i(s,x)dxds\nonumber\\ 
		&\le G_f^{b}(t)-G_f^{b}(0) +C\int_0^t \|b(s)\|_{L^2_x}\left\|(I-P)f(s)\right\|_{L^2_{x,v}}ds+\int_0^t \|b(s)\|_{L^2_x}\|(I-P_\gamma)f(s)\|_{L^2_{\gamma_+}}ds \nonumber\\
		&\quad+C\int_0^t\|b(s)\|_{L^2_x}\|\Gamma(f,f)(s)\|_{L^2_{x,v}}ds + \epsilon \int_0^t\|a(s)\|_{L_x^2}^2ds + \epsilon \int_0^t\|c(s)\|^2_{L_x^2}ds\nonumber\\
		& \quad + C(\epsilon)\int_0^t \|(I-P)f(s)\|_{L_{x,v}^2}^2ds \nonumber\\
		&\le G_f^{b}(t)-G_f^{b}(0) +\epsilon\int_0^t \|b(s)\|_{L^2_x}^2ds + \epsilon \int_0^t\|a(s)\|_{L_x^2}^2ds + \epsilon\int_0^t \|c(s)\|_{L^2_x}^2ds\nonumber\\
		&\quad+C(\epsilon)\int_0^t\|(I-P_\gamma)f(s)\|_{L^2_{\gamma_+}}^2ds  +C(\epsilon)\int_0^t\|(I-P)f(s)\|_{L_{x,v}^2}^2ds\nonumber\\
		& \quad +C(\epsilon) \int_0^t\| \Gamma(f,f)(s)\|_{L^2_{x,v}}^2ds,
	\end{align}
	where we have used the Young's inequality and $|G_f^b(s)| \lesssim \|f(s)\|_{L^2_{x,v}}$.\\
	Combining the estimate \eqref{bo2} with the estimate \eqref{bo1} for $\partial_{x_i}\partial_{x_j}\left(\Delta_x^{-1}b_j\right)b_i$, for all $i \not= j$,
	\begin{align} \label{bo3}
		&\int_0^t \int_{\Omega}\left(\partial_{jj}\Delta_x^{-1}b_i\right)(s,x)b_i(s,x)dxds\nonumber\\
		&\le G_f^{b}(t)-G_f^{b}(0) +\epsilon\int_0^t \|b(s)\|_{L^2_x}^2ds + \epsilon \int_0^t\|a(s)\|_{L_x^2}^2ds + \left(\epsilon+C(\epsilon)\right)\int_0^t \|c(s)\|_{L^2_x}^2ds\nonumber\\
		&\quad +C(\epsilon)\int_0^t\|(I-P)f(s)\|_{L_{x,v}^2}^2ds+C(\epsilon)\int_0^t\|(I-P_\gamma)f(s)\|_{L^2_{\gamma_+}}^2ds\nonumber\\
		&\quad  +C(\epsilon) \int_0^t\| \Gamma(f,f)(s)\|_{L^2_{x,v}}^2ds.
	\end{align}
	From the estimates for $\partial_{jj}(\Delta_x^{-1}b_j)b_j$ and $\partial_{jj}\left(\Delta_x^{-1}b_i\right)b_i$, summing over $j=1,2,3$,
	\begin{equation} \label{co2}
	\begin{aligned}
		\int_0^t \|b(s)\|_{L_x^2}^2ds &\le G_f^{b}(t)-G_f^{b}(0)+\epsilon\int_0^t \|b(s)\|_{L^2_x}^2ds + \epsilon \int_0^t\|a(s)\|_{L_x^2}^2ds \\ 
		& \quad + \left(\epsilon+C(\epsilon)\right)\int_0^t \|c(s)\|_{L^2_x}^2ds +C(\epsilon)\int_0^t\|(I-P)f(s)\|_{L_{x,v}^2}^2ds\\
		& \quad+C(\epsilon)\int_0^t\|(I-P_\gamma)f(s)\|_{L^2_{\gamma_+}}^2ds +C(\epsilon) \int_0^t\| \Gamma(f,f)(s)\|_{L^2_{x,v}}^2ds .
	\end{aligned}
	\end{equation}
	\newline
	$\mathbf{(Estimate\ for \ \nabla_x\partial_t\phi_b^i)}$\\
	We consider the weak formulation \eqref{wf} over $[t,t+\epsilon]$. We choose the test function
	\begin{align*}
		\psi = \psi(x,v) = \phi(x)v_i\mu^{\frac{1}{2}}(v),
	\end{align*}
	where $\phi(x)$ depends only on $x$. Note that
	\begin{align*}
		\int_{\mathbb{R}^3}v_iv_j\mu(v)dv = \delta_{ij},\quad \int_{\mathbb{R}^3}v_iv_j\frac{|v|^2-3}{\sqrt{6}}\mu(v)dv = \sqrt{6}\delta_{ij}.
	\end{align*}
	For $\bar{\mathcal{W}}_1$ and $\bar{\mathcal{W}}_2$, we have
	\begin{align} \label{L2B15}
		\bar{\mathcal{W}}_1 = \int_{\Omega}\phi(x)b_i(t+\epsilon,x)dx, \quad \bar{\mathcal{W}}_2 = \int_{\Omega} \phi(x)b_i(t,x)dx.
	\end{align}
	For $\bar{\mathcal{W}}_3$, we decompose $h=Pf+(I-P)f$ to get
	\begin{align} \label{L2B16}
		\bar{\mathcal{W}}_3&= \sqrt{6} \int_t^{t+\epsilon}\int_{\Omega} \partial_{x_i} \phi(x) \left[(a(s,x)+c(s,x)\right]dxds\nonumber\\
		&\quad + \int_{t}^{t+\epsilon} \int_{\Omega \times \mathbb{R}^3} v \cdot \nabla_x \phi(x)v_i\mu^{1/2}(v)\left[\left(I-P\right)f(s)\right]dxdvds,
	\end{align}
	where we have used the oddness of integration in $v$.\\
	For fixed $t>0$, define $\phi(x) = \Phi_b^i(x)$ with
	\begin{align*}
		\begin{cases}
			-\Delta_x \Phi_b^i(x) = \partial_t b_i(t,x)\\
			\Phi_b^i|_{\partial \Omega} = 0
		\end{cases}.
	\end{align*}
	Then we have for fixed t,
	\begin{align} \label{L2B17}
		\Phi_b^i (x) = -\Delta_x ^{-1}\partial_t b_i(t,x) = \partial_t \phi_b^i(t,x).
	\end{align}
	For $\bar{\mathcal{W}}_4$, from the fact $\Phi_b^i|_{\partial \Omega} = 0$, we get
	\begin{align} \label{L2B18}
		\bar{\mathcal{W}}_4=0.
	\end{align} 
	For $\bar{\mathcal{W}}_5$, $\bar{\mathcal{W}}_6$, and $\bar{\mathcal{W}}_7$, it holds that
	\begin{align} \label{L2B19}
		\bar{\mathcal{W}}_5 = 0, \quad \bar{\mathcal{W}}_6 =0, \quad \bar{\mathcal{W}}_7=0.
	\end{align}
	Gathering \eqref{L2B15}, \eqref{L2B16}, \eqref{L2B18}, \eqref{L2B19} and taking the difference quotient in \eqref{wf}, for all $t \ge 0$,
	\begin{align*}
		\int_{\Omega} \phi(x) \partial_t b_i(t,x) dx &\lesssim \left(\|a(t)\|_{L_x^2} + \|c(t)\|_{L_x^2} \right)\|\nabla_x \phi\|_{L_x^2} +\|(I-P)f(t)\|_{L^2_{x,v}}^2\|\nabla_x \phi\|_{L_x^2}\\
		& \lesssim \left(\|a(t)\|_{L_x^2} + \|c(t)\|_{L_x^2} +\|(I-P)f(t)\|_{L^2_{x,v}}^2\right)\|\phi\|_{L_x^2}.
	\end{align*}
	where we have used the Young's inequality and the Poincar\'e inequality. This implies that
	\begin{align*}
		\|\partial_t b_i(t)\|_{(H_x^1)^*} \lesssim \|a(t)\|_{L_x^2} + \|c(t)\|_{L_x^2} +\|(I-P)f(t)\|_{L^2_{x,v}}^2.
	\end{align*}
	From the above equality and \eqref{L2B17}, for all $t \ge 0$,
	\begin{align*}
		\|\nabla_x \partial_t \phi_b^i(t)\|_{L^2_x}^2 &= \int_\Omega \left|\nabla_x \Phi_b^i(x)\right|^2 dx= -\int_\Omega \Phi_b^i(x) \left(\Delta_x \Phi_b^i(x)\right)dx= \int_\Omega \phi(x) \partial_t b_i(t,x) dx\\
		& \lesssim \epsilon \left\{\|\nabla_x \Phi_b^i\|_{L^2_x}^2+\|\Phi_b^i\|_{L^2_x}^2\right\} + \|a(t)\|_{L^2_x}^2+ \|c(t)\|_{L^2_x}^2+\|(I-P)f(t)\|_{L^2_{x,v}}^2,
	\end{align*}
	where we have used the integration by parts and the Young's inequality. We use the Poincar\'e inequality to obtain
	\begin{align*}
		\|\nabla_x \partial_t \phi_b^i(t)\|_{L^2_x}^2 \lesssim \epsilon \|\nabla_x \Phi_b^i\|_{L^2_x}^2 + \|a(t)\|_{L^2_x}^2+\|c(t)\|_{L^2_x}^2+\|(I-P)f(t)\|_{L^2_{x,v}}^2.
	\end{align*}
	Choosing sufficiently small $\epsilon >0$, we have for all $t \ge 0$,
	\begin{align} \label{L2B30}
		\|\nabla_x \partial_t \phi_b^i(t)\|_{L^2_x} \lesssim \|a(t)\|_{L^2_x}+\|c(t)\|_{L^2_x}+\|(I-P)f(t)\|_{L^2_{x,v}}^2.
	\end{align}
	\newline
	$\mathbf{(Estimate \ for \ a(t,x))}$\\
	Similar to the estimate for $c(t,x)$, we consider the weak formulation \eqref{wwf} for the equation \eqref{FPBER}. Since $\int_\Omega a(t,x)dx = \int_{\Omega \times \mathbb{R}^3} \mu^{1/2}(v)f(t,x,v)dxdv = 0$ by the mass conservation, we choose the test function
	\begin{align*}
		\psi = \psi_a(t,x,v) = \left(|v|^2-\beta_a\right)\mu^{1/2}(v)v\cdot\nabla_x \phi_a(t,x),
	\end{align*}
	where
	\begin{align*}
			\begin{cases}
			-\Delta_x \phi_a(t,x) = a(t,x)\\
			\frac{\partial}{\partial n}\phi_a |_{\partial \Omega} = 0
		\end{cases}
	\end{align*}
	and $\beta_a>0$ is chosen such that
	\begin{align*}
		\int_{\mathbb{R}^3}\left(|v|^2-\beta_a\right)\frac{|v|^2-3}{\sqrt{6}}v_i^2\mu(v)dv = 0 \quad \text{for all }i=1,2,3.
	\end{align*}
	From the standard elliptic estimate, we get 
	\begin{align} \label{ellipticest3}
		\|\phi_a(t)\|_{H_x^2} \lesssim \|a(t)\|_{L_x^2}
	\end{align}
	for all $t\ge 0$.\\
	For $\bar{\mathcal{W}}_2$ and $\bar{\mathcal{W}}_3$, we use the H\"{o}lder inequality and the elliptic estimate \eqref{ellipticest3} to obtain
	\begin{align} \label{L2A1}
		|\mathcal{W}_2| \lesssim \|f(0)\|_{L^2_{x,v}}^2, \quad |\mathcal{W}_3| \lesssim \|f(t)\|_{L^2_{x,v}}^2.
	\end{align}
	For $\mathcal{W}_4$ and $\mathcal{W}_5$, using the elliptic estimate \eqref{ellipticest3}, we get
	\begin{align} \label{L2A2}
		|\mathcal{W}_4|\lesssim \int_0^t \|a(s)\|_{L^2_x}\left\|(I-P)f(s)\right\|_{L^2_{x,v}}ds, \quad |\mathcal{W}_5|\lesssim \int_0^t \|a(s)\|_{L^2_x}\left\|(I-P)f(s)\right\|_{L^2_{x,v}}ds.
	\end{align}
	For $\mathcal{W}_6$, we decompose it into two terms for $P_\gamma f$ and $(I-P_\gamma)f$, and use the trace theorem to obtain
	\begin{align} \label{L2A3}
		|\mathcal{W}_6|\lesssim \int_0^t \|a(s)\|_{L^2_x}\|(I-P_\gamma)f(s)\|_{L^2_{\gamma_+}}ds.
	\end{align}
	For $\mathcal{W}_8$, we use the Cauchy-Schwartz inequality to derive
	\begin{align} \label{L2A4}
		|\mathcal{W}_8| \lesssim \int_0^t \|a(t)\|_{L^2_{x,v}}\|\Gamma(f,f)(s)\|_{L^2_{x,v}}ds.
	\end{align}
	For $\mathcal{W}_1$, from the construction for $\beta_a$ and the oddness of integration in $v$, we deduce
	\begin{align} \label{L2A5}
		\mathcal{W}_1
		& = \sum_{i=1}^3 \int_0^t \int_{\Omega \times\mathbb{R}^3} \left(|v|^2-\beta_a\right)\mu(v)v_i^2\partial_{ii}\phi_a(s,x) a(s,x)dxdvds\nonumber\\
		& = 5\int_0^t\int_{\Omega}\left(-\Delta_x\phi_a(s,x)\right)a(s,x)dxds\nonumber\\
		& = 5\int_0^t \left\|a(s)\right\|_{L_x^2}^2ds
	\end{align}
	For $\mathcal{W}_7$, we decompose $f=Pf + (I-P)f$ to get
	\begin{align*}
		\mathcal{W}_7 
		&= -5 \sum_{i=1}^3 \int_0^t \int_{\Omega}\partial_t\partial_{x_i} \phi_a(s,x)b_{i}(s,x)dxds \nonumber \\
		&\quad +\sum_{i=1}^3\int_0^t\int_{\Omega  \times\mathbb{R}^3} \left(|v|^2-\beta_a\right)\mu^{1/2}(v)v_i\left[\partial_t\partial_{x_i}\phi_a\right]\left[(I-P)f(s)\right]dxdvds,
	\end{align*}
	where we have used the construction for $\beta_a$ and the oddness of integration in $v$. Lately, we will demonstrate the estimate  \eqref{L2A30} for $\nabla_x \partial_t \phi_a$ as following :
	\begin{align*}
		\left\|\nabla_x \partial_t \phi_a(t)\right\|_{L^2_x} \lesssim \|b(t)\|_{L_x^2}.
	\end{align*}
	By the estimate \eqref{L2A30} for $\nabla_x \partial_t \phi_a$, we get
	\begin{align} \label{L2A6}
		|\mathcal{W}_7|
		& \lesssim \int_0^t \left\|\nabla_x \partial_t \phi_a(s)\right\|_{L^2_x}\left(\|(I-P)f(s)\|_{L_{x,v}^2}+\|b(s)\|_{L^2_x}\right)ds\nonumber\\
		& \lesssim \int_0^t \|b(s)\|_{L_x^2} \left(\|(I-P)f(s)\|_{L_{x,v}^2}+\|b(s)\|_{L^2_x}\right)ds\nonumber\\
		& \lesssim  \int_0^t\|b(s)\|_{L_x^2}^2ds+\int_0^t\|(I-P)f(s)\|_{L_{x,v}^2}^2ds ,
	\end{align}
	where we have used the Young's inequality. Gathering \eqref{L2A1}, \eqref{L2A2}, \eqref{L2A3}, \eqref{L2A4}, \eqref{L2A5}, and \eqref{L2A6}, we obtain
	\begin{align*}
		\int_0^t\|a(s)\|_{L_x^2}^2ds &\le G_f^a(t)-G_f^a(0) +C\int_0^t \|a(s)\|_{L^2_x}\left\|(I-P)f(s)\right\|_{L^2_{x,v}}ds\\
		&\quad+C\int_0^t \|a(s)\|_{L^2_x}\|(I-P_\gamma)f(s)\|_{L^2_{\gamma_+}}ds +C\int_0^t \|a(t)\|_{L^2_{x,v}}\|\Gamma(f,f)(s)\|_{L^2_{x,v}}ds\\
		&\quad + C \int_0^t\|b(s)\|_{L_x^2}^2ds+C\int_0^t\|(I-P)f(s)\|_{L_{x,v}^2}^2ds \\
		&\le  G_f^a(t)-G_f^a(0) +\epsilon\int_0^t \|a(s)\|_{L^2_x}^2ds + C \int_0^t\|b(s)\|_{L_x^2}^2ds\\
		&\quad +C(\epsilon)\int_0^t\|(I-P)f(s)\|_{L_{x,v}^2}^2ds+C(\epsilon)\int_0^t\|(I-P_\gamma)f(s)\|_{L^2_{\gamma_+}}^2ds\\
		&\quad +C(\epsilon)\int_0^t\|\Gamma(f,f)(s)\|_{L^2_{x,v}}^2ds,
	\end{align*}
	where we have used the Young's inequality and $G_f^a(s) \lesssim \|f(s)\|_{L^2_{x,v}}^2$.\\
	Thus choosing sufficiently small $\epsilon >0$, we conclude 
	\begin{align} \label{co3}
		\int_0^t\|a(s)\|_{L_x^2}^2ds &\lesssim G_f^a(t)-G_f^a(0) + \int_0^t\|b(s)\|_{L_x^2}^2ds+\int_0^t\|(I-P)f(s)\|_{L_{x,v}^2}^2ds\nonumber\\
		&\quad +\int_0^t\|(I-P_\gamma)f(s)\|_{L^2_{\gamma_+}}^2ds+\int_0^t\|\Gamma(f,f)(s)\|_{L^2_{x,v}}^2ds.
	\end{align}
	\newline
	$\mathbf{(Estimate \ for \ \nabla_x\partial_t\phi_a)}$\\
	We consider the weak formulation \eqref{wf} over $[t,t+\epsilon]$.
	We choose the test function
	\begin{align*}
		\psi = \psi(x,v) = \phi(x) \mu^{\frac{1}{2}}(v),
	\end{align*}
	where $\phi(x)$ depends only on $x$.\\
	For $\bar{\mathcal{W}}_1$ and $\bar{\mathcal{W}}_2$, we have
	\begin{align} \label{L2A15}
		\bar{\mathcal{W}}_1 = \int_{\Omega}\phi(x)a(t+\epsilon,x)dx, \quad \bar{\mathcal{W}}_2 = \int_{\Omega} \phi(x)a(t,x)dx.
	\end{align}
	For $\bar{\mathcal{W}}_3$, we decompose $f=Pf+\left(I-P\right)f$ to get
	\begin{align} \label{L2A17}
		\bar{\mathcal{W}}_3 = \sum_{i=1}^3	\int_t^{t+\epsilon}\left(\int_{\Omega }\partial_{x_i} \phi(x) b_i(s,x)dx\right)ds,
	\end{align}
	where we have used the oddness of integration in $v$ and $\int_{\mathbb{R}^3}v_i\mu^{1/2}(v)[\left(I-P\right)f(s)]dv=0$.\\
	For $\bar{\mathcal{W}}_5$, $\bar{\mathcal{W}}_6$, and $\bar{\mathcal{W}}_7$, it holds that
	\begin{align} \label{L2A16}
		\bar{\mathcal{W}}_5 = 0, \quad \bar{\mathcal{W}}_6 =0, \quad \bar{\mathcal{W}}_7=0.
	\end{align}
	For $\bar{\mathcal{W}}_4$, we have
	\begin{align*}
		\int_{\gamma} \psi f \{n(x) \cdot v\}dS(x)dv=\int_{\gamma_+} \psi \left[\left(I-P_\gamma\right)f\right]\{n(x) \cdot v\}dS(x)dv + \int_{\gamma} \psi \left(P_\gamma f\right)\{n(x) \cdot v\}dS(x)dv.
	\end{align*}
	By the oddness of integration in $v$, we get
	\begin{align*}
		\int_{\gamma} \psi \left(P_\gamma f\right)\{n(x) \cdot v\}dS(x)dv &= \sum_{i=1}^3 \int_{\partial \Omega} \phi(x)z(t,x)n_i(x)\left(\int_{\mathbb{R}^3}v_i\mu(v)dv\right)dS(x)=0,
	\end{align*}
	where $z(t,x) =c_\mu \int_{n(x) \cdot v' >0}f(x,v')\mu^{\frac{1}{2}}(v')\{n(x)\cdot v'\}dv'$. From the fact $(I-P_\gamma) \perp P_\gamma$, we obtain
	\begin{align*}
		\int_{\gamma_+} \psi \left[\left(I-P_\gamma\right)f\right]\{n(x) \cdot v\}dS(x)dv
		&=\int_{\partial \Omega} \phi(x) \left(\int_{n(x) \cdot v >0} \mu^{1/2}(v) \left[\left(I-P_\gamma\right)f\right]\{n(x) \cdot v\}dv\right)dS(x)\\
		&=0.
	\end{align*}
	This yields
	\begin{align} \label{L2A18}
		\bar{\mathcal{W}}_4=0.
	\end{align}
	Gathering \eqref{L2A15}, \eqref{L2A17}, \eqref{L2A16}, \eqref{L2A18}, and taking the difference quotient in \eqref{wf}, for all $t \ge 0$,
	\begin{align*}
		\int_{\Omega } \phi(x) \partial_t a(t,x) dx &= \int_{\Omega }b(t,x)\cdot \nabla_x \phi(x) dx\lesssim \|b(t)\|_{L_x^2} \|\nabla_x \phi\|_{L_x^2}\lesssim \|b(t)\|_{L_x^2}  \|\phi\|_{H_x^1}.
	\end{align*}
	where we have used the Poincar\'e inequality. Thus we get
	\begin{align*}
		\|\partial_t a(t)\|_{(H_x^1)^*} \lesssim \|b(t)\|_{L_x^2},
	\end{align*}
	where $(H_x^1)^*$ is dual space of $H_x^1$.\\
	Since $\int_{\Omega} \partial_t a(t,x)dx=0$ by the mass conservation, for fixed $t>0$, define $\phi(x) = \Phi_a(x)$ with
	\begin{align*}
		\begin{cases}
			-\Delta_x \Phi_a(x) = \partial_t a(t,x)\\
			\frac{\partial}{\partial n}\Phi_a|_{\partial \Omega} = 0.
		\end{cases}
	\end{align*}
	Then we have for fixed t,
	\begin{align*}
		\Phi_a (x) = -\Delta_x^{-1}\partial_t a(t,x) = \partial_t \phi_a(t,x).
	\end{align*}
	Hence we use the standard elliptic estimate to obtain
	\begin{align} \label{L2A30}
		\|\nabla_x \partial_t \phi_a(t)\|_{L^2_x} &= \|\nabla_x \Delta_x^{-1}\partial_t a(t)\|_{L^2_x}= \|\Delta_x^{-1}\partial_t a(t)\|_{H_x^1}=\|\Phi_a(t)\|_{H_x^1} \nonumber\\
		& \lesssim \|\partial_t a(t)\|_{(H_x^1)^*}\nonumber\\
		&\lesssim \|b(t)\|_{L^2_x} +\|g(t)\|_{L^2_{x,v}}.
	\end{align}
	\newline
	$\mathbf{(Conclusion)}$\\
	From \eqref{co1}, \eqref{co2}, and \eqref{co3}, for $\eta,\delta >0$, we have
	\begin{align*}
		&\int_0^t\left(\|a(s)\|_{L_x^2}^2+\eta\|b(s)\|_{L_x^2}^2+\delta\|c(s)\|_{L_x^2}^2\right)ds\\
		&\le G_f(t)-G_f(0) + \eta \epsilon_b\int_0^t\|a(s)\|_{L_x^2}^2ds + \left(\delta \epsilon_c + C_a\right)\int_0^t\|b(s)\|_{L_x^2}^2ds + \eta \left(\epsilon_b+C(\epsilon_b)\right)\int_0^t\|c(s)\|_{L_x^2}^2ds\\
		&\quad +C\int_0^t\left\|(I-P)f(s)\right\|_{L_{x,v}^2}^2ds + C\int_0^t \|(I-P_\gamma)f(s)\|_{L^2_{\gamma_+}}^2ds+C\int_0^t \|\Gamma(f,f)(s)\|_{L^2_{x,v}}^2ds,
	\end{align*}
	where $G_f(s ) \lesssim \|f(s)\|^2_{L^2_{x,v}}$.
	First of all, we choose large $\eta>0$ such that $\eta > C_a$, and then choose small $\epsilon_b >0$ such that $1>\eta \epsilon_b$. Next, we choose large $\delta >0$ such that $\delta > \eta\left(\epsilon_b + C(\epsilon_b)\right)$. Lastly, we choose small $\epsilon_c>0$ such that $\eta > C_a + \delta \epsilon_c$. Therefore we conclude 
	\begin{align*}
		\int_0^t\left\|Pf(s)\right\|^2_{L^2_{x,v}}ds &\lesssim G_f(t) - G_f(0) +\int_0^t\left[\left\|\left(I-P\right)f(s)\right\|^2_{L^2_{x,v}}+\left\|\left(I-P_\gamma\right)f(s)\right\|^2_{L_{\gamma+}^2} \right] ds\\
		&\quad +\int_0^t \|\Gamma(f,f)(s)\|_{L^2_{x,v}}^2ds.
	\end{align*}

\end{proof}

\noindent{\bf Data availability:} No data was used for the research described in the article.
\newline

\noindent{\bf Conflict of interest:} The authors declare that they have no conflict of interest.\newline

\noindent{\bf Acknowledgement}
J. Kim and D. Lee are supported by the National Research Foundation of Korea(NRF) grant funded by the Korea government(MSIT)(No.RS-2023-00212304 and No.RS-2023-00219980). Dingqun Deng is supported by JSPS KAKENHI Grant Number JP25K23329.

\bibliographystyle{abbrv}
\bibliography{Lpframe.bib}
\end{document}